\documentclass[tikz, sepfignums]{nmd/nmd-article}
\pdfoutput=1 

\usepackage{xparse}
\usepackage{tikz-cd}
\usepackage{enumitem}

\usetikzlibrary{arrows, shapes}
\tikzset{
  commutative diagrams/.cd,
  arrow style=tikz,
  diagrams={>=stealth', line width=0.7pt},
  every label/.append style = {font = \small}
}

\definecolor{locus}{HTML}{8B008B}
\definecolor{simplealex}{HTML}{48D1CC}
\definecolor{galoisgeom}{HTML}{ADF802}
\definecolor{otherparabolic}{gray}{0.2}
\definecolor{axesgray}{gray}{0.5}
\definecolor{Lline}{HTML}{00468B}
\definecolor{SUEL}{cmyk}{0.7283, 0.1433, 0.744, 0.0134}

\tikzset{%
  TEL figure/.style={%
      nmdstd,
  line width=1pt,
  every node/.style={color=black}
  },
  coor grid/.style={%
    color=black!20,
    dashed,
    line width=0.7pt
  },
  L line/.style={%
    every path/.style={color=Lline, line width=0.8pt},
    every circle/.style={radius=0.02},
  },
  parabolic/.style={%
    fill=otherparabolic,
    line width=0.5pt,
    every circle/.style={radius=0.03},
  }
}

\renewcommand{\phi}{\varphi}
\providecommand{\SU}{}
\renewcommand{\SU}{{\mathrm{SU}_2}}
\newcommand{\SUoneone}{{\mathrm{SU}_{1,1}}}
\newcommand{\SLR}{{\SL{2}{\R}}}
\newcommand{\SLRpm}{{\mathrm{SL}_2^\pm(\R)}}
\newcommand{\PSLRpm}{{\mathrm{PSL}_2^\pm(\R)}}
\newcommand{\tSLR}{{\widetilde{\SL{2}{\R}}}}
\newcommand{\PSLR}{{\PSL{2}{\R}}}

\newcommand{\SLC}{{\SL{2}{\C}}}
\newcommand{\PSLC}{{\PSL{2}{\C}}}

\newcommand{\PO}{{\mathrm{PO}}}
\renewcommand{\O}{{\mathrm{O}}}

\newcommand{\irr}{{\mathrm{irr}}}
\newcommand{\red}{{\mathrm{red}}}
\newcommand{\cirr}{{c,\,\mathrm{irr}}}
\newcommand{\cred}{{c,\,\mathrm{red}}}

\definecolor{jakecomment}{rgb}{0.067, 0.412, 0.067}

\DeclareMathOperator{\Sim}{Sim}

\DeclareMathOperator{\mIm}{Im}

\DeclareMathOperator{\sign}{sign}
\DeclareMathOperator{\mRe}{Re}
\DeclareMathOperator{\im}{im}
\DeclareMathOperator{\interior}{int}
\DeclareMathOperator{\Diff}{Diff}

\newcommand{\RSU}{R_{\SU}}
\newcommand{\RSLR}{R_{\SLR}}
\newcommand{\XSU}{X_{\SU}}
\newcommand{\XSLR}{X_{\SLR}}

\newcommand{\XC}{X_{\C}}
\newcommand{\RC}{R_{\C}}

\newcommand{\bv}{\mathbf{v}}

\newcommand{\bw}{\mathbf{w}}
\newcommand{\bc}{\mathbf{c}}

\newcommand{\bA}{\boldsymbol{A}}
\newcommand{\bg}{{\boldsymbol{g}}}
\newcommand{\bal}{\mathrm{bal}}

\newcommand{\eh}{\htil}
\newcommand{\elll}{{\mathrm{ell}}}
\newcommand{\hypp}{{\mathrm{hyp}}}

\newcommand{\EK}{{M_K}}
\newcommand{\sslash}{\mathbin{/\mkern-5mu/}}
\newcommand{\XA}{\widetilde{X}}
\newcommand{\chiep}{\chi_{[\epsilon]}}
\newcommand{\Mhe}{\Mhat_{[\epsilon]}}
\DeclareMathOperator{\mtr}{Tr}
\DeclareMathOperator{\codim}{codim}

\renewcommand{\Ltil}{\kerntilde{4}{L}{2}\vphantom{L}}
\newcommand{\Ltilde}{\Ltil}
\newcommand{\Sirr}{{S_\irr}}
\newcommand{\SSL}{{S_\tSLR}}
\newcommand{\SLO}{S_{LO}}

\newcommand{\ELG}{\mathit{EL}_{\Gtil}}
\newcommand{\OTEL}[1]{i^*\big(\Etil(#1)\big)}
\newcommand{\Ka}{\EK(\alpha)}

\newcommand{\IsoE}{{\Isom^+(\E^2)}}
\newcommand{\SimE}{{\Sim^+(\E^2)}}
\newcommand{\SimEall}{{\Sim(\E^2)}}
\newcommand{\IsoH}{{\Isom^+(\H^2)}}
\newcommand{\IsoS}{{\Isom^+(S^2)}}
\newcommand{\IsoHall}{{\Isom(\H^2)}}
\newcommand{\IsoSall}{{\Isom(S^2)}}
\newcommand{\CnE}[1][n]{{\cC_{#1}(\E^2)}}
\newcommand{\CnsignE}[1][n]{{\cC_{#1}^\pm(\E^2)}}
\newcommand{\CnH}[1][n]{{\cC_{#1}(\H^2)}}
\newcommand{\CnsignH}[1][n]{{\cC_{#1}^\pm(\H^2)}}
\newcommand{\CnS}[1][n]{{\cC_{#1}(S^2)}}
\newcommand{\CnP}[1][n]{{\cC_{#1}\left(P^2(\R)\right)}}
\newcommand{\CnSdelta}[1][n]{{\cC_{#1}^\delta(S^2)}}

\newcommand{\Ydelta}{\cY^\delta}
\newcommand{\quo}[2]{\rightquom{#1}{#2}{2pt}{\big}}
\newcommand{\pmone}{{\{\pm1\}}}
\newcommand{\pmonen}{{\{\pm1\}^n}}
\newcommand{\sigmaarrow}{{\xrightarrow{\, \sigma \,}}}

\newcommand{\cRbar}{{\kernoverline{3}{\cR}{4}}\vphantom{\cR}}
\providecommand{\Qbar}{}
\renewcommand{\Qbar}{{\kernoverline{2.5}{Q}{0}}\vphantom{Q}}
\newcommand{\One}{{\mathbb{1}}}
\newcommand{\Utcn}{{(U_t^c)^n}}
\newcommand{\id}{\mathrm{id}}
\newcommand{\rhotilde}{\widetilde{\rho}}

\newcommand{\picred}{\pitil_c^{\,\red}}
\newcommand{\cev}[1]{\reflectbox{\ensuremath{\vec{\reflectbox{\ensuremath{#1}}}}}}
\newcommand{\Rrev}{\cev{\R}}
\newcommand{\Xhatirr}{\Xhat\vphantom{X}^\irr}

\newcommand{\ut}[1][t]{\mathfrak{u}_{#1}}

\newcommand{\ds}{\mathit{ds}}
\newcommand{\df}{\mathit{df}}
\newcommand{\dx}{\mathit{dx}}
\newcommand{\dy}{\mathit{dy}}
\newcommand{\reg}{\mathrm{reg}}

\renewcommand{\MCG}{\mathrm{MCG}}

\newcommand{\jhatstar}{\jhat\vphantom{j}^*}
\newcommand{\lambdabarstar}{\lambdabar\vphantom{\lambda}^*}
\newcommand{\itilstar}{\itil\vphantom{i}^*}
\newcommand{\phihatstar}{\phihat\vphantom{\phi}^*}
\newcommand{\betahatprime}{\betahat\vphantom{\beta}'}

\DeclareMathOperator{\trans}{trans}
\newcommand{\cMtil}{{\kerntilde{5}{\cM}{-1}}}

\DeclareDocumentCommand{\derivat}{ o m m }{%
  \IfNoValueTF{#1}{%
    {\left. {#2} \, \right|_{#3}}%
  }{%
    {{#2} \, #1|_{#3}}%
  }
}
    
\newcommand{\dbydt}{{\frac{d}{\mathit{dt}}}}
\newcommand{\pairLbetaL}{{\langle \, \cL^c, \ \beta^* \cL^c \, \rangle_{\cX^c_\bal(S_{2m})}}}
\newcommand{\pairbetaLL}{{\langle \, \beta^* \cL^c, \ \cL^c \, \rangle_{\cX^c_\bal(S_{2m})}}}

\newcommand{\pairLbetaLatZ}{{\langle \, \cL^c, \ \beta^* \cL^c \,
    \rangle \big\vert_{Z}}}
\DeclareDocumentCommand{\pairLbraidLatZ}{ o m }{%
  \IfNoValueTF{#1}{%
    {{\big\langle \, \cL^c, \ #2 (\cL^c) \, \big\rangle \big\vert_{Z}}}
  }{%
    {{\big\langle \, #1 \cL^c, \ #2 (\cL^c) \, \big\rangle \big\vert_{Z}}}
  }
}
\newcommand{\pairlocal}[4][]{{#1\langle \, #2, \, #3 #1\rangle #1|_{#4}}}

\DeclareDocumentCommand{\pairLbraidLatZp}{ o m }{%
  \IfNoValueTF{#1}{%
    {{\big\langle \, \cL^c, \ #2 (\cL^c) \, \big\rangle \big\vert_{Z'}}}
  }{%
    {{\big\langle \, #1 \cL^c, \ #2 (\cL^c) \, \big\rangle \big\vert_{Z'}}}
  }
}

\newcommand{\pairLbetaLSU}{{\big\langle \, \cL^c_+, \ \beta^* (\cL^c_+) \, \big\rangle_{\cX^c_+(S_{2m})}}}
\newcommand{\pairLbetaLSLR}{{\big\langle \, \cL^c_-, \ \beta^*(\cL^c_-) \, \big\rangle_{\cX^c_-(S_{2m})}}}

\newcommand{\PiPresentation}[1]{\rightquom{F_{2m}}{\pair{#1}}{0.3em}{\big}}

\declaretheorem[nmd, style=plain, sibling=subsection, name={L-space Conjecture}]{LspaceConj}
\declaretheorem[nmd, style=plain, sibling=subsection, name={Enhanced Riley Conjecture}]{EnhancedRiley}

\makeatletter
\renewcommand*\l@subsection{\@dottedtocline{2}{1.5em}{3.2em}}
\renewcommand{\@pnumwidth}{1.7em} 
\renewcommand{\@tocrmarg}{2.7em}
\makeatother

\title{A unified Casson-Lin invariant for the real forms of SL(2)}

\author{Nathan M. Dunfield}
\givenname{Nathan}
\surname{Dunfield}
\address{ Dept.~of Math., MC-382 \\
          University of Illinois \\
          1409 W. Green St. \\
          Urbana, IL 61801 \\ 
          USA
}
\email{nathan@dunfield.info}
\urladdr{https://dunfield.info}

\author{Jacob Rasmussen}
\givenname{Jacob}
\surname{Rasmussen}
\address{ Dept.~of Math., MC-382 \\
          University of Illinois \\
          1409 W. Green St. \\
          Urbana, IL 61801 \\ 
          USA
}
\email{rasmusj@illinois.edu}
\urladdr{https://rasmusj.web.illinois.edu}

\arxivreference{}
\arxivpassword{}

\begin{document}

\begin{abstract}
  We introduce a unified framework for counting representations of
  knot groups into $\SU$ and $\SLR$.  For a knot $K$ in the 3-sphere,
  Lin and others showed that a Casson-style count of $\SU$
  representations with fixed meridional holonomy recovers the
  signature function of $K$. For knots whose complement contains no
  closed essential surface, we show there is an analogous count for
  $\SLR$ representations.  We then prove the $\SLR$ count is
  determined by the $\SU$ count and a single integer $h(K)$, allowing
  us to show the existence of various $\SLR$ representations using only
  elementary topological hypotheses.

  Combined with the translation extension locus of Culler-Dunfield, we
  use this to prove left-orderability of many 3-manifold groups
  obtained by cyclic branched covers and Dehn fillings on broad
  classes of knots.  We give further applications to the existence of
  real parabolic representations, including a generalization of the
  Riley Conjecture (proved by Gordon) to alternating knots.  These
  invariants exhibit some intriguing patterns that deserve
  explanation, and we include many open questions.
  
  The close connection between $\SU$ and $\SLR$ comes from viewing their
  representations as the real points of the appropriate $\SLC$
  character variety.  While such real loci are typically highly
  singular at the reducible characters that are common to both $\SU$
  and $\SLR$, in the relevant situations, we show how to resolve these
  real algebraic sets into smooth manifolds. We construct these
  resolutions using the geometric transition $S^2 \to \E^2 \to \H^2$,
  studied from the perspective of projective geometry, and they allow
  us to pass between Casson-Lin counts of $\SU$ and $\SLR$
  representations unimpeded.  
  \end{abstract}
\maketitle

\newpage

\tableofcontents

\newpage

\section{Introduction}

In analogy with Casson's invariant for homology 3-spheres, Lin
introduced an invariant of knots in $S^3$ in \cite{Lin1992}. In
essence, Lin's invariant is a signed count of (conjugacy classes of)
irreducible representations \(\rho \maps \pi_1(S^3-K) \to \SU\) where
\(\tr \rho(\mu) = 0\) for a meridian \(\mu\) of \(K\). This was
generalized in \cite{Herald1997a, HeusenerKroll1998} to the
\emph{Casson-Lin invariant} \(h_\SU^c(K) \in \Z\), which counts
irreducible representations \(\rho \maps \pi_1(S^3-K) \to \SU\) with
\(\tr \rho(\mu) = c\) for a fixed value of \(c \in [-2,2]\). Here, one
needs to exclude those $c$ corresponding to roots of the Alexander
polynomial \(\Delta_K(t)\) on the unit circle, specifically avoiding
\(D_K := \setdef{a + 1/a}{\mbox{$a \in \C$ and $\Delta_K(a^2) = 0$}}.\)

Here, we introduce a similar invariant \(h^c_\SLR(K)\), which
counts representations \(\rho \maps \pi_1(S^3-K) \to \SLR\) with
\(\tr\rho(\mu) = c\). The fact that \(\SLR\) is noncompact introduces
difficulties that we sidestep, at least for this introduction, by
requiring that $K$ is \emph{small}, that is, \(S^3-K\) contains no
closed essential surface; in particular, this implies $K$ is prime.
Following the standard approach, we show that:

\begin{theorem}
  \label{thm: hSLR intro}
  If \(K\) is a small knot in \(S^3\), there is an integer-valued
  invariant \(h^c_\SLR(K)\) for each \(c \in [-2,2] \setminus
  D_K\). If \(h^c_\SLR(K) \neq 0\), then there is an irreducible
  representation \(\rho\maps \pi_1(S^3-K) \to \SLR\) with
  \(\tr \rho(\mu) = c\).
\end{theorem}
Fixing the knot \(K\), we can view \(h^c_\SLR(K)\) and \(h^c_\SU(K)\)
as functions of \(c\) which are constant on the components of
$[-2, 2] \setminus D_K$. The central goal of this paper is to show
that these two functions are related:
\begin{theorem}
  \label{Thm:MainTheorem}
  If \(K\) is a small knot in \(S^3\), there is an integer \(h(K)\)
  such that \(h(K) = h^c_\SU(K) + h^c_\SLR(K)\) for all
  \(c \in [-2,2] \setminus D_K\).
\end{theorem}
Herald \cite{Herald1997a} and Heusener-Kroll \cite{HeusenerKroll1998}
showed that the $\SU$ Casson-Lin invariant is determined by the
Levine-Tristram signature function $\sigma_K \maps S^1 \to \Z$,
specifically
\begin{equation}
  \label{eq: intro SU is sig}
  h^c_\SU(K) = -\frac{1}{2}\sigma_K(\omega^2)
  \mtext{where \(c = \omega+\omegabar\) is in $[-2, 2] \setminus D_K$.}
\end{equation}
Hence by Theorem~\ref{Thm:MainTheorem}, our new \(\SLR\) Casson-Lin
invariant is determined by the signature function and the single
integer \(h(K)\). As we outline in Section~\ref{sec: intro
  apps}, the connection of $h^c_\SLR(K)$ to the classical signature
function is actually key to its usefulness. Specifically, it will
allow us to prove the existence of $\SLR$ representations from
elementary topological hypotheses.

For many knots, the invariant \(h(K)\) is determined by the
Trotter-Murasugi signature $\sigma(K) := \sigma_K(-1)$; specifically,
for alternating knots and Montesinos knots we show
\(h(K) = -\frac{1}{2}\sigma(K)\) in Corollary~\ref{cor: h for alt} and
Proposition~\ref{prop: Montesinos}. In contrast, for the positive
torus knot $K = T(p, q)$, we have $h(k) = g(K) = \frac{(p-1)(q-1)}{2}$
by Corollary~\ref{Cor:T(p,q)}, which is greater than
\(-\frac{1}{2} \sigma_K\) unless \((p,q) = (2,2n+1), (3,4)\), or
\((3,5)\). Many similar examples can be found by considering knots
which have lens space surgeries.

\subsection{Applications}
\label{sec: intro apps}

Our principal motivation for studying \(h(K)\) and \(h^c_\SLR(K)\) is
to prove the existence of irreducible representations to \(\SLR\) from
\(\pi_1(\EK)\), where $\EK = S^3 \setminus \nu(K)$ is the knot
exterior, as well as the fundamental groups of manifolds constructed
from \(K\) via branched coverings or Dehn surgery.  Similar questions
for representations to \(\SU\) have been extensively studied in the
context of instanton Floer homology, for example in
\cite{AkbulutMcCarthy1990,
  CollinSteer1999,KronheimerMrowka2004,BaldwinSivek2019}.
 
In recent years, interest in representations to \(\SLR\) has been
raised by the L\hyp space conjecture of Boyer, Gordon, and Watson
\cite{BGW2013}, which predicts that a prime \3-manifold \(Y\) is a
Heegaard Floer L-space if and only if \(\pi_1(Y)\) is not
left-orderable. By results in \cite{BRW2005}, the existence of
certain representations \(\pi_1(Y) \to \SLR\) provides one of our most
effective criteria for proving that \(\pi_1(Y)\) is left-orderable,
see the overview in \cite[\S 1.5]{CullerDunfield2018}.
 
A representation \(\rho\maps \pi_1(\EK) \to \SLR\) is called
\emph{elliptic}, \emph{parabolic}, or \emph{hyperbolic} according to
whether \(\rho(\mu)\) is an elliptic, parabolic, or hyperbolic element
of \(\SLR\). Our first application is a criterion for the existence of
elliptic representations analogous to results of Herald and
Heusener-Kroll for \(\SU\):
\begin{corollary}
  \label{cor: ell rep SLR}
  If \(K\) is a small knot with \(\sigma_K(\omega)\)
  nonconstant, then \(\pi_1(M_K)\) admits an irreducible elliptic
  representation to \(\SLR\).
\end{corollary}
The point here is that if \(\sigma_K(\omega)\) is nonconstant, so is
$h^c_\SU(K)$ by (\ref{eq: intro SU is sig}), and then also
$h^c_\SLR(K)$ by Theorem~\ref{Thm:MainTheorem}; hence
$h^c_\SLR(K) \neq 0$ for some open set of $c$, giving the needed
representation by Theorem~\ref{thm: hSLR intro}. The complement of the
figure-eight knot has no such representations, so the condition that
\(\sigma_K(\omega)\) is nonconstant in Corollary~\ref{cor: ell rep
  SLR} is necessary.  For context, recall for nontrivial $K$ there is
always an irreducible $\SU$ representation by the deep results of
\cite{KronheimerMrowka2004}. Reid has asked if the same is true for
$\SLR$, and Corollary~\ref{cor: ell rep SLR} is perhaps the strongest
general result in that direction.

Parabolic representations of \(\pi_1(M_K)\) are of particular
interest.  Since \(h^{\pm2}_\SU(K) = 0\) by Theorem~\ref{thm:
  hSL+hSU}, the invariant \(h(K)\) can be interpreted as a signed count of
parabolic \(\SLR\) representations of \(\pi_1(M_K)\), see
Corollary~\ref{Cor:parabolics}.  Combined with Corollary~\ref{cor: h
  for alt} and Proposition~\ref{prop: Montesinos}, we have:
\begin{theorem}
  \label{thm: alternating riley}
  If \(K\) is a small knot with \(\sigma(K) \neq 0\) that is either
  alternating or Montesinos, then \(\pi_1(M_K)\) admits an irreducible
  parabolic representation to \(\SLR\).  There are at least
  $\abs{\sigma(K)}$ conjugacy classes of such
  representations when counted with Casson-Lin multiplicities.
\end{theorem}
When \(K\) is a 2-bridge knot, the fact that \(K\) should have at
least \(\sigma(K) \) conjugacy classes of parabolic representations,
without any hypothesis on multiplicities, was conjectured fifty years
ago by Riley \cite{Riley1972} and recently proved by Gordon
\cite{Gordon2017}.  Thus Theorem~\ref{thm: alternating riley} can be
naturally viewed as an extension of Riley's conjecture to alternating
and Montesinos knots.  In Section~\ref{sec: almost riley} we give a
new proof of Gordon's result and suggest a refinement in the form of
Conjecture~\ref{conj: enhanced riley}.

Our most satisfactory results for left-orderability apply to cyclic
branched covers. If \(\Sigma_n(K)\) is the \(n\)-fold cyclic branched
cover of \(K\), we prove:
\begin{theorem}
  \label{thm: branched intro}
  If \(K\) is a small knot with \(\sigma_K(\omega)\) nonconstant,
  then \(\pi_1\big(\Sigma_n(K)\big)\) is left-orderable for all
  sufficiently large \(n\).
\end{theorem}
Here, an explicit lower bound on $n$ can be computed easily from the
roots of $\Delta_K(t)$, see Theorem~\ref{thm: branched LO} and also
Remark~\ref{rem: better branched bound}.  For small knots, this
answers a question raised by Boileau, Boyer, and Gordon who proved a
similar theorem to the effect that large cyclic branched covers of
quasipositive knots are not \(L\)-spaces in
\cite{BoileauBoyerGordon2019a}.

For Dehn surgery, our results are more complicated to state, but can
be described succinctly when \(K\) is a 2-bridge knot. If
\(\alpha = \frac{p}{q} \in \Q\), we write \(M_K(\alpha)\) for the
result of Dehn surgery with slope \(p \mu + q\lambda\) on \(K\).  Then
Theorem~\ref{thm: 2-bridge LO} gives:
\begin{theorem}\label{thm: 2-bridge intro}
  If \(K\) is a 2-bridge knot with \(\sigma(K)\neq 0\), the
  \(\pi_1\big(M_K(\alpha)\big)\) is left-orderable either for all
  \(\alpha \in (-\infty,1)\) or for all \(\alpha \in (-1,\infty)\).
\end{theorem}
In the special case of double-twist knots, see the prior work of
\cite{KhoiTeragaitoTran2021, Gao2022}; for a family of 2-bridge knots
not covered by Theorem~\ref{thm: 2-bridge LO}, see \cite{Le2021}.

If \(K\) is a 2-bridge knot other than the unknot or the trefoil,
every nontrivial surgery on \(K\) results in a non-L-space, and thus
should give a manifold whose fundamental group is left-orderable.
Hence Theorem~\ref{thm: 2-bridge intro} has a less satisfactory
conclusion with respect to the L-space conjecture than
Theorem~\ref{thm: branched intro}, which can be used to settle it for
all but finitely many manifolds in its family.  This discrepancy
primarily represents a failure of the technique of using elliptic
\(\SLR\) representations to prove left-orderability, rather than a gap
in our understanding of $\SLR$ representations for these knots. If
\(K\) is a 2-bridge knot with \(\sigma(K)>0\), we expect that the
interval of left-orderable filling slopes can be improved to
\((-\infty,\, \sigma(K) - 1)\) using elliptic \(\SLR\)
representations, but no further.

A similar phenomenon is evident for knots with lens space
surgeries. By \cite{OzsvathSzaboLens}, if \(K\) has a positive lens
space surgery, the set of L-space filling slopes of \(K\) is the
interval \([2g(K)-1, \infty]\).  In Section~\ref{sec: lens space sur},
and specifically Theorem~\ref{thm:lens arcs}, we use elliptic \(\SLR\)
representations to prove that for many such knots, the set of
left-orderable filling slopes contains an interval of the form
\((-\infty, k)\), where \(k\) is a positive integer strictly less than
\(2g(K)-1\). The L-space conjecture predicts that the fillings in the
interval \([k,2g(K)-1)\) should be left-orderable as well, but it
seems unlikely that we will be able to construct left-orderings there
using representations to \(\SLR\).

Our applications to left-orderability are detailed in
Section~\ref{sec: LO apps}.  We make heavy use of the translation
extension locus of \cite{CullerDunfield2018}, which organizes
representations of $\pi_1(S^3 \setminus K)$ to the universal covering
group \(\tSLR\) of $\SLR$. Indeed, that paper motivated much of our
study here.  Our key advance over \cite{CullerDunfield2018} is that
$h^c_\SLR(K)$ allows us to give concrete lower bounds on the sizes of
the intervals of left-orderable branched covers/Dehn fillings, whereas
in \cite{CullerDunfield2018} only the existence of a nonempty open
interval is proved. In complete generality, for Dehn surgery such
intervals can be quite small, see Figure~\ref{fig: K12a380}, meaning
one cannot hope to improve e.g.~Theorem~1.2 of
\cite{CullerDunfield2018} using only $\SLR$ techniques without
additional hypotheses.

Using this perspective, in Section~\ref{sec: extended Lin} we refine
the count $h(K) = h^2_\SLR(K)$ of parabolic $\SLR$ representations
into the \emph{extended Lin invariant} \(\eh(K)\) in $\Z[t^{\pm 1}]$
by taking into account the ``longitudinal heights'' of lifted
representations \(\pi_1(S^3-K) \to \tSLR\). This invariant exhibits
some interesting formal similarities to the Seiberg-Witten invariant
of the knot complement.  We suspect these should be explained by the
work of Haydys \cite{Haydys2020} relating solutions to the 2-spinor
Seiberg-Witten equations to representations into \(\SLR\).  We show in
Theorem~\ref{thm: extended Lin to LO} that if $\eh(K)$ has positive
degree, then one gets left-orderability for large ranges of Dehn
fillings.

\begin{figure}
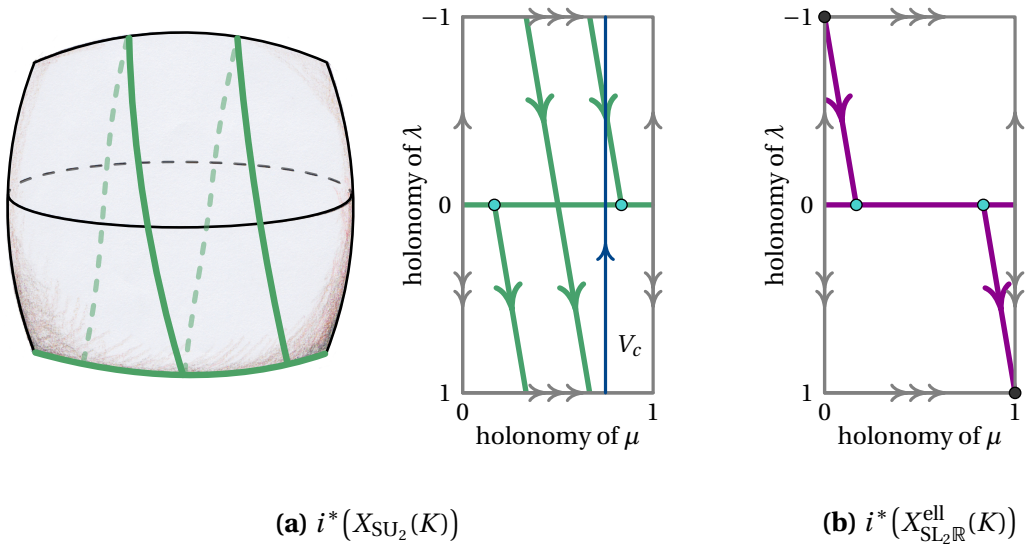

  \begin{center}
    \input figures/trefoil_SU_EL
  \end{center}
  \caption{ At left in (a) is the pillowcase orbifold
    $X_\SU(\partial M_K)$ containing the image of $X_\SU(M_K)$ under
    the restriction map, where $K$ is the positive trefoil knot.
    Details are given in Section~\ref{sec: SLR and T2}, but the
    horizontal and vertical coordinates on the pillowcase are the
    holonomies of $\mu$ and $\lambda$, where $\lambda$ is the Seifert
    longitude. The blue dots correspond to the roots of $\Delta_K(t)$
    on the unit circle.  At right in (b) is the corresponding picture
    for $X_\SLR^\elll(M_K)$.  }
  \label{Fig:CVs for trefoil}
\end{figure}

\subsection{The character variety and the translation extension locus}

We now outline the geometrical interpretations of $h^c_\SU(K)$ and
$h^c_\SLR(K)$, which are best understood in terms of character
varieties of representations to $\SU$ and $\SLR$; throughout, see
Section~\ref{sec: intro char var} for precise definitions and
technical background.

Recall that the $\SU$ character variety of $M_K$ is
\[
  X_\SU(M_K) := \Hom(\pi_1(M_K), \SU)/\sim
\]
where two representations \(\rho\) and \(\rho'\) are equivalent if
they define the same character, i.e.  \(\tr \rho(x) = \tr \rho'(x)\)
for all \(x \in \pi_1(M_K)\); for $\SU$, this is the same as saying
$\rho$ and $\rho'$ are conjugate. Restriction gives a natural map
\[
  i^* \maps X_\SU(M_K) \to X_\SU(\partial M_K) \cong X_\SU(T^2).
\]
The character variety \(X_\SU(\partial M_K)\) is the pillowcase
orbifold shown in Figure~\ref{Fig:CVs for trefoil}(a). The character
variety \( X_\SU(M_K)\) can be divided into  the part $X^\red_\SU(M_K)$
coming from reducible representations, whose image is always the arc
at the bottom of the pillowcase, and the part \(X^\irr_\SU(M_K)\)
coming from irreducible representations, which has expected dimension
\(1\). If \(X^\irr_\SU(M_K)\) is transversally cut out, it carries a
natural orientation, and \(h_\SU^c(K)\) can be interpreted
\cite{HeusenerKroll1998} as the intersection number between
\(i^*\big(X_\SU^\irr(M_K)\big)\) and the vertical curve
\[
  V_c = \setdef{[\rho]}{\mbox{$\rho\maps \pi_1(\partial M_K) \to \SU$ %
      and $\tr \rho(\mu) = c$}},
\]
see Figure~\ref{Fig:CVs for trefoil}(a) where $h^c_\SU(K) = 1$ for the
$V_c$ shown.  This orientation can also be understood in terms of
Reidemeister torsion of cohomology with certain local coefficients
\cite{Dubois2006}.

The situation for \(\SLR\) is simliar in many respects, but there are
a few key differences. The part of \(X_\SLR(T^2)\) where $\pi_1(T^2)$
acts by elliptic elements is again a pillowcase, but now
\(X_\SLR(T^2)\) also contains four noncompact hyperbolic components
which intersect the elliptic component at the orbifold points of the
pillowcase (see Section~\ref{sec: SLR and T2} and especially
Figure~\ref{Fig:CV of T^2}).  In this paper, we focus on the part
$X^\elll_\SLR(M_K)$ of $X_\SLR(M_K)$ coming from representations whose
restriction to $\pi_1(\partial M_K)$ is elliptic.  This gives a very
similar picture to the $\SU$ case, see Figure~\ref{Fig:CVs for
  trefoil}(b).  Provided \(X^\irr_\SLR(M_K)\) is transversally cut out,
we will again be able to view \(h_\SLR^c(K)\) as the intersection
number of \(i^*(X_\SLR^\irr(M_K))\) with $V_c$.

\subsection{Sketch of the definition}

As reinterpreted in \cite{Heusener2003}, Lin's construction begins by
taking a bridge diagram of \(K\) with \(n\) maxima. Splitting along
the middle \(S^2\), we can decompose \(\EK = H_1 \cup_{S_{2n}} H_2\),
where \(S_{2n}\) denotes the sphere with \(2n\) punctures, and \(H_1\)
and \(H_2\) are handlebodies (see Figure~\ref{fig:plat}). If \(\mu_i\)
is a meridional loop about the \(i\)th puncture, we define
\(X^c_\SU(S_{2n})\) to be the variety of characters where all $\mu_i$
have trace $c$.  While \(X^c_\SU(S_{2n})\) is singular at the
reducible characters (a finite set), the irreducible characters
\(X^\cirr_\SU(S_{2n})\) form a smooth open stratum of dimension
\(4n-6\). Each \(L_i = X^{c,\irr}_\SU(H_i)\) is a half-dimensional
subvariety of \(X^{c,\irr}_\SU(S_{2n})\).  For $c \notin D_K$, the
intersection of \(L_1\) and \(L_2\) is compact; the invariant
\(h^c_\SU(K)\) is then defined to be the intersection number of
\(L_1\) and \(L_2\) in \(X^\cirr_\SU(S_{2n})\).  (One must of course
also show that \(h^c_\SU(K)\) does not depend on the particular bridge
diagram.)

Turning to $\SLR$, the difficulty of repeating this construction is
that $\SLR$ is not compact.  However, by requiring that $\EK$ is small
(or more generally \emph{real representation small}, see
Section~\ref{sec: knot reps}), the only possible noncompactness of
$L_1 \cap L_2$ in the smooth stratum $X^\cirr_\SLR(S_{2n})$ is 
the kind already present for $\SU$, namely intersections running out
to some \emph{reducible} character of $\pi_1(\EK)$.  This is prevented
by requiring $c \notin D_K$, just as in the $\SU$ case, allowing us to
define $h^c_\SLR(K)$ for such $c$ and so prove Theorem~\ref{thm: hSLR
  intro}.

\begin{figure}
  \begin{center}
    \begin{tikzpicture}[
  scale=2.6,
  nmdstd,
  line width=1.2pt,
  every node/.style={color=black},
  every circle/.style={radius=0.03}]

  \def\basicxscale{2.5}
  \def\basicyscale{1}
  
  \coordinate (X) at (\basicxscale, 0);
  \coordinate (Y) at (0, \basicyscale);
  \coordinate (mY) at (0, -\basicyscale);

  \coordinate (U) at (${0.1666667}*(X)$);
  \coordinate (V) at (${1-0.1666667}*(X)$);

  \coordinate (A) at (0, 1);
  \coordinate (B) at (U);
  \coordinate (C) at ($2*(U) - (Y)$);
  \coordinate (D) at ($2*(U) + (Y)$);
  \coordinate (DE) at ($3*(U)$);
  \coordinate (E) at ($4*(U) - (Y)$);
  \coordinate (F) at ($4*(U) + (Y)$);
  \coordinate (G) at ($5*(U)$);
  \coordinate (H) at ($(X) - (Y)$);
  \begin{scope}[line width=2pt, color=SUEL, line cap=butt]
    \draw (0, 0) -- (X);
    \draw[mid arrow=0.57] (B) -- (C);
    \draw[mid arrow=0.57] (D) -- (DE);
    \draw[mid arrow=0.57] (DE) -- (E);
    \draw[mid arrow=0.57]  (F) -- (G);
    \draw[mid arrow=0.57, color=locus] (A) -- (B);
    \draw[mid arrow=0.57, color=locus] (G) -- (H);
    \draw[color=locus, dashed] (0, 0) -- (X);
  \end{scope}

  \begin{scope}[color=axesgray]
    \draw (A) rectangle (H);
    \draw (mY) -- (H)
    node[pos=0, below] {$0$}
    node[pos=0.5, below=0.3] {holonomy of $\mu$}
    node[pos=1, below] {$1$};
    
    \draw (mY) -- (Y);
    \node[left] at (Y) {$-1$};
    \node[left] at (0, 0) {$0$};
    \node[rotate=90] at (-0.25, 0) {holonomy of $\lambda$};
    \node[left] at (mY) {$1$};     

    \draw[->] (0, 0) -- ($0.5*(Y)$);
    \draw[->] (X) -- +($0.5*(Y)$);

    \draw[->>] (0, 0) -- ($-0.5*(Y)$);
    \draw[->>] (X) -- +($-0.5*(Y)$);

    \draw[->>>] (Y) -- +($0.55*(X)$);
    \draw[->>>] (mY) -- +($0.55*(X)$);
  \end{scope}

  \node[below left] at (U) {$\frac{1}{6}$};
  \node[above right] at (V) {$\frac{5}{6}$};

  \begin{scope}[dashed, color=simplealex, line width=1]
  \draw ($(U) + 1.2*(Y)$) -- ($(U) - 2.4*(Y)$);
  \draw ($(V) + 1.2*(Y)$) -- ($(V) - 2.4*(Y)$);
  \end{scope}

  \draw[fill=otherparabolic, line width=0.5pt] (A) circle;
  \draw[fill=otherparabolic, line width=0.5pt] (H) circle;
  \draw[fill=simplealex, line width=0.5pt] (B) circle;
  \draw[fill=simplealex, line width=0.5pt] (G) circle;





  \draw[color=Lline, mid arrow=0.7]
  ($0.28*(X) - (Y)$) -- +($2*(Y)$)
  node[right, pos=0.7] {$V_c$};

  \coordinate (L1) at ($0.5*(U)$);
  \coordinate (L2) at ($0.5*(X)$);
  \coordinate (L3) at ($(V)!0.5!(X)$);
  \coordinate (L0) at  ($-1.8*(Y)$);
  
  \node[left] at (L0) {$\sigma_K$};
  \node at ($(L0)+(L1)$) {$0$};
  \node at ($(L0)+(L2)$) {$-2$};
  \node at ($(L0)+(L3)$) {$0$};
    
  \coordinate (L0) at  ($-2.0*(Y)$);
  
  \node[left] at (L0) {$h_\SU^c$};
  \node at ($(L0)+(L1)$) {$0$};
  \node at ($(L0)+(L2)$) {$1$};
  \node at ($(L0)+(L3)$) {$0$};

  \coordinate (L0) at  ($-2.2*(Y)$);

  \node[left] at (L0) {$h_\SLR^c$};
  \node at ($(L0)+(L1)$) {$1$};
  \node at ($(L0)+(L2)$) {$0$};

  \node at ($(L0)+(L3)$) {$1$};

\end{tikzpicture}
  \end{center}
  \caption{For the positive trefoil knot $K$, this figure shows
    $i^*\big(X_\SU(M_K)\big)$ in green and
    $i^*\big(X^\elll_\SLR(M_K)\big)$ in purple on a single pillowcase;
    compare Figure~\ref{Fig:CVs for trefoil}.  Excluding the locus of
    reducible characters, which is common to both, they come together
    only at the points corresponding to roots of $\Delta_K$ on the
    unit circle.  With respect to the orientations indicated, one has
    $h^c_\SU(K) = \pair{ i^*\big(X^\irr_\SU(M_K)\big),\ V_c}$ and
    similarly for $h^c_\SLR(K)$.  Hence $h(K) = 1$ for this knot.}
  \label{fig: trefoil unified} 
\end{figure}
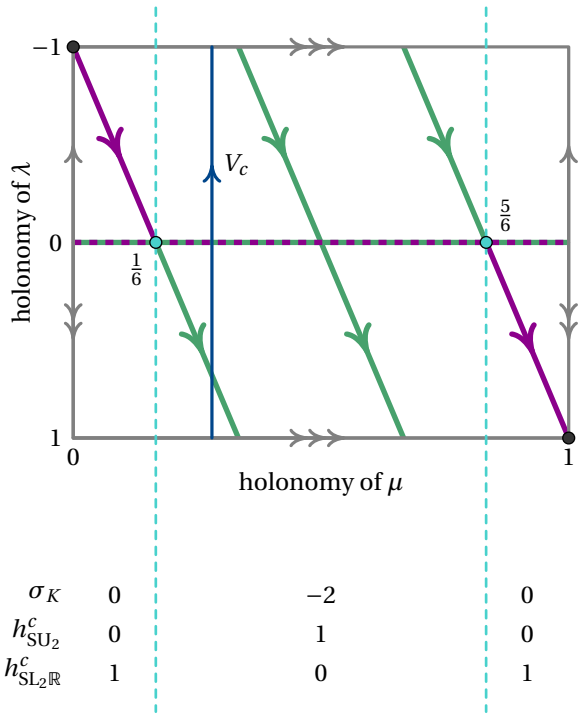

\subsection{Motivation for Theorem~\ref{Thm:MainTheorem}}
\label{sec: motivation}
  
To motivate Theorem~\ref{Thm:MainTheorem}, recall from
\cite{FrohmanKlassen1991} that if $c \in D_K$ corresponds to a simple
root of $\Delta_K(t)$ on the unit circle, then one gets a 1-parameter
family of representations $\rho_t \maps (-1, 1) \to \SLC$ where
$\rho_t$ is an irreducible $\SU$ representation for $t >0$, an
irreducible $\SLR$ representation for $t < 0$, and $\rho_0$ is a
reducible representation with image in $\SU \cap \SLR = S^1$ with
$\tr\big(\rho_0(\mu)\big) = c$.  For simplicity, assume the parameter
$t$ actually corresponds to the trace of $\mu$, say
$\tr\big(\rho_t(\mu)\big) = c + t$.  Then $\rho_t$ will contribute to
$h^{c + t}_\SU(K)$ for $t > 0$ and to $h^{c + t}_\SLR(K)$ for $t < 0$;
part of Theorem~\ref{Thm:MainTheorem} is that these contributions are
the same and hence $h^{c + t}_\SU(K) + h^{c + t}_\SLR(K)$ is unchanged
as $t$ crosses $0$.  Figure~\ref{fig: trefoil unified} visualizes
this for the trefoil by putting the $\SU$ and $\SLR$ pictures from
Figure~\ref{Fig:CV of T^2} together on the same pillowcase.

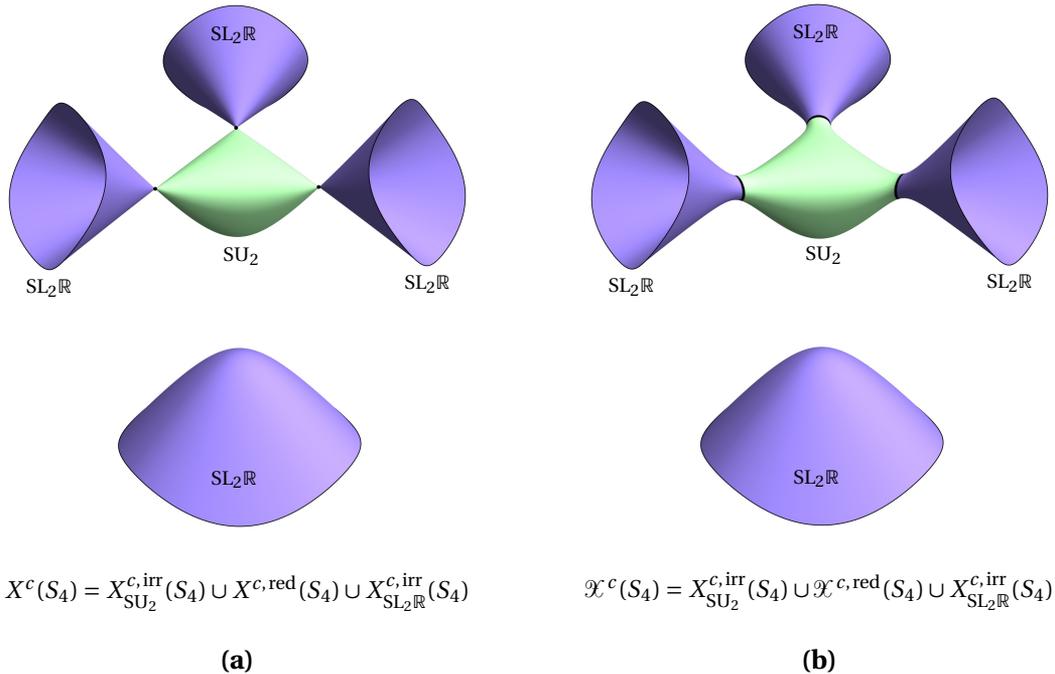
\begin{figure}
  \begin{center}
    \begin{tikzoverlay}[width=6cm]{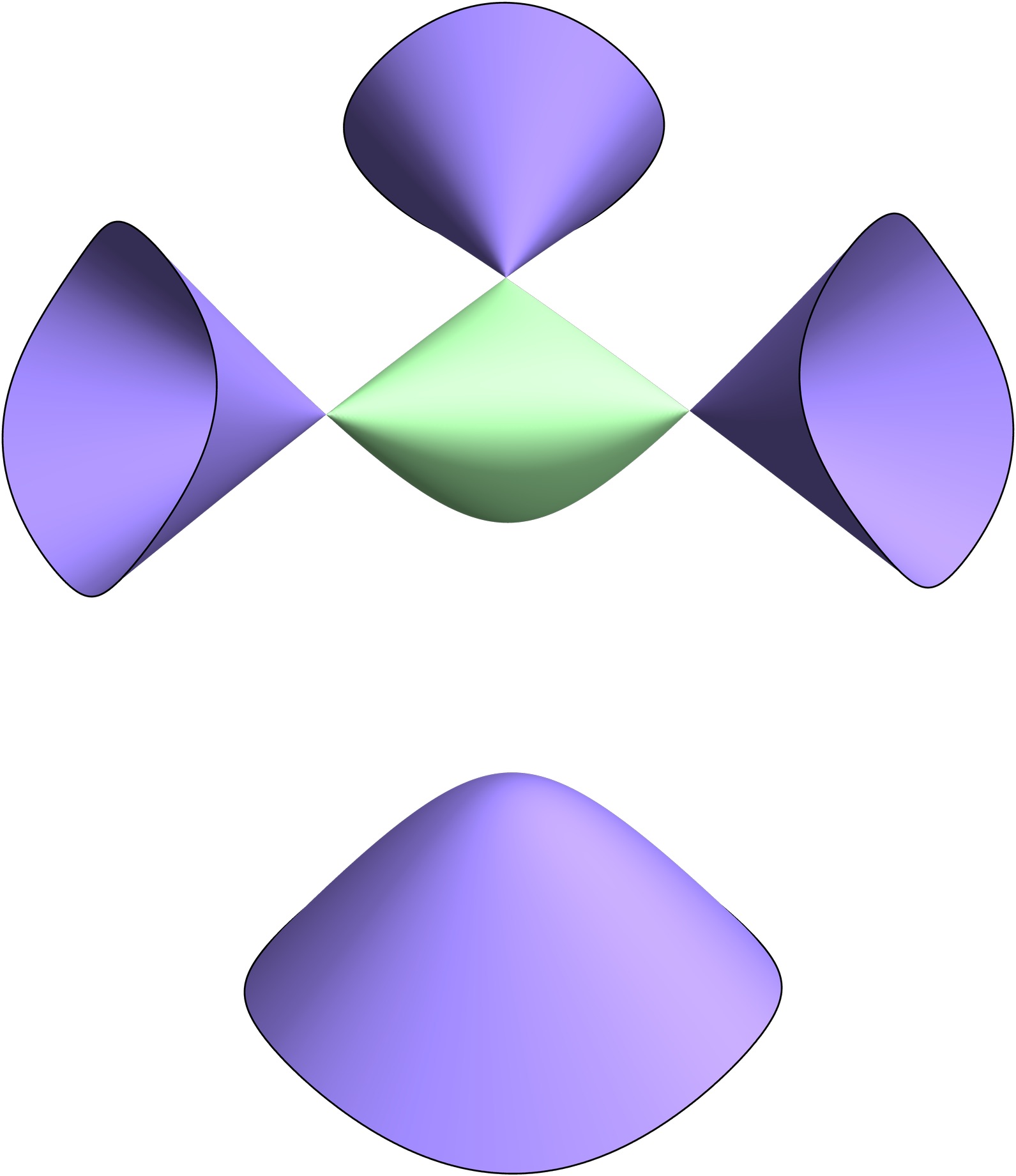}
      \begin{scope}[font=\scriptsize, every circle/.style={radius=0.3}]
        \node[below=8] at (50,0) {\footnotesize $X^c(S_4) = X_\SU^\cirr(S_4) \cup
          X^\cred(S_4) \cup X_\SLR^\cirr(S_4)$};
        \node[below=25] at (50,0) {\small \textbf{(a)}};
        \node[] at (51,60.0) {$\SU$};
        \node[below] at (8.9,57.1) {$\SLR$};
        \node[below] at (91.5,57.9) {$\SLR$};
        \node[] at (49.2,109.0) {$\SLR$};
        \node[] at (49.3,10.6) {$\SLR$};
        \filldraw (32.14,74.9) circle;
        \filldraw (49.8,88.4) circle;
        \filldraw (67.83,75.3) circle;       
      \end{scope}
    \end{tikzoverlay}
    \hspace{1cm}
    \begin{tikzoverlay}[width=6cm]{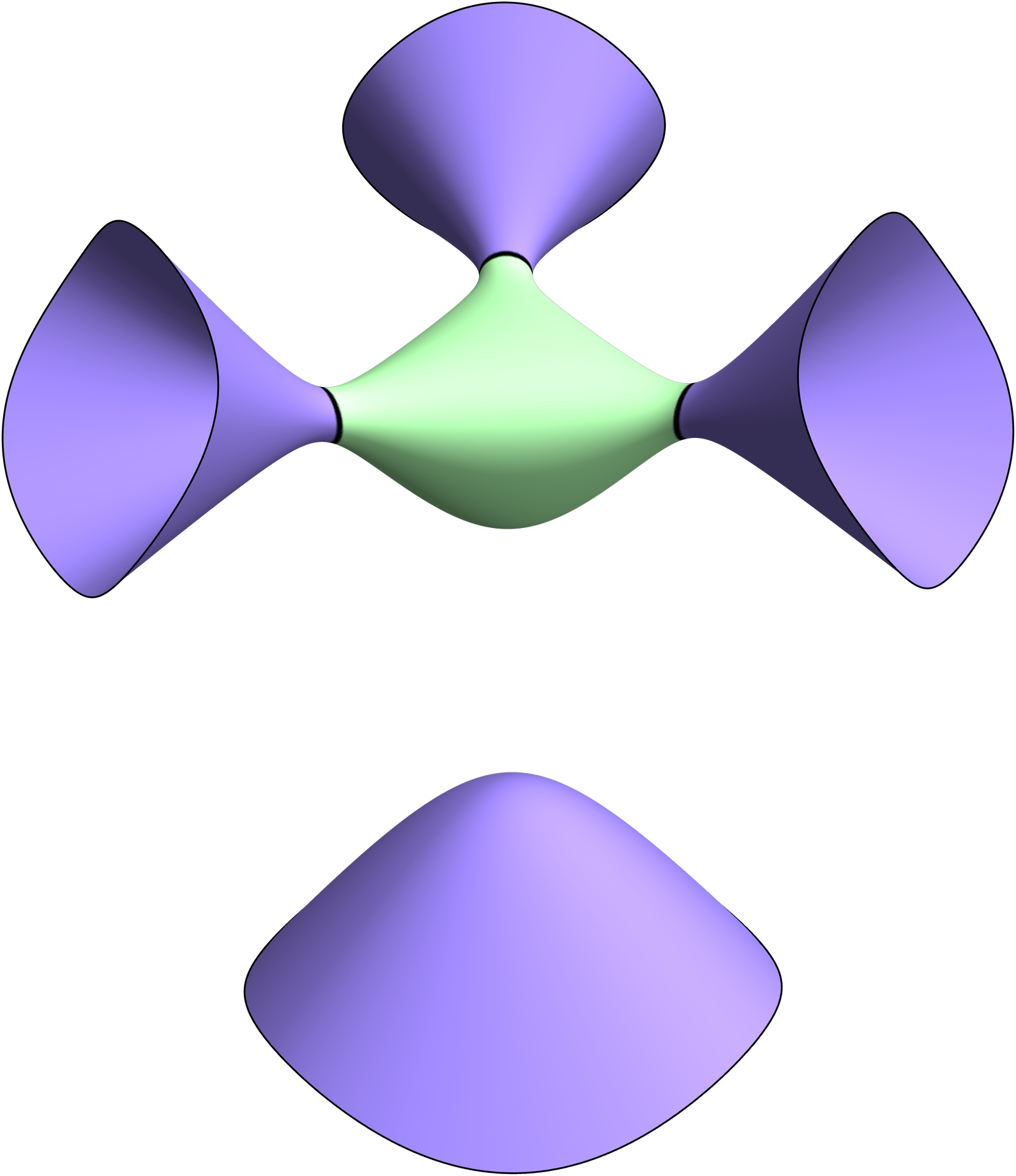}
      \begin{scope}[font=\scriptsize]
        \node[below=8] at (50,0) {\footnotesize $\cX^c(S_4) = X_\SU^\cirr(S_4) \cup
          \cX^\cred(S_4) \cup X_\SLR^\cirr(S_4)$};
        \node[below=25] at (50,0) {\small \textbf{(b)}};
        \node[] at (51,60.0) {$\SU$};
        \node[below] at (8.9,57.1) {$\SLR$};
        \node[below] at (91.5,57.9) {$\SLR$};
        \node[] at (49.2,109.0) {$\SLR$};
        \node[] at (49.3,10.6) {$\SLR$};
      \end{scope}
    \end{tikzoverlay}
  \end{center}
  \caption{At left in (a) is $X^c(S_4)$ for
    $c = 2 \cos \frac{2 \pi}{5}$, drawn using the equations from
    \cite{BenedettoGoldman1999}.  Topologically, $X_\SU(S_4)$ is a
    2-sphere, whereas $X_\SLR(S_4)$ has four distinct components, each
    of which is a plane.  The intersection
    $X^c_\SU(S_4) \cap X^c_\SLR(S_4) = X^\red(S_4)$ is just three
    points, indicated by small dots. The bottom component of
    $X_\SLR(S_4)$ corresponds to the Teichm\"uller space of hyperbolic
    structures on the orbifold with underlying space $S^2$ and four
    points labeled $\Z/5$.  At right in (b) is the resolution
    $\cX^c(S_n)$ where each of the three points in $X^\red(S_4)$ has
    been replaced by a circle.  }
  \label{fig: X(S_4)}
\end{figure}

More generally, Theorem~\ref{Thm:MainTheorem} can be viewed as a
statement about representations of $\pi_1(\EK)$ transitioning from
$\SU$ to $\SLR$ and vice versa as the trace of $\mu$ varies; any such
transition must happen at a reducible representation corresponding to
$c \in D_K$ (see e.g.~Section~\ref{sec: alex and the reds}).  In
particular, the starting point for Theorem~\ref{Thm:MainTheorem} is
the observation that
\[
  X^c(S_{2n}) := X^c_\SU(S_{2n}) \cup X^c_\SLR(S_{2n})
\]
is precisely the real points of the $\SLC$ character variety
\(X^c_\SLC(S_{2n})\).  Moreover, the intersection of
\(X^c_\SU(S_{2n})\) and \(X^c_\SLR(S_{2n})\) is the set of reducible
characters, which is finite. The key difficulty is that the reducible
characters are highly singular points of $X^c(S_{2n})$, as shown in
Figure~\ref{fig: X(S_4)}(a) for the case $n=2$.  Thus it is very
unclear how to ``track'' the relevant intersection numbers used to
define $h^c_\SU(K)$ and $h^c_\SLR(K)$ as one moves through the
reducible locus.

To prove Theorem~\ref{Thm:MainTheorem}, we will resolve
\(X^c(S_{2n})\) to produce a {smooth manifold} \(\cX^c(S_{2n})\)
containing half-dimensional smooth submanifolds \(\cL_1, \cL_2\) which
are resolutions of the $X^c_\SU(H_i) \cup X^c_\SLR(H_i)$.  (See
Figure~\ref{fig: X(S_4)}(b) for a picture of $\cX^c(S_{2n})$ when
$n=2$.)  For any $c \in [-2, 2]$, we then define an invariant
\(h^c(K)\) as the intersection number between \(\cL_1\) and \(\cL_2\).
When $c \notin D_K$, the intersection \(\cL_1 \cap \cL_2\) will
naturally subdivide into $\SU$ and $\SLR$ parts, showing
\(h^c(K) = h^c_\SU(K) + h^c_\SLR(K)\).  We then show there is a
manifold $\cX(S_{2n})$ with a submersion
$\tr \maps \cX(S_{2n}) \to [-2, 2]$ so that
$\tr^{-1}(c) = \cX^c(S_{2n})$.  Moreover, $\cX(S_{2n})$ will have
submanifolds $\cL_i$ so that $\cL_i \cap \cX^c(S_{2n}) = \cL_i^c$.
Then continuity of intersection numbers, in the form of
Theorem~\ref{Thm:InvarianceOfIntersection}, ensures that $h^c(K)$ is
independent of $c$.  This will prove the existence of $h(K)$, at least
for a fixed bridge diagram.

\subsection{Resolutions of real points of character varieties}
\label{sec: resolve intro}

Next, we discuss the general setting for building \(\cX^c(S_{2n})\),
which we call a \emph{resolution} of $X^c(S_{2n})$.  For ease of
exposition, we first discuss the corresponding resolution in the case
of the free group $F_n = \pair{s_1, \ldots, s_n}$ and restrict to
$c \in (-2, 2)$.  Let $X^c_\SLC(F_n)$ be the character variety of
$\SLC$ representations $\rho$ of $F_n$ where all $\rho(s_i)$ have
trace $c$.  Its real points $X^c(F_n)$ are again the union of
$X^c_\SU(F_n)$ and $X^c_\SLR(F_n)$, which meet in exactly
$X^\cred(F_n)$.  (Here, $X^c(F_n)$ is a real algebraic set with
$X^c_\SU(F_n)$ and $X^c_\SLR(F_n)$ real semialgebraic subsets.)  It
turns out both $X^\cirr_\SU(F_n)$ and $X^\cirr_\SLR(F_n)$ are smooth
manifolds of dimension $2n - 3$, whereas $X^\cred(F_n)$ is just
$2^{n-1}$ points. Our resolution of $X^c(F_n)$ from
Theorem~\ref{thm:def of pi} has the following properties, many of
which can be visualized by comparing parts (a) and (b) of
Figure~\ref{fig: X(S_4)}:

\begin{enumerate}
\item The resolution is a smooth manifold $\cX^c(F_n)$ equipped with a
  smooth surjection $\pi \maps \cX^c(F_n) \to X(F_n)$.  We take
  $\cX^\cred(F_n)$ and $\cX^\cirr(F_n)$ to be the preimages under
  $\pi$ of $X^\cred(F_n)$ and $X^\cirr(F_n)$ respectively.

\item The subset $\cX^\cred(F_n)$ is a smooth submanifold of $\cX^c(F_n)$ of
  (real) codimension~1 and $\cX^\cirr(F_n)$ is an open submanifold.

\item The map $\pi$ restricts to a diffeomorphism between the open
  sets $\cX^{\cirr}(F_n)$ and $X^{\cirr}(F_n)$.  Each connected
  component of $\cX^\cirr(F_n)$ maps diffeomorphically under $\pi$ to a
  connected component of either $X^{\cirr}_\SU(F_n)$ or
  $X^{\cirr}_\SLR(F_n)$.

\item \label{item: U0} The subset $\cX^\cred(F_n)$ can be naturally
  identified with the ``character variety'' of representations to
  \(
  U_0 = \setdef{\mysmallmatrix{a}{b}{0}{\abar}}{%
    \mbox{$a, b \in \C$ with $\abs{a} = 1$}} \)
  that are not conjugate to a representation with diagonal image.
  
\item \label{item: funct} For an automorphism $\sigma$ of $F_n$ coming
  from the braid group, consider the induced automorphism $\sigma^*$
  of $X^c(F_n)$.  There is a unique diffeomorphism $\sigmatil^*$ of
  $\cX^c(F_n)$ which is compatible with $\sigma^*$ in that sense that
  $\sigma^* \circ \pi = \pi \circ \sigmatil^*$.

\end{enumerate}

In (\ref{item: U0}), we are replacing each character
$\chi_0 \in X^\red(F_n)$ with something made out of reducible $\SLC$
representations with that character. However, we should point out that
$U_0$ is not conjugate into either $\SU$ or $\SLR$.  Indeed, any
representation $\rho$ from $F_n$ into $\SU$ or $\SLR$ with character
$\chi_0$ is will in fact be diagonalizable in $\SLC$.  Still, the
appearance of $U_0$ here is very natural from another vantage point:
back in the setting of a knot $K$, de Rham showed the roots of
$\Delta_K$ on the unit circle characterize the $c$ where there are
representations $\rho \maps \pi_1(M_K) \to U_0$ with
$\tr \rho(\mu) = c$ that are not diagonalizable.

\subsection{Geometric transition via projective geometry}
The work of Frohman-Klassen \cite{FrohmanKlassen1991} mentioned in
Section~\ref{sec: motivation} starts from the fact that, after modding
out by $\pm I$, the groups $\SU$, $U_0$, and $\SLR$ are the
orientation preserving isometry groups of the round sphere $S^2$, the
Euclidean plane $\E^2$, and the hyperbolic plane $\H^2$.  Their
construction of the arc of characters, half in $X_\SU(K)$ and half in
$X_\SLR(K)$, uses a 1-parameter family of metrics on the plane that go
from positive curvature to zero curvature to negative curvature.  When
the root of $\Delta_K$ is simple, they show how a nondiagonalizable
representation $\rho \maps \pi_1(M_K) \to \Isom^+(\E^2)$ can be
deformed into isometries of these nearby metrics.

To build the resolution $\cX^c(F_n)$, we follow the lead of
\cite{FrohmanKlassen1991} and consider the geometric transition
$S^2 \to \E^2 \to \H^2$, which we study via projective geometry using
the perspective of \cite[\S 2.1]{CooperDancigerWienhard2014}.  Because
we are restricting to representations of $F_n$ where the generators go
to elliptic elements of a fixed trace, a relevant representation into
$\SU$, $U_0$, or $\SLR$ is largely determined by the fixed points of
the generators in their action on $S^2$, $\E^2$, or $\H^2$.  This
leads us to study configuration spaces of $n$ points in each of these
three geometries in Section~\ref{sec:config}, where the points are
taken modulo isometry (for $S^2$ and $\H^2$) or similarity (for
$\E^2$).  Here, we require that not all $n$ points are the same, and
denote the resulting configuration spaces by $\CnS$, $\CnE$, and
$\CnH$ which have dimensions $2n - 3$, $2n - 4$, and $2n - 3$
respectively.

We now sketch a natural way of combining these configuration spaces
together into a single smooth manifold.  To begin, we consider $\CnH$,
which has two ends: one where all the points coalesce and the other where
their diameter goes to infinity.  We focus on the former end and leave
the other alone. Given an element in $\CnH$, we can rescale the metric
on $\H^2$ so that the diameter of the set of points is exactly 1 in
the new metric. The closer the original points are in $\H^2$, the
flatter the rescaled metric is. When the original points are very
close, the new space is nearly isometric to $\E^2$ on the scale of the
points in the final configuration.  This makes it plausible that we
should compactify this end of $\CnH$ by adding a copy of $\CnE$ at
infinity to produce a manifold with boundary.  Looking now at the end
of $\CnS$ where all the points come together, the exact same story
applies to suggest that we should also compactify this end by adding a
copy of $\CnE$. We could then glue our two compactifications together
to get a nice manifold structure on $\CnH \cup \CnE \cup \CnS$.

To resolve $\cX^c(F_n)$, the correct object to look at involves
``signed points'' to account for the rotation directions of the
generators at their fixed points. We assemble the relevant
configuration spaces into a smooth manifold $\cY$ in
Theorem~\ref{thm:consY}, and then use it to build $\cX^c(F_n)$ in
Theorem~\ref{thm:def of pi}.  The generalization to $X(S_n)$ and to
allowing $c$ to vary are too involved for this introduction but are
tackled in Sections~\ref{sec:punc sphere} and~\ref{sec:varies}.

\subsection{Open problems}

Our work here raises many questions, both general and specific.
General problems include:

\begin{enumerate}

\item Can the definition of $h^c_\SLR(K)$ be extended to all knots,
  not just those which are real representation small?  Naively, could
  one simply ignore any non-compact components of the intersection?
 Alternatively, could the sheaf-theoretic
 perspective of Abouzaid and Manolescu \cite{AbouzaidManolescu2020}
  in the case of $\SLC$ be adapted to our real-algebraic setting?
  Given its connection to the signature, one expects that
  $h(K_1 \# K_2) = h(K_1) + h(K_2)$.  However, the connected sum of
  nontrivial knots is typically not real representation small, meaning
  that we cannot make sense of $h(K_1 \# K_2)$.  More generally, it
  would be nice to be able to define and compute \(h\) for satellite
  knots.
  
\item Can we count hyperbolic representations with fixed trace and so
  extend the definition of $h^c_\SLR(K)$ to all
  \(c \in \R \setminus D_K\)?  Assuming this can be done, can
  $h^c_\SLR(K)$ jump as we move through $c \in D_K$ for $\abs{c} > 2$?
  If it does not, we would get interesting consequences about ideal
  points of $X_\SLR(K)$ in terms of $h(K)$.  A weak version of such an
  argument is provided by Corollary~\ref{cor: really ideal}, which
  counts the number of such points modulo 2.

\item Can our theory be extended to knots in other closed
  \3-manifolds?  One approach would be to use bridge diagrams in
  Heegaard splittings and \((g,n)\) knots.  Here, the bridge
  presentations we use are $(0, n)$ knots and the doubly-pointed
  Heegaard diagrams from \cite{OzsvathSzabo2004} are the $(g, 1)$
  case.  The first step would be to understand how to resolve the real
  points of the relevant character variety of a surface of genus $g$
  with $2n$ punctures. Our approach here uses that representations of
  $\pi_1(S_{2n})$ have image generated by elliptic elements, which is
  no longer the case in higher genus. This more general perspective
  might also help us understand the behavior of \(h\) for Berge knots and
  other \((1,1)\) knots. 
  
\item Is it possible to give a more intrinsic definition of
  \(h^c_\SLR \) by working in the space of all \(\SLR\) connections,
  as in Herald's treatment \cite{Herald1997a} of the \(\SU\)
  Casson-Lin invariant? If this could be done, it might provide a way
  to extend $h$ to knots in general 3-manifolds. What would a proof
  of Theorem~\ref{Thm:MainTheorem} look like in this context?
  
\item The translation extension locus of \cite{CullerDunfield2018}
  shares some interesting similarities with the moduli space of
  solutions to the Seiberg-Witten equations on the knot complement
  equipped with a cylindrical end. In this setting, the extended Lin
  invariant \(\eh(K)\) is analogous to the Seiberg-Witten invariant.
  In \cite{Haydys2020}, Haydys relates solutions to the 2-spinor
  Seiberg-Witten equations on a manifold \(Y\) to the space of
  \(\SLR\) connections on \(Y\). Can these similarities be explained
  by this? Conjectures~\ref{conj:L-space torsion} and \ref{conj:
    enhanced riley}, which concern \(\eh(K)\) for Berge knots and
  2-bridge knots, were originally formulated with this analogy in
  mind.
  

\item The proof of Theorem~\ref{thm: 2-bridge intro} relies on the
  fact that if \(K\) is 2-bridge, then every parabolic
  \(\rho\maps \pi_1(\EK) \to \SLR\) lifts to
  \(\rhotil: \pi_1(\EK) \to \tSLR\) where \( \rho(\lambda)\)
  has nonzero translation number. Are there other interesting classes
  of knots for which this statement holds? This would enable one to
  prove results similar to Theorem~\ref{thm: 2-bridge LO} for such
  knots. It is conceivable that alternating or Montesinos knots have
  this property, see Remark~\ref{rem: extend LO}.
\end{enumerate}

Among the specific problems, we point out
Conjectures~\ref{conj:L-space torsion} and~\ref{conj: enhanced riley}
which describe the expected structure of the enhanced Lin invariant
for Berge knots and 2-bridge knots.  Other specific problems include
those of Section~\ref{sec: odd punct}, Remark~\ref{rem: long at
  roots}, Question~\ref{qu: real ideal}, Remark~\ref{rem: extend LO},
Remark~\ref{rem: lens space conj}, and
Question~\ref{question:interval}.

\subsection{Plan of the paper}

The first half of the paper is devoted to the definition of \(h\) and
the proof of Theorem~\ref{Thm:MainTheorem}. Sections~\ref{sec: smooth
  background}--\ref{sec: tSLR basics} contain background material.
Section~\ref{sec: smooth background} gives our conventions for smooth
manifolds and orientations, and collects some basic results about
intersection numbers that are needed for the main theorems.
Section~\ref{sec: intro char var} gives background on representations
into $\SLC$ and its subgroups $\SU$ and $\SLR$, their associated
character varieties, and various standard or easy lemmas.
Section~\ref{sec: tSLR basics} discusses representations into $\tSLR$
and shows there is a ``character variety'' in this context
(Theorem~\ref{thm: tSLRchar}).

Sections~\ref{sec:config}--\ref{sec:join configs} discuss
configurations of points in $S^2$, $\E^2$, and $\H^2$, and how they
fit together to form a space $\cY$ in Theorem~\ref{thm:consY}.
Section \ref{sec:map to X} then uses the space $\cY$ to build the
resolution of $X^c(F_n, S)$ for $c \in (-2, 2)$ in
Theorem~\ref{thm:def of pi}.  Section \ref{sec:punc sphere} builds the
corresponding resolution for the punctured sphere $S_{2n}$ in
Theorem~\ref{thm:cXSn}.  Section~\ref{sec:varies} contains
Theorems~\ref{thm:sX(F_n)} and~\ref{thm:sX(S_n)} which build the
``parameterized'' resolutions where $c$ is allowed to vary.
Section~\ref{Sec:Orientations} studies orientations on what has
been constructed so far.  Section~\ref{Sec:Definition of h} shows the
functoriality of $\cX^c$ with respect to automorphisms coming from the
braid group, and uses this to define $h^c_\SU$, $h^c_\SLR$,
and $h^c$ for the plat closure of a braid.  For a fixed braid,
Theorem~\ref{thm: h indep of c} shows that $h^c$ is independent of $c$
and Theorem~\ref{thm: hSL+hSU} proves that $h = h^c_\SU + h^c_\SLR$.
Section~\ref{sec: invariance} shows that $h^c_\SU$, $h^c_\SLR$, and
$h^c$ do not depend on choice of plat diagram, completing the proofs of
Theorems~\ref{thm: hSLR intro} and \ref{Thm:MainTheorem}.

The second half of the paper is focused on the applications outlined
in Section~\ref{sec: intro apps}.  Section~\ref{sec: ext Lin} refines
what $h$ says about parabolic representations and begins to study how
it relates to $\tSLR$ representations and the translation extension
locus of \cite{CullerDunfield2018}; this leads to
Definition~\ref{Def:extended invariant} for the extended Lin invariant
$\eh(K)$.  Section~\ref{sec: properties} establishes further
properties of $h(K)$ and $\eh(K)$, including their behavior under
mirroring.  Therein, Section~\ref{sec: rel with signature} pins down
the connection between our version of $h^c_\SU(K)$ and that in
\cite{Heusener2003}, formally proving equation (\ref{eq: intro SU is
  sig}) as Corollary~\ref{Cor: hSU=sig}.  Section~\ref{sec: LO apps}
computes $h(K)$ for alternating knots and gives the applications to
left-orderability.  Finally, Section~\ref{sec:computations} computes
\(h\) for Montesinos knots and torus knots, and gives some conjectures
about the extended Lin invariant for Berge knots and 2-bridge knots.

\subsection{Acknowledgements}

Dunfield was partially supported by the US National Science Foundation
under grants DMS-1510204, DMS-1811156, and DMS-2303572, as well as by a
fellowship from the Simons Foundation (673067, Dunfield).  Rasmussen
was partially supported by EPSRC grant EP/M000648/1 and by the US
National Science Foundation under grant DMS 2405302.  We thank Chris
Herald for helpful correspondence related to Lemma~\ref{lem: herald
  power}, and Steve Boyer, Marc Culler, Cameron Gordon, Sean Lawton,
and Liam Watson for helpful conversations. We also thank the referee for their
careful reading of this paper and helpful comments.  This research was
conducted in part at Cambridge University, the Newton Institute, the
Universit\'e du Qu\'ebec \`a Montr\'eal, and the Institute for
Advanced Study.

\section{Background on orientations and smooth topology}
\label{sec: smooth background}

\subsection{Manifolds}

Our terminology for smooth manifolds is as follows. A \emph{smooth
  embedding} $S \to M$ of smooth manifolds $S$ and $M$ is a smooth
immersion that is a homeomorphism onto its image.  The image of a
smooth embedding is called a \emph{smooth submanifold}. When $S$ has
boundary, we insist $\partial S$ is contained in $\partial M$ with $S$
transverse to $\partial M$; we refer to the situation when
$\partial S$ is contained in the interior of $M$ as an \emph{embedded
  submanifold with boundary}. A smooth submanifold $S$ is
\emph{closed} when it is a closed subset of $M$; this is equivalent to
the inclusion $S \to M$ being a proper map.

We will need the following basic fact starting in Section~\ref{sec:punc
  sphere}:

\begin{lemma}\label{lem:quotient is manifold}
  Suppose $\pi \maps M \to N$ is a surjective submersion of smooth
  manifolds without boundary.  If $S \subset M$ is a closed
  submanifold that is a union of fibers of $\pi$, then $\pi(S)$ is a
  closed submanifold of $N$ of the same codimension as $S$.
\end{lemma}

\begin{proof}
By Proposition 4.28 of \cite{Lee2013}, the map $\pi$ is a topological
quotient map, and hence $\pi(S)$ is closed, and so it will be a closed
submanifold if it is a submanifold at all.  Set $m = \dim M$ and
$n = \dim N$, and let $k$ be the codimension of $S$. To see that
$\pi(S)$ is an embedded submanifold of the claimed codimension, we
apply the local slice criterion of Theorem~5.8 of \cite{Lee2013}: it
is enough to show that given $s_0 \in S$, we can find a chart $U$ on
$N$ containing $\pi(s_0)$ so that $\pi(S) \cap U$ is a local slice of
the form $\R^{n - k} \times (0, \cdots, 0)$ in $U \cong \R^n$.

Choose local charts \(V \cong \R^m\) on $M$ and \(U \cong \R^n\) on
$N$ so that $s_0 = 0$ in $V$, the image $\pi(V)$ is $U$, and the map
$\pi:V \to U$ is projection onto the first \(n\) coordinates. Let
\(S' = S \cap V\). Since \(\pi(V) = U\) and \(S\) is a union of fibers
of \(\pi\), we have \(\pi(S) \cap U = \pi(S')\).  As $S$ has
codimension $k$, we can shrink our charts so that additionally there
is a submersion $f \maps \R^m \to \R^k$ with $S' = f^{-1}(0)$.  Since
$S'$ contains $0 \times \R^{m - n}$, it follows that the last $m - n$
columns of $D_0f$ are $0$, and hence the initial $n$ columns of $D_0f$
have rank $k$.  Thus, after zooming in if necessary, we have that the
restriction of $f$ to $\R^n \times (0, \ldots, 0)$ is also a
submersion.  In particular,
$H = S' \cap \big(\R^n \times (0, \ldots, 0)\big)$ is a smooth
embedded submanifold of $\R^n \times (0, \ldots, 0)$, and, as $\pi$ is
coordinate projection, we have $S' = H \times \R^{m - n}$.  As $\pi$
is a local diffeomorphism onto $U$ when restricted to
$\R^n \times (0, \ldots, 0)$, we have that $\pi(S')$ is contained in a
local slice near $\pi(s_0)$ as needed.
\end{proof}

\subsection{Orientations}
\label{sec:orientshortex}

A short exact sequence of oriented vector spaces
\[0 \to A \xrightarrow{i} B \xrightarrow{\pi} C \to 0\] is
\emph{compatible} when given a positively oriented basis
\(\langle a_1, \ldots a_n \rangle \) of \(A\) and a positively
oriented basis \(\langle \pi(b_1), \ldots, \pi(b_m) \rangle \) of
\(C\), then the basis
\(\langle b_1, \ldots, b_m, i(a_1), \ldots, i(a_n)\rangle\) for \(B\)
is also positively oriented.  For such a sequence, orientations on any
two of \(A,B,\) and \(C\) determine a unique orientation on the third
that makes the sequence compatible. With this convention, the standard
orientation on the direct sum \(A\oplus B\) is the one compatible with
the sequence \(0 \to B \to A \oplus B \to A \to 0 \).

When $E \xrightarrow{\pi} B$ is a smooth submersion of oriented
manifolds, the fibers $\pi^{-1}(b)$ can be \emph{compatibly
  oriented} by the requirement that for all $e \in E$ the sequence
\[
  0 \to T_e
  F_e \xrightarrow{d i} T_e E \xrightarrow{d \pi} T_{b} B \to 0
\]
is compatibly oriented, where $b = \pi(e)$ and $F_e = \pi^{-1}(b)$.
In the special case of a smooth fiber bundle of oriented manifolds
$F \to E \xrightarrow{\pi} B$ with structure group $\Diff^+ \! (F)$,
orientations on any two of $F, E$, and $B$ determine a unique
orientation on the third that makes the orientations compatible with
respect to $\pi$. Finally, if \(M\subset M'\) is an embedding of
smooth oriented manifolds, the normal bundle \(\nu_{M'/M}\) is
oriented by the requirement that the short exact sequence
\[
  0 \to TM \to TM' \to \nu_{M'/M} \to 0
\]
is compatibly oriented. 

We orient the boundary of a manifold via the ``outwards normal first''
convention of Chapter 15 of \cite{Lee2013}.  Thus if
$M = \setdef{x \in \R^n}{x_1 \leq 0}$, and $\pair{e_1, e_2, \ldots, e_n}$
gives the preferred orientation of $M$ then an oriented basis for
$\partial M$ is $\pair{e_2,\ldots,e_n}$.

\begin{remark}\label{rem: orient vs Heusener}
  Our orientation conventions for both short exact sequences and
  $\partial M$ differ from \cite{Heusener2003}: we use
  $(\mbox{base}) \oplus (\mbox{fiber})$ and ``outwards normal first''
  and whereas he uses $(\mbox{fiber}) \oplus (\mbox{base})$ and
  ``inwards normal last''.
\end{remark}

The rest of this section is not used until Section~\ref{Sec:Definition
  of h} or later, so you will want to skip ahead to Section~\ref{sec:
  intro char var} at first reading. If \(V\) is an oriented vector
space, we write \(-V\) for the vector space with the opposite
orientation.

\begin{lemma}
  \label{lem: 3 x 3 with orient}
  Suppose the following diagram of oriented finite-dimensional vector
  spaces commutes:
  \[
    \begin{tikzcd}
      V_{11} \arrow{r} \arrow{d} & V_{12} \arrow{r} \arrow{d} & V_{13}\arrow{d} \\
      V_{21} \arrow{d} \arrow{r} & V_{22} \arrow{d} \arrow{r} & V_{23} \arrow{d} \\ 
      V_{31} \arrow{r}           & V_{32} \arrow{r}           & V_{33}
  \end{tikzcd}
  \]
  If all of the rows and the leftmost two of the columns are short
  exact, then the rightmost column is short exact as well. If
  additionally they are compatibly oriented, the rightmost column
  is also compatibly oriented if and only if $(\dim V_{13}) (\dim V_{31})$
  is even.
\end{lemma}

\begin{proof}
The rightmost column is short exact by the $(3 \times 3)$-lemma,
so it remains to puzzle over the orientations.  With our conventions,
the middle term can be written as an oriented direct sum in two ways:
\begin{align*}
  V_{22} &= V_{23} \oplus V_{21} =
            V_{23} \oplus \big( V_{31} \oplus V_{11}\big) \\
         &= V_{32} \oplus V_{12} =
            \big( V_{33} \oplus V_{31}\big) \oplus
           \big( V_{13} \oplus V_{11}\big) \\
         &= (-1)^{(\dim V_{13}) (\dim V_{31})}
            V_{33} \oplus V_{13} \oplus V_{31} \oplus V_{11}.
\end{align*}
Cancelling the $V_{31} \oplus V_{11}$ from the ends of first and last lines
gives
\[
  V_{23} = (-1)^{(\dim V_{13}) (\dim V_{31})} V_{33} \oplus V_{13}
\]
proving the lemma.
\end{proof}

\begin{corollary}
  \label{Cor:Normal bundles}
  Suppose \(M \subset M'\) and \(X \subset X'\) are inclusions of
  smooth oriented manifolds and that \(f:M' \to X'\) is a smooth
  submersion that restricts to a smooth submersion on \(M\). Let
  \(N' = f^{-1}(X)\), and \(N = N' \cap M\). The submersions
  \(f:M'\to X'\) and \(f:M \to X'\) determine orientations on \(N'\)
  and \(N\), respectively, and we have
  \((\nu_{M'/M})|_N \cong (-1)^{e}\nu_{N'/N}\) where
  $e = (\codim_{N'}N)(\codim_{X'}X)$.
\end{corollary}

\begin{proof} 
We have a commuting diagram of vector bundles on $N$:
\[
  \begin{tikzcd}
    TN \arrow{r} \arrow{d} & TN' \arrow{r}\arrow{d} & \nu_{N'/N} \arrow{d}\\
    TM \arrow{d} \arrow{r}& TM' \arrow{d} \arrow{r}& \nu_{M'/M} \arrow{d} \\
    \nu_{X'/X} \arrow{r} & \nu_{X'/X} \arrow{r} & 0
  \end{tikzcd}
\]
to which Lemma~\ref{lem: 3 x 3 with orient} applies, proving our claim.
\end{proof}

A very similar argument can be used to prove:

\begin{corollary} \label{Cor:Normal Quotient} Suppose that
  \(Y\subset Y'\) is an embedding of smooth oriented manifolds.
  Suppose \(G\) is an oriented connected Lie group acting freely and
  properly on \(Y'\) where this action leaves $Y$ invariant. Then
  \(X = Y/G\) and \(X' = Y'/G\) inherit orientations from the
  orientations on \(G\), \(Y\), and \(Y'\). The action of \(G\) on
  \(Y\) extends to an action on \(\nu_{Y'/Y}\), giving
  \(\nu_{X'/X} \cong (\nu_{Y'/Y})/G\) as oriented vector bundles.
\end{corollary}

\subsection{Intersection numbers}

The invariants of Theorems~\ref{thm: hSLR intro}
and~\ref{Thm:MainTheorem} are defined as intersection numbers of
certain pairs of submanifolds.  In this final subsection, we record
the basic facts we will use starting in Section~\ref{Sec:Definition of
  h}; they are well-known, but we could not find a good reference.

\begin{theorem}\label{thm:intersection}
  Let $M$ be a smooth oriented $2n$-manifold without boundary,
  possibly noncompact.  Suppose $A$ and $B$ are oriented closed
  submanifolds of $M$ of dimension $n$.  If $A \cap B$ is compact,
  then there exists a compactly supported ambient isotopy of $A$ to
  $A'$ so that $A'$ and $B$ intersect transversely. The
  algebraic intersection number of $A'$ and $B$ is independent of
  the choice of such isotopy, and will be denoted $\pair{A, B}_M$.
\end{theorem}

\begin{proof}
First, we construct the claimed ambient isotopy.  Let
$f \maps A \to M$ be the inclusion map. As $f$ is proper, we can
choose a compact $n$-dimensional submanifold $V \subset A$ with
boundary containing $f^{-1}(B)$ in its interior. Then $f$ is trivially
transverse to $B$ on the closed set $A \setminus \interior(V)$.  By
Proposition~4.5.7 of \cite{Wall2016}, there is an arbitrarily small
perturbation of $f$ to a smooth $f'$ transverse to $B$ where $f = f'$
outside $\interior(V)$. By Proposition 4.4.4 of \cite{Wall2016}, we
can arrange that $f'|_V$ is also an embedding that is isotopic to
$f|_V$ via an isotopy that is constant on $\partial V$. (Here,
you have to go back to the proof of Corollary 2.2.5 of \cite{Wall2016}
to see that the isotopy is constant where $f$ and $f'$ agree.)
Applying Theorem 2.4.2 of \cite{Wall2016} to the isotopy from $f$ to
$f'$, we get a compactly supported ambient isotopy of $f$ to $f'$ as
required.

Because our isotopies are compactly supported, the proof that
$\pair{A, B}_M$ is well-defined is essentially the same as when $M$ is
compact.  Specifically, consider an isotopy of embeddings
$F \maps A \times I \to M$ where $F_0$ and $F_1$ are transverse to $B$
and that is constant outside a compact set $K \subset A$ containing
$F^{-1}(B)$. We can use Proposition~4.5.7 of \cite{Wall2016} to
perturb $F$ to a map $F'$ that agrees with $F$ on both
$A \times \{0, 1\}$ and $A \setminus K$ and is transverse to $B$.
Then $(F')^{-1}(B)$ is a closed 1-dimensional submanifold of
$A \times I$ contained in the compact set $K \times I$.  Therefore
$(F')^{-1}(B)$ is a finite union of oriented circles and arcs with
endpoints in $A \times \{0, 1\}$, and the usual argument from
e.g.~\cite[Chapter 3]{GuilleminPollack1974} shows that
$\pair{F_0(A), B}_M = \pair{F_1(A), B}_M$ as needed.
\end{proof}

\begin{theorem}
\label{Thm:InvarianceOfIntersection}
  Let $M$ be a smooth oriented $(2n + 1)$-manifold with boundary,
  possibly noncompact, with $\pi \maps M \to I$ a submersion where
  $\partial M = \pi^{-1}(\partial I)$.  For $t \in I$, let $M_t$ be
  the closed submanifold $\pi^{-1}(t)$.  Suppose that $A$ and $B$ are
  closed oriented submanifolds of $M$ of dimension $n + 1$, that both
  $\pi|_A$ and $\pi|_B$ are submersions, and that $A \cap B$ is compact.
  Consider the $n$-dimensional submanifolds $A_t = A \cap M_t$ and
  $B_t = B \cap M_t$.  Then the intersection number
  $\pair{A_t, B_t}_{M_t}$ is independent of $t \in I$.
\end{theorem}

\begin{proof}
It suffices to show that
$\pair{A_0, B_0}_{M_0} = \pair{A_1, B_1}_{M_1}$ where $I = [0, 1]$,
since any case can be reduced to this one by replacing $M$ with
$\pi^{-1}\big([a, b]\big)$ and reparametrizing $[a, b]$ by
$[0, 1]$.

Applying Theorem~\ref{thm:intersection} to $A_0$ and $B_0$ in $M_0$,
we get a compactly supported isotopy of $M_0$ that takes $A_0$ to an
$A'_0 \subset M_0$ that is transverse to $B_0$.  We also have an
analogous isotopy of $M_1$ taking $A_1$ to some $A'_1$ transverse to
$B_1$.  Since $\partial M = M_0 \cup M_1$, we can use the collar
structure near $\partial M$ to extend these isotopies to a compactly
supported isotopy of $M$ that takes $A$ to some $A'$ so that $A'_0$
and $A'_1$ are as previously constructed.  As embedded submanifolds,
both $A'$ and $B$ are transverse to $\partial M$ and hence the
inclusion map $f' \maps A' \to M$ is transverse to $B$ near
$\partial A' = A'_0 \cup A'_1$.  Now $f'$ is proper and $A' \cap B$ is
compact, so we can choose an open set $U \subset A'$ with compact
closure containing $(f')^{-1}(B)$.  As $f'$ is transverse to $B$ on
$\partial A' \cup (A' \setminus U)$, by Proposition~4.5.7 of
\cite{Wall2016} we can perturb $f'$ to $f''$ which is transverse to
$B$ without changing the values on
$\partial A' \cup (A' \setminus U)$.  Thus $(f'')^{-1}(B)$ is a closed
submanifold of $A'$ of codimension $n$ which is contained inside the
compact set $\Ubar$, that is, a finite union of circles and arcs with
endpoints in $\partial A'$.  This \1-manifold comes with an orientation
which we can use in the standard way to conclude that
$\pair{A_0', B_0}_{M_0} = \pair{A'_1, B_1}_{M_1}$ as needed.
\end{proof}

In the setting of Theorem~\ref{thm:intersection}, suppose that, in
addition to being compact, the intersection \(A\cap B\) has only
finitely many connected components.  
Now connected components are always closed, and since $A \cap B$ has
only finitely many such, each connected component \(Z\) of $A \cap B$
is also open in \(A\cap B\).  Thus we can find an open \(U \subset M\)
with \(U \cap A \cap B = Z\). In this situation, there is a local
intersection number
\(\langle A,B \rangle |_Z := \langle A \cap U, B \cap U \rangle_U \)
that is independent of the choice of $U$.  Note that
\begin{equation}\label{eq: component sum}
  \langle A, B \rangle = \sum_{Z} \langle A, B \rangle|_Z
\end{equation}
where the sum runs over the connected components \(Z\) of $A \cap
B$. 

Now suppose that we are in the situation of
Theorem~\ref{Thm:InvarianceOfIntersection}, where moreover
\(A \cap B\) has finitely many connected components.  If $W$ is a
connected component of \(A\cap B\), pick an open neighborhood $U$ of
it with $U \cap A \cap B = W$.  We then have a local intersection
number for each $t \in I$ given by
\(\langle A_t, B_t \rangle |_W := \langle A_t \cap U, B_t \cap U
\rangle_{U_t} \) which is independent of $U$.  Applying
Theorem~\ref{Thm:InvarianceOfIntersection} to
$(U, A \cap U, B \cap U)$ and using (\ref{eq: component sum}) gives:

\begin{corollary}
\label{Cor:component sum}
The intersection number \(\langle A_t, B_t \rangle|_W\) is independent
of \(t\).  If each $A_t \cap B_t$ has finitely many connected
components, then for all $t$ 
\[
  \langle A_t, B_t \rangle|_W = \sum_{Z\subset W} \langle A_t, B_t
  \rangle|_Z,
\]
where the sum is over the connected components of \(A_t\cap B_t\)
contained  in \(W\).

\end{corollary}

Finally, recall that smooth submanifolds $X$ and $Y$ of $M$ intersect
\emph{cleanly} when $X \cap Y$ is a smooth submanifold with
$T_x(X \cap Y) = T_x X \cap T_x Y$ for all $x \in X \cap Y$.  This is
a weaker notion than intersecting transversely; for example, all pairs
of affine subspaces in $\R^n$ intersect cleanly.

\begin{lemma}
  \label{Lem:Submanifold Intersection}
  Let \(M\) be a smooth oriented submanifold of \(M'\) with
  \(A',B' \subset M'\) smooth oriented submanifolds, each of which
  intersects cleanly with $M$.  Suppose \(A'\cap B' \subset M\) and
  let \(A= A'\cap M\), \(B = B'\cap M \).  Assume further that there
  are oriented bundles \(V_A, V_B\) on \(M\) such that
  \(\nu_{A'/A} \cong V_A|_{A}\), \(\nu_{B'/B} \cong V_B|_{B}\), and
  \(\nu_{M'/M} = V_A \oplus V_B\) as oriented vector bundles.  If
  \(Z\) is a compact component of \(A\cap B = A' \cap B'\) which is an
  open subset of \(A\cap B\), then
  \[
    \big(\pair{ A',B'}_{M'}\big)\big|_{Z} =
    (-1)^{(\dim A)(\dim V_B)}\big(\pair{A,B}_M\big)\big|_{Z}.\]
\end{lemma}

\begin{proof}
If \(A\) and \(B\) are transverse at \(x\) as submanifolds of \(M\),
then \(A'\) and \(B'\) are transverse at \(x\) as submanifolds of
\(M'\), and we are just comparing the oriented vector space
\(T_x A' \oplus T_x B' = \big( V_A|_x \oplus T_xA \big) \oplus \big(
V_B|_x \oplus T_xB \big)\) with
\(\nu_{M'/M} \oplus T_xM = \big(V_A|_x \oplus V_B|_x \big) \oplus
\big( T_xA\oplus T_xB \big) \).
 
In general, the assumption that \(A'\) and \(B'\) intersect cleanly
with \(M\) means that, after restricting to a tubular neighborhood of
\(M\), we can assume that \(M' = \nu_{M'/M}\), \(A' = V_A|_A\), and
\(B' = V_B|_B\). Suppose \(f_t\) is an isotopy of \(A\) such that
\(f_1(A)\) is transverse to \(B\). By choosing a connection on
\(V_A\), we can can extend \(f_t\) to an isotopy \(F_t\) of \(A'\)
with the property that \(F_t(V_A|_x) = V_A|_{f_t(x)}\). Then
\(F_1(A')\) is transverse to \(B'\) and
\(F_1(A')\cap B' = f_1(A) \cap B\), so the statement follows from the
transverse case.
\end{proof}

\section{Background on character varieties}
\label{sec: intro char var}

In this section, we establish notation for the character varieties we
will consider and recall some of their basic properties; for further
background, see \cite[\S 1.4]{CullerShalen1983}.  Beyond their mention
in the introduction, the contents of this section are not used until
Section~\ref{sec:map to X}. Let $\Gamma$ be a finitely-generated
group.  The representation variety $\RC(\Gamma)$ is
$\Hom(\Gamma, \SLC)$ viewed as an affine algebraic set over \(\C\).
The group \(\SLC\) acts on \(\RC(\Gamma)\) by conjugation, and the
character variety \(\XC(\Gamma)\) is the geometric invariant theory
quotient of \(\RC(\Gamma)\) by this action.  For each $g \in \Gamma$,
there is a regular function $\tr_g \maps \XC(\Gamma) \to \C$ defined
by $[\rho] \mapsto \tr(\rho(g))$.  The character variety $\XC(\Gamma)$
is also an affine algebraic set over \(\C\); concretely, there is a
finite set of $g_i \in \Gamma$ such that the $\tr_{g_i}$ give global
coordinates for $\XC(\Gamma)$.

We will be interested in representations where a distinguished finite
subset $S = \{s_1, \ldots s_n\}$ of $\Gamma$ share a common conjugacy
class in $\SLC$.  Here, a typical pair $(\Gamma, S)$ to keep in mind
is when $\Gamma = \pair{s_1, \ldots s_n}$ is the free group generated
by $S$.

Define the \emph{\(\SLC\)-representation variety of the pair}
\((\Gamma, S)\) to be
\[
  \RC(\Gamma, S) := \setdef{\rho \in \RC(\Gamma)}{
    \mbox{All $\rho(s_i)$ and $\rho(s_j)$ are conjugate in $\SLC$}}.
\]
Despite the name, $\RC(\Gamma, S)$ is rarely an affine algebraic set,
but rather is constructible, a term which we now recall.  If $W$ and
$V$ are affine algebraic sets (i.e.~are Zariski closed in $\C^m$),
their difference $W - V$ is a \emph{quasi-affine algebraic set}.  A
\emph{constructible set} is then a finite union of quasi-affine sets.
For context, recall that the image of an affine algebraic set under a
polynomial map is generally only constructible, but the image of a
constructible set under such a map is still constructible.  We define
$\XC(\Gamma, S)$ to be the image of
$\RC(\Gamma, S) \subset \RC(\Gamma, \emptyset) = \RC(\Gamma)$ under
the natural map $\tau \maps \RC(\Gamma) \to \XC(\Gamma)$.  Hence
$\XC(\Gamma, S)$ is a constructible subset of the locus
$\{\tr_{s_i} = \tr_{s_j}\}$ of $\XC(\Gamma)$.
 
A representation \(\rho \in \RC(\Gamma)\) is \emph{reducible} when
$\rho(\Gamma)$ leaves invariant a line in $\C^2$; otherwise $\rho$ is
\emph{irreducible}.  The corresponding subsets of $\RC(\Gamma, S)$ are
denoted \(\RC^\red(\Gamma,S)\) and \(\RC^\irr(\Gamma,S)\)
respectively.  These two sets have disjoint images under \(\tau\) and
are denoted by \(\XC^\red(\Gamma,S)\) and \(\XC^\irr(\Gamma,S)\).  The
subsets \(\RC^\red(\Gamma,S)\) and \(\XC^\red(\Gamma,S)\) are Zariski
closed in $\RC(\Gamma,S)$ and $\XC(\Gamma, S)$, respectively, and
hence \(\RC^\irr(\Gamma,S)\) and \(\XC^\irr(\Gamma,S)\) are Zariski
open.  For $\chi$ in $\XC^\irr(\Gamma, S)$, all representations in
$\tau^{-1}(\chi)$ are conjugate, whereas this need not be the case for
$\chi$ in $\XC^\red(\Gamma, S)$.

We will need the following general fact about character varieties.
\begin{lemma}
\label{Lem:PropInc}
If \(\phi \maps \Gamma_1 \to \Gamma_2\) is an epimorphism, then the induced map
\(\phi^* \maps \XC(\Gamma_2) \to \XC(\Gamma_1)\) given by $[\rho] \mapsto [\rho
\circ \phi]$ is proper in the classical topology.
\end{lemma}

\begin{proof}
  First note that for any $\Gamma$, a closed subset
  $C \subset \XC(\Gamma)$ is compact if and only if $\tr_g(C)$ is
  bounded for all $g \in \Gamma$; this is because the functions
  $\tr_g$ give coordinates on $\XC(\Gamma)$.  So suppose
  $K \subset \XC(\Gamma_1)$ is compact and set $L = (\phi^*)^{-1}(K)$ in
  $\XC(\Gamma_2)$.  Given $g_2 \in \Gamma_2$, choose $g_1 \in \Gamma_1$
  with $g_2 = \phi(g_1)$.  Then $\tr_{g_2}(L) = \tr_{g_1}(K)$ and
  hence the former is bounded as the latter is.  Thus $L$ is compact
  and $\phi^*$ is proper.
\end{proof}

\subsection{Real characters}
\label{sec:realchar}

Our focus throughout this paper is on the real characters of
\(\Gamma\), that is, the subset $X(\Gamma, S)$ of \(\XC(\Gamma,S)\)
where $\tr_g (\chi) \in \R$ for all $g \in \Gamma$.  The algebraic set
\(\XC(\Gamma, \emptyset)\) is defined over $\Q$ \cite[\S
III.1]{MorganShalen1984}, and the same is true for $\XC(\Gamma, S)$;
hence, complex conjugation preserves \(\XC(\Gamma,S)\) with fixed
point set the real locus $X(\Gamma, S)$. In particular, the subset
$X(\Gamma, S)$ is a \emph{real semialgebraic set}, i.e., it is
defined by a system of real polynomial equations and inequalities.
Note that \(\tau^{-1}\big(X(\Gamma,S)\big)\) is usually considerably
larger than the real locus of \(\RC(\Gamma,S)\). However, we will show
that the points in \(X(\Gamma,S)\) all come from representations into
the subgroups $\SU$ and $\SLR$ of $\SLC$.  For any subgroup $H$ of
$\SLC$, let $R_H(\Gamma, S)$ be the subset of $R_\C(\Gamma, S)$
consisting of $\rho$ where $\rho(\Gamma) \subset H$, and let
$X_H(\Gamma, S)$ be the image of $R_H(\Gamma, S)$ under $\tau$.

\begin{proposition}
  \label{prop: really basic}
  The set \(X(\Gamma,S)\) is
  \(X_\SU(\Gamma,S) \cup X_\SLR(\Gamma,S)\). Moreover, for
  \[
    D = \SU \cap \SLR = \left\{\mysmallmatrix{\cos t}{-\sin t}{\sin
        t}{\cos t} \right\} \cong S^1
  \]
  we have
  \(X_\SU(\Gamma,S) \cap X_\SLR(\Gamma,S) = X^\red(\Gamma,S)=
  X_D(\Gamma, S)\).  Consequently, $X^\irr(\Gamma, S)$ is the disjoint
  union of $X^\irr_\SU(\Gamma,S)$ and $X^\irr_\SLR(\Gamma, S)$.
\end{proposition}
\begin{proof}
By Proposition III.1.1 of \cite{MorganShalen1984}, we have that
\(X(\Gamma, \emptyset) = X_\SU(\Gamma, \emptyset) \cup X_\SLR(\Gamma,
\emptyset)\).  Moreover, by Lemma 2.10 of \cite{CullerDunfield2018},
we have
\(X_\SU(\Gamma, \emptyset) \cap X_\SLR(\Gamma, \emptyset) =
X^\red(\Gamma, \emptyset)= X_D(\Gamma, \emptyset)\).  To transfer
these facts to \(X(\Gamma, S)\), it suffices to show
\[
  X_H(\Gamma, S) = X_H(\Gamma, \emptyset) \cap X(\Gamma, S)
  \mtext{for $H$ each of $D$, $\SU$, and $\SLR$.}
\]
As the definitions immediately give
$X_H(\Gamma, S) \subset X_H(\Gamma, \emptyset) \cap X(\Gamma, S)$, we
need only show the other containment.

Given $\chi \in X_D(\Gamma, \emptyset) \cap X(\Gamma, S)$, choose
$\rho \in R_D(\Gamma, \emptyset)$ with $\tau(\rho) = \chi$.
Now, on $X(\Gamma, S)$ we have $\tr_{s_i} = \tr_{s_j}$ for all
$s_i, s_j \in S$.  Moreover, elements of $D$ are conjugate in $\SLC$
if and only if they have the same trace.  Hence $\rho$ must be in
$R_D(\Gamma, S)$ and hence $\chi$ is also in $X_D(\Gamma, S)$ as
needed.

The identical argument works for $\SU$, but not $\SLR$ as there are
two conjugacy classes of elements with trace $2$ and $-2$.  However,
since $X^\red_\SLR(\Gamma, \emptyset) = X_D(\Gamma, \emptyset)$, we
need only consider $\chi \in X^\irr_\SLR(\Gamma, \emptyset)$, where
every fiber of $\tau$ in $R_\C(\Gamma, \emptyset)$ consists of
conjugate representations. In particular, if
$\chi \in X^\irr_\SLR(\Gamma, \emptyset) \cap X(\Gamma, S)$ comes from
$\rho \in R^\irr_\SLR(\Gamma, \emptyset)$, it follows that
$\rho \in R^\irr_\SLR(\Gamma, S)$ since some $\rho'$ with the same
character is in $\RC(\Gamma, S)$.  This completes the proof.
\end{proof}

Both $X_\SU(\Gamma,S)$ and $X_\SLR(\Gamma,S)$ are real semialgebraic
sets since they are images of the real semialgebraic subsets
$R_\SU(\Gamma, S)$ and $R_\SLR(\Gamma, S)$ of $\RC(\Gamma, S)$ under
the polynomial map $\tau$ that is defined over $\Q$.  Note that while
$X_\SU(\Gamma, \emptyset)$ can be identified with the topological
quotient $R_\SU(\Gamma, \emptyset)/\SU$, the analogous statement for
$\SLR$ is false.  This is for two reasons: first, because there are
non-conjugate elements of \(\SLR\) with the same trace, and second,
because the normalizer of $\SLR$ in $\SLC$ is larger than just
$\SLR$. To be precise, the normalizer is the full stabilizer of the
copy of $\H^2 \subset \H^3$ that is preserved by $\SLR$, and hence is
a double-cover of the disconnected group $\Isom(\H^2)$. We use
$\SLRpm$ to denote the normalizer of $\SLR$ in $\SLC$.

Note also that 
\[
  R_\SU(\Gamma, S) = \setdef{\rho \in R_\SU(\Gamma)}{\mbox{All
      $\rho(s_i)$ and $\rho(s_j)$ are conjugate in $\SU$}}
\]
and
\[
 R_\SLR(\Gamma, S) = \setdef{\rho \in R_\SLR(\Gamma)}{\mbox{All
      $\rho(s_i)$ and $\rho(s_j)$ are conjugate in $\SLRpm$}}.
\]
Recalling that for $\chi$ in $\XC^\irr(\Gamma, S)$, all
representations in $\tau^{-1}(\chi) \subset \RC(\Gamma, S)$ are
conjugate, it is easy to show:
\begin{lemma}\label{lem: SLRpm quo}
  As topological spaces, $X_\SLR^\irr(\Gamma, S)$ is the quotient of
  $R_\SLR^\irr(\Gamma, S)$ by the conjugation action of $\SLRpm$.
\end{lemma}

\subsection{The free group}
\label{subsec:X(F_n)}

We now establish some basic facts about the character variety of the
free group. We start by considering \(\XC(F_n,\emptyset)\) for the
free group $F_n$ on generators $S = \{s_1,\ldots, s_n\}$.  By
Proposition~5.8 of \cite{HeusenerPorti2004}, for \(n>2\), the smooth
locus of \(\XC(F_n, \emptyset)\) is precisely
\(\XC^\irr(F_n,\emptyset)\); for $n \leq 2$, the entirety of
$\XC(F_n, \emptyset)$ is smooth.  We define
\[X^*_\C(F_n) = \setdef{[\rho] \in X^\irr_\SLC(F_n,\emptyset)}{
    \mbox{$\rho(s_i) \neq \pm I$ for all $i$}}
\]
which is a Zariski-open subset of \(\XC^\irr(F_n,\emptyset) \) as one has
\[\XC^*(F_n) = \setdef{\chi \in \XC^\irr(F_n,\emptyset)}{
    \mbox{for all $i$ there exists $j$ with $\chi\big([s_i, s_j]\big) \neq 2$}}.
\]
Since $\chi \in  \XC^\red(F_n,\emptyset)$ satisfies $\chi\big([s_i,
s_j]\big) = 2$ for all $i, j$, we have
\[
\XC^*(F_n) = \XC(F_n, \emptyset) \setminus \bigcup_{i = 1}^n
\setdef{\chi \in \XC(F_n, \emptyset)}{\mbox{$\chi\big([s_i, s_j]\big)
    = 2$ for all $j$}}. 
\]
In particular, $\XC^*(F_n)$ is a quasi-affine algebraic set. 

Consider
the map \(\tr_S: \XC^*(F_n) \to \C^n\) given by
\(\tr_S(\chi) = \big(\chi(s_1), \ldots, \chi(s_n)\big)\).  For
\(c\in \C\), we take $\bc = (c, \ldots, c) \in \C^n$ and define
\[
  \XC^{c,\irr}(F_n,S) := \tr^{-1}_S(\bc) \subset \XC^*(F_n).
\] Observe that
\(\displaystyle \XC^\irr(F_n,S) = \bigcup_{c \in \C} \XC^\cirr(F_n,S)
\).
 
\begin{lemma}\label{lem:c slice smooth}   
  The sets \(\XC^{c,\irr}(F_n,S) \) and \(\XC^\irr(F_n,S)\) are
  quasi-affine and smooth.  The map
  \(\tr \maps \XC^\irr(F_n,S) \to \C\) given by
  $\chi \mapsto \chi(s_1)$ is a submersion, with
  $\tr^{-1}(c) = \XC^{c,\irr}(F_n,S)$.
\end{lemma}
\begin{proof}
For clarity, we denote the map $\XC^\irr(F_n,S) \to \C$ in the
statement as $\tr_{s_1}$ to distinguish it from
$\tr_S \maps \XC^*(F_n) \to \C^n$ and also \(\tr: \SLC \to \C\).
Set
\[
  \RC^*(F_n) = \setdef{\rho\maps F_n \to \SLC}{
    \mbox{$\rho(s_i) \neq \pm I$ for all $i$}} \cong  (\SLC
  \setminus \{\pm I\})^n,
\]
so that \(\tau\big(\RC^*(F_n)\big) = \XC^*(F_n)\).  The only critical
points of the map \(\tr: \SLC \to \C\) are \(\pm I\), so
\(\tr_S \circ \tau: \RC^*(F_n) \to \C^n\) is a surjective submersion.
It follows that $\tr_S$ is a surjective submersion.  Now
\(\XC^\irr(F_n,S) = \tr_S^{-1}(\Delta)\), where
\(\Delta = \{\bc \, | \, c \in \C\} \subset \C^n\), so
$\XC^\irr(F_n,S)$ is a smooth quasi-affine algebraic set.  
The restriction $\tr_S \maps \XC^\irr(F_n,S) \to \Delta$ is also a
surjective submersion whose composition with any coordinate projection
$\Delta \to \C$ is $\tr_{s_1} \maps \XC^\irr(F_n, S) \to \C$; hence
the latter map is a surjective submersion and
$\XC^\cirr(F_n, S) = \tr_{s_1}^{-1}(c)$ is also smooth and
quasi-affine.
\end{proof}  

For \(c \in \R\), we define \(X^\cirr (F_n,S)\) and \(X^\irr(F_n,S)\)
to be the real loci of \(\XC^\cirr (F_n,S)\) and \(\XC^\irr(F_n,S)\),
respectively. Since they are the real loci of smooth algebraic sets,
we see that:
\begin{corollary}
\label{cor: real part smooth} 
The sets \(X^\cirr (F_n,S)\) and \(X^\irr(F_n,S)\) are smooth.
\end{corollary}

 Using Proposition~\ref{prop: really basic}, we have
\[
  X^\cirr(F_n,S) = X^\cirr_\SU (F_n,S) \coprod X^\cirr_\SLR (F_n,S)
  \mtext{and}
  X^\irr(F_n,S) = \bigcup_{c \in \R} X^\cirr(F_n,S),
\]
as well as the following consequence of
Lemma~\ref{lem:c slice smooth}:
\begin{corollary}
  \label{cor:tr is submersion}
  The map \(\tr: X^\irr(F_n,S) \to \R\) is a submersion.
\end{corollary}
We define $\RC^\cirr(F_n, S)$ inside $\RC^\irr(F_n, S)$ to be the
preimage of $\XC^\cirr(F_n, S)$ under $\tau$ and use analogous
notation for the other target groups.

\begin{lemma}\label{lem:charsubmer}
  The maps
  \(\tau \maps R^\cirr(F_n, S) \to X^\cirr(F_n, S)\) and
  \(\tau \maps R^\irr(F_n, S) \to X^\irr(F_n, S)\) are
  submersions. 
\end{lemma}

\begin{proof}
As $X^\irr(F_n, S)$ is the disjoint union of $X^\irr_\SU(F_n, S)$ and
$X^\irr_\SLR(F_n, S)$ by Proposition~\ref{prop: really basic}, this
breaks up into two cases, corresponding to $\SU$ and $\SLR$.  We do
$\SLR$; the other case is identical.

By Lemma~2.11 of \cite{CullerDunfield2018}, for any group $\Gamma$ the
map
$\tau \maps R_\SLR(\Gamma, \emptyset) \to X_\SLR(\Gamma, \emptyset)$
has the weak-path lifting property: given
$f \maps I \to X_\SLR(\Gamma, \emptyset)$ there exists
$\ftil \maps I \to X_\SLR(\Gamma, \emptyset)$ with
$f = \tau \circ \ftil$.  By that same lemma, if $f(0)$ is in
$X_\SLR^\irr(\Gamma, \emptyset)$ then $\ftil(0)$ can be required to be
any representation in $\tau^{-1}\big(f(0)\big)$.  To show
\(\tau \maps R^\cirr_\SLR(F_n, S) \to X^\cirr_\SLR(F_n, S)\) is a
submersion, consider $\rho \in R^\cirr_\SLR(F_n, S)$ and $v$ any
tangent vector to $X^\cirr_\SLR(F_n, S)$ at $\tau(\rho)$.  Choose a
smooth path $f \maps I \to X^\cirr_\SLR(F_n, S)$ with
$f(0) = \tau(\rho)$ and $f'(0) = v$.  Let
$\ftil \maps I \to R_\SLR(F_n, \emptyset)$ be a lift with
$\ftil(0) = \rho$.  As $R^\cirr_\SLR(F_n, S)$ is the preimage under
$\tau$ of $X^\cirr_\SLR(F_n, S)$, it follows that $\ftil(I)$ is
contained in $R^\cirr_\SLR(F_n, S)$ and so $\ftil'(0)$ is a tangent
vector to $R^\cirr_\SLR(F_n, S)$ that maps to $v$ under $D\tau$,
proving $\tau$ is a submersion.  The proof that
\(\tau \maps R^\irr_\SLR(F_n, S) \to X^\irr_\SLR(F_n, S)\) is a
submersion is identical.
\end{proof}

\begin{lemma}
  \label{lem: empty char vars}
  When $c$ is not in $(-2, 2)$, the set $X^\cirr_\SU(F_n, S)$ is
  empty.  For $c = \pm 2$, every $[\rho] \in X^\cirr_\SLR(F_n, S)$ has
  $\rho(s_i) \neq \pm I$ for all $i$.  The same holds for any
  pair $(\Gamma, S)$ where $\Gamma$ is generated by conjugates of
  elements of $S$.
\end{lemma}

\begin{proof}
The first claim is immediate for $\abs{c} > 2$ since $\SU$ has no
elements with that trace.  For $\abs{c} = 2$, if
$[\rho] \in X^c_\SU(F_n, S)$ the only possibilities for $\rho(s_i)$
are $I$ or $-I$ depending on the sign of $c$, which means the image of
$\rho$ is contained in $\{\pm I\}$ and so $[\rho]$ is reducible; hence
$X^{\pm 2, \irr}_\SU(F_n, S)$ is empty.  Finally, suppose $c = \pm 2$
and $[\rho] \in X^\cirr_\SLR(F_n, S)$ with
$\rho \in R^\irr_\SLR(F_n, S)$.  As all $\rho(s_i)$ are conjugate by
definition, if any $\rho(s_i) = \pm I$ then all are, which is
impossible as $\rho$ is irreducible.
\end{proof}

We will also need the following standard lemma:

\begin{lemma}\label{lem: PSLR cent}
  If an irreducible $\rho \maps \Gamma \to \PSLR$ is centralized by
  $A \in \PSLR$ then $A = I$.
\end{lemma}

\begin{proof}
We provide the proof to caution that the conclusion does not hold when
we allow $A$ to be in the larger group $\PSLC$.  Following Section 3
of \cite{HeusenerPorti2004}, a representation
$\rhohat \maps \Gamma \to \PSLC$ is $\Ad$-reducible if it leaves
invariant either a point in $P^1(\C)$ \emph{or} a geodesic in $\H^3$.
Moreover, $\rho$ is $\Ad$-irreducible exactly when its stabilizer
under the conjugation action of $\PSLC$ is trivial by
Proposition~3.16(ii) of \cite{HeusenerPorti2004}.  So if
$\rho \maps \Gamma \to \PSLR$ is Ad-irreducible, then the lemma holds.
If instead $\rho \maps \Gamma \to \PSLR$ is irreducible but
Ad-reducible, then it is a metabelian representation into the
stabilizer of some geodesic $L$ in $\H^3$, which is a copy of
$\C^\times \rtimes (\Z/2)$.  Since any finite representation into
$\PSLR$ is reducible, $\rho(\Gamma)$ contains elements that translate
along $L$.  This forces $L$ to be in $\H^2$.  As $\rho$ is
irreducible, $\rho(\Gamma)$ also contains elliptic elements of order
two that interchange the two endpoints of $L$.  Thus there is a unique
nontrivial centralizer of $\rho(\Gamma)$, namely the elliptic of order
two with fixed point set $L$.  As that element is not in $\PSLR$, we
are done.
\end{proof}

\subsection{Reducibles}
\label{subsec:reducibles}
Suppose \(\rho \in R^\red_\C(F_n,S)\), so \(\rho(F_n)\) preserves a
line in \(\C^2\). Conjugating by an element of \(\SLC\), we may assume
this line is spanned by $e_1$. Hence the elements of
\(\rho(\Gamma)\) are upper triangular matrices, say
\(\rho(s_j) = \mysmallmatrix{\lambda_j}{*}{0}{\lambda_j^{-1}}\).
Since \(\rho\) has the same character as an abelian representation,
namely the one which sends $s_j$ to
\(\mysmallmatrix{\lambda_j}{0}{0}{\lambda_j^{-1}}\), we see
\(X^\red_\C(F_n,S)\) is the image of the subset
\(\RC^{\mathrm{ab}}(F_n,S)\) of abelian representations.

Now suppose that \([\rho]\) is a real character, moveover one in
$X^c(F_n, S)$, so for each generator we have
\(\tr \rho(s_j) = \lambda_j+ \lambda_j^{-1}= c \in \R\).  If we
further assume that \(c \in (-2,2)\), then \(\lambda_j \in S^1\).
Letting \(\alpha = \cos^{-1}(c/2)\), we can write
\(\lambda_j = e^{i\epsilon_j \alpha}\), where \(\epsilon_j = \pm 1\).
Hence \([\rho]\) is determined by the vector
\(\epsilon = (\epsilon_1, \ldots, \epsilon_n) \in \{\pm 1\}^n\). It is
easy to see that the representations determined by \(\epsilon\) and
\(\epsilon'\) have the same character if and only if
\(\epsilon = \pm \epsilon '\).  Hence, the set \(X^{c, \red}(F_n,S)\)
is in bijection with \(\{\pm 1\}^n/ \pm \One\), where
$\One = (1, \ldots, 1)$. We denote the reducible character
associated to \([\epsilon] \in \{\pm 1\}^n/ \pm \One\) by
\(\chiep^c\), or just \(\chiep\) when \(c\) is clear from context.

\subsection{Knot complements}
\label{sec: knot reps}

Suppose $M$ is the exterior of a knot \(K\) in $S^3$, and
$\mu \in \pi_1(\partial M)$ is a meridian.  We will use $X(K)$ to
denote $X(M, \{\mu\})$; as $S = \{\mu\}$ has only one element,
$X(M, \{\mu\}) = X(M)$, and so the $X(K)$ notation is just encoding
that $\tr \maps X(K) \to \R$ is $\tr_\mu$. Similarly, we define
$X^c(K) := X^c(M, \{\mu\})$.  The knot $K$ is \emph{small} when $M$
does not contain a closed essential surface. Our main results will
apply to small knots as well as to the following much larger class.
Specifically, we say $K$ is \emph{real representation small} if the
preimage of $[-2, 2]$ under $\tr \maps X(K) \to \R$ is compact.  As
the nomenclature suggests:

\begin{lemma}\label{lem:small to few chars}
  If \(K\) is small then it is real representation small. Moreover,
  for each $c \in \R$ the set $X^c(K) = \tr^{-1}(c)$ is finite.
\end{lemma}
We prove this below using the Culler-Shalen machinery.

\begin{proof}[Proof of Lemma~\ref{lem:small to few chars}]
We need to prove $E := \tr_\mu^{-1}\big([-2, 2]\big)$ in $X(M)$ is
compact.  Consider the complex affine algebraic set \(X_\C(M)\).  By
Section 2.4 of \cite{CCGLS}, as $M$ is small, every irreducible
component of \(X_\C(M)\) has complex dimension 1.  As there are
finitely many such components $C$, it suffices to prove each
$E \cap C$ is compact.  Following \cite[\S 1.3]{CullerShalen1983}, let
$\Ctil$ be a smooth projective curve with a birational isomorphism
$\iota \maps \Ctil \to C$.  In particular, $\Ctil$ is a compact
Riemann surface.  The points where $\iota$ is regular (i.e.~defined)
will be denoted $\Ctil_{\reg}$; the finitely many points in
$\Ctil \setminus \Ctil_\reg$ are called \emph{ideal points}.  The
trace map \(\tr_\mu\maps C \to \C\) induces a \emph{regular} map
\(\tr_\mu \maps \Ctil \to P^1(\C)\).  

We next claim $\Etil := \iota^{-1}\left(E \cap C\right)$ is closed in
$\Ctil$, where throughout this proof open and closed refer to the
classical rather than Zariski topology. We know $E$ is closed in
$X(M)$ which is closed in $X_\C(M)$; thus, $E \cap C$ is closed in
$C$.  Hence $\Etil$ is closed in $\Ctil_\reg$, so it suffices to show
that no ideal point $x$ is a limit point of $\Etil$.  If $x$ were such
an ideal point, then $\tr_\mu(x) \in [-2, 2]$ by continuity of
$\tr_\mu \maps \Ctil \to P^1(\C)$.  By Chapter~1 of
\cite{CullerGordonLueckeShalen1987}, there is an essential surface $S$
associated to \(x\), and as $\tr_\mu(x) \neq \infty$, it follows that
$S$ is either closed or has boundary some number of parallel copies of
$\mu$.  In the latter case, Theorem~2.0.3 of
\cite{CullerGordonLueckeShalen1987} implies that $M$ also contains a
closed essential surface.  In either case, this contradicts that $M$
is small, so $\Etil$ is closed in $\Ctil$.  As $\Ctil$ is compact, we
have $\Etil$ is compact and, by definition, it is contained in
$\Ctil_\reg$.  As $\iota(\Ctil_\reg)$ is all of $C$, we have
$E \cap C = \iota(\Etil)$ is the continuous image of a compact set and
hence compact as desired.  So $K$ is real representation small.

The final claim that each $X^c(K)$ is finite follows provided each
\(\tr_\mu \maps \Ctil \to P^1(\C)\) is nonconstant.  This is the case
since otherwise we would again have an ideal point where $\tr_\mu(x)
\neq \infty$, violating smallness.
\end{proof}

\subsection{Reducible representations and the Alexander polynomial}
\label{sec: alex and the reds}

As noted in Section~\ref{subsec:reducibles}, the subset $X_\C^\red(K)$
of reducible characters is the same as the characters of diagonal
representations.  As the latter factor through $H_1(M; \Z) \cong \Z$,
we can parameterize $X_\C^\red(K)$ by $\tr_\mu(\chi) \in \C$; we will
use $\chi^c$ to denote the reducible character with $\tr_\mu = c$.

For the Alexander polynomial $\Delta_K(t)$ of $K$ define
\[
  D_K: = \setdef{a + 1/a}{\mbox{$a \in \C$ with $\Delta_K(a^2) = 0$}}.
\]
Note that \(\Delta_K(\pm 1)\) is always an odd integer, so $D_K$ is
finite with $0$ and $\pm 2$ not in $D_K$.  By de Rahm, the set $D_K$
is exactly the $c \in \C$ where there is a reducible representation
$\rho$ to $\SL{2}{\C}$ with \emph{nonabelian image} and character
$\chi^c$, see \cite[\S 6]{CCGLS} or
\cite[Proposition~14.6]{BurdeZieschangHeusener2014}.  Any reducible
character that is the limit of irreducible ones must come from $D_K$:

\begin{lemma}[{\cite[Prop.~6.2]{CCGLS}}]
  \label{lemma: deform red}
  If a reducible character $\chi^c$ is in the Zariski closure of
  $X^\irr_\C(K)$, then $c \in D_K$.
\end{lemma}

\section{Representations to the universal covering group}

\label{sec: tSLR basics}

Next, we discuss the basics of representations into the Lie group
\(\tSLR\), the universal cover of \(\SLR\); throughout, see
\cite[\S~3]{CullerDunfield2018} for background and the basic facts
that we use, noting that $\tSLR = \PSLRtilde$.  This
material will not be used until Section~\ref{sec: trans ex locus}.

The Lie group \(\tSLR\) is the universal central extension of \(\SLR\)
by \(\pi_1(\SLR) \cong \Z\):
\begin{equation}\label{eq: tSLR}
  1 \to \Z \to \tSLR \to \SLR \to 1
\end{equation}
Note the $\Z$ in (\ref{eq: tSLR}) is actually index two in the center
$\pair{c} \cong \Z$ in $\tSLR$, with $c \mapsto -I$ under
the map $\tSLR \to \SLR$   and $\tSLR \, / \! \pair{c} = \PSLR$.

For a finitely generated group $\Gamma$, the obstruction to lifting
\(\rho \maps \Gamma \to \SLR\) to a representation
\(\widetilde{\rho} \maps \Gamma \to \tSLR\) is the Euler class
\(e(\rho) \in H^2(\Gamma;\Z)\). When \(e(\rho)=0\), the cohomology
group \(H^1(\Gamma;\Z)\) acts freely on the set of lifts: viewing
$\alpha \in H^1(\Gamma; \Z)$ as a homomorphism from $\Gamma$ to the
$\Z$ in (\ref{eq: tSLR}), we set
$(\alpha \cdot \rhotil)(g) = \alpha(g) \rhotil(g)$ for $g \in \Gamma$,
which is another representation as $\alpha(\Gamma)$ is central in
$\tSLR$.  Consider the set of liftable $\SLR$ representations:
\[
  R_\SLR(\Gamma)_0 = \setdef{\rho \in R_\SLR(\Gamma)}{e(\rho) = 0}
\]
which is a union of connected components of $R_\SLR(\Gamma)$ by
\cite[\S 3.3]{CullerDunfield2018}.  The group $\tSLR$ is not linear and
so has a real analytic rather than real algebraic structure. Still,
\(R_\tSLR(\Gamma)\) makes sense as a real analytic variety.  The
cohomology group $H^1(\Gamma; \Z)$ acts freely and properly
discontinuously on \(R_\tSLR(\Gamma)\) with quotient
\(R_\SLR(\Gamma)_0\); in particular, \(R_\tSLR(\Gamma) \to
R_\SLR(\Gamma)_0\) is a covering map.

Before introducing the analog of the character variety for $\tSLR$, we
take a brief detour to study
\(\XA_\SLR(\Gamma) = R_\SLR(\Gamma)\sslash \sim\), where \(\sslash\)
denotes the geometric invariant theory quotient and \(\sim\) is the
equivalence relation generated by conjugation by elements of $\SLR$.
Equivalently, $\XA_\SLR(\Gamma)$ is the polystable quotient of
$R_\SLR(\Gamma)$ under the action of $\SLR$ by conjugation, as
discussed in Section~\ref{sec: pain} below.  This differs from
$X_\SLR(\Gamma)$ which is $R_\SLR(\Gamma) \sslash \SLRpm$ as discussed
in Section~\ref{sec:realchar}.  In particular, the conjugation action
of \(\mysmallmatrix{i}{0}{0}{-i}\) on $R_\SLR(\Gamma)$ descends to a
$\Z/2$ action that quotients $\XA_\SLR(\Gamma)$ down to
$X_\SLR(\Gamma)$; thus the natural map
$\XA_\SLR(\Gamma) \to X_\SLR(\Gamma)$ is a (possibly branched) 2-fold
cover.  As with $X_\SLR(\Gamma)$, the space $\XA_\SLR(\Gamma)$ is a
real semialgebraic set \cite[Theorem~7.6]{RichardsonSlodowy1990}, and
thus Hausdorff and locally contractible by
\cite[Theorem~5.43]{BasuPollackRoy2006}.  The following will be needed
in Section~\ref{sec: ext Lin}:

\begin{lemma}
  \label{lem: SLR branching}
  The map $\XA_\SLR(\Gamma) \to X_\SLR(\Gamma)$ can be branched only
  at reducible characters. In particular, each
  $\chi \in X_\SLR^\irr(\Gamma)$ has two distinct preimages in
  $\XA_\SLR^\irr(\Gamma)$.
\end{lemma}

\begin{proof} Lest you think the claim obvious, we point out that it
fails for $\PSLR$ (consider an infinite dihedral group in $\PSLR$
generated by two distinct elements of order 2).  Suppose
$\rho \in R_\SLR(\Gamma)$ represents $\chi$.  The claim is equivalent
to showing that the $\SLRpm$ orbit of $\rho$ has two connected
components (each of which is then a distinct $\SLR$ orbit, giving two
distinct points in $\XA_\SLR^\irr(\Gamma)$).  Since the stabilizer of
an irreducible representation in $\SLC$ is just $\{\pm I\}$ (e.g.~use
the setup of the proof of Lemma 1.5.1 of \cite{CullerShalen1983}), we
see $\SLRpm \cdot \rho$ is disconnected by the orbit-stabilizer
theorem.
\end{proof}

We will define $\Xtil_\tSLR(\Gamma)$ to be the polystable quotient of
$R_\tSLR(\Gamma)$ under the action of $\tSLR$ by conjugation, see
Section~\ref{sec: pain} and especially Corollary~\ref{cor: X of Gtil}
for details.  The main result of this section tells us that while
$\tSLR$ is not an algebraic group, it still has a very reasonable
``character variety'' in the form of $\Xtil_\tSLR(\Gamma)$:
\begin{theorem}
  \label{thm: tSLRchar}
  The image $\XA_\SLR(\Gamma)_0$ of $R_\SLR(\Gamma)_0$ in
  $\XA_\SLR(\Gamma)$ is a union of connected components.  The map
  $R_\tSLR(\Gamma) \to R_\SLR(\Gamma)_0$ induces a regular covering
  map $\Xtil_\tSLR(\Gamma) \to \Xtil_\SLR(\Gamma)_0$, where the
  covering group is $H^1(\Gamma; \Z)$.
\end{theorem}
We point out that the final conclusion of this theorem does not hold
for the central extension $\Z/2 \to \SLC \to \PSLC$; the map
$X_\SLC(\Gamma) \to X_\PSLC(\Gamma)$ can have branching, at
irreducible representations no less \cite[\S 4]{HeusenerPorti2004}.
The proof of Theorem~\ref{thm: tSLRchar} is unconnected with the rest
of this paper, and so the trusting reader can skip ahead to
Section~\ref{sec:config}.

\subsection{Polystable quotients}
\label{sec: pain}

Following Section 7.2 of \cite{RichardsonSlodowy1990}, consider the
following construction.  Let $H$ be a locally compact group acting on
a locally compact space $V$.  This action is \emph{polystable} when
for all $v$ in $V$ there is a unique closed $H$-orbit in the closure
of $H \cdot v$.  For a polystable action, let $V^*$ be those points in
$V$ whose $H$-orbits are closed, which is called the set of
\emph{polystable points}.  Now define the set $V \sslash H$ to be
$V^*/H$.  Let $\pi \maps V \to V \sslash H$ send $v$ to $[x]$ where
$H \cdot x$ is the unique closed orbit in the closure of $H \cdot v$.
Finally, give $V \sslash H$ the quotient topology induced by $\pi$;
this space is called the \emph{polystable quotient}. As the notation
suggests, the space $V \sslash H$ is also the geometric invariant
theory quotient in many situations, including when $H$ is a reductive
real or complex linear algebraic group acting on an affine variety
\cite{RichardsonSlodowy1990}.  In particular, this gives an alternate
perspective on the construction of $R_\C(\Gamma) \to X_\C(\Gamma)$.
Since the polystable quotient is purely a topological construction, we
can try to apply it to the conjugation action of $\tSLR$ on
$R_\tSLR(\Gamma)$ where there is no geometric invariant theory
quotient to be had; we succeed in doing so in Corollary~\ref{cor: X of
  Gtil}.

For ease of notation, set $G = \SLR$, $\Gtil = \tSLR$, and
$\Gbar = \PSLR = G/Z(G) = \Gtil/Z(\Gtil)$.  We are interested in the
conjugation action of $\Gbar$ on $R_G(\Gamma)$ and $R_\Gtil(\Gamma)$.
There are five types of orbits $\Gbar \cdot \rho$ for
$\rho \in R_G(\Gamma)$:
\begin{enumerate}
\item When the image of $\rho$ is in $Z(G)$, then
  $\Gbar \cdot \rho = \{\rho\}$ which is closed.
  
\item When $\rho$ is irreducible, the orbit $\Gbar \cdot \rho$ is
  closed and homeomorphic to $\Gbar$ itself (compare Lemma~\ref{lem:
    PSLR cent}).
  
\item When $\rho$ is a reducible and has a unique fixed point in
  $\H^2$, then the orbit $\Gbar \cdot \rho$ is closed and homeomorphic
  to $\H^2$.
  
\item \label{item: diag} When $\rho$ is reducible and has exactly two
  fixed points in $P^1(\R)$, then $\rho$ is conjugate to a
  representation into
  $\setdef{\mysmallmatrix{a}{0}{0}{1/{a}}}{a \in \R^\times}$. The
  orbit $\Gbar \cdot \rho$ is again closed and homeomorphic to
  $S^1 \times \R$.
  
\item \label{item: nonclosed}
  When $\rho$ is reducible and has a unique fixed point in $P^1(\R)$,
  then the orbit $\Gbar \cdot \rho$ is not closed.  There is a unique
  orbit of type (\ref{item: diag}) in the closure of
  $\Gbar \cdot \rho$: if $\rho$ is upper triangular, that orbit is
  that of its ``diagonal part''.  Indeed, if $\rho$ is upper
  triangular, define $\rho_t$ for $t \in [0, 1]$ by
  $\rho_t(g) = \mysmallmatrix{a}{b t^2}{0}{a^{-1}}$ when
  $\rho(g) = \mysmallmatrix{a}{b}{0}{a^{-1}}$.  This is a smooth path
  in $R_G(\Gamma)$ where $\rho_0$ is diagonal and $\rho_t$ is
  conjugate to $\rho$ for $t > 0$ by $\mysmallmatrix{t^{-1}}{0}{0}{t}$.
\end{enumerate}

For an action of a group $H$ on a space $V$, for $v \in V$ define
\[
  [v]_H = \setdefm{\big}{w \in V}%
  {\mbox{the closures of $H\cdot v$ and $H \cdot w$ intersect}}.
\]
When $H$ acts polystably on $V$ with induced map
$\pi \maps V \to V \sslash H$, then $[v]_H = \pi^{-1}\big(\pi(v)\big)$
for each $v$ in $V$.  The basic facts about the action of $\Gbar$ on
$R_G(\Gamma)$ and its quotient are then:

\begin{lemma}
  \label{lem: polystable}
  The action of $\Gbar$ on $R_G(\Gamma)$ is polystable.  The fibers of
  $\tau \maps R_G(\Gamma) \to \XA_G(\Gamma) = R_G(\Gamma)\sslash
  \Gbar$ are closed and path connected. Finally
  $\tau^{-1}\big(\tau(\rho)\big) = [\rho]_\Gbar$ for each
  $\rho \in R_G(\Gamma)$.
\end{lemma}

\begin{proof}
The first claim can be deduced with some effort from the concrete
description of the orbits, but it is also a general fact about actions
of real reductive groups \cite[Theorem~7.3.1]{RichardsonSlodowy1990}.
As noted above, $\XA_G(\Gamma)$ is Hausdorff, so points are closed and
hence so are their preimages under the continuous map $\tau$.  That the
fibers are path connected follows from the discussion in (\ref{item:
  nonclosed}).  Finally, as previously mentioned, the last statement
holds for any polystable action.
\end{proof}

We turn now to studying $\pi \maps R_\Gtil(\Gamma) \to R_G(\Gamma)_0$,
which we recall is a regular cover with deck group $H^1(\Gamma; \Z)$.

\begin{lemma}
  \label{lem: orbits upstairs}
  For each $\rho \in R_G(\Gamma)_0$, there is an open set $U$
  containing $[\rho]_\Gbar$ where the following holds: For every lift
  $\rhotil \in R_\Gtil(\Gamma)$ of $\rho$ there is an open
  $\Util \subset \pi^{-1}(U)$ containing $[\rhotil]_\Gbar$ where
  $\pi|_\Util$ is a homeomorphism onto $U$.  Moreover,
  $[\rhotil]_\Gbar$ is closed in $R_\Gtil(\Gamma)$ and $\pi$ restricts
  to a homeomorphism $[\rhotil]_\Gbar \to [\rho]_\Gbar$.
\end{lemma}

\begin{proof}
Consider the continuous translation number function
$\trans \maps \Gtil \to \R$ discussed in Section 3.1 of
\cite{CullerDunfield2018}, which is a conjugacy invariant. Let
$s_1, \ldots, s_d$ in $\Gamma$ generate
$H_1(\Gamma; \Z)/(\mathrm{torsion})$.  Define
$T \maps R_\Gtil(\Gamma) \to \R^d$ by
$T(\rho)_i = \trans\big(\rho(s_i)\big)$, which is continuous and
constant on $\Gbar$ orbits.  Let $\alpha \in H^1(\Gamma; \Z)$ act on
$x \in \R^d$ by $(\alpha \cdot x)_i = x_i + 2 \alpha(s_i)$.  If $c$ is
the preferred generator of $Z(\Gtil)$, one has
$\trans(\gtil c^k) = \trans(\gtil) + k$; from this it follows that $T$
is $H^1(\Gamma; \Z)$ equivariant (recall here that the $\Z$ in
(\ref{eq: tSLR}) is $\pairm{}{c^2}$ not $\pairm{}{c}$).

For now, fix any $\rhotil$ in $\pi^{-1}(\rho)$.  Let $\Vtil$ be
the open box in $\R^d$ about $T(\rhotil)$ that is the product of
intervals $\left(t_i - 1/3, t_i + 1/3\right)$ where
$t_i = T(\rhotil)_i$.  Let $\Ctil$ be the closed box which is the
product of intervals $\left[t_i - 1/6, t_i + 1/6\right]$. Consider
$\Util = T^{-1}(\Vtil)$ and $\Dtil = T^{-1}(\Ctil)$ which are open and
closed respectively, with $\rhotil \in \Dtil \subset \Util$;
moreover, both $\Dtil$ and $\Util$ are invariant under the $\Gbar$ action.
As $T(\rhotil) = T(\Gbar \cdot \rhotil) =
T\big([\rhotil]_\Gbar\big)$, we have that $[\rhotil]_\Gbar$ and
its closure are contained in $\Dtil$.

For all nonzero $\alpha \in H^1(\Gamma; \Z)$, we have
$\Vtil \cap \alpha \cdot \Vtil$ is empty since $\alpha$ moves some
coordinate of $\R^d$ by at least two.  As $T$ is equivariant, it
follows that $\Util \cap \alpha \cdot \Util$ is empty. This implies
$\pi|_\Util$ is injective.  As $\Util$ is open and $\pi$ is an open
map, so \(\pi|_\Util\) is a homeomorphism onto the open set
$U = \pi(\Util)$.  To show that $\pi(\Dtil)$ is closed in
$R_G(\Gamma)$, we note that its complement is the image under the open
map $\pi$ of the open set
$T^{-1}\big(\R^d \setminus \coprod \setdefm{\big}{\alpha \cdot
  \Ctil}{\alpha \in H^1(\Gamma; \Z)}\big)$.

As $\pi$ is $\Gbar$ equivariant, the subsets $\pi(\Dtil)$ and $U$ are
also invariant under the $\Gbar$ action.  So $\Gbar \cdot \rho$ and
indeed its closure is contained in $\pi(\Dtil)$.  If $\rho_1$ is any
representation where the closure of $\Gbar \cdot \rho_1$ meets the
closure of $\Gbar \cdot \rho$, then $\Gbar \cdot \rho_1$ must meet $U$
and hence be contained in it.  Thus $[\rho]_\Gbar$ is contained in
$U$, and of course $\pi\big([\rhotil]_\Gbar \big)$ is contained in
$[\rho]_\Gbar$; in fact, since $\pi|_\Util$ is a $\Gbar$ equivariant
homeomorphism, they are equal.  So $[\rhotil]_\Gbar$ is closed in
$\Util$ since $[\rho]_\Gbar$ is closed in $R_G(\Gamma)_0$ by
Lemma~\ref{lem: polystable}.  Since it is also a subset of $\Dtil$
which is closed in $R_\Gtil(\Gamma)$, it follows that
$[\rhotil]_\Gbar$ is closed in $R_\Gtil(\Gamma)$ as needed.

This establishes that $U$ has the properties claimed in the statement
of the lemma for the particular lift $\rhotil$ that we fixed near the
start.  The $H^1(\Gamma; \Z)$ action on $R_\Gtil(\Gamma)$ shows that
the same $U$ also works for any other lift $\rhotil'$, completing the
proof.
\end{proof}

It follows from Lemma~\ref{lem: orbits upstairs} that \(\pi\) maps the
closure of \(\Gbar \cdot \rhotil\) to the closure of
\(\Gbar \cdot \rho\) homeomorphically and G-equivariantly. Since the
action of \(\Gbar\) on \(R_G(\Gamma)\) is polystable, the action of
\(\Gbar\) on \(R_\Gtil(\Gamma)\) is polystable as well. Hence we can
define $\XA_\Gtil(\Gamma) = R_\Gtil(\Gamma) \sslash \Gbar.$

\begin{corollary}
  \label{cor: X of Gtil}
  The fibers of $R_\Gtil(\Gamma) \to \XA_\Gtil(\Gamma)$ are closed and
  the map $\pi \maps R_\Gtil(\Gamma) \to R_G(\Gamma)$ induces a
  continuous map $\pibar \maps \XA_\Gtil(\Gamma) \to \XA_G(\Gamma)$
  that makes the obvious diagram commute.  Similarly, for each
  $\gamma \in \Gamma$ the function,
  $\trans_\gamma \maps R_\Gtil(\Gamma) \to \R$ given by
  $\trans_\gamma(\rho) := \trans\big(\rho(\gamma)\big)$ descends to a
  continuous map $\trans_\gamma \maps \XA_\Gtil(\Gamma) \to \R$.
\end{corollary}

\begin{proof}
The first claim is immediate from Lemma~\ref{lem: orbits upstairs},
since each $[\rhotil]_\Gbar$ is closed.  For the second claim, we need
to check that each fiber of $R_\Gtil(\Gamma) \to \XA_\Gtil(\Gamma)$
has constant image under the composition
$R_\Gtil(\Gamma) \to R_G(\Gamma) \to \XA_G(\Gamma)$. This also follows
from Lemma~\ref{lem: orbits upstairs}, since
\(\pi([\rhotil]_\Gbar) = [\rho]_\Gbar\) and
$\XA_G(\Gamma) = R_G(\Gamma) \sslash \Gbar$.  Finally, the last claim
follows since $\trans_\gamma$ is constant on each fiber
$[\rhotil]_\Gbar$.
\end{proof}

Consider the maximal compact subgroup $\Kbar = \mathrm{PSO}_2$ in
$\Gbar$.  From \cite{RichardsonSlodowy1990}, there is a closed subset
$\cM$ of $R_G(\Gamma)$ such that:
\begin{enumerate}
\item
  \label{item: closed orbits}
  The closed $\Gbar$ orbits in $R_G(\Gamma)$ are exactly those
  that meet $\cM$.

\item
  \label{item: orbits in M}
  For $\rho \in \cM$, the intersection of $\Gbar \cdot \rho$ with
  $\cM$ is $\Kbar \cdot \rho$.

\item
  \label{item: def retract}
  There is a continuous $\Kbar$ equivariant deformation retraction
  $\phi \maps R_G(\Gamma) \times [0, 1] \to R_G(\Gamma)$ of
  $R_G(\Gamma)$ onto $\cM$ that is \emph{along the orbits of} $\Gbar$,
  that is, for all $\rho \in R_G(\Gamma)$ the set
  $\setdef{\phi(\rho, t)}{t \in [0, 1)} \subset \Gbar \cdot \rho$ and
  $\phi(\rho, 1)$ is in the closure of $\Gbar \cdot \rho$.

\item \label{item: M mod K} The map
  $\cM/\Kbar \to R_G(\Gamma) \sslash \Gbar$ induced by the inclusion
  $\cM \to R_G(\Gamma)$ is a homeomorphism. (That this map
  is a continuous bijection follows from (\ref{item: closed orbits}) and
  (\ref{item: orbits in M}).)
\end{enumerate}
We call a closed $\cM$ satisfying (\ref{item: closed
  orbits})--(\ref{item: M mod K}) a \emph{Kempf-Ness} subset
\cite{KempfNess1979}. Their usefulness lies in that (\ref{item: M mod
  K}) allows us to understand the topology of
$R_G(\Gamma) \sslash \Gbar$ in terms of an ordinary topological
quotient by a \emph{compact} group.  To prove Theorem~\ref{thm:
  tSLRchar}, we will need:

\begin{theorem}
  \label{thm: Kempf-Ness}
  Let $\cM_0 = \cM \cap R_G(\Gamma)_0$ and set $\cMtil =
  \pi^{-1}(\cM_0)$ inside $R_\Gtil(\Gamma)$.  Then $\cMtil$ is a
  Kempf-Ness subset for the action of $\Gbar$ on $R_\Gtil(\Gamma)$.
\end{theorem}

\begin{proof}
We check the properties in turn, with $\rho$ denoting a
representation in $R_G(\Gamma)_0$ and $\rhotil$ one of its preimages
in $R_\Gtil(\Gamma)$.  Throughout, we use Lemma~\ref{lem: orbits
  upstairs} without further comment.

For (\ref{item: closed orbits}), to start we know $[\rho]_\Gbar$
and $[\rhotil]_\Gbar$ are closed in $R_G(\Gamma)$ and
$R_\Gtil(\Gamma)$ respectively, and that the restriction
$\pi \maps [\rhotil]_\Gbar \to [\rho]_\Gbar$ is a homeomorphism.  Thus
$\Gbar \cdot \rhotil$ is closed if and only if $\Gbar \cdot \rho$ is
closed, and of course $\Gbar \cdot \rhotil$ meets $\cMtil$ if and only
if $\Gbar \cdot \rho$ meets $\cM_0$.  Hence property (\ref{item:
  closed orbits}) for $R_G(\Gamma)$ gives the same for
$R_\Gtil(\Gamma)$.  Property (\ref{item: orbits in M}) follows since
$\pi$ restricts to a $\Gbar$ equivariant bijection from
$\Gbar \cdot \rhotil$ to $\Gbar \cdot \rho$.

For (\ref{item: def retract}), recall $R_G(\Gamma)$ is a real
semialgebraic set and hence locally contractible; in particular, every
connected component is path connected.  Also, $R_G(\Gamma)_0$ is a
union of connected components of $R_G(\Gamma)$.  If $\phi$ is the
deformation retraction for $R_G(\Gamma)$, consider
\[
  \begin{tikzcd}
    R_\Gtil(\Gamma) \times [0, 1] \arrow[dashed]{r}{\phitil}
    \arrow{d}[left]{\pi \times \id}
    & R_\Gtil(\Gamma) \arrow{d}{\pi} \\
    R_G(\Gamma)_0 \times [0, 1] \arrow{r}{\phi} & R_G(\Gamma)_0
  \end{tikzcd}
\]
Here the vertical maps are covering maps and $\phi$ is a homotopy
equivalence.  Now $R_\Gtil(\Gamma) \times [0, 1]$ is locally path
connected since $R_G(\Gamma) \times [0, 1]$ is, so basic covering
space theory gives us a lift $\phitil$ of
$\phi \circ (\pi \times \id)$ whose restriction $\phitil_0$ is the
identity map on $R_\Gtil(\Gamma)$.  Since $\pi$ is bijective on each
$\Gbar$ orbit, the claims that $\phitil$ is $\Kbar$ equivariant and is
along the orbits of $\Gbar$ follow from the corresponding properties
of $\phi$.

To see that $\phi$ is a deformation retract onto $\cMtil$, suppose
$\rhotil \in \cMtil$. Note the path $\phitil(\rhotil, t)$ for
$t \in [0, 1]$ is a lift of $\phi(\rho, t)$ which is the constant path
at $\rho$.  Since $\pi^{-1}(\rho)$ is discrete, it follows that
$\phitil(\rhotil, t)$ must be the constant path at $\rhotil$.  In
particular, each $\phi_t$ restricts to the identity on $\cMtil$.  As
$\phitil_1\big(R_\Gtil(\Gamma)\big)$ is contained in
$\pi^{-1}\big(\im(\phi_1)\big) = \cMtil$, we have shown $\phitil$
satisfies (\ref{item: def retract}).

For (\ref{item: M mod K}), we use the argument from Section 9.6 of
\cite{RichardsonSlodowy1990}.  By (\ref{item: closed orbits}) and
(\ref{item: orbits in M}), the inclusion $\cMtil \to R_\Gtil(\Gamma)$
induces a continuous bijection
$\cMtil/\Kbar \to R_\Gtil(\Gamma) \sslash \Gbar$.  If $\pi_\Kbar$ is the
quotient map $\cMtil \to \cMtil/\Kbar$, we claim the composite
$\pi_\Kbar \circ \phitil_1 \maps R_\Gtil(\Gamma) \to \cMtil/\Kbar$ is
constant on each fiber of $\pi$, i.e.~on each $[\rhotil]_\Gbar$.  To
see this, note that since $\phitil$ is along $\Gbar$ orbits, 
$\phitil_1(\rhotil) \in [\rhotil]_\Gbar$, so if $\rhotil \in
\cMtil$ we have $\phi_1^{-1}(\Kbar \cdot \rhotil) = [\rhotil]_\Gbar$.  In
particular, $\pi_\Kbar \circ \phitil_1$ descends to the map
shown at right:
\[
 \cMtil/\Kbar \to R_\Gtil(\Gamma) \sslash \Gbar \to \cMtil/\Kbar.
\]
The definitions give that the composition above is the identity
on $\cMtil/\Kbar$.  Both maps are continuous, and the leftmost one is a
bijection, so  both are homeomorphisms. This completes the
proof of the theorem.
\end{proof}

\begin{proof}[Proof of Theorem~\ref{thm: tSLRchar}]
First, we show that $\XA_G(\Gamma)_0$ is a union of connected
components.  By Lemma~\ref{lem: polystable}, the fibers of
$\tau \maps R_G(\Gamma) \to \XA_G(\Gamma)$ are connected. Since
$R_G(\Gamma)_0$ is a union of components, this implies that
$\tau^{-1}\big(\XA_G(\Gamma)_0\big) = R_G(\Gamma)_0$, so like
$R_G(\Gamma)_0$ the subset $\XA_G(\Gamma)_0$ is both open and closed,
i.e.~a union of components.

To show $\XA_\Gtil(\Gamma) \to \XA_G(\Gamma)_0$ is a covering map,
consider Kempf-Ness subsets $\cM$ and $\cMtil$ as in Theorem~\ref{thm:
  Kempf-Ness} and set $\cM_0 = \cM \cap R_G(\Gamma)_0$.  We have the
commutative diagram:
\[
  \begin{tikzcd}
    \cMtil/\Kbar \arrow{r}{} \arrow{d}{} & \XA_\Gtil(\Gamma) \arrow{d}{} \\
    \cM_0/\Kbar \arrow{r}{} &
    \XA_G(\Gamma)_0
  \end{tikzcd}
\]
where the horizontal maps are homeomorphisms and vertical ones are
induced by $\pi \maps R_\Gtil(\Gamma) \to R_G(\Gamma)_0$.  Moreover,
the top homeomorphism is equivariant with respect to the
$H^1(\Gamma; \Z)$ action.  Now $\pi \maps \cMtil \to \cMtil_0$ is a
$\Kbar$ equivariant covering map with covering group
$H^1(\Gamma; \Z)$.  As $\Kbar$ is compact, this implies that
$\cMtil/\Kbar \to \cM_0/\Kbar$ is also a $H^1(\Gamma; \Z)$ cover,
which proves the theorem since the diagram commutes.
\end{proof}

\section{Geometric point configurations}
\label{sec:config}

Having dispensed with the background, we move to the main body of the
paper. Our first task is to construct a resolution of the character
variety $X^c(F_n, S)$, in the sense of Section~\ref{sec: resolve
  intro}. We do this in Sections~\ref{sec:config} to \ref{sec:map to
  X}.  The goal of this section is to state Theorem~\ref{thm:consY},
which describes a smooth manifold $\cY$ built out of configurations of
$n > 1$ points in all three \2-dimensional geometries: $\E^2$, $S^2$,
and $\H^2$. We prove this theorem in Section~\ref{sec:join
  configs}.  In Section~\ref{sec:map to X}, we show that $\cY$ is
the desired resolution of $X^c(F_n, S)$.

\subsection{Conventions and models} \label{sec:models}

Throughout this section, the integer $n > 1$.  We fix for the whole
paper the following preferred models for the three \2-dimensional
geometries $\E^2$, $S^2$, and $\H^2$.  First, we define the Euclidean
plane $\E^2$ as $\R^2$ with the Riemannian metric associated to the
standard quadratic form $x^2 + y^2$ oriented so that $e_1, e_2$ is a
postive basis.  Second, the sphere $S^2$ is the points in $\R^3$ where
$x^2 + y^2 + z^2 = 1$ with the Riemannian metric determined by
restricting the quadratic form $x^2 + y^2 + z^2$ to the tangent bundle
$TS^2$; we orient $S^2$ so that $e_1, e_2$ is a positive basis for
$T_{e_3}S^2$.  Finally, the hyperbolic plane $\H^2$ is the upper sheet
of the hyperboloid $x^2 + y^2 - z^2 = -1$ with the Riemannian metric
coming from $x^2 + y^2 - z^2$, where $e_1, e_2$ is a positive basis
for $T_{e_3}\H^2$.

\subsection{The Euclidean case}

Consider the action of the orientation preserving isometry group
$\IsoE$ on ordered $n$-tuples of points in $\E^2$, that is, the
diagonal action of $\IsoE$ on $\left(\E^2\right)^n$.  This action is
free except when all the points are the same, in which case the
stabilizer is $S^1$.  Restricting to
$\cE_n = \left(\E^2\right)^n \setminus \Delta$ where
$\Delta = \setdef{(p, p, \ldots, p)}{p \in \E^2}$, we have a free
action of $\IsoE$. As this action is also proper, the quotient
$\cE_n/\IsoE$ is a smooth manifold.  The manifold $\cE_n/\IsoE$ is
noncompact with two ends: one corresponding to all the points
coalescing and the other to the diameter of the set of points going to
infinity.

Let $\SimE$ denote the group of orientation preserving similarities of $\E^2$,
which is generated by $\IsoE$ and dilations about any given point.
While the action of $\SimE$ on $\E^2$ is not proper (point stabilizers
are noncompact), its action on $\cE_n$ is proper, roughly because the
amount of dilation is detected in the change in the distance between
any pair of distinct points. So we get a manifold
\[
  \CnE = \quo{\cE_n}{\SimE}
\]
whose topology is as follows.
\begin{lemma}
  \label{lem:CnE is Pn-2}
  The manifold $\CnE$ is diffeomorphic to $P^{n-2}(\C)$.
\end{lemma}
\begin{proof}
We will typically denote a point in $\cE_n$ by
$v = (v_1,v_2, \ldots,v_n)$ where each $v_i \in \E^2$ is regarded as a
vector in $\R^2$. Using translations and dilations, we see that the
subset
\[
\cA = \setdef{v \in \cE_n}{\mbox{$v_1 = 0$ and $\sum_{i = 2}^n \abs{v_i}^2 = 1$}}
\]
meets every $\SimE$ orbit, and indeed that $\CnE$ is the quotient of
$\cA$ under the $S^1$ action of rotations about $0$. Now $\cA$ is just
$S^{2n - 3} \subset (\R^2)^{n - 1} = \C^{n-1}$, and so we see that the
quotient under the circle action is exactly $P^{n-2}(\C)$ as claimed. 
\end{proof}

\subsection{Signed points} We will also need to work with $n$
ordered points where each one has an associated sign in $\{\pm 1\}$.
The set of such configurations is $\cE_n^\pm = \cE_n \times \pmonen$
where if $(v, \epsilon) \in \cE_n^\pm$ then $\epsilon_i$ is the sign
associated to the point $v_i$.  We define an action of the full group
of similarities $\SimEall$ on $\cE_n^\pm$ so that an orientation
preserving element leaves the signs unchanged, but an orientation
reversing element multiplies them all by $-1$.  The action of
$\SimEall$ on $\cE_n^\pm$ is again free and proper, and we set
\[
  \CnsignE = \quo{\cE^\pm_n}{\SimEall}
\]
Note that the signs are important to make this action free: any $v$ in
$\cE$ where all $v_i$ lie on a common line $L$ is fixed by reflection
in $L$.  

There is an action of $\pmonen$ on $\CnsignE$ by changing the signs of
the points; more formally, the action of $\pmonen$ on $\cE_n^\pm$ by
multiplication on the second factor commutes with action of $\SimEall$,
and so descends to an action on $\CnsignE$.  Using this action, you
can show that $\CnsignE$ is the disjoint union of $2^{n-1}$ copies of
$\CnE$; a subtlety here is that while $(-1, \ldots, -1)$ does not
permute the connected components of $\CnsignE$, it does act
nontrivially on each of them when $n \geq 3$.

\subsection{The hyperbolic case}  
\label{Subsec:hyperbolic}

For the hyperbolic plane $\H^2$, our focus will be on 
\[
  \CnH = \quo{\cH_n}{\IsoH} \mtext{where} \cH_n = (\H^2)^n
  \setminus \Delta
\]
which, as in the Euclidean case, is a smooth manifold with two ends. We
will show below that $\CnH$ is diffeomorphic to
$\CnE \times \R \cong P^{n-2}(\C) \times \R$, but for now just note
that $\dim \CnH = 2 n - 3$.  There
is also a signed variant
\[
  \CnsignH= \quo{\cH_n \times \{\pm 1\}^n}{\Isom(\H^2)}
\] 
defined as above, which again has an action of $\pmonen$ and consists
of $2^{n-1}$ copies of $\CnH$.

\subsection{The spherical case}
\label{sec:spherical}

Turning now to the round sphere $S^2$, we must remove a little more to
get a free action since a pair of antipodal points has nontrivial
stabilizer.  Specifically, we have a free action of $\IsoS$ on  
\[
  \cS_n = \setdef{v \in \left(S^2\right)^n}{\mbox{$v_i \notin \{-v_j,
      v_j\} $ for at least one $i,j$}}
\]
which is automatically a proper action since $\IsoS \cong \SO_3$ is
compact.  Hence
\[
  \CnS = \quo{\cS_n}{\IsoS}
\]
is a smooth manifold of dimension $2n-3$. This time, it has $2^{n-1}$
ends corresponding to the different ways the points can coalesce to a
pair of antipodal points; concretely, each $v_i$ for $i > 1$ can
approach either $v_1$ or $-v_1$. We will show in
Theorem~\ref{thm:localmodel} that each end is diffeomorphic to
$\CnE \times \R$.

The action of $\pmonen$ on $\cS_n$ by coordinatewise multiplication
descends to an action on $\CnS$.  This action is faithful but not
free; for example, when $n = 2$ the entire group fixes $[(e_1,
e_2)]$. We will not use this directly, but the quotient orbifold
$\CnS/\pmonen$ is $\CnP$, that is, configurations of $n$ points in
$P^2(\R)$, not all the same, modulo
$\Isom\left(P^2(\R)\right) \cong \PO_3 \cong \SO_3$.

\subsection{Orientations}\label{sec:configorient}

We will eventually use these configuration spaces to construct an
invariant via an algebraic intersection number, so we next give them
preferred orientations. These will be derived from the orientations of
$\E^2$, $\H^2$, and $S^2$ specified in Section~\ref{sec:models}, and
we use freely the orientation conventions of
Section~\ref{sec:orientshortex}.  To start, we need to orient
$\SimEall$, $\IsoHall$, and $\IsoSall$.  We orient the tangent bundle
of an oriented surface $X$ so that if $(x_1, x_2)$ are oriented
coordinates on $X$ then
$\big\langle x_1, x_2, \frac{\partial}{\partial x_1},
\frac{\partial}{\partial x_2}\big\rangle$ are oriented coordinates on
$TX$.  We orient the unit tangent bundle $UTX$ by viewing it as the
boundary of the submanifold $\setdef{v \in TX}{\abs{v} \leq 1}$;
equivalently, if $\pair{v_1, v_2}$ is an oriented orthonormal basis of
$T_p X$, then we orient $T_{(p, v_1)}(UTX) = T_pX \oplus T_{v_1} S^1$
by $\pair{v_1, v_2, v_3}$ where
$v_3 = \frac{d}{dt}\big\vert_{t=0} \cos(t) v_1 + \sin(t) v_2 $.

If we fix a tangent vector $v_0 \in T\E^2$, then the map
$\SimE \to T\E^2$ given by $g \mapsto g \cdot v_0$ is a
diffeomorphism, and we orient $\SimE$ so that this map is orientation
preserving.  We then orient all of $\SimEall$ by insisting that the
orientation be left-invariant.  We orient $\IsoHall$ by the
corresponding identification of $\IsoH$ with $UT\H^2$ and orient
$\IsoSall$ analogously.

For the unsigned configuration spaces, we first give $\cE_n$ the
product orientation, and then orient $\CnE$ by requiring that the fiber bundle
$\SimE \to \cE_n \to \CnE$ be compatibly oriented.  Analogously, we
orient $\CnH$ and $\CnS$ via $\IsoH \to \cH_n \to \CnH$ and
$\IsoS \to \cS_n \to \CnS$ respectively.

The signed case is handled as follows, focusing on $\E^2$ for ease of
notation.  Orient $\E^2 \times \pmone$ so that $\E^2 \times \{+1\}$
has the preferred orientation of $\E^2$ and $\E^2 \times \{-1\}$ has
its reverse.  Consider the homomorphism
$\sign \maps \SimEall \to \pmone$ that is $1$ on orientation
preserving similarities and $-1$ on orientation reversing ones.
We
make $g \in \SimEall$ act on $\E^2 \times \{\pm 1\}$ by
$g \cdot (v, s) = (g(v), \sign(g) s)$ and note that this action is
always orientation preserving.  Now orient $\cE_n^\pm = \cE_n \times \pmonen$ via
the product orientation on $\left(\E^2 \times \pmone \right)^n$.  The
action of $\SimEall$ on $\cE_n^\pm$ preserves this
orientation, and the orientation induced on $\CnsignE$ from $\SimEall
\to \cE_n^\pm \to \CnsignE $ will be our
preferred one there.

Each of $\CnsignH$, $\CnsignE$, and $\CnS$ is acted on by $\pmonen$,
and it will be important to understand how these actions interact with
our preferred orientations.  In all cases, an $\epsilon \in \pmonen$
preserves orientation if the product of the $\epsilon_i$ is $+1$ and
reverses it otherwise.

\subsection{Putting the geometries together} 
\label{subsec:join geometries}

We now sketch a natural way of combining these configuration spaces
for $\H^2$, $\E^2$, and $S^2$ together into a single smooth manifold.
To begin, we consider the end of $\CnH$ where all the points coalesce.
Given an element in $\CnH$, we can rescale the metric on $\H^2$ so
that the diameter of the set of points is exactly 1 in the new
metric. The closer the original points are in $\H^2$, the flatter the
rescaled metric is. When the original points are very close, the new
space is nearly isometric to $\E^2$ on the scale of the points in the
final configuration.  This makes it plausible that we should
compactify this end of $\CnH$ by adding a copy of $\CnE$ at infinity
to produce a manifold with boundary.  Looking now at the particular
end of $\CnS$ where all the points come together, the exact same story
applies to suggest that we should also compactify this end by adding a
copy of $\CnE$. We could then glue our two compactifications together
to get a nice manifold structure on $\CnH \cup \CnE \cup \CnS$.

As discussed, the manifold $\CnS$ has $2^{n - 1}$ ends, and we can
bring the others into our picture as follows.  Regard an element of
$\CnS$ as $n$ pairs of antipodal points where one point in each pair
is labeled $+1$ and the other $-1$; the action of $\pmonen$ on $\CnS$
now simply swaps the labels as appropriate. An end of $\CnS$ can now
be labeled by an $\epsilon$ in $\pmonen$ corresponding to the pattern
of signs on the points that are coalescing; in fact, each end has two
such labels, $\epsilon$ and $-\epsilon$, as there are two clusters of
points we can focus on.  The points in these two clusters have
opposite signs and their positions differ by the antipodal map, which
is orientation reversing.  Thus to compactify all the ends of $\CnS$,
we should consider \emph{signed} points in $\E^2$, up to the
equivalence used to define $\CnsignE$.  Hence we should be able to
make $\CnsignE \cup \CnS$ into a compact manifold with boundary to
which we can attach a copy of $\CnsignH$.  

That this whole sketch can be made precise is the content of:
\begin{theorem}\label{thm:consY}
  There is a smooth structure on $\CnsignH \cup \CnsignE \cup \CnS$
  making it into an oriented manifold $\cY$ of dimension $2n - 3$
  that is compatible with the previously defined smooth structures on
  $\CnsignH$, $\CnsignE$, and $\CnS$; the orientation of $\cY$ agrees
  with that of$\CnS$ but is opposite of that on $\CnsignH$. The action
  of $\pmonen$ on the set $\cY$ is in fact an action by
  diffeomorphisms.  The subset $\CnsignH \cup \CnsignE$ is a closed
  submanifold with boundary equal to $\CnsignE$, and it is
  diffeomorphic to $2^{n-1}$ copies of
  $P^{n-2}(\C) \times [0, \infty)$.
\end{theorem}
One consequence of Theorem~\ref{thm:consY} is that $\cY$ is 
diffeomorphic to $\CnS$ itself.  We will prove Theorem~\ref{thm:consY}
over the next two sections.

\section{Groups and quadratic forms}
\label{sec:groups and forms}

Our framework for building the smooth structure on $\cY$ will be the
geometric transition $\H^2 \to \E^2 \to S^2$, studied via projective
geometry from the perspective of
\cite[\S 2.1]{CooperDancigerWienhard2014}.  To describe this
transition, we use a family of quadratic forms and their automorphism
groups.

\subsection{The quadratic forms \(B_t\)}
\label{subsec:B_t}
On $\R^3$ with coordinates $(x,y,z)$, consider the family of quadratic
forms:
\[
B_t = t (x^2 + y^2) + z^2 \mtext{where $t \in \R$.}
\]
We will use $B_t(\cdotspaced, \cdotspaced)$ for the associated
bilinear forms.  Let $\Qbar_t$ be the quadric surface where $B_t = 1$,
which is an ellipsoid when $t > 0$, two planes when $t = 0$, and a
hyperboloid of two sheets when $t < 0$. We define \(Q_t\) to be the
subset of \(\Qbar_t\) where \(z>0\).

The group $\R_{>0}$ acts on $\R^3$ by dilating the $x$ and $y$
coordinates:
\begin{equation}\label{eq:dilate}
s\cdot (x, y, z) = (s x, s y, z).
\end{equation}
Note this action satisfies $B_t(s \cdot v) = B_{s^2 t}(v)$ for
$v \in \R^3$ and hence $v \mapsto s \cdot v$ gives an isometry
$(\R^3, B_t) \to (\R^3, B_{t/s^2})$. Thus
$s \cdot \Qbar_t = \Qbar_{t/s^2}$.  Consequently, each
$(\Qbar_t, B_t)$ with $t > 0$ is isometric to
$S^2 = (B_{1}, \Qbar_{1})$ and each $(\Qbar_t, B_t)$ with $t < 0$ is
isometric to $(\Qbar_{-1}, B_{-1})$, which is two copies of $\H^2$ with
the negative of its metric from Section~\ref{sec:models}. We
give the two planes making up $\Qbar_{0}$ the Euclidean metric
associated to the quadratic form $x^2 + y^2$. We will always orient
$\Qbar_t$ by the convention that  $u_1, u_2$ is a positive basis
for $T_v \Qbar_t$ if and only if $v, u_1, u_2$ is a positive basis for
$\R^3$; in particular, $e_1, e_2$ is a positive basis for
$T_{e_3} \Qbar_t$ and  $e_1, -e_2$ is a positive basis for
$T_{-e_3} \Qbar_t$.

\subsection{The groups \(G_t\)}
\label{subsec:Gt}
For $t \neq 0$, let $G_t = \SO(B_t) $ be the subgroup of $\SL{3}{\R}$
that preserves $B_t$.  For $t=0$, we define \(G_0\) to be the subgroup
of $\SL{3}{\R}$ that preserves $B_0$ and acts on $\Qbar_{0}$ by an
isometry; concretely,
\begin{equation}\label{eq:G0}
  G_{0} = \setdef{
  \left(\begin{array}{cc} 
          A & b \\
          0  & \det A
        \end{array}\right)\in \SL{3}{\R}}{ 
      A \in \O(2)}.
\end{equation}
From our definition of the orientation on $\Qbar_t$, since
$G_t \leq \SL{3}{\R}$, we see that $G_t$ acts on $\Qbar_t$ preserving
orientation.  For \(t>0\), the group $G_t \cong G_1 \cong \IsoS$ is
connected. However, for \(t\leq 0\), the group \(G_t\) has two
connected components, corresponding to whether an element preserves or
interchanges the two components of \(\Qbar_t\).

If $t < 0$, the identity component of \(G_t\) is
\[
  \Isom^+(Q_t) \cong \Isom^+(Q_{-1}) = \Isom^+(\H^2).
\]
We can regard the full isometry group $\Isom(Q_t)$ as the subgroup of
the full orthogonal group $\O(B_t)$ that preserves $Q_t$. There is an
isomorphism \(\Isom(Q_t) \to G_t\) given by \(g \mapsto (\det g) g\),
so $G_t \cong G_{-1}\cong \Isom(\H^2)$. Note that the action of
\(G_{-1}\) on \(\Qbar_{-1}\) corresponds to the action of
$\Isom(\H^2)$ on \(\H^2 \times \{\pm 1\}\) from
Section~\ref{Subsec:hyperbolic}, with $Q_t$ playing the role of
$\H^2 \times \{1\}$.  A similar argument shows that
\(G_0 \cong \Isom(\E^2) \).

When \(t=0\), the groups \(\R_{>0}\) and \(G_0\) both act on
\(\Qbar_0\); their actions can be combined to give the action of a
larger group \(\Ghat_0 \cong \Sim(\E^2)\). More precisely, define
\(\Ghat_0\) to be the subgroup of \(\SL{3}{\R}\) generated by \(G_0\)
and matrices of the form \( \mysmallmatrix{sI}{0}{0}{1}\) for $s >
0$. The action of \(\Ghat_0\) preserves \(\Qbar_0\). Restriction to
\(\Qbar_0\) defines a homomorphism from \(\Ghat_0\) to
\( \Sim(\E^2)\), which you can check is an isomorphism.

\subsection{The groups \(U_t\)}\label{sec:SL2Csubs}

The groups $G_t$ are closely related to certain subgroups of $\SLC$
which we now describe.  For \(t\in \R\), let
\[
  U_t = \setdef{\ \twobytwomatrix{a}{b}{-t\bbar}{\abar}}{
     \mbox{$a, b \in \C$ with $\abs{a}^2 +t\abs{b}^2 = 1$}}.
\]
For \(t\neq 0\), the group $U_t$ is the subgroup of $\SLC$ that
preserves the Hermitian form given by
$J_t = \mysmallmatrix{t}{0}{0}{1}$, as can be seen by solving
$A^* J_t A = J_t$ in the form $A^* J_t = J_t A^{-1}$.  In particular,
$U_1 = \SU$ and $U_{-1} = \SUoneone \cong \SLR$. More generally, if
$T_s = \mysmallmatrix{s}{0}{0}{1}$, you can easily check that
$T_s U_t T_s^{-1} = U_{t/s^2}$, so $U_t \cong U_1 = \SU$ for $t > 0$
and $U_t \cong U_{-1} = \SUoneone$ for $t < 0$.

Let \(U_t'\) be the normalizer of \(U_t\) in \(\SLC\).  When
$t \neq 0$, we define
$\gamma_t = \mysmallmatrix{0}{t^{-1/2}}{-t^{1/2}}{0}$, which is in
$U_t'$ with conjugation by $\gamma_t$ inducing the automorphism of
$U_t$ that is complex conjugation of the matrix entries. For \(t>0\),
we have \(U_t' = U_t\), and, for \(t=0\), the group \(U_0'\) consists
of upper triangular matrices with determinant \(1\). For \(t<0\), the
subgroup \(U_t'\) is conjugate to the subgroup \(\SLRpm\) from
Section~\ref{sec:realchar}; it is generated by \(U_t\) and
\(\gamma_t\).

The Lie algebra $\ut$ of $U_t$ is easy to compute: when $t \neq 0$,
the relation $A^* J_t = J_t A^{-1}$ turns into 
$X^* J_t = -J_t X$; combining this with the trace 0 condition defining
$\mathfrak{sl}_2\C$ gives
\begin{equation}\label{eq:Ut}
  \ut = \setdef{\ i \twobytwomatrix{z}{-x-iy}{-t(x-iy)}{-z}}{
     \mbox{$x,y,z \in \R$}}. 
\end{equation}
It is easy to check that this formula for $\ut$ is also valid when
$t = 0$.  Let \(\varphi_t\maps \R^3 \to \ut\) be given by
$$ \varphi_t(x,y,z) = i \twobytwomatrix{z}{-x-iy}{-t(x-iy)}{-z}. $$
From now on,
we identify \(\ut\) with $\R^3$ via \(\varphi_t\), for example getting
a linear action of \(G_t\) on \(\ut\) by requiring
\(g \cdot \phi_t(v) = \phi_t(g \cdot v)\).  With this identification,
a routine calculation shows that the Killing form on \(\ut\) is a
multiple of \(B_t\), specifically $-8 B_t$.  As any Lie algebra
automorphism in $\Aut(\ut)$ preserves the Killing form, we see
$\Aut(\ut) \leq \O(B_t)$.

As $U'_t$ normalizes $U_t$, we get a homomorphism
$\psi_t \maps U_t' \to \Aut(\ut)$ by the adjoint action
\(\psi_t(A) = \Ad_A\).  A key relationship between the groups $G_t$
and $U_t$ is:
\begin{lemma}
 \label{Lem:Gt=Aut Ut}
 For \(t\neq 0 \), the action of $G_t$ on $\ut$ gives $G_t = \Aut(\ut)$ as
 subgroups of $\GL{}(\ut)$.  Moreover,
 $\Aut(\ut) = \Aut(U_t) = \psi_t(U'_t) \cong U_t'/{\pm I}$.
\end{lemma} 

\begin{proof}

First, note \(\ut\) carries a natural orientation determined by the
ordered basis \(v , w, [v, w]\), where \(v\) and \(w\) are any
linearly independent elements of \(\ut\).  Now $\Aut(\ut)$
preserves this orientation and the Killing form, giving
\(\Aut(\ut) \leq \SO(B_t) = G_t\).

Now \(U_t'\) acts on \(U_t\) by conjugation, giving
$\psitil_t \maps U_t' \to \Aut(U_t)$ whose kernel is \(\pm I\).  The
group \(U_t\) is connected, so the derivative at \(I\) gives an
injective homomorphism \(\Aut(U_t) \to \Aut(\ut)\).  Note that the
$\psi_t$ defined previously is the composition of $\psitil_t$ and
\(\Aut(U_t) \to \Aut(\ut)\), so the kernel of $\psi_t$ is also
\(\pm I\).

So we have $\psi_t(U_t) \cong U_t/\{\pm I\}$ is a connected Lie
subgroup of dimension 3 inside $G_t$, which forces $\psi_t(U_t)$ to be
the identity subgroup of $G_t$.  If $t > 0$, then $U_t' = U_t$ and
$G_t$ are all connected, so $\psi_t(U_t') = G_t$.  If $t < 0$, recall
that for \(u \in \ut\), we have \(\Ad_{\gamma_t}(u) = \ubar\),
implying 
\(\psi_t(\gamma_t) = \left(\begin{smallmatrix} -1 & 0 & 0 \\ 0 & 1 & 0
    \\ 0 & 0 & -1 \end{smallmatrix}\right)\) is in \(G_t\), and again
$\psi_t(U_t') = G_t$.  All the claims of the lemma now follow since
$\psi_t(U'_t) \leq \Aut(U_t) \leq \Aut(\ut) \leq G_t$ must all
actually be equalities.
\end{proof}

The analogue of Lemma~\ref{Lem:Gt=Aut Ut} when $t = 0$ is:

\begin{lemma}
  \label{Lem:G0=Aut U0}
  As subgroups of $\GL{}{(\ut[0])}$, we have $\Aut(\ut[0]) = \Ghat_0$.
  Moreover, $\Aut(\ut[0]) = \Aut(U_0)$.  Finally,
  $\psi_0(U'_0) \cong U_0'/{\pm I}$ and
  $\psi_0(U_0) \cong U_0/{\pm I}$ are the identity components of $\Ghat_0$
  and $G_0$ respectively.
\end{lemma}

\begin{proof}
Write \(\ut[0] = L \oplus \langle h \rangle\), where
\(L= \setdef{ \mysmallmatrix{0}{\alpha}{0}{0}}{\alpha \in \C}\), which
we view as a 1-dimensional complex vector space, and
\(h = \mysmallmatrix{i}{0}{0}{-i}\).  We have \([v, w]=0\) if
\(v, w \in L\), and \([h, v] = 2i v\) for \(v \in L\).  Observe that
\(g \in \Aut(\ut[0])\) must preserve \(L\), and so is represented by a
matrix of the form \(\mysmallmatrix{A}{*}{0}{a}\). In order for $g$ to
respect the bracket, we must have \([ah,A(v)] = A(2i v)\), or
equivalently \(iaA(v) = A(i v)\). As a self-map of $L$, we know $A$ is
\(\R\)-linear, so applying this relation twice we see that
\(a^2=1\). When \(a=1\), the map \(A\) is \(\C\)-linear; when
\(a=-1\), it is conjugate linear. In either case, $A = s \cdot \Abar$
where $s \in \R_{> 0}$ and $\Abar \in \O(2)$. Hence
\(\Aut(\ut[0]) \subset \Ghat_0\).

For the converse, observe the adjoint action gives a homomorphism
\(\psitil_0\maps U'_0 \to \Aut(U_0)\) whose kernel is $\pm I$. As
before, $\Aut(U_0) \to \Aut(\ut[0])$ is injective as $U_0$ is
connected, and so $\psi_0$ also has kernel $\pm I$.  Hence
$\psi_0(U'_0) = U_0'/\{\pm I\}$ is a connected Lie group of dimension
4 inside $\Ghat_0$.  As $\dim \Ghat_0 = 4$, this means
$\psi_0(U'_0) = U_0'/\{\pm I\}$ is the identify component of
$\Ghat_0$.  The elements of $\Aut(U_0)$ and $\Aut(\ut[0])$ induced by
entrywise complex conjugation correspond to
\(C = \left(\begin{smallmatrix} -1 & 0 & 0 \\ 0 & 1 & 0 \\ 0 & 0
    &-1 \end{smallmatrix}\right)\) in $\Ghat_0$.  As $\Ghat_0$ is
generated by its identity component and $C$, we have
$\Aut(U_0) = \Aut(\ut[0]) = \Ghat_0$ as desired.
\end{proof}

\subsection{Further properties of $G_t$ and $U_t$}
\label{sec: props Gt and Ut}

Combining Lemmas~\ref{Lem:Gt=Aut Ut} and \ref{Lem:G0=Aut U0}, we will
view $G_t$ as a subgroup of $\Aut(U_t)$, and hence $G_t$ acts on $U_t$.
Indeed, for $t \neq 0$, the homomorphism \(\psi_t: U_t' \to G_t\) is
surjective, so concretely the action of \(G_t\) on \(U_t\) is given by
\(\psi_t(A) \cdot B = ABA^{-1}\).  For all $t$, the element
\(C = \left(\begin{smallmatrix} -1 & 0 & 0 \\ 0 & 1 & 0 \\ 0 & 0
    &-1 \end{smallmatrix}\right)\) acts on $\ut$ as complex
conjugation, and $C = \psi_t(\gamma_t)$ for $t \neq 0$. The action also satisfies
\[
  \exp(g \cdot u) = g \cdot \exp u \ \text{and} \ \tr (g \cdot B) = \tr B,
\]
which may be easily checked for \(g \in \im \psi_t(U_t')\) or
\(g = C\), and so hold for all \(g \in G_t\).

Recall that  $T_s U_t T_s^{-1} = U_{t/s^2}$, where $T_s =
\mysmallmatrix{s}{0}{0}{1}$.   The adjoint action of
$T_s$ on $\mathfrak{sl}_2{\C}$ thus takes $\ut$ to
$\mathfrak{u}_{t/s^2}$; under our identification of the latter two
with $\R^3$, this is precisely the dilation action defined in
(\ref{eq:dilate}).  In particular, it carries $\Qbar_t \subset \ut$ to
$\Qbar_{t/s^2} \subset \mathfrak{u}_{t/s^2}$ and $B_t$ to
$B_{t/s^2}$.  Moreover, the following commutes for each $g \in U_t$
\[
\begin{tikzcd}[column sep=large]
  \Qbar_t \arrow{r}{\Ad_{T_s}} \arrow{d}[left]{\psi_t(g)} & 
  \Qbar_{t/s^2} \arrow{d}{\psi_{t/s^2}(T_s g T_s^{-1})} \\
  \Qbar_t  \arrow{r}{\Ad_{T_s}} &  \Qbar_{t/s^2}
\end{tikzcd}
\]
since both compositions send $v \in
\ut$ to $T_s g v g^{-1} T_s^{-1}$.

\begin{remark}\label{rem:actsame}
  We will always orient $\ut$ by taking $(x,y,z)$ from
  (\ref{eq:Ut}) to be oriented coordinates.  This is consistent with
  our earlier conventions in the following sense: if we orient $G_t$
  as isometries of $\Qbar_t$ using the convention of
  Section~\ref{sec:configorient}, then the homomorphism $U_t \to G_t$
  is orientation preserving.  We leave the detailed check of this to
  you, but note that as the orientations of $\Qbar_t$, $G_t$, and
  $U_t$ all vary continuously in $t$, it suffices to check this for
  $t = 0$.
\end{remark}

\subsection{A \(\C^\times\) action}
\label{sec:Cstaracts}
The action of \(\R_{>0}\) on \(\R^3\) extends to an action of
\(\C^\times\) as follows. Writing \(\R^3 = \C \times \R\), where
\((x,y,z)\) is identified with \((x+iy,z)\), we define
\(u\cdot (w,z) = (uw,z)\) for $u \in \C^\times$. An easy calculation
shows that
\begin{equation}
\label{eq:u action}
\varphi_{t/|u|^2}(u \cdot v) = T_u \varphi_t(v) T_{u}^{-1} \mtext{where $T_u = \mysmallmatrix{u}{0}{0}{1}.$}
\end{equation} 

The \(\C^\times\) action combines the actions of \(\R_{>0}\) and the
adjoint action of the subgroup of \(U_t\) consisting of diagonal
matrices, in the sense that if \(\iota\maps S^1 \to U_t\) is given by
\(\iota(u) = \mysmallmatrix{\sqrt{u}}{0}{0}{\sqrt{u}^{-1}}\), then
\( \varphi_t(u \cdot v) = \iota(u) \cdot \varphi_t(v).\)

\subsection{The exponential map}
 
Suppose \(A \in U_t\), and let \(c= \tr A\) . A quick look at the
definitions of \(U_t\) and \(\ut\) shows that we can write
\(A = \frac{c}{2}I + u\), where \(u \in \ut\). Moreover, if
\(g \in G_t\), we have \(g \cdot A = \frac{c}{2} I + g \cdot u\),
where \(G_t\) acts on \(U_t\) in the right-hand side of the equality
and acts on \(\ut\) in the left-hand side; this is easily checked
when \(g \in \im \psi_t\) and when \(g= C\), and so holds for all
elements of \(G_t\).

\begin{lemma} 
  \label{lem:exp formula}
  If \(u \in \Qbar_t \subset \ut \), then
  \(\exp (\alpha u) = (\cos \alpha) I + (\sin \alpha) u\).
  Geometrically, the element $\exp(\alpha u)$ acts on $\Qbar_t$ as the
  isometry in $G_t$ that fixes $u$ and rotates about it anticlockwise
  by angle $2 \alpha$ with respect to the orientation of $\Qbar_t$
  fixed in Section~\ref{subsec:B_t}.
\end{lemma}
\begin{proof}
We first prove the claims for $h = (0, 0, 1) \in \Qbar_t$, which
corresponds to $\mysmallmatrix{i}{0}{0}{-i} \in \ut$ under
$\varphi_t$.  First, we have
\( \exp(\alpha h) = \mysmallmatrix{e^{i \alpha}}{0}{0}{e^{-i \alpha}}
= (\cos \alpha) I + (\sin\alpha) h \) as claimed.  A straightforward
calculation shows
\[
  \psi_t\big(\exp(\alpha h)\big) =
  \left(\begin{array}{ccc}
          \cos(2\alpha) & -\sin(2\alpha) & 0 \\
          \sin(2\alpha) & \hphantom{-}\cos(2\alpha) & 0 \\
               0       &        0      & 1
        \end{array}\right)
\]
which verifies the geometric claim.

For a general \(u \in \Qbar_t\), choose \(g \in G_t\) with 
\(g \cdot h = u\). Then 
\begin{equation*}
  \exp(\alpha u) = \exp( g \cdot \alpha h) =  g \cdot \exp(\alpha h) =
  (\cos \alpha) I + (\sin \alpha) g \cdot h = (\cos \alpha) I + (\sin \alpha) u.
\end{equation*}
We also have
$\psi_t\big(\exp(\alpha u)\big) = \psi_t\big(g \cdot \exp(\alpha
h)\big) = g  \psi_t\big(\exp(\alpha h)\big)  g^{-1}$ and, as
noted in Section~\ref{subsec:Gt}, the action of $G_t$ on $\Qbar_t$
preserves orientation.  Therefore $\exp(\alpha u)$ acts as a
$2\alpha$-anticlockwise rotation about $u$ as desired.
\end{proof}

\begin{corollary}
  \label{cor:exp Qt}
  If \(c \in (-2,2)\) and \(c = 2 \cos \alpha\), the map
  \(r\maps \Qbar_t \to \tr^{-1}(c)\cap U_t\) given by
  \(r(u)= \exp(\alpha u)\) is a diffeomorphism.
\end{corollary}
\begin{proof} It is evident from Lemma~\ref{lem:exp formula} that
\(\tr\big(\exp(\alpha u)\big) = 2 \cos \alpha\), so
\(\im r \subset \tr^{-1}(c)\). The hypothesis \(c \in (-2,2)\) implies
that \(\sin \alpha \neq 0\), so the same lemma also shows that \(r\) and
its derivative are injective. If \(\tr A = c\), we can write
\(A = (\cos \alpha) I + (\sin \alpha ) u\) for some \(u \in \ut\); the
condition that \(\det A = 1\) implies that \( u \in \Qbar_t\), so
\(r\) is surjective.
\end{proof}

\section{Joining configuration spaces}
\label{sec:join configs}

In this section, we construct the smooth structure on
\(\CnS \cup \CnsignE \cup \CnsignH\) promised in Section
\ref{subsec:join geometries}. Consider the set
\[
\cRbar_n = \setdef{(t, v_1, \ldots, v_n) \in \R \times (\R^3)^n}
                  {\mbox{$v_i \in \Qbar_t$ for all $i$;
                      $v_i \notin \{v_j,  -v_j\}$
                    for some $i, j$ if $n > 1$}}
\]
By calculating the Jacobian matrix of its $n$ defining equations, you
can check $\cRbar_n$ is a smooth $2n + 1$ dimensional submanifold of
$\R \times (\R^3)^n \cong \R^{3n + 1}$. A typical element of
$\cRbar_n$ will be denoted by $v$ with constituent vectors $v_i$ and
$t$ value denoted $t(v)$.  Given \(w\in \R^3\) which does not lie on
the \(z\)-axis, there is a unique \(t\) such that
\(w \in Q_t\). It follows that for \(n>1\), any \(v \in \cRbar_n\) is
uniquely determined by the vectors \((v_1, \ldots, v_n)\).

Recall that \(\R_{>0}\) acts on \(\R^3\) by dilating the \(x\) and
\(y\) coordinates: \(s \cdot (x,y,z) = (sx,sy,z)\). Since
\(s\cdot \Qbar_t = \Qbar_{t/s^2}\), the group \(\R_{>0}\) acts on
\(\cRbar_n\) by
\[s\cdot (t,v_1, \ldots, v_n) = (t/s^2,s\cdot v_1, \ldots, s \cdot
  v_n).\] For $t \in \R$, we set
$\cRbar_n^{t} = \setdefm{\big}{v \in \cRbar_n}{t(v) = t}$. The group
\(G_t\) acts on the slice \(\cRbar_n^{t}\) via
\(g \cdot (t, v_1, \ldots, v_n) = (t, g\cdot v_1, \ldots, g \cdot
v_n)\). The action of \(\R_{>0}\) intertwines these actions, in the
sense that if \(g \in G_t\) and we define
\((s\cdot g)(v) = s \cdot (g\cdot (s^{-1} \cdot v))\) for
$v \in \Qbar_{t/s^2}$, then \( s \cdot g \in G_{t/s^2}\).

We define the space \(\cY\) to be the quotient of \(\cRbar_n\) by the
simultaneous actions of \(\R_{>0}\) and all the \(G_t\), where $\cY$ has
the quotient topology.  We let \(\cY_-, \cY_0\), and \( \cY_+\)
be the images of \(\cRbar_n^{<0}, \cRbar_n^{0}\), and
\(\cRbar_n^{>0}\) in \(\cY\).  The main result of this section is:

\begin{theorem} 
\label{thm:varconsY}
The space \(\cY\) can be given the structure of a smooth
\((2n-3)\)\hyp dimensional manifold such that the quotient map
\(\sigma\maps \cRbar_n \to \cY\) is a smooth submersion.  The subsets
\(\cY_-\) and \(\cY_+\) are open submanifolds of \(\cY\), while
\(\cY_0\) is a closed codimension-1 submanifold. They are
diffeomorphic to \(\CnsignH\), \(\CnS\), and \(\CnsignE\)
respectively.  Finally, the action of $\pmone^n$ on $\cY$ is by
diffeomorphisms.
\end{theorem}

We will need the following technical tool in our eventual proof of
Theorem~\ref{thm:varconsY}.

\begin{lemma}\label{lem:normalize}
  Suppose $(t, v) \in \cRbar_1$ with $v = (x, y, z)$ and $z > 0$.  Let $L$
  be the unique shortest geodesic segment in $Q_t$ joining $v$ to
  $e_3 = (0, 0, 1)$. Then 
  \[\arraycolsep=4pt\def\arraystretch{1.5}
    W(t, v) = \left(\begin{array}{ccc}
       1 - \frac{t x^2}{1 + z} & -\frac{t x y}{1 + z} & x \\
       -\frac{t x y}{1 + z}   & 1-\frac{t y^2}{1 + z} & y \\
       -t x & -t y & z 
    \end{array}\right)
  \]
  is an element of $G_t$ taking $e_3$ to $v$.  When $t \leq 0$, it is
  a translation whose axis contains $L$; when $t > 0$, it is a
  rotation with angle $\theta < \pi$ whose invariant equator
  contains $L$.  Finally, the inverse of $W(t, v)$ is
  $W\big(t, (-x, -y, z)\big)$.
\end{lemma}
\begin{proof}
When $t = 0$, we have $z = 1$ and the formula for $W(t, v)$ greatly
simplifies making it clear that $W(t, v)$ is the element of $G_0$ that
translates $e_3$ to $v$.  So we henceforth assume $t \neq 0$.  A
tedious but straightforward calculation using that $Q_t(v) = 1$ shows
that the columns of $W(t, v)$ have the same $B_t$ inner products as
the standard basis $\{e_1, e_2, e_3\}$ for $\R^3$.  Hence $W(t, v)$
preserves $B_t$, and moreover it is in $G_t$ as we can connect
$W(t, v)$ to the identity matrix by moving $v$ to $e_3$ along $L$.
Now $W(t, v)$ takes $e_3$ to $v$, and we can see the claimed geometric
description by noting that $W(t, v)$ fixes $(-y, x, 0)$, which is a
basis for the $B_t$-orthogonal complement to the plane spanned by
$\{e_3, v\}$.  Noting that $e_3$, $v$, and $(-x, -y, z)$ all lie on
the geodesic containing $L$, the fact that $W(t, v)$ and
$W\big(t, (-x, -y, z)\big)$ are inverses follows from their geometric
characterizations.
\end{proof}

\subsection{A local model} 

Rather than tackle Theorem~\ref{thm:varconsY} all at once, we start
with a more specialized statement that exhibits all the key issues.
Given \(0<\delta<\pi/2\), let
\(\cV_n = \setdefm{\big}{v \in \cRbar_n}{\mbox{all $v_i \in Q_t$}}\),
and consider
\[
  \cR_n= \setdef{v \in \cV_n}{\mbox{$B_t(v_1, v_i) > \cos(\delta)$ for all $i$}}.
\]
Note that \(\cR_n\) is an open subset of $\cRbar_n$.  When $t \leq 0$,
the condition that $B_t(v_1, v_i) > \cos(\delta)$ is automatically
satisfied: applying an element of $G_t$ we can assume that $v_1 = e_3$
and then $B_t(v_1, v_i)$ is just the $z$-coordinate of $v_i \in Q_t$,
which is at least $1$.  When $t > 0$, the condition on $B_t(v_1, v_i)$
is equivalent to the geometric distance in $Q_{t}$ between $v_1$ and
$v_i$ being less than $\delta$.  Since the dilation action of
$\R_{> 0}$ on $\cRbar_n$ preserves the $B_t$, the set $\cR_n$ is
invariant under this action.

Let \(\Ydelta = \sigma(\cR_n)\) where $\sigma \maps \cRbar_n \to \cY$
is the quotient map, and take \(\Ydelta_-\),
\(\Ydelta_0\), and \(\Ydelta_+\) to be its intersection with the
subsets $\cY_-$, $\cY_0$, and $\cY_+$ of $\cY$. We will show:

\begin{theorem}\label{thm:localmodel}
  For any $n > 1$ and $0 < \delta < \pi/2$, there is a smooth
  structure on \(\Ydelta\) making it diffeomorphic to
  $P^{n-2}(\C) \times \R$. With respect to this smooth structure, the
  projection \(\sigma\maps \cR_n \to \Ydelta\) is a submersion. The subsets
  \(\Ydelta_-\), \(\Ydelta_0\), and \(\Ydelta_+\) are identified with
  $P^{n-2}(\C) \times (-\infty, 0)$, $P^{n-2}(\C) \times \{0\}$, and
  $P^{n-2}(\C) \times (0, \infty)$, respectively.
\end{theorem}
  
  The smooth structure on $\Ydelta$ is uniquely characterized by the
requirement that $\sigma$ is a submersion; this is a
consequence of the following basic fact that we will use repeatedly:
if $\pi \maps M \to \Mbar$ is a surjective submersion and
$f \maps M \to N$ is any smooth map that is constant on each fiber of
$\pi$, then the induced map $\fbar \maps \Mbar \to N$ is smooth
\cite[Theorem 4.30]{Lee2013}.

To prove Theorem~\ref{thm:localmodel}, we focus on a subset
$\cR_n''$ of $\cR_n$ which still surjects $\Ydelta$ but where
the fibers $\cR_n'' \sigmaarrow \Ydelta$ are simple enough that we can
push the smooth structure forward. To this end, consider the
(degenerate) quadratic form $H$ on $\R^3$ given by $x^2 + y^2$ and
define
\[
\cR_n' = \setdef{v \in \cR_n}{v_1 = e_3} \mtext{and} \cR_n'' = \setdefm{\bigg}{v \in \cR_n'}{\sum_{i=2}^n H(v_i) = 1}
\]
Unrolling  all the definitions, if we set $v_i = (x_i, y_i, z_i)$
then equivalently
\begin{equation}\label{eq:Rmess}
\cR_n'' = \setdef{(t, v_1, \ldots, v_n) \in \R^{3 n + 1}}{
  \substack{
    \mbox{$t (x_i^2 + y_i^2) + z_i^2 = 1$ and 
      $z_i > \cos(\delta)$ for all $i$}\\ 
    \mbox{$v_1 = e_3$ and $\sum_{i = 2}^n x_i^2 + y_i^2 = 1$}
  }
}
\end{equation}
It is straightforward to check that $\cR_n'$ and $\cR_n''$ are smooth
submanifolds of $\cR_n$.

\begin{lemma}
 \label{lem:2submersions}
  There are submersions $\cR_n \xrightarrow{\, p_1\, } \cR'_n$ and
   $\cR'_n  \xrightarrow{\, p_2\, } \cR''_n$ that commute with
   $\cR_n \xrightarrow{\, \sigma\, } \Ydelta$.
\end{lemma}

\begin{proof}
  Let us start with $p_2 \maps \cR'_n \to \cR''_n$, which we define via the
  $\R_{>0}$ action by
\[
 v \mapsto \frac{1}{\sqrt{\sum H(v_i)}} \cdot v
\]
The sets $s \cdot \cR''_n$ as $s$ varies are a family of local
sections to $p_2$ which pass through every point of $\cR'_n$, showing
that $p_2$ is a submersion.

Our map $p_1 \maps \cR_n \to \cR'_n$ is defined using
Lemma~\ref{lem:normalize} to move $v_1$ to $e_3$ in a consistent way:
given $v \in \cR_n$ with $v_i = (x_i, y_i, z_i)$, define 
\[
  p_1(v) = W\big(t(v), v_1\big)^{-1} \cdot v =  W\big(t(v), (-x_1,
  -y_1, z_1)\big) \cdot v 
\]
which is a smooth function by the formula for $W$ in
Lemma~\ref{lem:normalize}.  To see this is a submersion at
$w \in \cR_n$, pick a path
$\gamma \maps \big(t(w) - \eta, t(w) + \eta\big) \to \cV_1$ for some
$\eta > 0$ with $t(\gamma(s)) = s$ for all $s$ and
$\gamma(t(w)) = w_1$. Then the map $\beta_w \maps \cR'_n \to \cR_n$
given by
\[
  \beta_w(v) = W\big(t(v),\  \gamma(t(v))\big) \cdot v
\]
gives a local section to $p_1$ whose image contains $w$.
\end{proof}

Note that both $\cR_n'$ and $\cR_n''$ are invariant under the
$S^1$ action of rotations about the $z$-axis. This \(S^1\) is a
subgroup of \(G_t\) for all \(t\), giving a continuous map
\[
  j \maps \cR''_n/S^1 \to \Ydelta \mtext{via the composition} \cR''_n
  \xrightarrow{\, p_3\, } \cR''_n/S^1 \to \sigma(\cR''_n)
  \hookrightarrow \Ydelta
\]
where $p_3$ is the topological quotient map $\cR''_n \to \cR''_n/S^1$.

\begin{lemma}\label{lem: j homeo}
The map \(j\) is a homeomorphism. 
\end{lemma}

\begin{proof}
This $S^1$ can be viewed as the subgroup of $G_t$ that stabilizes
$e_3$, so any $g \in G_t \setminus S^1$ moves
$\setdef{v \in \cR_n'}{t(v) = t}$ to something disjoint from $\cR_n'$.
Additionally, any $s \neq 1$ in $\R_{>0}$ takes $\cR_n''$ completely
off itself since $H(s \cdot v) = s^2 H(v)$.  It follows that \(j\) is
injective. Lemma~\ref{lem:2submersions} implies that \(j\) is
surjective, so it remains to check that \(j^{-1}\) is continuous.

Given \(U \subset \cR_n''/S^1\) open, we must check that \(j(U)\) is
open in \(\Ydelta\), i.e. that \(\sigma^{-1}(j(U))\) is open in
\(\cR_n\). This too follows from Lemma~\ref{lem:2submersions}, since
\(\sigma^{-1}(j(U)) = p_1^{-1}(p_2^{-1}(p_3^{-1}(U)))\), and \(p_1\),
\(p_2\), and \(p_3\) are all continuous.
\end{proof}

\begin{proof}[Proof of Theorem~\ref{thm:localmodel}]
The \(S^1\) action on \(\cR_n''\) is free and proper, so
\(\cR_n''/S^1 = \Ydelta\) inherits a smooth structure from $\cR_n''$,
where we have identified \(\cR_n''/S^1\) with \(\Ydelta\) via $j$
from Lemma~\ref{lem: j homeo}.  In addition,
$\sigma\maps \cR_n \to \Ydelta$ is a submersion, since
\(\sigma = p_3 \circ p_2 \circ p_1\) and each of the individual
factors is a submersion.

It remains to identify the topology of $\Ydelta$ and its various
pieces.  Consider the smooth map
$\kappa \maps \cR_n'' \to \R \times S^{2 n - 3}$ defined by
$\kappa(v) = (t(v), x_2, y_2, \ldots, x_n, y_n)$, where as usual
$v_i = (x_i, y_i, z_i)$.  From the description of $\cR_n''$ in
(\ref{eq:Rmess}), we see $\kappa$ is injective with image the open set
\begin{equation}\label{eq: Tn}
  \cT_n = \setdef{ (t, x_2, y_2, \ldots, x_n, y_n) \in \R \times
    S^{2n-3}} {t < \frac{\sin^2(\delta)}{\max_i\left(x_i^2 +
        y_i^2\right)}}
\end{equation}
and moreover that $\kappa$ has a smooth inverse on $\cT_n$.
Therefore, the map $\kappa$ is a diffeomorphism onto $\cT_n$. There is
a natural action of $S^1$ on $\R \times S^{2n-3}$ by viewing the
latter as a subset of $\R \times \C^{n-1}$ and acting diagonally on
the $\C$ factors. As this action is compatible with the $S^1$ action
on $\cR_n''$, we have an embedding of $\Ydelta = \cR_n''/S^1$ into
$\R \times P^{n-2}(\C)$.  Since $\cT_n$ is simply the region lying
below the graph of a piecewise smooth function $S^{2n - 3} \to \R$, we
see that $\Ydelta$ is indeed diffeomorphic to $\R \times P^{n-2}(\C)$.
The remaining conclusions of the theorem are now easily checked.
\end{proof}

\begin{proof}[Proof of Theorem~\ref{thm:varconsY}]
For \( 0 < \delta < \pi/2\), let
\( \cR_n^*= \setdefm{\big}{v \in \cRbar_n}%
{\mbox{$B_t(v_1, v_i) > \cos(\delta)$ for all $i$}}\), which contains
$\cR_n$ as an open subset.  A \(v \in \cRbar_n\) with $t(v) > 0$ is in
\(\cR^*_n\) if and only if the geometric distance in $Q_{t(v)}$ from
\(v_1\) to \(v_i\) is less than \(\delta\) for all \(i\). When
\(t(v) \leq 0 \), a \(v \in \cRbar_n\) is in \(\cR^*_n\) exactly when
the \(v_i\) all belong to the same component of \(\Qbar_{t(v)}\).  As
\(\cR^*_n\) is an open subset of $\cRbar_n$ which is invariant under the
actions of \(\R_{>0}\) and the \(G_t\), \(\sigma(\cR_n^*)\) is an
open subset of \(\cY\).  Note that 
\(\sigma(\cR_n) \subset \sigma(\cR_n^*)\); we claim that these two
sets are equal.  Indeed, given \(v \in \cR_n^*\), we can find
\(g \in G_{t(v)}\) with \(g \cdot v_1 = e_3\). Since
\( \delta < \pi/2\), the element \(g\cdot v\) is in \(\cR_n\), and so
\(\sigma(\cR_n^*) \subset \sigma(\cR_n)\). Hence
\(\sigma(\cR_n^*) = \sigma(\cR_n) = \Ydelta\) is open in \(\cY\).
 
The group \(\{\pm 1\}^n\) acts on \(\cRbar_n\), where the action is
given by \((\epsilon \cdot v)_i = \epsilon_i v_i\). This action
commutes with the actions of \(G_t\) and \(\R_{>0}\). Let
\(\One\) be the identity element of \(\{\pm 1\}^n\). The
subgroup \(\pm \One\) preserves \(\cRbar^*_n\), and if
\( \epsilon \neq \pm \One\), then we have
\( \epsilon \cdot \cR^*_n \cap \cR_n^* = \emptyset\).  For
\(\chi \in \{\pm 1\}^n / \pm \One\), let
\( \cR_n^\chi = \chi \cdot \cR_n^*\). The sets \(\cR_n^\chi\) are
disjoint open subsets of \(\cRbar_n\) which are invariant under the
actions of \(\R_{>0}\) and the \(G_t\). Together with
\(\cRbar_n^{>0}\), they form an open cover of \(\cRbar_n\).

Let \(\cY^\chi = \sigma(\cR_n^\chi)\), which is an open subset of
\(\cY\).  By Theorem~\ref{thm:localmodel}, \(\cY^\chi\) can be given the
structure of a smooth manifold diffeomorphic to
\(P^{n-2}(\C) \times \R\) for which the map
\(\sigma\maps \cR_n^\chi \to \cY^\chi\) is a submersion.

The set \(\cRbar_n^{>0}\) is invariant under the actions of
\(\R_{>0}\) and \(G_t\), so \(\cY_+ = \sigma \big(\cRbar_n^{>0}\big)\) is an open
subset of \(\cY\). The set \(\cRbar_n^1\) is a global slice for the
action of \(\R_{>0}\) on \(\cRbar_n^{>0}\), so
\[\cY_+ \cong \cRbar_n^1/G_1 \cong \cS_n/\Isom^+(S^2) =
  \CnS\] can be given the structure of a smooth manifold. The map
\(\sigma\maps \cRbar_n^{>0} \to \CnS\) is the composition of the map
\(p_4\maps \cRbar_n^{>0} \to \cRbar_n^{1}\) given by
\(p_4(v) = \sqrt{t(v)} \cdot v\) with the projection
\( p_5\maps \cRbar_n^1 \to \CnS\). Both \(p_4\) and \(p_5\) are
submersions, so \(\sigma\maps \cRbar_n^{>0} \to \CnS\) is a submersion as
well.

In summary, we have shown that the \(\cY^\chi\), together with
\(\cY_+\), form an open cover of \(\cY\), each of whose elements is a
smooth manifold. To show that \(\cY\) is a smooth manifold, we must
show that the transition functions relating different sets in our
cover are smooth and that \(\cY\) is Hausdorff. For the first,
consider the overlap $\cY^\chi_+ = \cY^\chi \cap \cY_+$.  Viewing
$\cY^\chi_+$ as a subset of $\cY^\chi$ and identifying $\cY_+$ with
$\CnS$, we get a transition function \(g_\chi\maps \cY^\chi_+ \to
\CnS\). As \(g_\chi\) is induced by the smooth map
\(\sigma\maps \cRbar_n^{>0} \to \CnS\), which is constant on the fibers of
the submersion \(\sigma\maps \cRbar_n^{>0} \cap \cR_n^*\to \cY^\chi_+\), we
know $g_\chi$ is smooth by \cite[Theorem 4.30]{Lee2013}. A similar
argument shows \(g_\chi^{-1}\) is smooth, making all transition functions
smooth as needed.
 
To show $\cY$ is Hausdorff, suppose \(x\) and \(y\) are distinct
points of \(\cY\). The open sets \(\cY^\chi\) and \(\cY_+\) forming
the cover are themselves Hausdorff, so if \(x\) and \(y\) both belong
to the same open set \(U\) of the cover, they can be separated by sets
\(U_x, U_y\) which are open in \(U\) and hence open in \(\cY\). On the
other hand, if \(x\in \cY^{\chi_1}\) and \(y \in \cY^{\chi_2}\) with
\(\chi_1 \neq \chi_2\), then \(x\) and \(y\) are separated by
\( \cY^{\chi_1}\) and \( \cY^{\chi_2}\). The only remaining case is
when \(x \in \cY^\chi \setminus \cY_+\) and
\(y \in \cY_+ \setminus \cY^\chi \). Since $x \notin \cY_+$, we have
$x \in \chi \cdot \cY^{\delta_0}$ for all $0 < \delta_0 < \pi/2$.
Write \(y = \sigma(v)\) for \(v \in \cRbar_n^1\), and let
\(\delta_0 = \frac{1}{2}\max(d(v_i,v_j))\). Then \(x\) and \(y\) are
separated by the open sets \(U_x = \chi \cdot \cY^{\delta_0}_n\) and
\(U_y = \setdef{\sigma(v)}{\mbox{$t(v)>0$ \, and \,
    $\max(d(v_i,v_j))>\delta_0$}}\).
 
With this smooth structure on \(\cY\), the group $\pmone^n$ acts by
diffeomorphisms since our open cover was constructed using that very
action and we already knew $\pmone^n$ acts smoothly on
$\cY_+ \cong \CnS$.  Moreover, the map
\(\sigma \maps \cRbar_n \to \cY\) is a smooth submersion, since its
restriction to each set in our open cover of \(\cRbar_n\) is a smooth
submersion.  It remains to check the descriptions of \(\cY_-,\cY_0\),
and \( \cY_+\) given in the statement. We have already shown that
\(\cY_+\) is a smooth submanifold of \(\cY\) diffeomorphic to
\(\CnS\). The same argument shows that \(\cY_-\) is diffeomorphic to
\(\CnsignH\).  Finally, Theorem~\ref{thm:localmodel} implies that
\(\cY_0\) is a closed submanifold of \(\cY\). The subspace \(\cY_0\)
is the quotient of \(\cR_{0}\) by the simultaneous actions of \(G_0\)
and \(\R_{>0}\). This is the quotient of
\[ \setdefm{\big}{v \in \Qbar_0^n}{v_i \neq \pm v_j \ \text{for some} \ i,
    j} \] by the action of the group \(\Ghat_0\), which is the space
\(\cE_n^\pm/\Sim(\E^2) = \CnsignE\).  So $\cY_0 \cong \CnsignE$,
as needed.
\end{proof}

Applying Theorem~\ref{thm:localmodel} to the open cover
$\cY_+ \cup \{\cY^\chi\}_{\chi \in \pmone^n/\pm \One}$ of $\cY$
in the preceding proof of Theorem~\ref{thm:varconsY}, we see that
$\cY$ is diffeomorphic to its subset $\cY_+$ and hence:
\begin{corollary}
 \label{cor:cYCnS}
 The manifold \(\cY\) is diffeomorphic to \(\CnS\). 
\end{corollary}
As we oriented $\CnS$ in Section~\ref{sec:configorient}, the
corollary shows that $\cY$ is orientable.  Moreover, with the
orientation on \(\CnH\) from Section~\ref{sec:configorient}:
\begin{proposition} 
\label{prop:orientation on Y}
If we orient \(\cY\) so that the induced orientation on \(\cY_+\)
matches that of \(\CnS\), then the induced orientation on
\(\cY_-\) is the reverse of that on \(\CnH\).
\end{proposition}

\begin{proof}
We begin with the case $n = 2$.  The maps $\CnH[2] \to \R_{>0}$ and
$\CnS[2] \to (0, \pi)$ that compute the distance between the two
points are diffeomorphisms which, given our conventions in
Section~\ref {sec:configorient}, are either both orientation
preserving or both orientation reversing (it is the former, but we do
not need this).  Recall the set $\cT_2$ from equation~(\ref{eq: Tn}).
Then we can identify $\Ydelta$ with the slice
$\setdef{(t, 1, 0)}{t < \sin^2(\delta)}$ of the $S^1$ action on
$\cT_2$.  For $t > 0$, increasing $t$ increases the distance between
the two points in $\CnS[2] \cap \Ydelta$. However, for $t < 0$,
increasing $t$ decreases the distance on $\CnH[2]$. Thus, in this case
$\CnS[2]$ and $\CnH[2]$ give opposite orientations to $\Ydelta$.

For general $n$, the moral is that for $t > 0$, we act by $\sqrt{t}$
to normalize a point of $\cR_n''$ into $\cV_n^1 = \cS_n$, whereas, for
$t < 0$, we act by $\sqrt{-t}$ to normalize a point of $\cR_n''$ into
$\cV_n^{-1} = \cH_n$; as the orientations on $\cV_n^1$ and
$\cV_n^{-1}$ match in a suitable sense, the orientation flip between
$\sqrt{t}$ and $\sqrt{-t}$ means that $\CnSdelta$ and $\CnH$ give
opposite orientations to $\Ydelta$.
  
This proves the proposition for the subset of $\CnsignH$ where all
points have the same sign.  As discussed at the end of
Section~\ref{sec:configorient}, the action of an element on $\pmone^n$
on $\CnS$ and $\CnsignH$ either preserves the orientation of both or
reverses it for both.  As the action of $\pmone^n$ is transitive on
the components of $\CnsignH$, this completes the proof.
\end{proof}

Finally, we can assemble all the pieces to give the result promised
in Section~\ref{sec:config}:

\begin{proof}[Proof of Theorem~\ref{thm:consY}]
All but the last claim are covered by Theorem~\ref{thm:varconsY} and
Proposition~\ref{prop:orientation on Y}; the last claim follows
from Theorem~\ref{thm:localmodel} and the $\pmone^n$ action.
\end{proof}
 
\section{The character variety of the free group}
\label{sec:map to X}

Fix $c \in (-2, 2)$ and consider the free group $F_n$ generated by
$S = \{ s_1, \ldots, s_n\}$. Recall from Sections~\ref{subsec:X(F_n)}
and \ref{subsec:reducibles} that the character variety
\(X^c(F_n,S) = X_\SLR^c(F_n,S) \cup X_\SU^c(F_n,S)\) is a smooth
manifold away from the reducible locus, which consists of \(2^{n-1}\)
reducible characters \(\chiep^c\) labeled by
\([\epsilon] \in \{\pm 1\}^n/\pm \One\). The main result of this
section is that the space \(\cY\) from Section~\ref{thm:varconsY} is a
resolution of \(X^c(F_n,S)\):
\begin{theorem}
\label{thm:def of pi}
For each \(c \in (-2,2)\), there is a proper smooth surjective map
\(\pi_c\maps \cY \to X^c(F_n,S)\). The function \(\pi_c\) maps
\(\cY_- \cong \CnsignH\) diffeomorphically onto
\(X^\cirr_\SLR(F_n,S)\) and maps \(\cY_+ \cong \CnS\) diffeomorphically
onto \(X^\cirr_\SU(F_n,S)\). It maps the component of
\(\cY_0\cong \CnsignE\) corresponding to
\([\epsilon] \in \{\pm 1\}^n/\pm \One\) to the reducible
character \(\chiep^c\).
\end{theorem}
To understand the statement that $\pi_c$ is smooth, note that while
only the dense open set $X^\cirr(F_n, S)$ is a smooth manifold for
$n > 2$, the semialgebraic set $X^{c}(F_n, S)$ can be viewed as a
subset of some $\R^m$ using suitable trace functions as
coordinates. While the same space \(\cY\) appears as the resolution
for all \(c \in (-2,2)\), when we want to emphasize the value of
\(c\), we write \(\cX^c(F_n,S)\) for \(\cY\). Furthermore, we set
\begin{equation}\label{eq: cXirr}
  \cX^\cirr(F_n,S) := \pi_c^{-1}\big(X^\cirr(F_n,S)\big) \cong \CnS \cup \CnsignH
\end{equation}
and note that $\cX^\cirr(F_n,S)$ is dense in $\cX(F_n, S)$ by the
last sentence in Theorem~\ref{thm:consY}.  Similarly, we define
$\cX^\cred(F_n, S) = \pi_c^{-1}\big(X^\cred(F_n, S)\big)\cong \CnsignE$.

We will define $\pi_c$ using a map $\pi_c'$ from the space $\cRbar_n$
of Section~\ref{sec:join configs}, whose quotient is $\cY$, to the
representation variety $\RC(F_n, S)$.  In doing so, we will make heavy
use of the subgroups $U_t$ of $\SLC$ from Section~\ref{sec:SL2Csubs}.
Recall the conjugacy class of $U_t$ in $\SLC$ depends only on whether
$t$ is positive, negative, or zero.  Thus much of the distinction
between the $U_t$ is meaningless at the level of character varieties.
However, this nicely varying family $U_t$ will be crucial in writing
down a smooth $\pi_c' \maps \cRbar_n \to \RC(F_n, S)$ even though
$U_1 = \SU$ and $U_{-1} \cong \SLR$ are radically different subgroups
of $\SLC$.

Turning now to the reducible locus, the points of $\cX^\cred(F_n, S)$
can be understood in terms of representations to $U_0$ as
follows. Recall from Sections~\ref{subsec:Gt} and
\ref{sec:SL2Csubs} that $G_0 \cong \Isom(\E^2)$ and $U_0 \to G_0$ is a
double-cover onto the identity component of $G_0$, which is
$\Isom^+(\E^2)$.  Recall also that $G_0$ is contained in
$\Ghat_0 \cong \Sim(\E^2)$ which is $\Aut(U_0)$ by
Lemma~\ref{Lem:G0=Aut U0}.  A representation $\rho$ to $U_0$ is
\emph{fully reducible} when it is conjugate to one with diagonal image.
Now set
\begin{equation}
  \label{eq: RcU0}
  R^c_{U_0}(F_n, S) = \setdef{\rho\maps F_n \to U_0}%
  {\mbox{all
      $\tr \rho(s_i) = c$ and $\rho$ not fully reducible}}.
\end{equation}
We define \(X^c_{U_0}(F_n, S)\) to be the quotient
\(R^c_{U_0}(F_n, S)/\Ghat_0\). Despite the notation, the set
\(X^c_{U_0}(F_n, S)\) is not a true character variety: it turns out to
be a manifold of dimension $2n - 4$, whereas the representations in
\( R^c_{U_0}(F_n, S)\) have only the finitely many distinct characters
$\chiep^c$. We will show:

\begin{theorem}
  \label{thm: red is U0}
  The group $\Ghat_0 \cong \Aut(U_0)$ acts freely and properly on
$R^c_{U_0}(F_n, S)$, so \(X^c_{U_0}(F_n, S)\) has a natural smooth structure, and 
  there is a diffeomorphism $\picred:
  \cX^\cred(F_n, S) \to X^c_{U_0}(F_n, S)$.  Under this
  identification, the restriction $\pi_c
  \maps \cX^\cred(F_n, S) \to X^\cred(F_n, S)$ corresponds to the
  map $X^c_{U_0}(F_n, S) \to X^\cred(F_n, S)$ taking $[\rho]$
  to its character.
\end{theorem}
Thus we can think of Theorem~\ref{thm:def of pi} as saying we can
resolve $X^c(F_n, S)$ by ``blowing up'' the finitely many reducible
representations into the space of ``characters'' of $U_0$
representations that are not fully reducible.

\subsection{Constructing $\pi_c$}
\label{sec: pi_c}

Given $v \in \Qbar_t$ and $\alpha \in (0, \pi)$, we regard $v$ as the
element $\varphi_t(v)$ of the Lie algebra $\ut$ and define
$A(\alpha, t, v) = \exp(\alpha \varphi_t(v)) \in U_t$; as per
Lemma~\ref{lem:exp formula}, the element $A(\alpha, t, v)$ acts on
$\Qbar_t$ as a $2\alpha$--anticlockwise rotation about $v$ with
\begin{equation}\label{eq:A}
  A(\alpha, t, v) = (\cos \alpha) I  +  (\sin \alpha) \varphi_t(v).
\end{equation}
Note that $\tr\big(A(\alpha, t, v)\big) = 2 \cos \alpha$ and
$A(\alpha, t, -v) = A(\alpha, t, v)^{-1}$.  For \(c \in (-2,2)\), let
\(\alpha = \cos^{-1}(c/2)\).  We define
$\pi_c' \maps \cRbar_n \to \RC^{c}(F_n, S)$ by sending
$v \in \cRbar_n$ to the representation where
\[
  s_i \mapsto A(\alpha, t(v), v_i).
\]
If \(t(v)\neq 0 \), the condition that \(v_i \neq \pm v_j\) for some
\(i,j\) implies that \(\pi_c'(v)\) is irreducible.  When $t(v) = 0$,
that condition instead implies that \(\pi_c'(v)\) is not fully reducible.

First, we study the image of $\pi'_c(\cRbar_n)$ under the
projection
\(\tau\maps \RC^c(F_n,S) \to \XC^c(F_n,S)\):

\begin{lemma}\label{lem: pi prime c}
  For each $t \neq 0$, the restriction of $\pi'_c$ to
  $\cRbar_n^t$ gives a diffeomorphism to $R^\cirr_{U_t}(F_n,S)$.
  Also, the restriction of $\pi'_c$ to $\cRbar_n^0$ gives a
  diffeomorphism to  $R^c_{U_0}(F_n,S)$ as defined in (\ref{eq: RcU0}).
  Finally, the subset \(\tau\big(\pi_c'(\cRbar_n)\big)\) is equal to
  \(X^c(F_n,S)\).
\end{lemma}

\begin{proof}
To start in on the initial claim where $t \neq 0$ is fixed, note that
for $v \in \cRbar_n^t$ each $A(\alpha, t, v_i)$ is in $U_t$ and so
$\pi'_c(v) \in R^\cirr_{U_t}(F_n,S)$.  Next, we show
\(\pi'_c\maps \cRbar_n^t \to R^\cirr_{U_t}(F_n,S)\) is onto.  By
Corollary~\ref{cor:exp Qt}, given $\rho \in R^\cirr_{U_t}(F_n,S)$,
there are unique $v_i \in \Qbar_t$ with
$\rho(s_i) = A(\alpha, t, v_i)$.  At least one $v_i$ is not $\pm v_1$
as otherwise all $\rho(s_i)$ would be $\rho(s_1)^{\pm 1}$, violating
the condition that $\rho$ is irreducible.  Hence the associated $v$ is in
$\cRbar_n^t$, establishing that
$\pi'_c\big(\cRbar_n^t\big) = R^\cirr_{U_t}(F_n,S)$.  By uniqueness of
the $v_i$, in fact \(\pi'_c\maps \cRbar_n^t \to R^\cirr_{U_t}(F_n,S)\)
is a bijection; Corollary~\ref{cor:exp Qt} implies it is also a
diffeomorphism since locally $\pi'_c$ on $\cRbar_n^t$ is
the restriction of the product of $n$ copies of the map $r$ in the
corollary.  This completes the proof of the first claim, and the
second claim covering the case $t = 0$ follows by the identical argument
with ``irreducible'' replaced with ``not fully reducible''.

Turning now to the final claim of the lemma, first recall that
\[
  X^c(F_n,S) = X_\SLR^\cirr(F_n,S) \cup X_\SU^\cirr(F_n,S) \cup  X^\cred(F_n, S).
\]
When \(t <0\), the subgroup \(U_{t} \leq \SLC\) is conjugate to \(\SLR\), so
\[
  \tau\big(\pi'_c(\cR^t_n)\big) = \tau\big(R^\cirr_{U_t}(F_n,S)\big)
  = \tau\big(R^\cirr_{\SLR}(F_n,S)\big)= X_\SLR^\cirr(F_n,S).
\]
For $t > 0$, the subgroup $U_t$ is conjugate to $\SU$ and
so $\tau\big(\pi'_c(\cR^t_n)\big) = X_\SU^\cirr(F_n, S)$. It thus
remains to show that $\tau\big(\pi'_c(\cR^0_n)\big) = X^\cred(F_n S)$.

In one direction, given $v \in \cRbar^0_n$, let
$\epsilon_i \in \pmone$ be the $z$-coordinate of $v_i$, and so
\[
  A(\alpha, 0, v_i) = \twobytwomatrix{e^{i \epsilon_i \alpha}}{*}{0}{e^{-i \epsilon_i \alpha}}.
\]
Hence the character of \(\pi_c(v)\) is the reducible character
\(\chiep^c\) from Section~\ref{subsec:reducibles}.  Conversely, given
$\chi_{[\epsilon]}^c \in X^\cred(F_n, S)$, take
$v_1 = (1, 0, \epsilon_1)$ and all other $v_i = (0, 0, \epsilon_i)$ to
get $v \in \cRbar_n^0$ with $\pi'_c(v) = \chi_{[\epsilon]}$. Combined,
we have  $\tau\big(\pi'_c(\cR^0_n)\big) = X^\cred(F_n S)$, completing
the proof of the lemma.
\end{proof}

\begin{lemma}\label{lem:diagram}
  There is a smooth map \(\pi_c\maps \cY \to X^c(F_n,S)\) making the
  following diagram commute:
  \begin{equation}\label{eq:diagram}
    \begin{tikzcd}
      \cRbar_n \arrow{r}{\pi'_c} \arrow{d}{\sigma} & 
      \im(\pi_c') \arrow{d}{\tau} \\
      \cY \arrow{r}{\pi_c} &  X^{c}(F_n, S)
    \end{tikzcd}
  \end{equation}
\end{lemma}

\begin{proof} To show that \(\pi_c\) exists as a function, we must
check that \(\tau \circ \pi_c'\) is invariant under the actions of
\(G_t\) and \(\R_{>0}\).  First, recall from
Section~\ref{sec:SL2Csubs} that $G_t \subset \Aut(U_t)$ and that the action of $G_t$ on $U_t$ satisfies
\(\exp(g\cdot v) = g \cdot \exp v\) and so
\begin{equation}
  \label{eq:UActionOnA}
  A(\alpha,t,g \cdot v) = g \cdot A(\alpha,t, v).
\end{equation}
Using this, we define an action of \(G_t\) on \(R^c_{U_t}(F_n, S)\) by
\( (g \cdot \rho)(s) = g \cdot \rho(s)\), where $g \cdot \rho$ is also
a representation of $F_n$ as $g$ acts as an element of $\Aut(U_t)$;
with this action, $\pi_c'$ is $G_t$ equivariant. Since
\(\tr(g \cdot B) = \tr B\) for all $B \in U_t$, we have
\(\tau (g \cdot \rho) = \tau(\rho)\).  Equation \eqref{eq:UActionOnA}
implies that
\(\tau\big(\pi_c'(g\cdot v)\big) = \tau\big(g \cdot \pi_c'(v)\big) =
\tau(\pi_c'(v))\), and thus $\tau \circ \pi'_c$ is invariant under
each $G_t$ action as desired.

Turning to the $\R_{>0}$ action, for the more general action of $u \in
\C^\times$, equation (\ref{eq:u action}) gives that
\begin{equation}
  \label{eq:CactionOnA}
  A(\alpha, t/\abs{u}^2, u \cdot v) = T_uA \, (\alpha,t,v) \, T_u^{-1}
  \mtext{where $T_u = \mysmallmatrix{u}{0}{0}{1}$.}
\end{equation}
So for $s \in \R_{>0}$, it follows
that  \(\pi_c'(s \cdot v) = T_s \cdot \pi_c'(v)\), where \(T_s\) acts
by conjugation, giving the claimed invariance
\(\tau\big(\pi_c'(s\cdot v)\big) = \tau\big(\pi_c'(v)\big)\).

To see $\pi_c$ is smooth, first note $\pi'_c$ is smooth by (\ref{eq:A})
and (\ref{eq:Ut}).  Moreover, the map $\tau \maps \RC(F_n, S) \to \XC(F_n, S)$ is
smooth since it corresponds to taking traces of matrices.  Hence
$\tau \circ \pi'_c$ in the diagram is smooth.  The map $\sigma$ is a
submersion by Theorem~\ref{thm:varconsY}, and so it follows that $\pi_c$ is
smooth as claimed.
\end{proof}

We can now prove the theorems stated at the beginning of this section.

\begin{proof}[Proof of Theorem~\ref{thm: red is U0}]
To connect $\cX^\cred(F_n, S)$ and $X^c_{U_0}(F_n, S)$, note that
Lemma~\ref{lem: pi prime c} gives a diffeomorphism
$\cRbar_n^0 \to R^c_{U_0}(F_n, S)$ by restricting $\pi'_c$.  Recall
that $\Ghat_0 \cong \Aut(U_0) \cong \Sim(\E^2)$ acts on both
$\cRbar_n^0$ and $R^c_{U_0}(F_n, S)$.  We now show $\pi'_c$ is
$\Ghat_0$ equivariant, using that the action of $\Ghat_0$
on $\cRbar_n^0$ is generated by those of $G_0$ and $\R_{> 0}$. In the
proof of Lemma~\ref{lem:diagram}, we showed
$\cRbar_n^0 \to R^c_{U_0}(F_n, S)$ is $G_0$ equivariant. The element
of $\Ghat_0$ corresponding to $s \in \R_{>0}$ is the image of
$\mysmallmatrix{\sqrt{s}}{0}{0}{1/\sqrt{s}}$ in $U'_0$ under $\psi_0$,
and these act equivariantly since conjugating by $T_u$ in
(\ref{eq:CactionOnA}) is the same as conjugating by
$\mysmallmatrix{\sqrt{u}}{0}{0}{1/\sqrt{u}}$.  Thus
$\pi'_c \maps \cRbar_n^0 \to R^c_{U_0}(F_n, S)$ is $\Ghat_0$
equivariant.

As
$\pi'_c \maps \cRbar_n^0 \to R^c_{U_0}(F_n, S)$ is a diffeomorphism,
it follows that the $\Ghat_0$ action on $R^c_{U_0}(F_n, S)$ is free
and proper.  Hence $X^c_{U_0}(F_n, S)$ is a smooth manifold
diffeomorphic to $\cRbar^0_n/\Ghat_0 = \cX^\cred(F_n, S)$ by the map
$\picred$ induced by $\pi'_c$.  The final claim of the theorem now
follows from commutativity of (\ref{eq:diagram}).
\end{proof}

\begin{proof}[Proof of Theorem~\ref{thm:def of pi}]
By Lemma~\ref{lem:diagram}, we have a smooth
$\pi_c \maps \cY \to X^c(F_n, S)$ which is surjective by
Lemma~\ref{lem: pi prime c}.  Momentarily deferring the proof that
$\pi_c$ is proper, we next show that $\pi_c$ gives diffeomorphisms
\(\cY_- \cong X^\cirr_{\SLR}(F_n,S)\) and
\(\cY_+ \cong X^\cirr_{\SU}(F_n,S)\). For \(t<0\) there is a
commutative square
\begin{equation}
  \label{eq: subsquare}
  \begin{tikzcd}
     \cRbar_n^{t} \arrow{r}{\pi'_c} \arrow{d}{\sigma} & 
     R^\cirr_{U_{t}}(F_n,S) \arrow{d}{\tau} \\
     \cY_- \arrow{r}{\pi_c} &  X^{\cirr}_\SLR(F_n, S)
   \end{tikzcd}
\end{equation}
which sits inside the commutative square (\ref{eq:diagram}). The map
\(\pi_c'\) is a diffeomorphism by Lemma~\ref{lem: pi prime c} and it
is equivariant with respect to the action of \(G_{t}\). The maps
\(\sigma\) and \(\tau\) are projections to their quotients by the
actions of \(G_t\) (the latter by Lemma~\ref{lem: SLRpm quo} and the
fact that \(U_t'\) is the conjugate of \(\SLRpm\) corresponding to
$U_t$), and both are smooth submersions. It follows that \(\pi_c\) is a
diffeomorphism. The argument for \(t>0\) is the same, and the case
\(t=0\) follows from Theorem~\ref{thm: red is U0} and the
identification \(\cX^{c, \red}(F_n,S) \cong \CnsignE\).
  
Finally, we must check that \(\pi_c\) is proper.  By
Theorem~\ref{thm:localmodel} and the discussion in the proof of
Lemma~\ref{lem: pi prime c}, near
$\cY_0 = \pi_c^{-1}\big(X^\red(F_n, S)\big)$ the map \(\pi_c\) is
modeled on the projection
\(P^{n-2}(\C) \times [-1,1] \to P^{n-2}(\C) \times [-1,1] \big/
P^{n-2}(\C) \times \{0\}\), which is proper as the domain is
compact. Away from $\cY_0$, the map $\pi_c$ is a
diffeomorphism. Combining these two pictures, we see $\pi_c$ is proper
as needed.
\end{proof}

\section{The punctured sphere}
\label{sec:punc sphere}

\subsection{Introduction}
\label{subsec: intro to punc sphere}
  
Let \(S_n\) be the sphere with \(n\) punctures. We use the presentation
$$
\pi_1(S_n) = \spandef{s_1, \ldots, s_n}{s_1s_2\cdots s_n = 1}
$$
where the \(s_i\) are a fixed set of loops representing the conjugacy
classes of the boundary circles of \(S_n\). As before, we take
$S = \{s_1, \ldots, s_n \}$, and consider the character variety
$X^c(S_n, S) := X(S_n, S) \cap X^c(F_n, S)$; as $S$ will be fixed
throughout this section, we omit it and simply use \(X^c(S_n)\) and
analogous notations.

The aim of this section is to construct a resolution of \(X^c(S_n)\).
As in the free group case, it turns out that \(X^{\cirr}(S_n)\) is a
smooth manifold (see Section~\ref{sec:irrchars} below), and
\(X^c(S_n)\) is usually singular along \(X^{\cred}(S_n)\).  To resolve
\(X^c(S_n)\), we make use of our resolution $\cX^c(F_n) = \cY$ from
Theorem~\ref{thm:def of pi}; to reduce clutter, we will use $\pi$
rather than $\pi_c$ for the map $\cX^c(F_n) \to X^c(F_n)$.  The main
result of this section is:

\begin{theorem}
  \label{thm:cXSn}
  For \(c \in (-2,2)\) and an even $n \geq 4$, there is a smooth
  closed codimension 3 submanifold \(\cX^c(S_n) \subset \cX^c(F_n)\)
  such that \(\pi\big(\cX^c(S_n)\big) = X^c(S_n)\) and
  \(\pi:\cX^c(S_n) \to X^c(S_n)\) is a diffeomorphism away from
  \(X^{c,\red}(S_n)\).  Moreover, the closure of $\cX^\cirr(S_n)$ in
  $\cX^c(F_n)$ is all of \, $\cX^c(S_n)$.
\end{theorem}
For odd $n$, the picture is slightly more complicated, and we do not
pin it down completely as it is not needed in the rest of the paper;
see Section~\ref{sec: odd punct} for details.  For $n < 4$, there is
no need to resolve $X^c(S_n)$ as you can check that it is empty when
$n = 1$ and a single point when $n$ is 2 or 3.

As with Theorem~\ref{thm: red is U0} for free groups, we can
understand $\cX^\cred(S_n) := \cX^c(S_n) \cap \cX^\cred(F_n)$ in terms
of certain $U_0$ representations.  Specifically, taking
$X^c_{U_0}(S_n)$ to be the subset of $X^c_{U_0}(F_n)$ corresponding to
representations of $\pi_1(S_n)$, we will see below that $\cX^c(S_n)$ is a
closed codimension \(1\) subset of $X^c_{U_0}(S_n)$.

For the rest of this section, we fix \(c = 2 \cos \alpha \in (-2,2)\)
and some $n \geq 4$; unless explicitly specified, $n$ can be either
even or odd.  Although \(\cX^c(F_n) \cong \cY\) for all
\(c\in(-2,2)\), the subset \(\cX^c(S_n)\), and indeed its topology,
will depend on \(c\).  To begin, we consider the map
\(F:\cRbar_n \to \SLC\) given by
$$
F(v) = 
\prod_{i=1}^n A(\alpha, t(v), v_i) =
\begin{pmatrix} F_{11} & F_{12} \\ F_{21}& F_{22} \end{pmatrix}.
$$
In other words, the map \(F\) is the composition of
\(\pi'_c: \cRbar_n \to R_\C^c(F_n)\) from Section~\ref{sec:map to X}
with the evaluation map $\rho \mapsto \rho(s_1 s_2 \cdots s_n)$.
Although the definition of \(F\) depends on the trace \(c\) and the
entries $F_{ij}$ are functions of $v$, this is not reflected in the
notation.  Note that $F^{-1}(I)$ is the set of $v \in \cRbar_n$ whose
representation $\pi'_c(v)$ is in $\RC^c(S_n)$.  Also, as all
\(A(\alpha,t(v),v_i) \in U_{t(v)}\), the matrix \(F(v)\) is in
\(U_{t(v)}\) as well.  For \(g \in G_{t(v)}\), equation
\eqref{eq:UActionOnA} gives that \(F(g \cdot v) = g \cdot F(v)\).
Similarly, for \(u \in \C^\times\), equation~\eqref{eq:CactionOnA}
implies \(F(u\cdot v) = T_u F(v) T_u^{-1}\).  Thus the subset
$F^{-1}(I)$ is invariant under the actions of \(G_t\) and
\(\C^\times\).  Naively, we might expect that the resolution
\(\cX^c(S_n)\) should be the quotient of \(F^{-1}(I)\) by these
actions, but this is incorrect: the part of \(F^{-1}(I)\) where
\(t=0\) is too big, as we now explain.

\subsection{Reducible characters}
\label{subsec: reduc char}

The reducible character \(\chiep^c \in X^c(F_n)\) corresponding to
\([\epsilon] \in \{\pm 1\}^n/\pm \One\) is in \(X^c(S_n)\)
precisely when \(\prod e^{i\epsilon_j \alpha} = 1\).  We distinguish
two kinds of reducible characters: the \emph{balanced} ones are those
where \(\sum \epsilon_j = 0\) and the rest are \emph{unbalanced}.
This distinction is relevant as every balanced
\(\chiep^c \in X^c(F_n)\) is in \(X^c(S_n)\) for all values of \(c\),
but an unbalanced \(\chiep^c\) is in \(X^c(S_n)\) for at most finitely
many \(c\). As \(c\) is fixed in this section, we will write
\(\chiep\) in place of \(\chiep^c\).

Let \(\Pi= \pi \circ \sigma \maps \cRbar_n \to X^c(F_n)\) be the
composition of maps from (\ref{eq:diagram}). The preimage
\(\Pi^{-1}(\chiep)\) for \(\chiep \in X^c(F_n)\) consists of points
of the form \((0,v_1, \ldots, v_n) \in \cRbar_n\), where
\(v_j = (x_j,y_j,z_j)\) with all \(z_j = \epsilon_j\) or all
\(z_j = - \epsilon_j\).  Setting \(\omega_j= y_j - ix_j\), equation
(\ref{eq:A}) gives
\[
  A(\alpha, 0, v_j) = \begin{pmatrix}
  e^{i z_j\alpha} & \omega_j \sin \alpha \\
  0 & e^{-i z_j\alpha} \end{pmatrix}
\]
For \(\chiep\) not in $X^c(S_n)$, it follows that \(\Pi^{-1}(\chiep)\)
is disjoint from $F^{-1}(I)$.  For \(\chiep\) that are in
\(X^c(S_n)\), then on \(\Pi^{-1}(\chiep)\) we have
$F(0,v_1, \ldots, v_n) = \mysmallmatrix{1}{F_{12}}{0}{1}$, so the
equation \(F(v)=I\) is a (real) codimension 2 condition on
\(\cRbar_n^0\), specifically $F_{12} = 0$.  (While it is not clear a
priori, we show \(F^{-1}(I) \cap \Pi^{-1}(\chiep)\) is nonempty for
\(\chiep\) in \(X^c(S_n)\) whenever $n$ is even in
Proposition~\ref{prop:preimage of chi}.)  In contrast, for $t \neq 0$,
the equation $F(v) = I$ turns out to be a codimension 3 condition on
\(\cRbar_n^t\), since $\dim U_t = 3$ (see
Proposition~\ref{prop:SnIrreducible} for details).  We now turn to
resolving this discrepancy by imposing an additional condition when
$t = 0$.

\subsection{Defining the resolution}
\label{subsec:define res}


Given \(\epsilon \in \{\pm1\}^n\), let 
\[
  Z_\epsilon = \setdef{v \in \cRbar_n^0}%
  {\mbox{$v_j = (x_j, y_j, \epsilon_j)$ for all $j$}}
\]
so \(\Pi^{-1}(\chiep) = Z_{\epsilon} \cup Z_{-\epsilon}\).  For
\(\chiep \in X^c(S_n)\), we define
\[
  V_\epsilon = \setdef{v \in \cRbar_n}%
  {\mbox{$\sign(z_j) = \epsilon_j$ for all $j$ and $\mRe(F_{11}(v)) > 0$}},
\] 
so \(V_\epsilon\) is an open neighborhood of \(Z_\epsilon\) in
\(\cRbar_n\).  Observe that the set \(V_\epsilon\) is preserved by the
\(\C^\times\) action on \(\cRbar_n\).

We have  local coordinates 
\((t, \omega_1, \ldots, \omega_n)\) on \(V_\epsilon \) where 
\[
  \omega_j = y_j - i x_j \quad \text{and} \quad z_j =
  \epsilon_j \big(1-t|\omega_j|^2\big)^{1/2}.
\]
When $t = 0$, note that an $\omega \in \C^n$ gives a point in
$Z_\epsilon$, that is, some $v_i \neq \pm v_j$, if and only if
$\omega$ is not in the complex subspace generated by $\epsilon$;
hence we can view $Z_\epsilon$ as $\C^n \setminus \pair{\epsilon}$.

It is easy to see that \(F_{11}\) is a real analytic function of
\((t, \omega_1, \ldots, \omega_n)\), and, as we observed above, 
\(F_{11}(0, \omega_1, \ldots, \omega_n) \equiv 1\).  It follows that
the function \((F_{11}-1)/t\) is well-defined and analytic on
\(V_\epsilon\). 
  The function
\(f_\epsilon:V_\epsilon \to \C \times \R\) defined by 
$$f_\epsilon = \big(F_{12}, \ \mIm  (F_{11} -1)/t\big)$$
will play a key role in what follows.

We define $\cRbar^c(S_n) \subset \cRbar_n $ to be the union of
\(F^{-1}(I) \cap \cRbar_n^{t \neq 0}\) with
\(f_{\epsilon}^{-1}(0,0) \subset V_\epsilon\) for all \(\epsilon\)
with $\chiep \in X^{c}(S_n)$, and begin by showing:
\begin{lemma}
  \label{lem:intersections}
  The set $\cRbar^c(S_n)$ is contained in $F^{-1}(I)$ and
  is equal to it inside $\cRbar_n^{t \neq 0}$.
\end{lemma}

\begin{proof}
Fix $\epsilon$ with $\chiep \in X^{c}(S_n)$.  For $v \in V_\epsilon$,
we have $f_\epsilon(v) = (0, 0)$ implies $F(v) = I$, so $f^{-1}(0, 0)$
is indeed contained in $F^{-1}(I)$.  It remains to show
\(f_\epsilon^{-1}(0,0) \cap V^{t \neq 0}_\epsilon = F^{-1}(I) \cap
V^{{t \neq 0}}_\epsilon\).  For \(v \in V^{t \neq 0}_\epsilon\), we
have \(f_\epsilon(v) = (0,0)\) if and only if \(F_{11}(v)\in \R_{>0}\)
and \(F_{12}(v)=0\). Since \(F(v) \in U_{t(v)}\), these conditions
hold if and only if \(F(v) = I\).
\end{proof}

The key to proving Theorem~\ref{thm:cXSn} will be:
\begin{proposition}
  \label{prop:RepV for S_n}
  The subset \(\cRbar^c(S_n)\) is a closed codimension
  3 submanifold of \(\cRbar_n\).  It is invariant under the
  \(\C^\times\) action on \(\cRbar_n\) and each \(G_t\) action on
  \(\cRbar_n^t\).
\end{proposition}
The proof divides into two cases depending on whether $t$ is zero and
takes the bulk of the rest of this section.  If
\(v \in \cRbar^c(S_n)\) has \(t(v) \neq 0\), then $\Pi(v)$ is in
$X^{\cirr}(S_n)$.  Following Lin and Heusener's work in the case of
$\SU$\hyp representations \cite{Lin1992,Heusener2003}, we will show in
Proposition~\ref{prop:SnIrreducible} that \(F^{-1}(I)\) is a
submanifold near \(v\).  Otherwise, if $v \in \cRbar^c(S_n)$ has
\(t(v)=0\), then \(v \in Z_\epsilon\) for some \(\epsilon\). In this
case, we will show in Proposition~\ref{prop:SnReducible} that
\(f_\epsilon\) is a submersion at \(v\).

\subsection{\(f_\epsilon\) is a submersion}
Building towards the proof of Proposition~\ref{prop:SnReducible}, we
begin by describing the behavior of \(f_\epsilon\) under the
\(\C^\times\) action.

\begin{lemma}
  \label{lem:Cstarf}
  If \(f_\epsilon(v) = (a,b)\), then
  \(f_\epsilon(u \cdot v) = (ua,|u|^2 b)\) for all $u \in \C^\times$.
\end{lemma}

\begin{proof}
  Recall that $u \cdot (t, v_1, \ldots, v_n) = \left(t/\abs{u}^2, u \cdot
  v_1, \ldots, u \cdot v_n\right)$.  Using (\ref{eq:CactionOnA}), we compute
\begin{align*}
  F\big(u \cdot(t, v_1, \ldots, v_n)\big) &=
      \prod_{i=1}^n A\big(\alpha, t/\abs{u}^2, u \cdot v_i\big)
      = \prod_{i=1}^n T_u A(\alpha, t,v_i) T_u^{-1} \\
      & =  \begin{pmatrix}u & 0 \\ 0 & 1 \end{pmatrix} 
           \begin{pmatrix} F_{11}(v) & F_{12}(v) \\ F_{21}(v) & F_{22}(v)
           \end{pmatrix}  
            \begin{pmatrix}u^{-1} & 0 \\ 0 & 1 \end{pmatrix}
      = \begin{pmatrix} F_{11}(v) & u F_{12}(v) \\  u^{-1}F_{21}(v)&
        F_{22}(v) \end{pmatrix}.
\end{align*}
So
\(f_\epsilon (u \cdot v) = \big(uF_{12}(v), \ \mIm(F_{11}(v)-1)\big/(|u|^{-2}
t)\big)\) as needed. 
\end{proof}

When \(t=0\), we can explicitly describe the function
\(f_\epsilon\). Recall that we have local coordinates
\((t, \omega):= (t,\omega_1,\ldots,\omega_n)\) on \(V_\epsilon\).
\begin{lemma}
  \label{lem: L and B are nice}
  There is a linear function \(L:\C^n \to \C\) and a  nondegenerate
  sesquilinear form  \(B:\C^n \times \C^n \to \C \) so that 
  \(f_\epsilon(0,\omega) =
  \big(L(\omega), \ \mIm B(\omega, \omega)\big)\).
\end{lemma}
We could guess that \(L\) is linear and \(B\) is sesquilinear from
Lemma~\ref{lem:Cstarf}, but we need the computation below to check
that \(B\) is nonsingular.
\begin{proof}
We  expand
\(A(\alpha,t,v_j)\) as a power series in \(t\):
\begin{align*}
  A(\alpha,t,v_j) & = \begin{pmatrix} \cos \alpha + i z_j \sin \alpha  & \omega_j \sin \alpha \\
			-t \omegabar_j \sin \alpha & \cos \alpha - i
                        z_j \sin \alpha \end{pmatrix} \\
                  & = 
                    \begin{pmatrix} \zeta_j & \omega_j \sin \alpha \\
                      0 & \zeta^{-1}_j \end{pmatrix} + 
                     t \sin \alpha \begin{pmatrix} -\frac{i \epsilon_j}{2} |\omega_j|^2 & 0 \\
			- \omegabar_j &  \frac{i\epsilon_j}{2} |\omega_j|^2  \end{pmatrix} + \cdots
\end{align*}
where \(\zeta_j = e^{i \epsilon_j  \alpha}\). 
Expanding   \(F \) to first order in \(t\), we find that 
\[
  F(t, \omega_1, \ldots, \omega_n) =  \prod_{j=1}^n  A(\alpha, t,
  \omega_j) = A_0 + (t \sin \alpha) A_1 + \ldots
\]
where 
\[
  A_0 = \begin{pmatrix} \xi &  \sin \alpha  \sum c_j \omega_j \\ 0 &
    \xi^{-1} \end{pmatrix}  \quad  \quad \xi = \prod_{l=1}^n \zeta_l
  \quad \text{and} \quad c_j = \prod_{1\leq l<j} \zeta_l \prod_{j < l
    \leq n} \zeta_l^{-1}.
\]
The hypothesis that \(\chiep \in X^c(S_n)\) implies that \(\xi=1\) and
so $c_j = \zeta_j \prod_{1 \leq l < j} \zeta_l^2$.
Taking \(L(\omega) = (\sin \alpha) \sum_{j=1}^n c_j \omega_j\),
we see that the first component of \(f\) has the desired form.
	
We further compute that the upper-left entry of the matrix
\((\sin \alpha) A_1\), which is $(F_{11} - 1)/t$ at $(0, \omega)$, is
given by
$$
B(\omega,\omega)=\sin \alpha \sum_{j=1}^n \frac{-i \epsilon_j}{2}
\zeta_j^{-1} |\omega_j|^2- \sin ^2 \alpha \sum_{k<j} b_{jk}
\omegabar_j\omega_k $$ where
\(b_{jk} = \displaystyle \prod_{1\leq l<k} \zeta_l \prod_{k<l<j}
\zeta_l ^{-1} \prod_{j<l\leq n} \zeta_l = \zeta_j \zeta_k \prod_{k
  \leq l \leq j}\zeta_l^{-2}\).  The form \(B\) is nonsingular, since
its corresponding matrix is lower-triangular with nonzero diagonal
entries.
\end{proof}

\begin{lemma}
  \label{lem:Hermitian}
  There is a Hermitian form $H$ on $\ker L$ such that  \(B(\omega,
  \omega) = i H(\omega,\omega)\) for \(\omega \in \ker L\).
\end{lemma}

\begin{proof}
From the proof of the last lemma, we have
\begin{equation}
  \label{eq:B as partial}
  B(\omega,\omega) = \lim_{t \to 0} \frac{F_{11}(t,\omega) -1}{t}
  = \derivat{\frac{\partial F_{11}}{\partial t}}{(0,\omega)}.
\end{equation}
For a fixed value of \(\omega\), consider the path
\(\gamma(t) = (t,F(t,\omega))\) in \(\R \times \SLC \). The image of
\(\gamma\) is contained in the submanifold
$\Util = \setdef{(t,g)}{g \in U_t} \subset \R \times \SLC$.
It is easy to see that
$T\Util|_{(0,I)} = (0,\mathfrak{u}_0) \oplus
\langle(1,\mathbf{0})\rangle.$
	
If \(\omega \in \ker L\), then \(F(0,\omega) = I\), so
$\gamma'(0) = \big(1,\frac{\partial F}{\partial
  t}\big|_{(0,\omega)}\big) \in T\Util\big|_{(0,I)}. $ Hence
\( \frac{\partial F}{\partial t}\big|_{(0,\omega)} \in
\mathfrak{u}_0\), which implies that
\( B(\omega, \omega) = \frac{\partial F_{11}}{\partial
  t}\big|_{(0,\omega)}\) is purely imaginary for all
\(\omega \in \ker L\). In other words, \(H := -i B|_{\ker L}\) is a
sesquilinear form with the property that \(H(\omega,\omega)\) is real
for all \(\omega \in \ker L\). Such a form must be Hermitian, proving
the lemma.
\end{proof}

For \(a \in \C\), the matrix
\( \tau_a = \mysmallmatrix{1}{ia/2}{0}{1} \in U_0\) acts on
\(Q_0 = \E^2\) by translation, with an easy computation showing that
$\psi_0(\tau_a) \cdot (0,\omega) = (0,\omega + a \epsilon)$ in our
local coordinates.  We next show that $L$ and $H$ are translation
invariant.

\begin{lemma}  
  \label{lem:Translationf}
  For all $a \in \C$ and $\omega \in \C^n$, one has
  \(L(\omega+a \epsilon) = L(\omega)\).  If \(\omega \in \ker L\),
  then
  \(H(\omega + a \epsilon, \omega + a \epsilon) = H(\omega, \omega)\).
\end{lemma}

\begin{proof}

Recall from Section~\ref{subsec: intro to punc sphere}
that
$F(g\cdot v) = g \cdot F(v)$ for \(g \in U_{t(v)}\).  Taking \(t=0\)
and \(g = \psi_0(\tau_a)\), we see that
\(F(0,\omega+a \epsilon) = \tau_a F(0,\omega) \tau_a^{-1}\).  As
\(F(0,\omega) = \mysmallmatrix{1}{F_{12}}{0}{1}\), it commutes with
\(\tau_a\). Thus \(F(0,\omega+a \epsilon) = F(0,\omega)\), giving the
first relation, since \(L(\omega) = F_{12}(0,\omega)\).
	
For the second relation, suppose that \(\omega \in \ker L\), so
\(F(0,\omega)=I\). By (\ref{eq:B as partial}), we must show that
\(\frac{\partial F_{11}}{\partial t}\) takes the same value at
$(0, \omega)$ and $(0, \omega + a \epsilon)$.  Since
\(F_{11}(0,\eta)\equiv 1\), the tangent space \(T\cRbar_n^0\) is
contained in \(\ker(dF_{11})\). Hence if
\(\gamma \maps \R \to \R \oplus \C^n\) is given by
\(\gamma(t) = (t,\omega(t))\), where \(\omega(0) = \omega\), then
\begin{equation} \label{eqn:partial=total}
  \derivat{\frac{\partial F_{11}}{\partial t}}{(0,\omega)} =
  \derivat{\frac{d}{dt} F_{11}\left(\gamma(t)\right)}{t=0}.
\end{equation}	
Choose a smooth family \(\tau_a(t) \in U_t\) satisfying
\(\tau_a(0) = \tau_a \) and consider the path
\(\gamma_a(t) = \psi_t\big(\tau_a(t)\big) \cdot (t,\omega) \).  Then
\begin{equation}
  \label{eq:deriv calc}
  \begin{split}
    \derivat{\dbydt F\left(\gamma_a(t)\right)}{t=0}
    & = \derivat{\dbydt F\left(\psi_t(\tau_a(t)) \cdot (t,\omega )\right)}{t=0} \\
    & = \derivat{\dbydt \tau_a(t) F(t,\omega ) \tau_a(t)^{-1}}{t=0}  \\
    & = \tau_a'(0)F(0,\omega) \tau_a(0)^{-1} -
    \tau_a(0) F(0,\omega) \tau_a(0)^{-1} \tau_a'(0)  \tau_a(0)^{-1} \\ 
    & \qquad + \tau_a(0)\left(\derivat{\dbydt F(t,\omega)}{t=0}\right) \tau_a(0)^{-1}  \\
    & =  \tau_a\left(\derivat{\frac{\partial F}{\partial t}}{(0,\omega)}
    \right) \tau_a^{-1}
  \end{split}
\end{equation}
where the last equality uses the fact that \(F(0,\omega) = I\). 
In the proof of Lemma~\ref{lem:Hermitian}, we saw that
\(\frac{\partial F}{\partial t}\big|_{(0,\omega)} \in \mathfrak{u}_0\)
for \(\omega \in \ker L\); in particular, it is upper
triangular. Hence the top-left entry of the matrix in the final line
of (\ref{eq:deriv calc}) is given by
\( \frac{\partial F_{11}}{\partial t}\big|_{(0,\omega)} \). By
equation~\eqref{eqn:partial=total}, the top-left entry in the initial
line of (\ref{eq:deriv calc}) is
\( \frac{\partial F_{11}}{\partial t}\big|_{(0,\omega+a\epsilon)}
\). Thus (\ref{eq:deriv calc}) gives the desired equality
\(H(\omega, \omega) = H(\omega + a \epsilon, \omega + a \epsilon)\).
\end{proof}

\begin{lemma}
  \label{Lem:NullSpace}
  The nullspace
  \(N = \setdef{\eta \in \ker L}{\mbox{$H(\eta, \omega) = 0$ for all
      $\omega \in \ker L$}}\) of $H$ is \(\langle \epsilon \rangle\).
\end{lemma}

\begin{proof}
Lemma \ref{lem:Translationf} gives $L(\epsilon) = L(0) = 0$ and so
\(\epsilon \in \ker L\). Moreover, that lemma implies
$H(\epsilon, \epsilon) = H(0, 0) = 0$ and also that 
\[
  H(\omega, \omega) = H(\omega + a \epsilon, \omega + a \epsilon) 
  = H(\omega, \omega) + 2 \mRe\left( a H( \epsilon, \omega)\right)
\]
for all $\omega \in \ker L$ and $a \in \C$.  Cancelling the
$H(\omega, \omega)$ terms, we get that
$\mRe\left( a H( \epsilon, \omega)\right) = 0$.  As this holds for all
$a \in \C$, we must have $H( \epsilon, \omega) = 0$ for each
$\omega \in \ker L$ and so $\epsilon \in N$.  For the converse, note
that \(H\) is the restriction of \(-i B\) to \(\ker L\). Since \(B\)
is nonsingular and \(\ker L\) has codimension 1 in \(\C^n\), the
nullspace \(N\) is at most \1-dimensional, proving
$N = \pair{\epsilon}$.
\end{proof}

\begin{proposition}
  \label{prop:SnReducible}
  If \(v \in Z_\epsilon\) satisfies \(f_\epsilon(v) = (0,0) \), then
  \(f_\epsilon|_{Z_\epsilon}\) is a submersion at \(v\).
\end{proposition}

\begin{proof}
Suppose \(v = (0, \omega_0)\) in local coordinates. We will identify
\(T_vZ_\epsilon\) with \(\C^n\).  As $f_\epsilon(v) = (0, 0)$, we have
$L(\omega_0) = 0$.  For all \(\omega \in \ker L\),
Lemma~\ref{lem:Hermitian} gives 
\(f_\epsilon(0, \omega) = (0, H(\omega, \omega))\).  Thus, for 
\(\eta \in \ker L\), we have
\[
  \derivat{\df_\epsilon}{v}(\eta)
  = \left(0,H(\omega_0,\eta)+H(\eta,\omega_0)\right)
  = \left(0, 2\mRe H(\omega_0,\eta)\right).
\]
So long as \(\omega_0\) is not in \(N\), there will be some
\(\eta \in \ker L\) for which $\df_\epsilon\big|_v(\eta) = (0,1)$.
Lemma \ref{Lem:NullSpace} gives that \(\omega_0 \in N\) if and only if
\(\omega_0 = a\epsilon\), which happens if and only if \(v_i=\pm v_j\)
for all \(i\) and \(j\). Since \(v \in Z_\epsilon \subset \cRbar_n\),
the latter is not the case, so \((0,1) \in \im \df_\epsilon |_v\). By
considering tangent vectors \(\eta \not \in \ker L\), we see that
\(\im \df_\epsilon |_v\) contains vectors of the form \((a,b)\) with
\(a \neq 0\). We conclude that \(f_\epsilon|_{Z_\epsilon}\) is a
submersion.
\end{proof}

We set \(M_\epsilon = f_\epsilon^{-1}(0,0) \cap Z_\epsilon \), which is a
smooth submanifold by Proposition~\ref{prop:SnReducible}.

\begin{lemma}
  \label{lem:M epsilon}
  The submanifold $M_\epsilon$ is in the closure of \ 
  $\cRbar^c(S_n) \setminus \cRbar_n^{0}$.  The \(\C^\times\)
  action leaves \(M_\epsilon\) invariant.  For $g \in G_0$, we have
  \(g \cdot M_\epsilon = M_\epsilon\) when \(g \) is in the identity
  component of $G_0$ and \(g \cdot M_\epsilon = M_{-\epsilon}\)
  otherwise.
\end{lemma}

\begin{proof}
Take $K_\epsilon$ to be the full $f_\epsilon^{-1}(0,0)$
in $V_\epsilon$, so that $M_\epsilon = K_\epsilon \cap Z_\epsilon$.
By Proposition~\ref{prop:SnReducible}, we have $K_\epsilon$ is a
smooth submanifold of $V_\epsilon$ near $M_\epsilon$, and moreover the
codimension of $K_\epsilon$ in $V_\epsilon$ is the same as that of
$M_\epsilon$ in $Z_\epsilon$, which is three. It follows that
$M_\epsilon$ is a codimension one submanifold of $K_\epsilon$, and
hence the closure of $K_\epsilon \setminus M_\epsilon$ includes
$M_\epsilon$.  This proves the first claim.

The statement for the $\C^\times$ action follows immediately from
Lemma~\ref{lem:Cstarf}. For the rest, first recall that the identity
component of \(G_0\) is isomorphic to \(\Isom^+(\E^2)\), and so is
generated by rotations about the origin, which are part of the
\(\C^\times\) action, and translations, corresponding to
\(\psi_0(\tau_a)\) for \(a \in \C\).  As Lemma~\ref{lem:Translationf}
implies that \(\psi_0(\tau_a) \cdot M_\epsilon = M_\epsilon\), it
follows \(g \cdot M_\epsilon = M_\epsilon\) for all \(g\) in the
identity component of \(G_0\) as claimed.

Knowing this, it suffices to check
\(g \cdot M_\epsilon = M_{-\epsilon}\) for a single \(g\) in the other
component of \(G_0\). We use the element \(C\) from
Section~\ref{sec: props Gt and Ut}. From the relation
\(F(C\cdot v) = C \cdot F(v)\), we see that if \(v \in Z_\epsilon\)
and \(f_\epsilon (v) = (a,b)\), then \(C \cdot v \in Z_{-\epsilon}\)
and \(f_{-\epsilon}(C \cdot v) = (\abar, - b)\). So
\( C \cdot M_{\epsilon} = M_{-\epsilon}\) as desired.
\end{proof}

Let \(\Mhe\) be the quotient of \(M_\epsilon \cup M_{-\epsilon}\) by
the simultaneous actions of \(\R_{>0}\) and \(G_0\), so
$\Mhe \subset \cY$.  Equivalently, it is the quotient of
\(M_\epsilon = \cRbar^c(S_n) \cap Z_\epsilon\) by the simultaneous
actions of \(\R_{>0}\) and the identity component of \(G_0\). Note
that \(\Mhe\) is contained in $\pi^{-1}(\chiep)$ which is
diffeomorphic to \( P^{n-2}(\C)\) by Lemma~\ref{lem:CnE is Pn-2}.  We
will later define $\cX^c(S_n)$ so that
\(\cX^c(S_n) \cap \pi^{-1}(\chiep) = \Mhe\). We conclude this
subsection by giving an explicit description of \(\Mhe\):

\begin{proposition}
  \label{prop:preimage of chi}
  Suppose $\chiep \in X^c(S_n)$.  There are linear subspaces
  \(V_1 \subset V_2 \subset \C^{n-1}\), where \(V_2\) has dimension
  \(n-2\), such that \(\Mhe\) can be identified with the boundary of a
  tubular neighborhood of the projective space \(P(V_1)\) in the
  larger projective space \(P(V_2)\), which itself sits in
  \(P^{n-2}(\C)\).  In particular, $\Mhe$ is compact.  Moreover, the
  manifold $\Mhe$ is nonempty, except possibly when $n$ is odd and
  $\epsilon = \pm\One$.  When \(\chiep\) is balanced (so \(n=2m\)),
  then \( \dim V_1 = m-1\).
\end{proposition}

\begin{proof}
Recall from Section~\ref{subsec:define res} that in our local
coordinates $Z_\epsilon = \C^n \setminus \pair{\epsilon}$.  The
submanifold
\(M_\epsilon= \setdef{\omega \in \C^n \setminus
  \pair{\epsilon}}{\mbox{$L(\omega) = 0$ and $H(\omega,\omega)=0$}}\)
is invariant under translation by any element of \(\pair{\epsilon}\).
We identify \(\C^{n-1}\) with \(\C^n/\langle \epsilon \rangle
\). Taking this quotient has the same effect as normalizing so that
the first point $\omega_1$ is $0$, which absorbs all the $G_0$ action
except the part contained in $\C^\times$.  Setting
\(V_2 = \ker L / \langle \epsilon \rangle \), the form \(H\) descends
to a Hermitian form \(H'\) on \(V_2\) which is nondegenerate by
Lemma~\ref{Lem:NullSpace}. Then \(\Mhe \) is the quotient of
\(M_\epsilon /\! \pair{\epsilon} = \setdef{\omega \in V_2 \setminus
  \{0\}}{H'(\omega,\omega) =0}\) by the action of \(\C^\times\).

With respect to an appropriate basis of \(V_2\), the form \(H'\) is
\[
  H' (\eta, \eta) = - |\eta_1|^2 - \cdots - |\eta_k|^2 +
  |\eta_{k+1}|^2 + \cdots + |\eta_{n-2}|^2
\]
where \(k\) is the index of \(H'\). Take
\(V_1 = \{(\eta_1,\eta_2, \ldots, \eta_k,0,\ldots,0)\}\); we leave it
as an exercise to see that \(\Mhe\) is the boundary of a tubular
neighborhood of \(P(V_1)\) in \(P(V_2)\).

Next, the manifold $\Mhe$ is nonempty provided $V_1$ is neither the
trivial subspace nor all of $V_2$; equivalently, we need to show $H'$
is not definite, i.e.~there exists a nonzero $\eta \in V_2$ with
$H'(\eta, \eta) = 0$.  (Recall from the Section~\ref{subsec: intro to
  punc sphere} that $n \geq 4$ so $\dim V_2 \geq 2$.) First, suppose
$[\epsilon]$ is not $[\One]$. By cyclically permuting the generators
$s_i$, we can assume $\epsilon_1 = 1$ and $\epsilon_2 = -1$.  In the
formulae used to prove Lemma~\ref{lem: L and B are nice}, we have
that $c_1 = c_2 = e^{i \alpha}$ and $b_{21} = 1$, and hence that
$\etatil = (1, -1, 0, \ldots 0)$ is in $\ker L$ with
$H(\etatil, \etatil) = 0$.  As $\etatil \notin \pair{\epsilon}$, this
gives the desired element $\eta \in V_2$ showing that $H'$ is
indefinite.

Second, suppose $\epsilon = \One$ and set $\zeta = e^{i \alpha}$.
Then each $\zeta_l = \zeta$, and so $\zeta^n = 1$ with
$\zeta \neq \pm 1$ since the trace $c \neq \pm 2$.  We then calculate
$c_j = \zeta^{2 j - 1}$ and $b_{jk} = \zeta^{2(k - j)}$.  So the
vectors $\nutil_j = \zeta^2 e_j - e_{j + 1}$ for $1 \leq j \leq n - 2$
together with $\epsilon$ form a basis for $\ker L$.  Setting
$H'' = (-2/\sin \alpha) H'$, consider the associated Hermitian matrix
$W$ with respect to the induced basis $\nu_1,\ldots,\nu_{n-2}$ of
$V_2$, that is, $H''(v, w) = w^* W v$.  A calculation finds the only
nonzero entries of $W$ are $ \zeta + \zeta^{-1} = 2 \cos \alpha$ along
the diagonal, $-\zeta$ immediately above the diagonal, and
$-\zeta^{-1}$ immediately below.  Using Lemma~\ref{lem: q mat det}
below with $q = \zeta$ and $m = n - 2$, we have
\[
  \det W =   \frac{\zeta^{n - 1} - \zeta^{-(n - 1)}}{\zeta - \zeta^{-1}}
  = \frac{\zeta^{-1} - \zeta}{\zeta - \zeta^{-1}} = -1 \mtext{as $\zeta^n
    = 1$ and $\zeta \neq \pm 1$.}
\]
As  $\det W < 0$,  the form
$H''$ must be indefinite when \(n\) is even.  Thus we have proved $\Mhe$ is nonempty for all
$\chiep \in X^c(S_n)$, except possibly when $n$ is odd and
$\epsilon = \pm \One$.
  
For the last claim, when \(\chiep\) is balanced, it is part of a
continuous family of reducible characters obtained by varying
\(\alpha\). The corresponding Hermitian forms are all nondegenerate
and form a continuous family, and so they have the same index. We
compute this index in the limit as \(\alpha \to 0\). From the formulae
in the proof of Lemma~\ref{lem: L and B are nice}, we have
\(\lim_{\alpha \to 0}(L/\sin \alpha) = L_0 \), where
\(L_0(w) = \sum_j \omega_j\), and
\(\lim_{\alpha \to 0}(B/\sin \alpha) = B_0\), where
\(B_0(w) = \sum_j -i\epsilon_j |\omega_j|^2/2\). Then
\(H_0 := -i B_0\) is a Hermitian form of index \(m = n/2\). Since
\(\langle\epsilon \rangle \) is contained in the nullspace of
\(H_0|_{\ker L_0}\), the form on \(\ker L_0/\langle \epsilon \rangle\)
induced by \(H_0\) has index \(m-1\) as needed.
\end{proof}

\begin{remark}
  \label{rem: annoying}
  When $n$ is odd, the form $H''$ with matrix $W$ can be (negative)
  definite.  Setting $\zeta_n = e^{2 \pi i/n}$, then for $n = 5$ the
  form is definite when $\zeta$ is either $\zeta_5^2$ or $\zeta_5^3$.
  More generally, experiment strongly suggests that $W$ is definite
  for $n$ odd if and only if
  $c = \zeta_n^{(n - 1)/2} + \zeta_n^{(n+1)/2} = \cos \frac{\pi (n -
    1)}{n}$.  Indeed, regardless of the parity of $n$, it seems that
  the number of negative eigenvalues of $H''$ for $\zeta_n^k$ with
  $k < n/2$ is exactly $2 k - 1$.
\end{remark}

Finally, here is the lemma needed for the proof of the previous
proposition:
\begin{lemma}
  \label{lem: q mat det}
  Let $J_m$ be the $m \times m$ matrix with entries in $\Z[q^{\pm 1}]$
  whose nonzero entries are $q + q^{-1}$ along the diagonal, $-q$
  immediately above the diagonal, and $-q^{-1}$ immediately below.
  Then
  $\det J_m = q^m + q^{m - 2} + \cdots + q^{-(m - 2)} + q^{-m} =
  \left(q^{m + 1} - q^{-(m + 1)}\right)\big/\left(q - q^{-1}\right)$.
\end{lemma}

\begin{proof}
We induct on $m$, with the base cases of $m = 1$ and $m = 2$ being easy
checks. In the illustrative case of $m = 4$, setting $\qbar =
q^{-1}$ we have:
\[
  J_4 =  \left(\begin{array}{cccc}
                 q + \qbar & -q & 0 & 0 \\
                 -\qbar & q + \qbar & -q & 0 \\
                 0 & -\qbar & q + \qbar & -q \\
                 0 & 0 & -\qbar & q + \qbar\\
               \end{array}\right)
\]
Thus, expanding on the first row and using the known values of
$\det(J_l)$ we get
\[
  \det J_m = \left(q + \qbar\right) \det(J_{m - 1}) - \det(J_{m - 2})
  =
  \frac{q^{m + 1} - \qbar\vphantom{q}^{m + 1}}{q - \qbar}
\]
as needed to complete the induction and hence the proof.
\end{proof}
 
\subsection{Irreducible characters}
\label{sec:irrchars}

Our goal in this subsection is to show that \(F^{-1}(I)\) meets
\(\cRbar_n^{t \neq 0}\) as a smooth codimension 3 submanifold.  The
argument follows the one given in \cite{Heusener2003}, but with the
group \(U_t\) in place of \(SU(2)\).  We start by briefly recalling
some properties of \(U_t\) from Section~\ref{sec:SL2Csubs}. Its Lie
algebra \(\ut\) consists of matrices of the form
\(\mysmallmatrix{ix}{w}{-t\wbar}{-ix}\), where \(x \in \R\) and
\(w \in \C\). Both $U_t$ and $\ut$ are contained in
\(V_t= \setdef{\mysmallmatrix{z}{w}{ -t \wbar} {\zbar}}%
{ z, w \in \C}\), which we regard as a vector space over $\R$ rather
than $\C$.  There is a symmetric bilinear form on \(V_t\) given by
\(\pair{u_1, u_2} = \tr( u_1 u_2)\); it is nondegenerate for
\(t \neq 0\). Note \(\ut = I^\perp\) with respect to this form, and
the orthogonal projection \(p:V_t \to \ut\) takes the tracefree part,
with $p(C) = C - \frac{\tr C}{2} I$.  The restriction of
\(\pair{\cdotspaced , \cdotspaced}\) to \(\ut\) is a scalar multiple
of our standard bilinear form from Section~\ref{subsec:B_t}, and hence
a scalar multiple of the Killing form.

Fix $t \neq 0$. We will show that the restriction of $F$ to
$\cRbar_n^t$ is a submersion.  This map can
be factored as a composition
\[
  \cRbar_n^t \xrightarrow{\bA} (U_t)^n \xrightarrow{m} U_t \mtext {where}
  \bA(t,v_1,\ldots,v_n) = \big(A(\alpha, t,v_1), \ldots, A(\alpha,t,v_n)\big)
\]
and \(m\) is the multiplication map.  Recall that $\cRbar_n^t$ is
defined by the open condition that $v_i \neq \pm v_j$ for some
$i \neq j$, and therefore each vector $v_k$ in a given
$v \in \cRbar_n^t$ can be moved independently in $\Qbar_t$ on a
small scale. Hence \(\im \, \bA\) will be an open subset of \(\Utcn\),
where \(U_t^c = \setdef{g \in U_t}{\tr g = c}\).  For \(g \in U_t\),
we identify \(T_gU_t\) with \(\ut\) via right multiplication, i.e.~the
map which sends \(u \in \ut\) to \(ug \in T_gU_t\).  
\begin{lemma}
  \label{Lem:Tangent}
  For $t \neq 0$ and any $c \in (-2, 2)$, the subset $U_t^c$ is a
  smooth submanifold of $U_t$, with \(T_gU_t^c = p(g)^\perp\) for
  each $g \in U_t^c$.
\end{lemma}

\begin{proof}
Identifying \(T_g U_t\) with \(\ut\) as above, for $u \in \ut$ we compute
\[
  d\tr_g (u) = \left. \frac{d}{\ds}\right\vert_{s = 0}\tr(e^{s u} g) =
  \tr\left( \left. \frac{d}{\ds}\right\vert_{s = 0} e^{s u}g \right)
  = \tr(u g) = \pair{u, g} = \pair{u, p(g)}.
\]
Since $g \neq \pm I$, we have $p(g) \neq 0$. As $t \neq 0$, the
bilinear form on $\ut$ is nondegenerate, so there is $u \in \ut$ with
$\pair{u, p(g)} \neq 0$.  Thus $\tr \maps U_t \to \R$ is a submersion
at $g$ and $T_g U_t^c = \ker( d\tr_g ) = p(g)^\perp$ as needed.
\end{proof}

\begin{lemma}
\label{Lem:AdAction}
Suppose \(\bg = (g_1,\ldots, g_n) \in \Utcn\).  If
$p(g_i)$ and $p(g_j)$ are linearly independent for some
\(i,j\), the map \(dm_\bg: T_\bg \Utcn \to \ut\) is surjective.
\end{lemma}
The proof is essentially the same as that of \cite[Lemma
3.1]{Heusener2003}, but we give it here for completeness.

\begin{proof}
If we identify \(T_\bg U_t^n\) with \(\ut^n\) via the map which sends
\((u_1,\ldots, u_n) \in \ut^n\) to
\((u_1g_1, \ldots, u_n g_n) \in T_\bg U_t^n\), the derivative
\(dm: \ut^n \to \ut \) at $\bg$ is given by
\[
  dm (u_1, \ldots, u_n)
  = \sum_{j=1}^n g_1 \cdots g_{j-1} u_j (g_1 \cdots g_{j-1})^{-1} 
  = \sum_{j=1}^n \Ad(g_1\cdots g_{j-1}) \cdot u_j.
\]
By Lemma~\ref{Lem:Tangent}, we see that \( dm \big(T_\bg \Utcn\big) \)
is spanned by the sum of the subspaces
\(\Ad(g_1 \cdots g_{j-1}) \left(p(g_j)^\perp\right)\).  Fix an index
\(j\) where \(p(g_j)\) and \(p(g_{j+1})\) are linearly independent, and
set $h = g_1 \cdots g_{j-1}$.  To see $dm$ is onto, as $\dim \ut = 3$,
it suffices to show that the 2-dimensional subspaces
\[
  \Ad(h)\left(p(g_j)^\perp\right) \mtext{and} \Ad(h g_j)\left(p(g_{j+1})^\perp\right)
\]
of $\ut$ are distinct.  Applying $\Ad(g_j^{-1} h^{-1})$, we
equivalently examine $\Ad(g_j^{-1})\left(p(g_j)^\perp\right)$ and
$p(g_{j+1})^\perp$.  As $\Ad$ preserves the bilinear form, we have
\[
  \Ad(g_j^{-1})\left(p(g_j)^\perp\right) =
  \left(\Ad(g_j^{-1})\left(p(g_j)\right)\right)^\perp = p(g_j)^\perp.
\]
As the form is nondegenerate with \(p(g_j)\) and \(p(g_{j+1}\)
linearly independent, it follows
\(p(g_j)^\perp \neq p(g_{j+1})^\perp\) as needed to show $m$ is a
submersion at $\bg$.
\end{proof}

We can now prove the main result of this subsection:

\begin{proposition}
  \label{prop:SnIrreducible}
  The intersection \(F^{-1}(I) \cap \cRbar_n^{t \neq 0}\) is a smooth
  closed codimension~3 submanifold of \(\cRbar_n^{t \neq 0}\).
\end{proposition}

\begin{proof}
Consider the submanifold
$\Util = \setdef{(t,g)}{g \in U_t} \subset \R \times \SLC $ , and
define a map \(\Ftil:\cRbar_n^{t \neq 0} \to \Util\) by
\(\Ftil(v)=(t(v),F(v))\). We claim that \(\Ftil\) is a submersion.
For $v \in \cRbar_n$, set \(g_i = A(\alpha, t, v_i)\).  By
Lemma~\ref{lem:exp formula}, we have $p(g_i) = (\cos \alpha) v_i$.  By
definition, a point \(v \in \cRbar_n\) has \(v_i \neq \pm v_j\) for
some \(j\), and hence \(p(g_i)\) and \(p(g_j)\) are linearly
independent.  Lemma~\ref{Lem:AdAction} thus shows that
\(d\Ftil(T_v\cRbar_n^t) = \{0\} \times T_g U_t \), where \(g =
F(v)\). By considering a path of the form
\(\gamma(t) = \left(t, \gamma_1(t), \ldots, \gamma_n(t)\right)\) in
\(\cRbar_n\), we see that \(\im(d\Ftil)\) contains a vector of the
form \((1,*)\). Hence \(\im(d\Ftil)\) has dimension 4 and \(\Ftil\) is
a submersion.

To conclude the proof, let \(\Itil = \setdef{(t,I)}{t \in
  \R}\). Now \(\Itil\) is a smooth closed codimension 3 submanifold of
\(\Util\), so
\(F^{-1}(I) \cap \cRbar_n^{t \neq 0} = \Ftil\vphantom{F}^{-1}(\Itil)\)
is a smooth codimension~3 submanifold of \(\cRbar_n^{t \neq 0}\).
\end{proof}

\subsection{Putting it all together} We can now prove the results
stated at the beginning of this section.

\begin{proof}[Proof of Proposition~\ref{prop:RepV for S_n}]
Recall from Section~\ref{subsec:define res} that \(\cRbar^c(S_n)\) is
the union of \(F^{-1}(I)\cap \cRbar_n^{t \neq 0}\) with all
\(f_\epsilon^{-1}(0,0)\), where the union is taken over those
\(\epsilon\) with \(\chiep \in X^c(S_n)\).  By
Lemma~\ref{lem:intersections}, we know
\(\cRbar^c(S_n) \subset F^{-1}(I)\) and
$\cRbar^c(S_n) \cap \cRbar_n^{t \neq 0} = F^{-1}(I) \cap \cRbar_n^{t
  \neq 0}$.  To show that \(\cRbar^c(S_n)\) is a codimension 3
submanifold of \(\cRbar_n\), we must check that for each
\(v \in \cRbar^c(S_n)\) there is an open subset
\(O_v \subset \cRbar_n\) containing \(v\) such that
\(\cRbar^c(S_n) \cap O_v\) is a codimension 3 submanifold of \(O_v\).
This follows from Proposition~\ref{prop:SnIrreducible} when
\(t(v)\neq 0 \) using $O_v = \cRbar_n^{t \neq 0}$ and from
Proposition~\ref{prop:SnReducible} when \(t(v) = 0\) taking $O_v$ to
be some $V_\epsilon$.

Turning now to showing that $\cRbar^c(S_n)$ is closed, the main worry
is that it could have a limit point in a component of $\cRbar_n^{0}$
corresponding to $\epsilon$ where $\chiep$ is not in $X^c(S_n)$. As
\(X^c(S_n)\) is a closed subset of \(X^c(F_n)\), for each \(\epsilon\)
with \(\chiep \not \in X^c(S_n)\), we choose an open set
\(W_\epsilon \subset X^c(F_n)\) disjoint from \(X^c(S_n)\).  Let
\(V_\epsilon = \Pi^{-1}(W_\epsilon)\), so
\(Z_\epsilon \subset V_\epsilon\) and
\(V_\epsilon \cap \cRbar^c(S_n) = \emptyset\).  For the other
$\epsilon$ with \(\chiep\) in \(X^c(S_n)\), we continue to define
\(V_\epsilon\) as in Section~\ref{subsec:define res}.  Then
\(\cRbar_n\) has an open cover consisting of \(\cRbar_n^{t \neq 0}\)
together with the \(V_\epsilon\)'s. If \(V\) is one of the open sets
in this cover, then \(\cRbar^c(S_n) \cap V\) is closed in \(V\). It
follows that \(\cRbar^c(S_n)\) is closed in \(\cRbar_n\).
	
Finally, regarding the group actions, as noted in Section~\ref{subsec:
  intro to punc sphere}, the subset $F^{-1}(I)$ is preserved by the
actions of \(\C^\times\) and the \(G_t\).  Those actions also preserve
$\cRbar_n^{t\neq 0}$, and hence also
$F^{-1}(I) \cap \cRbar_n^{t \neq 0} = \cRbar^c(S_n) \cap \cRbar_n^{t
  \neq 0}$.  The rest, namely $\cRbar^c(S_n) \cap \cRbar_n^{0}$, is
the union of the submanifolds $M_\epsilon$ from Lemma~\ref{lem:M
  epsilon}, and that lemma gives the needed invariance.
\end{proof}

\begin{proof}[Proof of Theorem~\ref{thm:cXSn}]
We define \(\cX^c(S_n) = \sigma\big(\cRbar^c(S_n)\big)\).  Combining
Lemma~\ref{lem:quotient is manifold}, Theorem~\ref{thm:varconsY}, and
Proposition~\ref{prop:RepV for S_n}, we see that \(\cX^c(S_n)\) is a
closed codimension 3 smooth submanifold of \(\cX^c(F_n)\).  To compute
\(\pi(\cX^c(S_n))\), first observe that it is
$\pi\big(\sigma\big(\cRbar^c(S_n)\big)\big)$ which is equal to
$\tau\big(\pi'_c\big(\cRbar^c(S_n)\big)\big)$ by
Lemma~\ref{lem:diagram}. By Lemma~\ref{lem:intersections}, the subset
$\cRbar^c(S_n)$ is contained in $F^{-1}(I)$, and, as discussed in
Section~\ref{subsec: intro to punc sphere}, the subset
$\pi_c'(F^{-1}(I))$ is contained in $R^c_\C(S_n)$.  Hence
$\pi(\cX^c(S_n)) \subset \tau(R^c_\C(S_n)) = X^c_\C(S_n)$.  As
$\pi(\cX^c(F_n)) = X^c(F_n)$, taking intersections yields
$\pi(\cX^c(S_n)) \subset X^c(S_n)$.  To show equality, we consider the
cases of $X^\cirr(S_n)$ and $X^\cred(S_n)$ separately.

Given $\chi \in X^\cirr(S_n)$, pick a corresponding representation
$\rho \maps \pi_1(S_n) \to U_t$ where $t = \pm 1$.  Viewing the domain
of $\rho$ as $F_n$, Lemma~\ref{lem: pi prime c} gives a unique
$v \in \cRbar_n^t$ with $\pi_c'(v) = \rho$.  As $\rho$ is in
$R^c_\C(S_n)$, we have $v \in F^{-1}(I)$ and so $v \in \cRbar^t(S_n)$.
Then $\pi(\sigma(v)) = \chi$ as needed.  If instead
$\chiep \in X^\cred(S_n)$, consider $\Mhe$ from
Proposition~\ref{prop:preimage of chi}.  By this proposition and the
discussion immediately before it, we have
$\Mhe = \cX^c(S_n) \cap \pi^{-1} (\chiep)$ and $\Mhe$ is nonempty
since we are requiring that $n$ is even. In particular, we have
$\chiep \in \pi(\cX^c(S_n))$.  So $\pi(\cX^c(S_n)) = X^c(S_n)$ as
claimed.

Next, the claim that the map \(\pi:\cX^c(S_n) \to X^c(S_n)\) is a
diffeomorphism away from \(X^{c,\red}(S_n)\) is immediate from the
analogue for $F_n$ in Theorem~\ref{thm:def of pi}.
Finally, to see that the closure of $\cX^\cirr(S_n)$ in $\cX^c(F_n)$
is $\cX^c(S_n)$, start by noting that $\cX^\cred(S_n) = \bigcup \Mhe$,
where the union is over $\chiep \in X^c(S_n)$.  The result now follows
from the first statement in Lemma~\ref{lem:M epsilon}.
\end{proof}

\subsection{Odd number of punctures}
\label{sec: odd punct}

In proving Theorem~\ref{thm:cXSn}, the only place the parity of $n$
was used was Proposition~\ref{prop:preimage of chi}, in the special
case that $\chiep \in X^c(S_n)$ for $\epsilon = \One$.  Thus
Theorem~\ref{thm:cXSn} still holds for odd $n$ provided we
exclude the finitely many $c$ of the form $2 \cos \frac{2 \pi k}{n}$
for $1 \leq k < n/2$.  For context, note that for a fixed odd $n$,
any reducible $\chiep \in X^c(S_n)$ is necessarily unbalanced, and
hence there are only finitely many $c \in (-2, 2)$ for which
$X^\cred(S_n)$ is nonempty.  Thus for most $c$, no resolution is
necessary as $X^c(S_n) = X^\cirr(S_n)$.

Moreover, given Remark~\ref{rem: annoying}, we conjecture that
Theorem~\ref{thm:cXSn} holds if we exclude the single value
$c = 2 \cos \frac{\pi(n - 1)}{n}$.  For that specific $c$, it turns
out that $X^\cirr_\SU(S_n)$ is empty (apply
\cite[Theorem~A]{Biswas1998} with $j = 0$), and experimental evidence
suggests $X^\cred(S_n) = \{\chi_{[\One]}\}$.  It seems plausible that
the corresponding $X^\cirr_\SLR(S_n)$ does not limit on
$X^\cred(S_n)$ and might even be empty; either way,
Theorem~\ref{thm:cXSn} would hold modulo the isolated point
$X^\cred(S_n)$ not being in the image of $\cX^c(S_n)$.

\section{Varying the trace}
\label{sec:varies}

Up to now, we have focused on character varieties where each of our
preferred generators have the same fixed trace \(c \in (-2, 2)\).  We
now consider the effect of varying $c$, starting with the
free group.  Recall from Section~\ref{subsec:X(F_n)} that
$\tr \maps X(F_n, S) \to \R$ is the map sending $\chi$ to $\chi(s_1)$,
or equivalently to any $\chi(s_i)$. We will show:

\begin{theorem}
  \label{thm:sX(F_n)}
  For each $n \geq 2$, there is a smooth \((2n-2)\)-dimensional
  manifold \(\cX(F_n,S)\) and a smooth map
  \(\pi\maps \cX(F_n,S) \to X(F_n,S)\) with the following properties.
  First, the map \(\mtr = \tr \circ \pi\maps \cX(F_n,S) \to \R\) is a
  submersion, with \(\mtr^{-1}(c) \cong \cX^c(F_n,S)\) for
  \(c \in (-2,2)\) and \(\mtr^{-1}(c) \cong X^{\cirr}_\SLR(F_n,S)\)
  otherwise. Second, the subset
  $\cX^\irr(F_n, S) := \pi^{-1}\big(X^\irr(F_n, S)\big)$ is dense in
  $\cX(F_n,S)$ and $\pi$ restricts to a diffeomorphism
  $\cX^\irr(F_n, S) \to X^\irr(F_n, S)$.
\end{theorem}

For the punctured sphere, the situation is more complicated. Recall
from Section~\ref{subsec: reduc char} that for balanced \(\epsilon\),
the reducible character \(\chiep^c\) is in \(X^c(S_n)\) for all \(c\),
but for unbalanced \(\epsilon\), the character \(\chiep^c\) is in
\(X^c(S_n)\) only at isolated \(c\). It turns out that we can fit the
resolutions of the balanced reducibles into a family, but not the
resolutions of the unbalanced reducibles. Our solution to this problem
is to excise the unbalanced reducibles from our moduli space. More
formally, we define
\[
  X_\bal(F_n,S) = X(F_n,S) \setminus \setdef{\chiep^c}%
  {\mbox{$\sum \epsilon_i \neq 0$}}
\]
and \(\cX_\bal(F_n,S) = \pi^{-1}\big(X_\bal(F_n,S)\big)\), so
\(\cX_\bal(F_n,S) \) is an open submanifold of
\(\cX(F_n,S)\). Likewise, we take
$X_\bal(S_n, S) = X_\bal(F_n, S) \cap X(S_n, S)$ and
$\cX_\bal(S_n,S) = \cX(S_n, S) \cap \cX_\bal(F_n,S)$, as well as 
$\cX^c_\bal(S_n,S) = \cX^c(S_n, S) \cap \cX_\bal(S_n,S)$, etc.
We will show:

\begin{theorem}
  \label{thm:sX(S_n)}
  For each even $n \geq 4$, there is a smooth closed
  codimension~\(3\) submanifold \(\cX(S_n,S)\) of \( \cX_\bal(F_n,S)\)
  satisfying
  \(X^\irr(S_n, S) \subset \pi\big(\cX(S_n,S)\big) \subset
  X_\bal(S_n,S)\) with the following properties.  The map
  \(\mtr\maps \cX(S_n,S)\to \R\) is a submersion, and
  \(\mtr^{-1}(c) = \cX^c_\bal(S_n,S)\) for \(c \in (-2,2)\); otherwise
  \(\mtr^{-1}(c) = X^\cirr_{\SLR}(S_n,S)\).  Also, the subset
  $\cX^\irr(S_n, S) := \pi^{-1}\big(X^\irr(S_n, S)\big)$ is dense in
  $\cX(S_n,S)$ and $\pi$ restricts to a diffeomorphism
  $\cX^\irr(S_n, S) \to X^\irr(S_n, S)$.
\end{theorem}
The generating set \(S\) will remain fixed throughout this section, so
we drop it from the notation, writing \(\cX(F_n)\) for \(\cX(F_n,S)\),
etc.

\subsection{The free group}
\label{sec:varying Fn}

The space \(\cX(F_n)\) will be constructed as the union of two sets
\(B\) and \(C\).  Let \(B= (-2,2) \times \cY \), and define
\(\pi_B \maps B \to X(F_n)\) by \(\pi_B(c, y) = \pi_c(y)\), where
\(\pi_c\) is the map constructed in Theorem~\ref{thm:def of pi}.  As
in Lemma~\ref{lem:diagram}, we have a commutative square
 \[
    \begin{tikzcd}
      (-2,2) \times \cRbar_n \arrow{r}{\pi'}
        \arrow{d}{\mathrm{id} \times \sigma } &
      \RC(F_n) \arrow{d}{\tau} \\
      (-2,2) \times \cY \arrow{r}{\pi_B} & X(F_n)
    \end{tikzcd}
\]
where the topmost map is defined by $\pi'(c,v) = \pi'_c(v)$, that is,
\[
  \pi'(c,v) (s_i) = A\left(\cos^{-1}(c/2), t(v), v_i\right),
\]
and so $\pi'$ is smooth by (\ref{eq:A}).  As $\tau$ is also smooth, so is
$\tau \circ \pi'$. As the lefthand arrow is a submersion, it follows
that \(\pi_B\) is smooth.

Next, we define $C = X^\irr(F_n)$, which is a smooth manifold by
Corollary~\ref{cor: real part smooth}.  Let
\(\pi_C \maps C \to X(F_n)\) be the inclusion.  We now define
\(\cX(F_n) = (B \coprod C)/ \sim\), where \(b \sim c\) if
\(\pi_B (b) =\pi_C(c)\).  We give $\cX(F_n)$ the quotient topology, so
the maps \(\pi_B\) and \(\pi_C\) combine to give a continuous map
\(\pi \maps \cX(F_n) \to X(F_n)\).  Finally, we let
\(i_B \maps B \to \cX(F_n)\) and \(i_C \maps C \to \cX(F_n)\) be the
maps induced by inclusion.  We define
\(B^\irr = \pi_B^{-1}\left(X^\irr(F_n)\right)\), which
Theorem~\ref{thm:def of pi} shows is
\[
  (-2, 2) \times \left(\cY \setminus \CnsignE\right) =
  (-2, 2) \times \left(\CnS \coprod \CnsignH \right).
\]
Next, notice that \(\im \pi_B \cap \im \pi_C \subset X^\irr(F_n) \) is
$\bigcup \setdef{X^\cirr(F_n)}{-2 < c < 2}$.  Since $\pi_C$ is
injective and $\pi_B$ is injective on
$\pi_B^{-1}(\im \pi_C) = B^\irr$, the maps \(i_B\) and \(i_C\)
are also injective. Since \(i_B\) and \(i_C\) are injective, we view \(B\) and \(C\) as
subsets of \(\cX(F_n)\).

\begin{proof}[Proof of Theorem~\ref{thm:sX(F_n)}]

We first describe the smooth structure on \(\cX(F_n)\). The pieces
\(B\) and \(C\) are both smooth manifolds, so we define
\(f:\cX(F_n) \to \R\) to be smooth if and only if \(f|_B\) and
\(f|_C\) are smooth. In order for this to make sense, we must check
that the transition function
\(i_C^{-1} \circ i_B = \pi_C^{-1} \circ \pi_B\) is a diffeomorphism
from $B^\irr$ onto its image.  As $\pi_C$ is the inclusion of the open
subset $C = X^\irr(F_n)$ into $X(F_n)$, this amounts to analyzing
$\pi_B$ restricted to $B^\irr$.  As noted, we have
$\pi_B(B^\irr) = \bigcup \setdef{X^\cirr(F_n)}{-2 < c < 2}$, and the
latter is open in $X(F_n)$.  Theorem~\ref{thm:def of pi} then gives
that $\pi_B|_{B^\irr}$ is injective. As
$\tr\big(\pi_B(c, y)\big) = c$, we have that
$d \pi_B(\frac{\partial}{\partial c})$ is not in the image of each
$T(\{c\} \times \cY)$ under $d \pi_B$, so Theorem~\ref{thm:def of pi}
also gives that $d \pi_B$ is an isomorphism at each point in $B^\irr$,
making $\pi_B|_{B^\irr}$ a diffeomorphism onto its image. We conclude
that \(\cX(F_n)\) admits a smooth structure with respect to which
\(i_B\) and \(i_C\) are open embeddings. The maps \(\pi|_B = \pi_B\)
and \(\pi|_C = \pi_C\) are smooth, so \(\pi\) is smooth.

Next, we check that \(\cX(F_n)\) is Hausdorff.  The pieces \(B\) and
\(C\) are Hausdorff open subsets which cover \(\cX(F_n)\), so it
suffices to show that if \(x \in B \setminus C\) and
\(y \in C \setminus B\), then \(x\) and \(y\) can be separated by open
sets in \(\cX(F_n)\).  Considering the map
\(\pi \maps \cX(F_n) \to X(F_n)\), we see $\pi(x) \neq \pi(y)$ since
$\pi(x)$ is reducible but $\pi(y)$ is irreducible. As \(X(F_n)\) is
Hausdorff and $\pi$ is continuous, we can separate $x$ from $y$ by the
preimages of open sets in \(X(F_n)\) separating $\pi(x)$ from
$\pi(y)$. Thus \(\cX(F_n)\) is Hausdorff.

Turning now to the map \(\mtr = \tr \circ \pi \maps \cX(F_n) \to \R\),
note that it satisfies \(\mtr(c,x) = c\) for \((c,x) \in B \) and is
given by the function $\tr$ on $C$. As \(\mtr|_B\) is thus a
submersion and \(\mtr|_C\) is also by Corollary~\ref{cor:tr is
  submersion}, we have that $\mtr$ is a submersion.  We leave the
statements about \(\mtr^{-1}(c)\) to you, but point out that for $c$
outside $(-2, 2)$ one has $X^\cirr(F_n, S) = X^\cirr_\SLR(F_n, S)$
since $X^\cirr_\SU(F_n, S)$ is empty by Lemma~\ref{lem: empty char
  vars}. Finally, $\pi$ maps $\cX^\irr(F_n) := C$ diffeomorphically
onto $X^\irr(F_n)$ by construction, and the density of $\cX^\irr(F_n)$
in $\cX(F_n)$ follows from the density of $\cX^\cirr(F_n)$ in
$\cX^c(F_n)$ for $c \in (-2, 2)$ which was noted after (\ref{eq:
  cXirr}).
\end{proof}

\subsection{Properness of projections}
\label{sec: prop of proj}

We now consider the smooth map \(\pi \maps \cX(F_n) \to X(F_n)\).  Its
image is
\[
  X_\circ(F_n):= \setdef{\chi \in X(F_n)}%
  {\chi \in X^{\irr}(F_n) \ \text{or} \ \tr(\chi) \in (-2, 2)}.
\]
We let \(\pi_\circ \maps \cX(F_n) \to X_\circ(F_n)\) be the map
obtained by restricting the range of \(\pi\).

\begin{proposition}
  \label{prop:proper 1}
  The map \(\pi_\circ:\cX(F_n) \to X_\circ(F_n)\) is proper.
\end{proposition}
The proof makes use of the following:
\begin{lemma}
  \label{lem:proper union}
  Suppose \(f \maps \cW \to W\) is continuous, where
  \(W = U_0 \cup U_1\) is a metric space and the \(U_i\) are open. If
  both restrictions \(\fhat_i \maps {f^{-1}(U_i)}\to U_i\) are proper,
  then so is \(f\).
\end{lemma}
  
\begin{proof} 
Let \(C_i= W-U_i\), so \(C_i\) is closed and
\(C_0\cap C_1 = \emptyset\).  We consider the maps
\(d_i\maps W \to [0, \infty) \) given by
\(d_i(x) = d(x,C_i) = \inf \setdef{d(x,y)}{y \in C_i}\), which are
continuous. As \(C_i\) is closed in \(W\), one has \(d_i(x) = 0 \) if
and only if \(x \in C_i\). Since the $C_i$ are disjoint,
\(d_0(x)+d_1(x)>0\) for all \(x \in W\). Hence if we define
\[g(x) = \frac{d_0(x)}{d_0(x)+d_1(x)},\] we have \(0\leq g(x)\leq 1\).
For each $i$, it follows that \(x \in C_i\) if and only if
\(g(x) = i\).
 
Now let \(K \subset W\) be compact, and let
\(K_0= g^{-1}\big([1/2, 1]\big)\cap K\), \(K_1= g^{-1}\big([0,
1/2]\big) \cap K \). Then \(K_i \subset U_i\) is a closed subset of a compact
set, hence compact, and \( K= K_0\cup K_1\). Since each \(\fhat_i\) is
proper,
\(f^{-1}(K) = \fhat\vphantom{f}^{-1}_1(K_1) \cup
\fhat\vphantom{f}^{-1}_2(K_2)\) is a union of compact sets, hence
compact. Thus \(f\) is proper.
  \end{proof}

\begin{proof}[Proof of Proposition~\ref{prop:proper 1}]
The sets \( U_1 := \tr^{-1}(-2,2)\) and \( U_2:= X^{\irr}(F_n)\) form
an open cover of \(X_\circ(F_n)\) with \(\pi_\circ^{-1}(U_1) = B\),
and \(\pi_\circ^{-1}(U_2) = C\). So by Lemma~\ref{lem:proper union},
it is enough to check that \(\pi_i:=\pihat_{\circ,i}\) is proper
for \(i=1,2\).  Now \(\pihat_{\circ, 2}:C \to U_2\) is proper as it is a
homeomorphism, and so it remains to consider $\pihat_{\circ, 1}$,
which we abbreviate to $\pi$.

To show that \(\pi \maps B\to U_1\) is proper, suppose
\(K\subset U_1\) is compact. By continuity of
$\tr \maps X_\circ(F_n) \to \R$, there is a closed interval
\([a,b] \subset (-2,2)\) such that \(\tr(K) \subset [a,b]\). Inside
$B$, consider \( Z= \mtr^{-1}\big([a,b]\big) = [a, b] \times \cY\),
and consider the smooth submanifold
$M:= \pi^{-1}\big(X^\red(F_n)\big) \cap Z$, which is
$[a, b] \times \CnsignE$ by Theorem~\ref{thm:varconsY}, and so compact
by Lemma~\ref{lem:CnE is Pn-2}.

Let \(D_1 \subset Z\) be a closed tubular neighborhood of \(M\), which
we identify with the product $M \times [-1, 1]$ (using
Theorem~\ref{thm:localmodel}), and let $D_2$ be $Z$ with
$M \times (-1/2, 1/2)$ removed. Then $D_1$ and $D_2$ are closed
subsets of \(B\), with \(\pi(D_1)\) and \(\pi(D_2)\) being closed
subsets of \(U_1\) which cover \(K\).  Now let
\(K_i = K \cap \pi(D_i)\).  Then \(K_1\) and \( K_2\) are closed
subsets of \(K\), hence compact, and \(K = K_1 \cup K_2\). The
preimage \(\pi^{-1}(K_1)\) is a closed subset of the compact set
\(D_1\), hence compact. The restriction of \(\pi\) to \(D_2\) is a
homeomorphism onto its image, so \(\pi^{-1}(K_2)\) is compact. Hence
\(\pi^{-1}(K) = \pi^{-1}(K_1) \cup \pi^{-1}(K_2)\) is compact. It
follows that \(\pi = \pihat_{\circ, 1}\) is proper and hence so is
$\pi_\circ$.
\end{proof}

\subsection{The punctured sphere}
\label{sec:varying Sn}

From now on, we assume $n$ is even and construct the submanifold
\(\cX(S_n) \subset \cX_\bal(F_n) \) from submanifolds \(B' \subset B\)
and \(C' \subset C\), where $B$ and $C$ were used to define
$\cX(F_n)$.  First, we define \(B'\), which is inside the open subset
\(B_\bal = B \cap \cX_\bal(F_n)\) of \(B\):

\begin{lemma}
\label{lem:B'}
  The set \(B' = \setdef{(c, y) \in B}{y \in \cX^c_\bal(S_n)}\)
  is a closed smooth codimension 3 submanifold of \(B_\bal\). The map
  \(\mtr \maps B' \to \R\) is a submersion. 
\end{lemma}

\begin{proof}
Recall from Section~\ref{subsec:define res} the subsets $Z_\epsilon$
and $\cRbar^c(S_n)$ of $\cRbar_n$.  We define
\[
  \cRbar_{n,\bal} = \cRbar_n \setminus \bigcup \setdefm{\big}{Z_\epsilon}%
  {\mbox{$\epsilon \in \pmonen$ and $\sum \epsilon \neq 0$}},
\]
so that \(\cRbar_{n,\bal} = \Pi_c^{-1}\big(X^c_\bal(F_n)\big)\) for
all \(c \in (-2,2)\).  Let
\(\cRbar(F_n) = (-2,2) \times \cRbar_{n,\bal} \), and define
\(\cRbar(S_n) = \setdefm{\big}{(c, v) \in \cRbar(F_n)}{v \in
  \cRbar^c(S_n)}\).
We will first show that \(\cRbar(S_n)\) is a closed smooth codimension
3 submanifold of \(\cRbar(F_n)\), and that the restriction of \(\mtr\)
to \(\cRbar(S_n)\) is a submersion.

Consider the function \(G \maps \cRbar(F_n) \to \SLC\) given by
\[
  G(c, v) = \prod_{i=1}^n A\big(\cos^{-1}(c/2),t(v),v_i\big),
\]
so that the restriction of \(G\) to the slice
\(c\times \cRbar_{n,\bal}\) is the function \(F\) from
Section~\ref{sec:punc sphere}.

Let \(\What = (-2,2) \times \cRbar_n^{t\neq 0}\), and, for each
balanced \(\epsilon\), let
\[
  \Vhat_\epsilon = \setdef{(c,v) \in \cRbar(F_n)}%
  {\mbox{$\sign z_j = \epsilon_j$ for all $j$ and
      $\mRe\big(G_{11}(c,v)\big) > 0$}}.
\]
Then \(\What\) together with the \(\Vhat_\epsilon\)'s form an open
cover of \(\cRbar(F_n)\).  Recall that
$\Util = \setdef{(t,g)}{g \in U_t} \subset \R \times \SLC $, and
define \(\Gtil\maps \What \to \Util\) by
\(\Gtil(c,v) = \big(t(v),G(c,v)\big)\).  Then
\(\cRbar(S_n) \cap \What = G^{-1}(I) =
\Gtil\vphantom{G}^{-1}(\Itil)\), where
$\Itil = \setdef{(t,I)}{t \in \R}$.  We showed in the proof of
Proposition~\ref{prop:SnIrreducible} that the restriction of \(\Gtil\)
to each slice \(c \times \cRbar_{n,\bal}\) is a submersion, so
\(\cRbar(S_n) \cap \What\) is a closed codimension 3 submanifold of
\(\What\), and \(\mtr\maps \cRbar(S_n) \cap \What \to \R\) is a
submersion.

As in Section~\ref{sec:punc sphere}, we have local coordinates
\((c,t,\omega_1, \ldots, \omega_n)\) on \(\Vhat_\epsilon\).  The
function \(G_{11}\) is real analytic on \(\Vhat_\epsilon\), and the
fact that \(\epsilon\) is balanced implies that
\(G_{11}(c,0,v)\equiv 1\). It follows that
\(g_\epsilon = \big(G_{12},\mIm (G_{11} - 1)/t\big)\) is real analytic
on \(\Vhat_\epsilon\). Its restriction to the slice
\(c\times V_\epsilon\) is the function \(f_\epsilon\) from
Section~\ref{sec:punc sphere}, so
\(g^{-1}(0,0) = \cRbar(S_n)\cap \Vhat_\epsilon\). We showed in
Section~\ref{sec:irrchars} that \(f_\epsilon\) is a submersion, and
hence $g_\epsilon$ is as well, so \(\cRbar(S_n)\cap \Vhat_\epsilon\)
is a smooth closed submanifold of \(\Vhat_\epsilon\).  Now
$\mtr \maps \Vhat_\epsilon \to \R$ is a submersion, being projection
onto the first coordinate $c$.  As $g_\epsilon$ is a submersion
restricted to each fiber of $\mtr$, as that is some $f_\epsilon$, it
follows that \(\mtr\maps \cRbar(S_n) \cap \Vhat_\epsilon \to \R\) is
also a submersion. This completes our proof that \(\cRbar(S_n)\) is
a smooth submanifold of $\cRbar(F_n)$ of codimension 3 and that
$\mtr \maps \cRbar(S_n) \to \R$ is a submersion.

Now \(B_\bal\) is the quotient of \(\cRbar(F_n)\) by the simultaneous
actions of the \(G_t\) and \(\R_{>0}\), and \(B'\) is the image of
\(\cRbar(S_n)\) under this quotient, so the statement of the lemma
follows from Lemma~\ref{lem:quotient is manifold}.
\end{proof}

\begin{lemma}
  \label{lem:C'}
  The subset \(C':=X^\irr(S_n)\) is a closed codimension 3 submanifold
  of \(C = X^\irr(F_n)\), and \(\tr \maps C' \to \R\) is a submersion.
\end{lemma}

\begin{proof}
Recall from Section~\ref{sec:realchar} that
\[
  X^\irr(F_n) = X^\irr_\SU(F_n) \coprod X^\irr_\SLR(F_n) =
  X^\irr_{U_1}(F_n) \coprod X^\irr_{U_{-1}}(F_n).
\]
Similarly
$X^\irr(S_n) = X^\irr_{U_1}(S_n) \coprod X^\irr_{U_{-1}}(S_n)$.  We
will show that the claims of the lemma hold for
\(X^\irr_{U_{-1}}(S_n)\); the proofs for \(X^\irr_{U_{1}}(S_n)\) are
identical.

Consider the map \(\mu \maps R^\irr_{U_{-1}}(F_n) \to U_{-1}\) given
by \(\mu(\rho) = \rho(s_1 \cdots s_n)\).  Applying
Lemmas~\ref{Lem:Tangent} and \ref{Lem:AdAction}, we see that the
restriction of \(\mu\) to each \(R^{\irr,c}(F_n)\) is a submersion.
It follows that \(\mu\) is a submersion, and so
\(R^\irr_{U_{-1}}(S_n) = \mu^{-1}(I)\) is a closed codimension 3
submanifold of \(R^\irr_{U_{-1}}(F_n)\).
 
The projection
\(\tau \maps R^\irr_{U_{-1}}(F_n) \to X^\irr_{U_{-1}}(F_n)\) is a
submersion by Lemma~\ref{lem:charsubmer}, and \(R^\irr(S_n)\) is a
union of fibers of \(\tau\).  Applying Lemma~\ref{lem:quotient is
  manifold}, we conclude that
\(X^\irr_{U_{-1}}(S_n) = \tau\big(R^\irr_{U_{-1}}(S_n)\big) \) is a
smooth closed submanifold of \( X^\irr_{U_{-1}}(F_n)\).

For the last claim, note that
\(\tr \maps R^\irr_{U_{-1}}(F_n) \to \R\) is a submersion by
Corollary~\ref{cor:tr is submersion} and Lemma~\ref{lem:charsubmer}.
The restriction of \(\mu\) to the fibers of \(\tr\) is a submersion,
so \(\tr \maps R^\irr_{U_{-1}}(S_n) \to \R\) is a submersion. This map
factors as the composition
\(R^\irr_{U_{-1}}(S_n) \xrightarrow{\tau} X^\irr_{U_{-1}}(S_n)
\xrightarrow{\tr} \R\), so \( \tr\maps X^\irr_{U_{-1}}(S_n) \to \R\)
is a submersion as well.
\end{proof}

\begin{proof}[Proof of Theorem~\ref{thm:sX(S_n)}]
Theorem~\ref{thm:cXSn} gives that
\(B'\cap C = C' \cap B = \tr^{-1}(-2,2) \cap X^\irr(S_n)\). We define
\(\cX(S_n) = B' \cup C'\), so \(\cX(S_n) \cap B = B'\) and
\(\cX(S_n) \cap C = C'\).  Since \(B\) and \(C\) form an open cover of
\(\cX(F_n)\), Lemmas~\ref{lem:B'} and \ref{lem:C'} give that
$\cX(S_n)$ is a closed submanifold of $\cX_\bal(S_n, S)$ of
codimension 3, and also $\pi\big(\cX(S_n)\big)$ is the union of
$\bigcup_{c \in (-2, 2)} X_\bal^c(S_n)$ with $X^\irr(S_n)$ which is a
subset of $X_\bal(S_n)$ as claimed.  The claims about $\mtr$ also
follow immediately from Lemmas~\ref{lem:B'} and \ref{lem:C'}. Finally,
the density of $\cX^\irr(S_n)$ in $\cX(S_n)$ follows from the last
sentence of Theorem~\ref{thm:cXSn}.
\end{proof}

\begin{remark}
  As discussed in Remark~\ref{sec: odd punct}, we do not determine
  the complete analogue of Theorem~\ref{thm:cXSn} when the number of
  punctures $n$ is odd.  However, the claims of
  Theorem~\ref{thm:sX(S_n)} are effectively trivial for odd $n$: every
  reducible is unbalanced, and hence
  \(X_\bal(S_n,S) = X^\irr(S_n, S)\), so one can simply define
  $\cX(S_n, S)$ as $ X^\irr(S_n, S)$.
\end{remark}

\section{Orientations}
\label{Sec:Orientations}

In this section, we orient the spaces constructed in the previous
sections: \(\cX^c(F_n)\), \(\cX(F_n)\), \(\cX^c(S_n)\), and
\(\cX(S_n)\). These spaces contain the \(\SU\) character varieties
used by Heusener in \cite{Heusener2003} as open subsets. We
compare our orientations with his, as this will be needed in
Section~\ref{sec: rel with signature}.

\subsection{The free group} By Corollary~\ref{cor:cYCnS} and
Theorem~\ref{thm:def of pi}, we know that
\(\cX^c(F_n) \cong X^\cirr_\SU(F_n) \cong \CnS\). We defined an
orientation on the latter space in Section~\ref{sec:configorient}.  To
compare this orientation with Heusener's, we explain how to orient
$\XSU^\irr(F_n)$ from the perspective of character varieties,
following \cite{Heusener2003}.

First, we orient $\SU = U_1$ via Remark~\ref{rem:actsame}.  For all
$c \in (-2, 2)$, we orient the conjugacy class
\(\RSU^c(F_1) = \tr^{-1}(c)\) so that the submersion
$\tr \maps \SU \setminus \{\pm I\} \to (-2, 2)$ is compatibly
oriented in the sense of Section~\ref{sec:orientshortex}, where
\((-2,2)\) has the \emph{opposite} of the standard orientation (this
seemingly odd choice will be explained in the proof of
Proposition~\ref{prop: orient agree SU2}).  Give
\(\RSU^c(F_n) = \big(\RSU^c(F_1)\big)^n\) the product orientation,
which then orients its open subset $R^\cirr(F_n)$.  Now orient
$\XSU^\cirr(F_n)$ by the requirement that
$\SO_3 \to \RSU^\cirr(F_n) \to \XSU^\cirr(F_n)$ is compatibly
oriented, where $\SO_3$ is oriented as $\SU/\{\pm I\}$.

Furthermore, we orient $\XSU^\irr(F_n)$ using
\[
\XSU^\cirr(F_n) \to \XSU^\irr(F_n) \xrightarrow{\tr} (-2, 2)
\]
where $(-2, 2)$ again has the reversed orientation, the map $\tr$ is a
submersion by Corollary~\ref{cor:tr is submersion}, and the
orientations on the fibers vary smoothly by construction.

\begin{proposition}
  \label{prop: orient agree SU2}
  The orientation on $\XSU^\cirr(F_n)$ is opposite to that of 
  \cite{Heusener2003}.  Also, the diffeomorphism
  $\pi_c \maps \CnS \to \XSU^\cirr(F_n)$ from
  Theorem~\ref{thm:def of pi} is orientation preserving where $\CnS$ is
  oriented as in Section~\ref{sec:configorient}.  Finally, the map
  $(-2, 2) \times \CnS \to \XSU^\irr(F_n)$ that sends
  $(c, p) \to \pi_c(p)$ is an orientation preserving diffeomorphism
  where $(-2, 2)$ has the reversed orientation.
\end{proposition}

\begin{proof}
For the first claim, our construction matches \cite{Heusener2003}
except there conjugacy classes in $\SU$ are parameterized by
$\alpha \in (0, \pi)$ where $c = 2 \cos(\alpha)$ and $(0, \pi)$ has
the usual orientation. However, this corresponds to the reversed
orientation we used on $(-2, 2)$.  Second, we must account for the
differences in orientation conventions noted in Remark~\ref{rem:
  orient vs Heusener}. The sequence used to orient \(R^c_\SU(F_1)\)
gives the same orientation as Huesener, since the fiber is even
dimensional.  However the second sequence gives the opposite of
Huesener's orientation, since both the fiber \(SO_3\) and base
\(X^\cirr_\SU(F_n)\) are odd-dimensional, proving the first claim.

For the second claim, the key is to show that the diffeomorphism of
$\Qbar_1 = S^2$ to $\RSU^c(F_1)$ that sends
$v \mapsto A(\alpha, 1, v)$ with $\alpha = \cos^{-1}(c/2)$ is
orientation preserving; this will suffice as the constructions of
orientations on $\CnS$ and $\XSU^\cirr(F_n)$ are then strictly
analogous.  Recall that
$ A(\alpha, t, v) = \exp(\alpha \varphi_t (v))$. Now the exponential
map restricts to an orientation-preserving diffeomorphism from
$\setdef{v \in \ut[1]}{0 < B_1(v) < \pi^2}$ to
$U_1 \setminus \{\pm I\}$ which sends
$\alpha \Qbar_1 = \setdef{v \in \ut[1]}{B_1(v) = \alpha^2}$ to
$\RSU^c(F_1)$.  Since the orientation of the sphere $\Qbar_1$ from
Section~\ref{subsec:B_t} is the same as orienting it with the outward
normals, this matches the orientation on $\RSU^c(F_1)$ which comes
from increasing $\alpha$.  So the map $\Qbar_1 \to \RSU^c(F_1)$ is
orientation-preserving as needed to prove the second claim.

  The final claim follows easily as $\XSU^\irr(F_n)$ was oriented using the
  submersion $\tr \maps \XSU^\irr(F_n) \to (-2,2)$ where $(-2, 2)$
  has the reversed orientation.
\end{proof}

To orient $\XSLR^\irr(F_n)$, we follow the same steps as for
$\SU$. First, fix the orientation on $\SLR$ coming from its
identification to $U_{-1} = \SUoneone$.  There is a submersion
$\tr \maps (\SLR \setminus \{\pm I\}) \to \Rrev$, where $\Rrev$
denotes $\R$ with the reversed orientation, and we orient
$\RSLR^c(F_1) := \tr^{-1}(c)$ as its fiber. By definition, the
singular points \(\pm I\) are not in \( \RSLR^{\pm 2}(F_1)\).  Unlike
for $\SU$, where each conjugacy class is a sphere, the topology of
$\RSLR^c(F_1)$ depends on $c$, but this makes no difference at this
stage.
Now orient $\RSLR^\cirr(F_n)$ as an open subset of the product
$\big(\RSLR^c(F_1)\big)^n$.  (When \(c = \pm 2\), recall that for
$\rho \in \RSLR^{c}(F_n)$, the \(\rho(s_i)\) are all conjugate in
\(\SLC\). Hence if $\rho \in \RSLR^{\pm 2, \irr}(F_n)$,
\(\rho(s_i) \neq \pm I\).)  Orient $\XSLR^\cirr(F_n)$ by the fiber
bundle
\[
  \PSLRpm \to \RSLR^\cirr(F_n) \to \XSLR^\cirr(F_n)
\]
where the orientation on $\PSLRpm = \SLRpm/\{\pm I\}$ comes from the
left-invariant extension to $\SLRpm$ of
the one on $\SLR$.  Finally, orient $\XSLR^\irr(F_n)$ by
$\tr \maps \XSLR^\irr(F_n) \to \Rrev$, which is a submersion by
Corollary~\ref{cor:tr is submersion}.

\begin{proposition}
  \label{prop: orient agree SL2R}
  For $c \in (-2, 2)$, the diffeomorphism
  $\pi_c \maps \CnsignH \to \XSLR^\cirr$ from
 Theorem~\ref{thm:def of pi} is orientation preserving where $\CnsignH$
  is oriented as in Section~\ref{sec:configorient}.
\end{proposition}
\begin{proof}
As in the proof of Proposition~\ref{prop: orient agree SU2}, the key
is to check that the diffeomorphism from $\Qbar_{-1}$ to
$\RSLR^c(F_1)$ that sends
$v \mapsto A(\alpha, 1, v) = \exp(\alpha\phi_{-1} (v))$, where
$\alpha = \cos^{-1}(c/2)$, is orientation preserving.  Recall from
Section~\ref{subsec:B_t} that both sheets of $\Qbar_{-1}$ are oriented
with respect to normals that point away from the origin in $\ut[-1]$.
For any $\beta \in (0, \pi)$, the exponential map gives a
diffeomorphism from $\beta \Qbar_{-1}$ to $\RSLR^{2 \cos \beta}(F_1)$;
here, $\beta \Qbar_{-1}$ denotes the dilation by ordinary scalar
multiplication of vectors, not the $\R_{>0}$ action from earlier.
Orient each $\beta \Qbar_{-1}$ as we did $\Qbar_{-1}$; hence the
dilation $\Qbar_{-1} \to \beta \Qbar_{-1}$ is orientation preserving.
Since $\RSLR^{2 \cos \beta}(F_1)$ was oriented using
$\tr \maps \RSLR(R_1) \to \Rrev$, i.e. with respect to increasing
$\beta$, and the map $\exp \maps \ut[-1] \to U_{-1}$ is orientation
preserving by definition, the diffeomorphism from $\beta \Qbar_{-1}$
to $\RSLR^{2 \cos \beta} (F_1)$ is orientation preserving.  Therefore
the map from $\Qbar_{-1}$ to $\RSLR^c(F_1)$ that sends
$v \mapsto A(\alpha, 1, v) = \exp(\alpha v)$ is the composition of two
orientation preserving maps and hence orientation preserving as
needed.
\end{proof}

\begin{theorem}\label{thm: res for F_n orient}
  The manifold \(\cX(F_n)\) is orientable, and we take its preferred
  orientation to be the one compatible with its open subset
  $\XSU^\irr(F_n)$.  Then the open subset $\XSLR^\irr(F_n)$ has the
  reverse of its standard orientation.  The submersion
  $\mtr \maps \cX(F_n) \to \Rrev$ gives a preferred orientation to
  each $\cX^c(F_n)$, which is the standard orientation on
  $\XSU^\cirr(F_n)$ but the reverse of it for $\XSLR^\cirr(F_n)$.
\end{theorem}

\begin{proof} Recall from Section~\ref{sec:varying Fn} that $\cX(F_n)$
is the union of open submanifolds $B = (-2, 2) \times \cY$ and
$C = X^\irr(F_n) = \XSU^\irr(F_n) \coprod \XSLR^\irr(F_n) $.  Both $B$
and $C$ are orientable, but the orientability of $\cX(F_n)$ is not
immediate as $B \cap C$ is disconnected.  Orient $B$ by giving
$(-2, 2)$ the reversed orientation and $\cY$ the orientation
compatible with the one  on \(\CnS\).  Now orient $C$ by giving
$\XSU^\irr(F_n)$ and $\XSLR^\irr(F_n)$ their standard and reversed
orientations respectively.  By Proposition~\ref{prop:orientation on
  Y}, the inclusions of $\CnS$ and $\CnsignH$ into $\cY$ are
orientation preserving and reversing respectively. Combining
Propositions~\ref{prop: orient agree SU2} and \ref{prop: orient agree
  SL2R} with the fact that the submersion
$\tr \maps X^\irr(F_n) \to \Rrev$ is compatible with the orientations
on $X^\cirr(F_n)$, we see that the restrictions of the two
orientations to \(B \cap C\) agree. This gives the desired global
orientation on $\cX(F_n)$, and the claims about $\cX^\cirr(F_n)$ are
immediate from the construction.
\end{proof}
 
\subsection{The punctured sphere}
To begin, we orient the manifold $\XSU^\cirr(S_n)$ as follows.  By
Lemma~\ref{Lem:AdAction}, the map $\mu \maps \RSU^\cirr(F_n) \to \SU$
where $\mu(\rho) = \rho(s_1 \cdots s_n)$ is a submersion, so we orient
$\mu^{-1}(I) = \RSU^\cirr(S_n)$ using our usual conventions. Now
orient $\XSU^\cirr(S_n)$ by
$\SO_3 \to \RSU^\cirr(S_n) \to \XSU^\cirr(S_n)$.  Finally orient
$\XSU^\irr(S_n)$ by $\tr \maps \XSU^\irr(S_n) \to \Rrev$, which is a
submersion by Lemma~\ref{lem:C'}, and the just constructed orientation
on the fibers $\XSU^\cirr(S_n)$, noting that said orientation varies
continuously in $c$.  Our standard orientations on $\XSLR^\cirr(S_n)$
and $\XSLR^\irr(S_n)$ are analogous.

\begin{remark}
  \label{rem: orient vs H for S_n}
  The orientations just constructed on \(\RSU^\cirr(S_n)\) and
  $\XSU^\cirr(S_n)$ are the opposites of those used in
  \cite{Heusener2003}. For $\RSU^\cirr(S_n)$, this is because
  $\RSU^\cirr(S_n)$ and $\SU$ are both odd-dimensional, so the
  convention difference discussed in Remark~\ref{rem: orient vs
    Heusener} matters.  This difference persists to $\XSU^\cirr(S_n)$
  since in both cases $\XSU^\cirr(S_n)$ is even-dimensional and hence
  the convention for orienting submersions is irrelevant.
\end{remark}

Unlike in the free group case, we can orient \(\cX^c(S_n)\) so that
both \(X^\cirr_\SU(S_n)\) and \(X^\cirr_\SLR(S_n)\) inherit their
standard orientations:
\begin{theorem}
  \label{Thm:Orientation for S}
  The manifold \(\cX^c(S_n)\) is orientable, and we take its preferred
  orientation to be the one compatible with its open subset
  \(X^\cirr_\SU(S_n)\).  The open subset \(X^\cirr_\SLR(S_n)\) then
  also has its standard orientation.
\end{theorem}
 
\begin{proof}
 We imitate the construction of the standard orientation on
 \(X^c_\SU(S_n)\) in the setting of our construction of $\cX^c(S_n)$
 from Section~\ref{sec:punc sphere}. We proceed in three steps. 
 
\subsection*{Step 1:  orient \(\cRbar_n\)}
Recall from Section~\ref{subsec:B_t} that each $\Qbar_t$ is oriented
so that $\pair{u_1, u_2}$ is a positive basis for $T_{v} (\Qbar_{t})$
exactly when $\pair{v, u_1, u_2}$ is a positive basis for $\R^3$. There is a
submersion \(\cRbar_n \to \R\) whose fibers are \(\Qbar^n_t\); we give
\(\Qbar^n_t\) the product orientation and use the fibration
\begin{equation}
\label{Eq:Orientation 1}
  \Qbar_t^n \to \cRbar_n \xrightarrow{t} \R 
\end{equation} 
to orient \(\cRbar_n\). 

 \subsection*{Step 2:  orient  \(\cRbar^{c}(S_n)\)}   Recall from
the proof of Proposition~\ref{prop:SnIrreducible} that
\(\cRbar^{c,t\neq 0}(S_n) := \Ftil^{-1}(\Itil)\), where
\(\Ftil\maps\cRbar_n^{t\neq 0} \to \Util\) is a submersion.  The normal
bundle \(\nu_{\Itil}\) is the bundle over \(\R\) whose fibers are
\(\ut \subset T\Util\); we orient each $\ut$ as described in
Remark~\ref{rem:actsame}, which gives smoothly varying orientations on
the fibers of \(\nu_{\Itil}\).

Over $\cRbar^{c,t \neq 0}(S_n)$, the short exact sequence
\begin{equation}
\label{Eq:Orientation 2}
 0  \to T \cRbar^{c,t \neq 0}(S_n) \to T \cRbar^{t\neq 0}_n \xrightarrow{d\Ftil} \nu_{\Itil} \to 0
\end{equation}
can be used to orient \( \cRbar^{c,t \neq 0}(S_n)\).  You might think
that we should extend this orientation to all of \(\cRbar^c(S_n)\),
but this is actually not possible.  To see this, consider the function
\(h\maps\Util \to \C \times \R \) given by
\(h(t,A) = ( A_{12},\mIm A_{11})\). We have
\(h^{-1}(0, 0) = \Itil \cup -\Itil\) where $-\Itil$ is everything of
the form $(t, -I)$; thus \(\nu_{\Itil}\) is one connected component of
\(h^*(T_0(\R\times \C))\). So the orientation on
\(\cRbar^{c, t \neq 0}(S_n)\) induced by \eqref{Eq:Orientation 2} is
the same as that induced by the submersion
\( \cRbar^{c, t\neq 0}(S_n) \to \cRbar^{t\neq 0}_n \xrightarrow{h
  \circ \Ftil} \R \times \C\).  But on the open set \(V_\epsilon\)
from Section~\ref{subsec: reduc char},
the subset \(\cRbar^c(S_n) \cap V_\epsilon\) is cut out by a
submersion \(f_\epsilon\maps V_\epsilon \to \C \times \R \), which
satisfies
\[
  f_\epsilon = ( F_{12},f_{11}) \quad \text{where} \quad h\circ \Ftil =
  (F_{12},tf_{11}).
\]
Consequently the induced orientation on \( f_{\epsilon}^{-1}(0,0)\) will
agree with the orientation induced by \eqref{Eq:Orientation 2} for
\(t>0\), but will be its opposite for \(t<0\) since $f_\epsilon$ differs
from $h\circ \Ftil$ by scaling the second coordinate of $\C \times \R$ by \(t\). 

To orient \(\cRbar^c(S_n)\), we give \(V_\epsilon\) the orientation
induced by our orientation on \(\cRbar_n\), and the orientation on
\(\C \times \R\) that pulls back to our standard orientation on
\(\ut\) under \(h^*\). Then the submersion \(f_\epsilon\) induces an
orientation on \(f_\epsilon^{-1}(0,0)\). We give
\( \Ftil^{-1}(\Itil) \cap \cRbar^{t>0}_n\) the orientation coming from
\eqref{Eq:Orientation 2}, and
\(\Ftil^{-1}(\Itil) \cap \cRbar^{t<0}_n\) the opposite orientation.
The discussion in the previous paragraph shows that these orientations
are all compatible and together give a global orientation on
\(\cRbar^c(S_n)\).

\subsection*{Step 3:  orient  \(\cX^{c}(S_n)\)} 
Now we can orient \(\cX^c(S_n)\).  We use the orientation on \(\ut\)
to orient \(G_t\) as in Remark~\ref{rem:actsame}, give \(\R_{>0}\) the
positive orientation, and give \(\R_{>0} \times G_t\) the product
orientation.  We will  orient
\(\cX^{c}(S_n)\) using the exact sequence
\begin{equation}
\label{Eq:Orientation 3}
 0\to T_e(\R_{>0} \times G_t) \to T_v \cRbar^{c}(S_n) \to T_p \cX^{c}(S_n)
 \to 0,
\end{equation}
where \(t=t(v)\), $e = (0, I)$, and $\sigma(v) = p$.  As the orientations on $T_v \cRbar^{c}(S_n)$ and
$T_e(\R_{>0} \times G_t)$ vary continuously, even as $t$ changes, the
main issue is that, for \(t<0\), the fibers $\R_{>0} \times G_t$ of
the submersion \(\sigma\maps \cRbar^c(S_n) \to \cX^c(S_n)\) are
disconnected.  So for $t < 0$, we must find \(v,v'\) belonging to two
different components of \(\sigma^{-1}(p)\) where the orientations on
\(T_p\cX^c(S_n)\) induced by the exact sequences
\begin{align*}
 & 0\to T_e(\R_{>0} \times G_t) \to T_v \cRbar^{c}(S_n) \to T_p \cX^{c}(S_n)  \to 0 \\
 & 0\to T_e(\R_{>0} \times G_t) \to T_{v'} \cRbar^{c}(S_n) \to T_p \cX^{c}(S_n)  \to 0
\end{align*}
are the same. 

To do this, consider
$ C = \left(\begin{smallmatrix} -1 & 0 & 0 \\ 0 & 1 & 0 \\ 0 & 0 &
    -1 \end{smallmatrix}\right)$, which is in every $G_t$.  Its action
on each $\ut$ is orientation preserving, and corresponds to complex
conjugation as noted in Section~\ref{sec: props Gt and Ut}.  We
successively deduce that \(C:\Qbar_t \to \Qbar_t\) and
\(C\maps\cRbar_n \to \cRbar_n\) are orientation preserving. Since
\(F(C\cdot v) = C \cdot F(v)\) as noted in Section~\ref{subsec: intro
  to punc sphere}, we see that \(C\maps\cRbar(S_n) \to \cRbar(S_n)\)
and is orientation preserving. Taking \(v' = C \cdot v\), we have a
commutative diagram
\[
  \begin{tikzcd}
0 \arrow{r}  &  T_e(\R_{>0} \times G_t) \arrow{r} \arrow{d}{\id \times \Ad_C} & T_v \cRbar^{c}(S_n) \arrow{r} \arrow{d}{dC} & T_p \cX^{c}(S_n)  \arrow{r} \arrow{d}{\id}& 0 \\
0 \arrow{r} &  T_e(\R_{>0} \times G_t) \arrow{r} & T_{v'} \cRbar^{c}(S_n) \arrow{r} & T_p \cX^{c}(S_n)  \arrow{r} & 0. 
\end{tikzcd}
\]
The two leftmost vertical maps are orientation preserving, so the
righthand vertical map is as well. In other words, the two induced
orientations on \(T_p \cX^{c}(S_n) \) are the same. We conclude that
\eqref{Eq:Orientation 3} determines an orientation on
\(\cX^c(S_n)\).

Having defined an orientation on $\cX^c(S_n)$, it remains to check
that the orientations it induces on its open subsets
\(\XSU^\cirr(S_n)\) and \(\XSLR^\cirr(S_n)\) agree with their
standard orientations. First, we consider the case \(t>0\). For the
standard orientation of \(\XSU^\cirr(S_n)\), we have the following
decomposition as oriented vector spaces:
\begin{equation}
  \label{Eq:Orientation 4}
  T \Qbar_1^n = T \big(\RSU^c(F_1)^n\big) = TU_1\oplus T
  \big(R^c_\SU(S_n)\big) =  TU_1\oplus \big[ T \big(X^c_\SU(S_n)\big)
  \oplus TU_1 \big]
\end{equation}
For the orientation of \(\cX^c(S_n)\), we have
\(T \cRbar_n^c = \R\oplus T\Qbar_t^n\), where the $\R$-factor corresponds
to the $t$ coordinate.  Again as oriented vector spaces, we have from
(\ref{Eq:Orientation 3}):
\begin{align*}
  T \cRbar_n^c &= TU_t \oplus T\cRbar_n^c(S_n) =
              TU_t \oplus \big(T \cX^c(S_n)
              \oplus T(\R_{>0} \times G_t) \big) \\
            &= - T(\R_{>0}) \oplus TU_t \oplus T \cX^c(S_n) \oplus TU_t
\end{align*}
since $\dim \cX^c(S_n) = 2n-6$ and $\dim U_t = 3$.  Now as $s \cdot
\Qbar_t = \Qbar_{t/s^2}$ and \(t>0\), increasing $s$ actually decreases $t$.  Hence
\[
  T \cRbar_n^c = - T(\R_{>0}) \oplus T \Qbar_t^n.
\]
Comparing, we conclude that $T \Qbar_t^n = TU_t \oplus T \cX^c(S_n)
\oplus TU_t$, and, taking $t = 1$, we see from (\ref{Eq:Orientation 4})
that $\cX^c(S_n)$ and $\XSU^\cirr(S_n)$ have matching orientations, as
desired.

The case of \(t<0\) is very similar, except that now
$T \cRbar_n^c = T(\R_{>0}) \oplus T \Qbar_t^n$ and 
\[
  T \cRbar_n^c = T(\R_{>0}) \oplus TU_t \oplus T \cX^c(S_n) \oplus TU_t
\]
since the orientation on this part of $\cX^c(S_n)$ comes from that of
\(\Ftil^{-1}(\Itil) \cap \cRbar^{t<0}_n\) which was given the opposite
orientation in Step 2.  The net effect is to insert a \(-\) sign
at one point and remove it from another, so we again find that the
induced orientation agrees with the standard one.
\end{proof}

\begin{corollary}
  \label{cor: orient cX(S_n)}
  The manifold \(\cX(S_n)\) is orientable, with preferred orientation
  that matches the standard ones on both $\XSU^\irr(S_n)$ and
  $\XSLR^\irr(S_n)$.  The submersion $\mtr \maps \cX(S_n) \to \Rrev$
  is compatible with the preferred orientations on
  $\cX_\bal^c(S_n) \subset \cX^c(S_n)$ from
  Theorem~\ref{Thm:Orientation for S}.
\end{corollary}

\begin{proof}
Recall from Section~\ref{sec:varying Sn} that
\(\cX(S_n) = B'\cup C'\). The map \(\mtr\maps B' \to (-2,2)\) is a
submersion by Lemma~\ref{lem:B'}, so there is a short exact sequence
\[
  0 \to T\cX^c_\bal(S_n) \to T B' \to T(-2,2)\to 0
\]
which we use to orient \(B'\); as usual, we give \((-2,2)\) the
reversed orientation. (Note \(\cX^c_\bal(S_n)\) is an open subset of
\(\cX^c(S_n)\), and inherits an orientation from it.)  We give both
components of \(C' = X^{\irr}_\SU(S_n) \amalg X^{\irr}_\SLR(S_n)\) our
standard orientations. The orientation on \(X^{\irr}_\SU(S_n)\) is
induced by the short exact sequence
\[0 \to X^{\cirr}_\SU(S_n) \to X^{\irr}_\SU(S_n) \to T(-2,2)\to 0\]
and similarly for \(X^{\irr}_\SLR(S_n)\). By
Theorem~\ref{Thm:Orientation for S}, the orientations on
\(\cX^c(S_n), X^{\cirr}_\SU(S_n)\), and \(X^{\irr}_\SLR(S_n)\) are
compatible. It follows that the orientations on \(B'\) and \(C'\) are
compatible and determine a global orientation on \(\cX(S_n)\) which
has the desired properties.
\end{proof}


\begin{figure}
  \centering
  \begin{tikzoverlay}[scale=0.95]{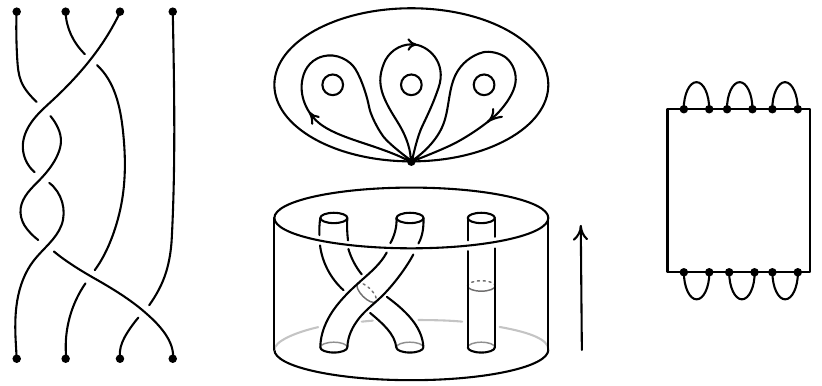}
  \begin{scope}[nmdstd]
    \node[below] at (2.1, 3.0) {$1$};
    \node[below] at (8.0, 3.0) {$2$};
    \node[below] at (14.7, 3.0) {$3$};
    \node[below] at (21.2 ,3.0) {$4$};
    \node[above] at (2.1, 46) {$1$};
    \node[above] at (8.0, 46) {$2$};
    \node[above] at (14.7,46) {$3$};
    \node[above] at (21.2,46) {$4$};
    \node[below, font=\small] at (11.35, -1)
         {$\sigma_2 \sigma_1^2 \sigma_2^{-1} \sigma_3^{-1}$};
    \node[above] at (39.6,40.2) {$s_1$};
    \node[above] at (50.1,41.8) {$s_2$};
    \node[above] at (60.3,40.4) {$s_3$};
    \node[right, font=\small] at (71.8,11.9) {$\phi_{\sigma_1}$};
    \node[below, font=\small] at (90.1,25.6) {$\betahat$};
  \end{scope}
\end{tikzoverlay}
  
  \caption{This figure shows our conventions for braids, including
    $\sigma_i$ versus $\sigma_i^{-1}$, that left-to-right in the braid
    word corresponds to  top-to-bottom in the picture, and that the
    generators $s_i$ of $\pi_1(S_n)$ correspond to clockwise loops.
    Moreover, the induced mapping class $\phi_\beta \in \MCG(D_n)$ of
    $\beta$ is the map that pushes curves from the bottom to top,
    hence $\phi_{\alpha \beta} = \phi_\alpha \circ \phi_\beta$ so that
    $B_{n} \to \MCG(D_n)$ is a homomorphism rather than an
    antihomorphism (here functions act on the left as usual); in
    particular, $\phi_{\sigma_1}(s_1) = s_2$ and
    $\phi_{\sigma_1}(s_2) = s_2^{-1} s_1 s_2$.}
  \label{fig:braids}
\end{figure}

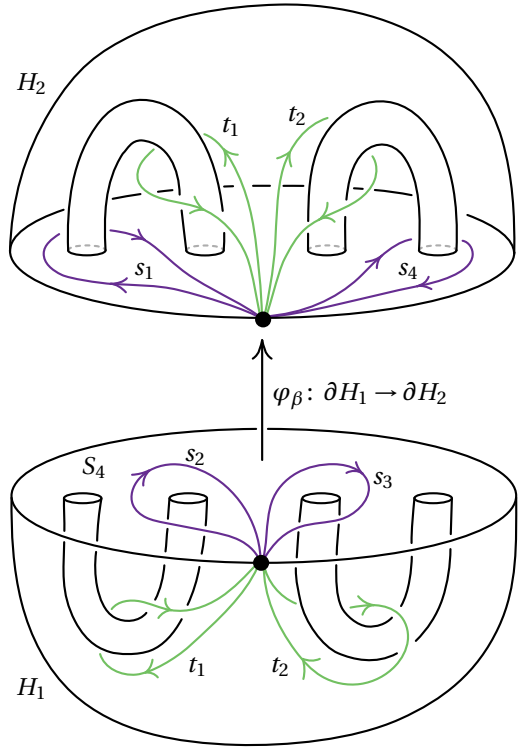
\begin{figure}
  \centering
  \begin{tikzoverlay}[scale=1.0]{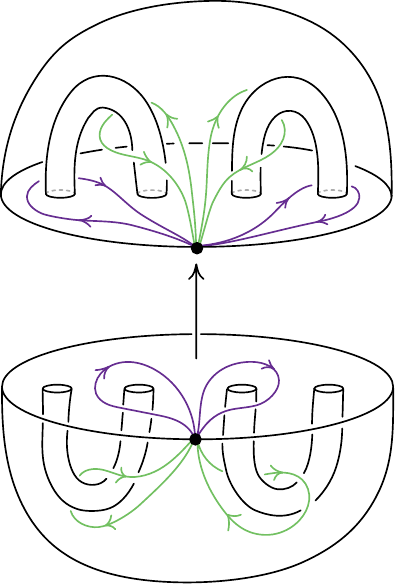}
    \begin{scope}[nmdstd]
      \node[] at (4.5,132.0) {$H_2$};
      \node[] at (4.5,11.6) {$H_1$};
      \node[right] at (50.0,69.8) {$\phi_\beta \maps \partial H_1 \to
        \partial H_2$};
      \node[above] at (26.7,91.1) {$s_1$};
      \node[above] at (78.5,91.1) {$s_4$};
      \node[below right] at (33,19.8) {$t_1$};
      \node[below left] at (57.3,19.8) {$t_2$};
      \node[left] at (40.8,57.6) {$s_2$};
      \node[] at (73.6,53) {$s_3$};
      \node[] at (16.8,55.6) {$S_4$};
      \node[above] at (43.6,119.5) {$t_1$};
      \node[above] at (56.6,121.7) {$t_2$};
    \end{scope}
  \end{tikzoverlay}
  \caption{This figure shows our conventions for the fundamental group
    of the exterior $\EK = H_1 \cup_{\phi_\beta} H_2$ of the plat
    closure $K = \betahat$.  Here the canonical copy of $S_{2m}$ is
    the one on the boundary of $H_1$.  For clarity, only half of the
    generators of $\pi_1(S_4)$ are drawn on each copy of $S_4$.}
  \label{fig:plat}
\end{figure}

\section{Definition of the invariant}
\label{Sec:Definition of h}

\subsection{Braids and plat closures}
\label{sec: braids and plats}

We are now ready to define the total Lin invariant of a knot
\(K \subset S^3\). Following Heusener \cite{Heusener2003}, we
represent knots as plat closures of braids.  Our conventions for
braids are given in Figure~\ref{fig:braids}, which also describes the
identification of the braid group \(B_{n}\) with \(\MCG( D_n)\), the
mapping class group of the \(n\)-punctured disk. In turn, there is a
homomorphism from \(\MCG( D_n)\) to \(\MCG( S_n)\) induced by capping
off the boundary to $D_n$ with a disk to form $S_n$. If
\(\beta \in B_n\), we abuse notation and write \(\beta\) for its image
$\phi_\beta$ in \(\MCG( S_n)\).  Recall a homomorphism
$\phi \maps \Gamma \to \Lambda$ of finitely generated groups has an
induced regular real algebraic map
$\phi^* \maps X(\Lambda) \to X(\Gamma)$ given by
$[\rho] \mapsto [\rho \circ \phi]$.  If $S$ and $T$ are finite
generating sets of $\Gamma$ and $\Lambda$ respectively where $\phi$
takes each $s_i \in S$ to a conjugate of some $t_j \in T$ or its
inverse, then we also have induced maps
$\phi^* \maps X^c(\Lambda, T) \to X^c(\Gamma, S)$ for each $c \in \R$.
Hence \(\beta_* \maps \pi_1(S_n) \to \pi_1(S_n)\) induces a map
\(\beta^*\maps X^c(S_n, S) \to X^c(S_n,S)\) for each $c \in \R$.

Given a braid \(\beta \in B_{2m}\), we can form its plat closure
\(\widehat{\beta}\) as shown in Figure~\ref{fig:braids}. Given a knot
\(K \subset S^3\), we choose \(\beta\) with \(\widehat{\beta} =
K\). As shown in Figure~\ref{fig:plat}, this presents the exterior
of \(K\) as the union of two genus \(m\) handlebodies glued together
along the \(2m\)--punctured sphere \(S_{2m}\):
\[
  \EK =  S^3 \setminus \nu(K) = H_1 \cup_{S_{2m}} H_2.
\]
Here, the surface $S_{2m}$ is identified with the boundary of $H_1$
and $x \in \partial H_1$ is glued to $\beta(x) \in \partial H_2$.  The
inclusion $S_{2m} \to \EK$ induces a surjection $\pi_1(S_{2m}) \to
\pi_1(\EK)$ and hence an injection of \(X^c(K) = X^c(\EK, \{\mu\})\) into
$X^c(S_{2m}, S)$.

Similarly, if \(T_k\) is the set of generators for \(\pi_1(H_k)\)
shown in Figure~\ref{fig:plat}, the inclusion
\(j_k\maps S_{2m} \to H_k\) induces an injection
\(j_k^*\maps X^c(H_k,T_k) \to X^c(S_{2m},S)\).  We let
\(L^c_k \subset X^c(S_{2m},S)\) be the image of \(j_k^*\). Then
\(X^c(K)\) can be identified with the intersection \(L^c_1\cap L^c_2\)
as in the original construction of Casson's invariant
\cite{AkbulutMcCarthy1990}.

More precisely, we index the standard generators
\(s_1, \ldots, s_{2m}\) of \(\pi_1(S_{2m})\) as shown in
Figures~\ref{fig:braids} and \ref{fig:plat} so that \(s_{2i-1}\) and
\(s_{2i}\) are meridians of the \(i\)th arc in the lower plat. The
generators shown in Figure~\ref{fig:plat} determine isomorphisms
\[
  \ell_k \maps \pi_1(H_k) \to F_m = \langle t_1, \ldots, t_m\rangle.
\]
If we define \(\jhat_k = \ell_k \circ j_{k*}\), then
\(\jhat_1(s_{2i-1}) = t_i \) and \( \jhat_1(s_{2i}) = t_{i}^{-1} \),
so \(\jhat_1\) induces a map
\(\jhatstar_1 \maps X^c(F_m,T) \to X^c(S_{2m},S)\) and
\( \im \jhatstar_1 = \im j_1^* = L_1^c\). We have
\(\jhat_2 = \jhat_1 \circ \beta_*\), so
\(L_2^c = \im \jhatstar_2 =\beta^* (L_1^c)\).  To simplify the notation,
we let \(L^c = L_1^c\), so
\(L_1^c \cap L_2^c = L^c \cap \beta^*(L^c)\), and write
\(i = \jhatstar_1\).  To summarize, we have established:

\begin{proposition}\label{prop: Casson style}
  For a knot $K = \betahat$ and each $c \in \R$, the inclusion
  $S_{2m} \to \EK$ gives an injection $X^c(K) \to X^c(S_{2m}, S)$ whose
  image is $L^c \cap \beta^* L^c$.
\end{proposition}

We would like to define the total Lin invariant \(h^c(K)\) to be the
algebraic intersection number \(\langle L^c, \beta^* L^c\rangle\), but
there are two problems with doing so.  The first is that
\( L^c \cap \beta^*(L^c)\) may not be compact, and the second is that
\(X^c(S_{2m},S)\) is not a manifold.  We more or less punt on the
first issue: we will impose hypotheses on \(K\) which prevent it from
happening. For the second problem, we restrict to $c \in (-2, 2)$ and
replace \(X^c(S_{2m},S)\) with the resolution \(\cX^c(S_{2m},S)\)
constructed in Section~\ref{sec:punc sphere}. This requires us to
construct smooth maps
\(i'\maps \cX^c(F_m,T) \to \cX^c_\bal(S_{2m},S)\) and
\(\beta^*\maps \cX^c(S_{2m},S) \to \cX^c(S_{2m},S)\) that are
``resolutions'' of the corresponding maps on character
varieties. Letting \(\cL^c = i'\big(\cX^c(F_m,T)\big)\), we will then
define \(h^c(K) = \langle \cL^c, \beta^*\cL^c \rangle\).  As described
in Theorem~\ref{Prop:Lcap beta L} below, the invariant \(h^c(K)\) can
be interpreted as a signed count of representations of \(\pi_1(\EK)\)
with \(\tr_\mu = c\) into certain double covers of
\(\Isom^+(S^2), \Isom^+(\E^2)\), and \(\Isom^+(\H^2)\).

\subsection{Induced maps of resolutions}
\label{sec: induced res}

For any of our resolutions $\cX \xrightarrow{\pi} X$ in
Theorems~\ref{thm:def of pi}, \ref{thm:cXSn}, \ref{thm:sX(F_n)}, and
\ref{thm:sX(S_n)}, define a decomposition of \(\cX\) into
$\cX_{+} \amalg \cX_0 \amalg \cX_{-}$ by
$\cX_{+} = \pi^{-1}(X^\irr_\SU)$, $\cX_0 = \pi^{-1}(X^\red)$, and
$\cX_{-} = \pi^{-1}(X^\irr_\SLR)$.  In particular,
$\cX^\irr = \cX_{+} \amalg \cX_{-}$ and $\cX^\red = \cX_0$.  
In each case $\pi$ gives diffeomorphisms $\cX_+ \to X^\irr_\SU$ and
$\cX_{-} \to X^\irr_\SLR$.  For $c \in (-2, 2)$, each $\cX^c$ is
contained in some $\cX^c(F_n) \cong \cY$ and this decomposition
corresponds to \(\cY = \cY_+ \amalg \cY_0 \amalg \cY_{-}\) from
Theorem~\ref{thm:varconsY}.  In turn, viewing $\cY$ as a quotient of
$\cRbar_n$, this corresponds to the regions $\cRbar_n^{t > 0}$,
$\cRbar_n^{t = 0}$, and $\cRbar_n^{t < 0}$ respectively.

For the case of $F_n$ with its standard generators
$S = \{s_1, \ldots, s_n\}$, Theorem~\ref{thm: red is U0} gives for
each $c \in (-2, 2)$ a preferred diffeomorphism $\picred$ from
$\cX_0^c(F_n, S)$ to $X^c_{U_0}(F_n, S)$.  Combined with the
identifications provided by $\pi$, for $c \in (-2, 2)$ we have a fixed
decomposition of sets (not spaces) given by
\[
\cX^c(F_n, S) = X^\cirr_{\SU}(F_n, S) \amalg X^c_{U_0}(F_n, S) \amalg
X^\cirr_{\SLR}(F_n, S).
\]
For \(c \not \in (-2,2)\), we define
$\cX^c(F_n, S) = X^\cirr_{\SLR}(F_n, S)$, which matches
Theorem~\ref{thm:sX(F_n)}.

Now suppose $T = \{t_1, \ldots, t_m\}$ is the standard generating set
for $F_m$ and $\phi \maps F_n \to F_m$ an epimorphism with each
$\phi(s_i)$ for $s_i \in S$ conjugate to some element of $T$ or its
inverse.  Because $\phi$ is onto,
$\phi^* \maps R_\C(F_m, T) \to R_\C(F_n, S)$ takes each irreducible
$\SU$ or $\SLR$ representation to another such representation and the
same for those $U_0$ representations that are not fully reducible.
Thus, we can take the induced maps on the character variety pieces of
the above decompositions of $\cX^c(F_m, T)$ and $\cX^c(F_n, S)$ to
construct a \emph{set-theoretic} map
$\phi^* \maps \cX^c(F_m, T) \to \cX^c(F_n, S)$ for each $c \in \R$.
These combine to give a set-theoretic map
$\phi^* \maps \cX(F_m, T) \to \cX(F_n, S)$.  By our construction and
Theorem~\ref{thm: red is U0}, we have a commutative diagram
\begin{equation}
  \label{eq:comm of induced}
  \begin{tikzcd}
    \cX(F_m,T) \arrow{r}{\phi^*} \arrow{d}{\pi} &
    \cX(F_n,S) \arrow{d}{\pi} \\
    X(F_m, T) \arrow{r}{\phi^*} & X(F_n, S)
  \end{tikzcd}
\end{equation}
For the $\phi$ relevant here, we show in
Lemma~\ref{Lem:BraidActionOnF_n} that the top $\phi^*$ is in fact
smooth.

We will also consider epimorphisms where $F_n$ is replaced by
$\pi_1(S_n)$ and/or $F_m$ is replaced by $\pi_1(S_m)$.  Because
$\cX^c_0(S_n, S)$ is a proper subset of $X^c_{U_0}(S_n, S)$ for
$c \in (-2, 2)$, there is something to check to see that there is a map
$\phi^* \maps \cX^c(S_m, T) \to \cX^c(S_n, S)$, but this will hold in
our cases of interest.  For the map $\cX(S_m, T) \to \cX(S_n, S)$, one
must also verify the correct behavior with regards to unbalanced
reducibles.

Note that our induced maps on resolutions inherit functoriality
from their component maps on character varieties: for example if
$\alpha \maps F_n \to F_m$ and $\beta \maps F_m \to F_\ell$ are
suitable epimorphisms then
$(\beta \circ \alpha)^* = \alpha^* \circ \beta^*$.

We can now state the analogue of Proposition~\ref{prop: Casson style}
for our resolved picture:
\begin{theorem}
  \label{Prop:Lcap beta L}
  Suppose  \(K = \betahat\) for $\beta \in B_{2m}$ and fix
  $c \in (-2, 2)$.  The induced map of the epimorphism
  $\jhat_1 \maps \pi_1(S_{2m}) \to F_m$ is a well-defined smooth
  proper embedding $i' \maps \cX^c(F_m, T) \to \cX^c(S_{2m}, S)$ whose
  image $\cL^c$ is thus a  smooth closed submanifold.  Moreover, the
  map induced by the action of $\beta$ on $\cX^c(S_{2m}, S)$ is a
  well-defined diffeomorphism
  $\beta^*:\cX^c(S_{2m}) \to \cX^c(S_{2m}) $.  At the level of
  character varieties, the inclusion $S_{2m} \to \EK$ induces an
  injection of
  \[
    \cX^c(K) : = X^\cirr_\SU(K)  \amalg X^c_{U_0}(K) \amalg
    X^\cirr_\SLR(K),
  \]
  into $\cX^c(S_{2m}, S)$ whose image is exactly
  $\cX^c(\beta) := \cL^c \cap \beta^* \cL^c$.  Finally, $\cL^c$ is
  contained in $\cX^c_\bal(S_{2m}, S)$.
\end{theorem}
The proof of Theorem~\ref{Prop:Lcap beta L} will be given in
Section~\ref{sec: proof of Lcap beta L}, after establishing the
smoothness of the induced maps on various resolutions.

\subsection{The braid group action}

The braid group $ B_{n}$ acts on the \(n\)-punctured disk \(D_n\).
With the conventions of Figure~\ref{fig:braids}, the induced action on
\(\pi_1(D_n) = F_n\) is given by
\begin{equation}
  \label{eq:braid act s_i}
  \sigma_i \cdot s_j = \begin{cases} s_{i+1} \quad & j=i \\
    s_{i+1}^{-1} s_i s_{i+1} \quad & j=i+1 \\ s_j \quad& j\neq i,
    i+1 \end{cases}
\end{equation}
The action of \(\beta \in B_{n}\) preserves our preferred generators
$S$ up to conjugacy, and thus induces a regular real algebraic map
\( \beta^* \maps X(F_n, S) \to X(F_n, S)\) as discussed above.  This
is an isomorphism of real algebraic sets since $(\beta^{-1})^*$ is the
inverse to $\beta^*$.  Turning to \(S_n\), as
\(\sigma_i \cdot (s_1s_2\ldots s_n) = s_1s_2 \ldots s_n\) from
(\ref{eq:braid act s_i}), the action of \(B_n\) on \(F_n\) descends to
an action on \(\pi_1(S_n)\).  Thus for \(\beta \in B_{n}\), we have an
induced automorphism \(\beta^*\maps X(S_n,S) \to X(S_n,S)\) which can
be viewed as the restriction of \(\beta^*\) on $X(F_n, S)$ to its
subset $X(S_n, S)$.  We next show the induced maps on resolutions are
as expected:

\begin{lemma}
  \label{Lem:BraidActionOnF_n}
  For each $\beta \in B_n$, the induced map $\beta^*$ on $\cX(F_n,S)$
  is a diffeomorphism that preserves the submanifolds $\cX^c(F_n, S)$,
  $\cX_\bal(F_n, S)$, $\cX^c(S_n, S)$, and $\cX(S_n, S)$, where for
  the last two we are assuming that $n$ is even.  Hence $\beta^*$ is
  well-defined on each of these five resolutions, and in each case it
  is a diffeomorphism.
\end{lemma}

\begin{proof}
As the induced maps in Section~\ref{sec: induced res} are functorial,
it suffices to prove the lemma for each generator $\sigma_i$ of
$B_n$. On the open submanifold $X^\irr_\SLR(F_n, S)$ of $\cX(F_n, S)$, the
map $\sigma^*_i$ is a real algebraic automorphism which you can easily
check preserves the claimed submanifolds.  So it suffices to understand
$\sigma_i^*$ on the subset $B = (-2, 2) \times \cY$ of $\cX(F_n, S)$,
which we accomplish by lifting to $(-2, 2) \times \cRbar_n$.  

Motivated by (\ref{eq:braid act s_i}), define
\(\Phibar_i \maps (-2, 2) \times \cRbar_n \to (-2, 2) \times
\cRbar_n\) as follows, 
\begin{equation}
  \label{Eq:Braid Action}
  \Phibar_i(c,t, v_1,\dots v_n) =
  (c,t,v_1,\ldots, v_{i-1},v_{i+1}, \rho(-v_{i+1})\cdot v_i,v_{i+2},\ldots,v_n),
\end{equation}
where
\(\rho(-v_{i+1}) = \psi_t\big(A(\cos^{-1}(c/2),t(v),-v_{i+1})\big).\)
This is a smooth map as \(A(\alpha,t,v)\) is a smooth function of
\((\alpha,t,v)\) and \(G_t\) acts smoothly on \(\Qbar_t\).  Its
inverse is given by
\[
  \Phibar_i^{-1}(c,t, v_1,\dots v_n) = (c,t, v_1,\ldots, v_{i-1},
  \rho(v_{i})\cdot v_{i+1},v_i,v_{i+2},\ldots,v_n)
\]
which is also smooth. Hence $\Phibar_i$ is a diffeomorphism on
$(-2, 2) \times \cRbar_n$.

Let $\pi' \maps (-2,2) \times \cRbar_n \to R_\C(F_n, S)$ be as in
Section~\ref{sec:varying Fn}. Using equations~\eqref{eq:UActionOnA}
and \eqref{eq:CactionOnA}, you can check that \(\Phibar_i\) is
equivariant with respect to the actions of \(\R_{>0}\) and
\(G_t\). Hence it descends to a diffeomorphism
\(\Phi_i \maps (-2,2) \times \cY \to (-2, 2) \times \cY \) which we
next show is $\sigma_i^*$.  Again using \eqref{eq:UActionOnA} and
\eqref{eq:CactionOnA}, we see that the following commutes:
\begin{equation}
  \label{eq: Phibar}
  \begin{tikzcd}
    (-2, 2) \times \cRbar_n \arrow{r}{\Phibar_i} \arrow{d}{\pi'} &
    (-2, 2) \times \cRbar_n\arrow{d}{\pi'} \\
    R_\C(F_n, S) \arrow{r}{\sigma_i^*} & R_\C(F_n, S)
  \end{tikzcd}
\end{equation}
To check $\Phi_i = \sigma_i^*$, we fix $c \in (-2, 2)$ and study the
restriction of $\Phi_i$ to each of $\cY_{+}$, $\cY_{0}$, and $\cY_{-}$
in turn.  In the case of $\cY_{-}$, which maps to
$X^\cirr_\SLR(F_n, S)$ under $\pi_c$, by restricting $\Phibar_i$ to
$\cRbar_n^{t}$ for some $t < 0$, we can combine (\ref{eq: subsquare})
and (\ref{eq: Phibar}) to see that $\Phi_i = \sigma_i^*$ on $\cY_{-}$.
The argument for $\cY_{+}$ is the same using $t > 0$ and
$X^\cirr_\SU(F_n,S)$.  Similarly, $\cY_{0}$ is handled by considering
$t = 0$ and $X^c_{U_0}(F_n,S)$; see also the proof of
Theorem~\ref{thm: red is U0}.  Thus we have shown
$\Phi_i = \sigma_i^*$ and hence $\sigma_i^*$ is a self-diffeomorphism
of $\cX(F_n, S)$.

It remains to show $\sigma_i^*$ preserves the listed submanifolds.
For $\cX^c(F_n, S)$, this is immediate from the definition of
$\sigma_i^*$.  The case of $\cX_\bal(F_n, S)$ follows from  its definition
at the start of Section~\ref{sec:varies} and (\ref{eq:braid act s_i}).
For $\cX^c(S_n, S)$, the only issue is whether $\sigma_i^*$ preserves
$\cX^c_0(S_n, S)$ inside $\cX^c_0(F_n, S)$, as the former is a proper
subset of $X^c_{U_0}(S_n, S)$.  By Theorem~\ref{thm:cXSn}, we have
that $\cX^c_0(S_n, S)$ is the intersection of the closure of
$\cX^\cirr(S_n, S)$ with $\cX^c_0(F_n, S)$; since $\sigma_i$ is a
diffeomorphism preserving $\cX^\cirr(S_n, S)$, it thus follows that it
preserves $\cX^c_0(S_n, S)$.  The final case of $\cX(S_n, S)$ is the
union of various intersections of the previous submanifolds as
described in Theorem~\ref{thm:sX(S_n)}, and hence also preserved.
\end{proof}

\begin{lemma}
  \label{Lem:BraidActionOrientation}
  Both \(\beta^*\maps \cX^c(F_n,S) \to \cX^c(F_n,S) \) and
  \(\beta^*\maps \cX^c(S_n,S) \to \cX^c(S_n,S)\) are orientation
  preserving.
\end{lemma}

\begin{proof}
It is enough to check this for the actions of the elementary braids
\(\sigma_i\). To do this, we go back and look at the short exact
sequences used to define the orientations in
Section~\ref{Sec:Orientations}. The restriction of \(\Phibar_i\) in
(\ref{Eq:Braid Action}) to $\{c\} \times \cRbar_n$ satisfies
\(t\big(\Phibar_i(v)\big) = t(v)\), so \(\Phibar_i \) acts
trivially on the \(T\R\) factor in equation~\eqref{Eq:Orientation
  1}. The action of \(\Phibar_i\) on the fiber \(T\Qbar_t^n\) is the
composition of a coordinate permutation
\((v_1,\ldots, v_n) \mapsto (v_1, \ldots v_{i+1},v_i,\ldots,v_n)\) and
a rotation on one factor:
\((v_1, \ldots, v_n) \mapsto (v_1,\ldots, v_i,\rho(-v_i)\cdot v_{i+1},
\ldots, v_n)\). The permutation is orientation preserving since
\(\Qbar_t\) is 2-dimensional, and the rotation is
orientation-preserving since it is homotopic to the identity.  It
follows that the action of \(\Phibar_i\) on \(\cRbar_n\) is
orientation preserving. Finally, \(\Phibar_i\) is equivariant with
respect to the actions of \(\R_{>0}\), and \(G_t\), as observed in the
proof of Lemma~\ref{Lem:BraidActionOnF_n}. Hence the action of
\(\sigma _i^*\) on \(\cX^c(F_n,S)\) is orientation preserving.

To see that the same statement holds for \(\cX^c(S_n,S)\), we must
also consider the exact sequence~\eqref{Eq:Orientation 2}. From the
definition of \(\Phibar_i\), a direct computation shows that
\(F \circ \Phibar_i = F\). Referring to
equation~\eqref{Eq:Orientation 2}, we thus see that the action of
\(\Phibar_i\) on \(\cRbar^{c,t \neq 0}(S_n,S)\) is orientation
preserving. As above, the action of \(\R_{>0}\) and \(G_t\) are
equivariant with respect to \(\Phibar_i\), so the action of
\(\sigma_i^*\) on \(\cX^{c,\irr}(S_n,S)\) is orientation preserving as
well. Finally, \(\cX^{c,\irr}(S_n,S)\) is a dense open subset of
\(\cX^c(S_n,S)\), so the action on \(\cX^c(S_n,S)\) is orientation
preserving.
\end{proof}

\subsection{The character variety of the handlebody}
\label{subsec:def of cL}

Recall the epimorphism $i = \jhatstar_1$ from $\pi_1(S_{2m})$ to $F_m$
coming from the inclusion of $S_{2m}$ into $\partial H_1$. Let
$\itil \maps F_{2m} \to F_m$ be the lift of $i$ defined by
$\itil(s_{2i - 1}) = t_i$ and by $\itil(s_{2i}) = t_i^{-1}$. By
Section~\ref{sec: induced res}, $\itil$ has an induced map
$\itilstar \maps \cX^c(F_m, T) \to \cX^c(F_{2m}, S)$ for each
$c \in \R$.

\begin{proposition}
  \label{Prop:Lift Of i}
  For each $c \in (-2, 2)$, the induced map $\itilstar$ is a smooth
  proper embedding with image contained in $\cX^c_\bal(S_{2m}, S)$.
  Thus $i$ induces a well-defined smooth proper embedding
  $i' \maps \cX^c(F_m, T) \to \cX^c_\bal(S_{2m}, S)$.
\end{proposition}

\begin{proof}
We first construct a map \(j'\maps \cRbar_m \to \cRbar_{2m}\) and show
it induces $\itilstar$ after we take quotients by the $\R_{>0}$ and
$G_t$ actions.  Specifically, define $j'$ by
\begin{equation}
  \label{eq: j'}
  j'(t,v_1, \ldots, v_m) = (t,v_1,-v_1,\ldots, v_m,-v_m).
\end{equation}
The map $j'$ is the restriction of a smooth embedding
\(\R \times (\R^3)^m \to \R \times (\R^3)^{2m}\) to $\cRbar_m$ and
$\cRbar_{2m}$, and so is also a smooth embedding. Since
\(A(\alpha,t, -v) = A(\alpha,t,v)^{-1}\), we see that the following
diagram commutes:
\begin{equation}\label{eq: j' comm}
  \begin{tikzcd}
    \cRbar_m \arrow{r}{j'} \arrow{d}{\pi'_c} &
    \cRbar_{2 m}  \arrow{d}{\pi'_c} \\
    R^c_\C(F_m, T) \arrow{r}{\itilstar} & R^c_\C(F_{2m}, S)
  \end{tikzcd}
\end{equation}
It is easy to see that \(j'\) is equivariant with respect to the
actions of \(\R_{>0}\) and \(G_t\), and so it descends to a smooth
embedding \(\cX^c(F_m,T) \to \cX^c(F_{2m},S)\) by
Lemma~\ref{lem:quotient is manifold}.  By restricting $j'$ to various
$\cRbar_m^t$ as in the proof of Lemma~\ref{Lem:BraidActionOnF_n}, you
can see that this map \(\cX^c(F_m,T) \to \cX^c(F_{2m},S)\) must be
$\itilstar$.  Hence $\itilstar$ is a smooth embedding.

To understand the image of $\itilstar$, first note that
for the map 
$F \maps \cRbar_{2m} \to \SLC$ from Section~\ref{subsec: intro to punc
  sphere} we have:
\[
  F\big(j'(v)\big) = \prod_{i=1}^m A(\alpha,t, v_i)A(\alpha,t, -v_i) = I.
\]
Thus
\(j'\big(\cRbar_m^{t\neq 0}\big) \subset F^{-1}(I) \cap \cRbar_{2m}^{t
  \neq 0}\). Since the closure of
\(F^{-1}(I) \cap \cRbar_{2m}^{t\neq 0}\) is \(\cRbar^c(S_{2m},S)\) by
Lemma~\ref{lem:M epsilon}, we have
\(j'(\cRbar_m) \subset \cRbar^c(S_{2m},S)\). Also, in the notation of
Section~\ref{subsec:define res}, if
\(v \in Z_\epsilon\) then \(j'(v) \in Z_{\epsilon'}\), where
\(\epsilon' = (\epsilon_1, -\epsilon _1, \epsilon_2, \ldots,
\epsilon_n, -\epsilon_n)\) is balanced.  Combining, we see that the
image of \(\itilstar\) is contained in \(\cX^c_\bal(S_n, S)\) as
claimed.

It remains to show that $\itilstar$ is proper.  Now
\begin{equation*}
  \begin{tikzcd}
    \cX^c(F_m,T) \arrow{r}{\itilstar} \arrow{d}{\pi} &
    \cX^c(F_{2m},S) \arrow{d}{\pi} \\
    X^{c}(F_m, T) \arrow{r}{\itilstar} & X^{c}(F_{2m}, S).
  \end{tikzcd}
\end{equation*}
commutes as (\ref{eq: j' comm}) does, the lefthand projection is
proper by Theorem~\ref{thm:def of pi}, and the bottom \(\itilstar\) is
proper by Lemma~\ref{Lem:PropInc}. It follows that the top
\(\itilstar\) is proper, completing the proof.
\end{proof}

\begin{proposition}
  \label{Prop:Lift Of i varied}
  The induced map $\itilstar \maps \cX(F_m, T) \to \cX(F_{2m}, S)$ is
  smooth, and its  image is contained in $\cX(S_{2m}, S)$.  Thus $i$ induces a
  well-defined map $i' \maps \cX(F_m, T) \to \cX(S_{2m}, S)$.
  Moreover, $i'$ is a smooth proper embedding.
\end{proposition}

\begin{proof}
Recall from Section~\ref{sec:varying Fn} that
\(\cX(F_m,T) = B \cup C\), where \(B \cong (-2, 2) \times \cY\) and
\(C \cong X^\irr (F_m,T)\) are open in \(\cX(F_m,T)\), as well as the
analogous decomposition of $\cX(S_{2m}, S)$ as $B' \cup C'$ from
Section~\ref{sec:varying Sn}.  The map $\itilstar$ restricted to $C$
has image contained in $C' \cong X^\irr(S_{2m}, S)$ by construction,
and is smooth as it is a regular real algebraic map.  Smoothness of
$\itilstar$ on $B$ follows as in the proof of
Proposition~\ref{Prop:Lift Of i} by using the map $(-2, 2)
\times \cRbar_m \to (-2, 2) \times \cRbar_{2m}$ sending $(c, v)$ to
$(c, j'(v))$, and this also shows that $\itilstar$ is an embedding.
The claim about the image of $\itilstar$ is an immediate
consequence of Proposition~\ref{Prop:Lift Of i}.

To show $i'$ is proper, recall from just before
Proposition~\ref{prop:proper 1} the definitions of $X_\circ(F_m, T)$
and $\pi_\circ \maps \cX(F_m, T) \to X_\circ(F_m, T)$ and the
analogous notions for $S_{2m}$.  Then 
\[
  \begin{tikzcd}
    \cX(F_m,T) \arrow{r}{i'} \arrow{d}{\pi_\circ} &
    \cX(S_{2m},S) \arrow{d}{\pi_\circ} \\
    X_\circ(F_m, T) \arrow{r}{i^*} & X_\circ(S_{2m}, S).
  \end{tikzcd}
\]
commutes by (\ref{eq:comm of induced}).  To see that $i'$ is proper,
note that the lefthand copy of \(\pi_\circ\) is proper by
Proposition~\ref{prop:proper 1} and \(i^*\) is proper by applying
Lemma~\ref{Lem:PropInc} to $i^* \maps X(F_m, T) \to X(S_{2m}, S)$ and
then analyzing its restriction.
\end{proof}

\subsection{Understanding the intersection}
\label{sec: proof of Lcap beta L}

We can now give:

\begin{proof}[Proof of Theorem~\ref{Prop:Lcap beta L}]

Proposition~\ref{Prop:Lift Of i} gives the claims about
$i' \maps \cX^c(F_m, T) \to \cX^c(S_{2m}, S)$ and $\cL^c$, including
the final claim that $\cL^c \subset \cX^c_\bal(S_{2m}, S)$.  The claim
about $\beta^*$ is immediate from Lemma~\ref{Lem:BraidActionOnF_n}.
By definition all the induced maps on our resolutions respect the
various decompositions $\cX = \cX_{+} \coprod \cX_0 \coprod \cX_{-}$
which in turn correspond to $\SU$, $U_0$, and $\SLR$ representations.
Applying the proof of Proposition~\ref{prop: Casson style} to each
piece of the induced decomposition
$\cL^c = \cL^c_+ \amalg \cL^{c}_{0} \amalg \cL^c_-$, we see that the
inclusion $S_{2m} \to \EK$ induces the promised bijection of $\cX^c(K)$
with $\cX^c(\beta)$.
\end{proof}

By Proposition~\ref{Prop:Lift Of i varied},
\(\cL := i'\big(\cX(F_m,T)\big)\) is a closed \((2m-2)\)-dimensional
submanifold of the \((4m-5)\)-dimensional manifold \(\cX(S_{2m},S)\).
Recall from Theorem~\ref{thm:sX(S_n)} the map
$\mtr \maps \cX(S_{2m},S) \to \R$, where $\mtr^{-1}(a)$ is
$\cX^a_\bal(S_{2m}, S)$ for $a \in (-2, 2)$ and
$X^{a,\irr}_\SLR(S_{2m}, S)$ otherwise.  For uniformity of notation,
we now define $\cX^a_\bal(S_{2m}, S) := \mtr^{-1}(a)$ for
$\abs{a} \geq 2$.  If \([a,b] \subset \R\) is an interval, we let
\(\cX^{[a,b]}(S_{2m},S) = \mtr^{-1}([a,b]) \subset \cX(S_{2m},S)\). It
is a (noncompact) manifold with boundary
\(\cX^a_\bal(S_{2m}, S) \amalg \cX^b_\bal(S_{2m}, S)\). We also set
\(\cL^{[a,b]} = \cL \cap \cX^{[a,b]}(S_{2m},S)\). The composition
\(\mtr \circ i'\) is the usual trace map on \(\cX(F_m,T)\), hence a
submersion. Thus \(\cL^{[a,b]}\) is a closed proper submanifold of
\(\cX^{[a,b]}(S_{2m},S)\).

We say \(\beta\in B_{2m}\) is real representation small if its plat
closure \(K = \betahat\) is real representation small in the sense of
Section~\ref{sec: knot reps}.

\begin{lemma}
  \label{lem: rrs gives compact}
  If \(\beta\) is real representation small, the intersections
  \(\cL^c \cap \beta^*\cL^c\) for $c \in [-2, 2]$ and
  \(\cL^{[-2,2]} \cap \beta^* \cL^{[-2,2]}\) are compact.
\end{lemma}

\begin{proof}
As \(\cL^c \cap \beta^*\cL^c\) is a closed subset of
\(\cL^{[-2,2]} \cap \beta^* \cL^{[-2,2]}\), it suffices to prove that
the latter is compact.  Setting $\cX(\beta) := \cL \cap \beta^* \cL$,
our goal is to prove that $\cX^{[-2, 2]}(\beta)$ is compact.
Throughout, we will view $\cX(\beta)$ as a subset of $\cX(F_{2m}, S)$
rather than $\cX(S_{2m}, S)$.

Recall from Section~\ref{sec: prop of proj} the surjection
$\pi_\circ \maps \cX(F_n, T) \to X_\circ(F_n, T)$.  Then (\ref{eq:comm
  of induced}) gives that $\pi_\circ(\cL) = L \cap X_\circ(F_{2m}, S)$
and $\pi_\circ(\beta^* \cL) = (\beta^* L) \cap X_\circ(F_{2m}, S)$.
Let $Z = \pi_\circ \big(\cX(\beta)\big)$.  By Proposition~\ref{prop:
  Casson style}, $Z \subset X(K)$, where we have identified $X(K)$
with a subset of $X(S_{2m}, S) \subset X(F_{2m}, S)$ via the map
induced by $S_{2m} \hookrightarrow \EK$.

Now \(Z\) is not all of $X(K)$, as $X(K)$ contains the line of
reducible representations $\chi^c$ for all $c \in \R$ by
Section~\ref{sec: alex and the reds}, while $Z $ is contained in
$ X_\circ(F_{2m}, S)$ which omits all reducibles with trace not in
\((-2,2)\).  However, $Z$ does contain
$X^\irr(K) = L^\irr \cap \beta^* L^\irr$.  Moreover, by
Theorem~\ref{Prop:Lcap beta L}:
\[
Z = X^\irr(K) \cup \setdefm{\big}{\chi^c}{\mbox{$c \in (-2, 2)$ and
    $X^c_{U_0}(K) \neq \emptyset$}}.
\]
The space $X^c_{U_0}(K)$ is nonempty exactly when there is a
representation $\rho \maps \pi_1(\EK) \to U_0$ with
$\tr_\mu(\rho) = c$ which is not fully reducible; this is equivalent
to having a reducible representation $\pi_1(\EK) \to \SL{2}{\C}$ with
nonabelian image and character $\chi^c$. Thus by Section~\ref{sec:
  alex and the reds}, we have
$Z = X^\irr(K) \cup \setdefm{\big}{\chi^c}{c \in D_K \cap (-2, 2)}$,
where recall $D_K$ is finite.  Equivalently, since $\pm 2 \notin D_K$
as $\Delta_K(\pm 1)$ is an odd integer,
$Z = X^\irr(K) \cup \setdefm{\big}{\chi^c}{c \in D_K \cap [-2, 2]}$.
Now Lemma~\ref{lemma: deform red} implies that any limit point of
$X^\irr(K)$ in $X(K)$ is $\chi^c$ for some $c \in D_K$ and hence
$Z^{[-2, 2]}$ is closed in $X^{[-2, 2]}(K)$.  As $K$ is real
representation small, $X^{[-2, 2]}(K)$ is compact by definition. Hence
so is its closed subset $Z^{[-2, 2]}$, which is equal to
$\pi_\circ\left(\cX^{[-2, 2]}(\beta)\right)$.

From the form of $\jhat_1 \maps \pi_1(S_{2m}) \to F_m$, any reducible
representation in $L$ is balanced.  Thus,
$Z^{[-2, 2]} \subset X(K) \subset L$ in fact lies in
$X_{\bal, \circ}(F_{2m}, S)$.  As
\(\pi_\circ:\cX(F_{2m}, S) \to X_\circ(F_{2m}, S)\) is proper by
Proposition~\ref{prop:proper 1}, we see that
$\pi_\circ^{-1}\left( Z^{[-2, 2]}\right)$ is both compact and a subset
of $\cX_\bal(F_{2m}, S)$.  Now $\cL$ and $\beta^* \cL$ are closed in
$\cX(S_{2m}, S)$ by Theorem~\ref{Prop:Lcap beta L}, so their
intersection $\cX(\beta)$ is closed in $\cX(S_{2m}, S)$.  As the
latter is closed in $\cX_\bal(F_{2m}, S)$ by
Theorem~\ref{thm:sX(S_n)}, $\cX(\beta)$ is  closed in
$\cX_\bal(F_{2m}, S)$.  So $\cX^{[-2, 2]}(\beta)$ is a closed
subset of the compact set $\pi_\circ^{-1}\left( Z^{[-2, 2]}\right)$,
and is hence compact.
\end{proof}

So that we can break up the intersection number defining $h^c(\beta)$
into ``local'' pieces to prove invariance in Section~\ref{sec:
  invariance}, we note the following: 

\begin{lemma}
  \label{lem: rrs gives finite}
  If \(\beta\) is real representation small, the intersections
  \(\cL^c \cap \beta^*\cL^c\) for $c \in [-2, 2]$ have finitely many
  connected components.
\end{lemma}

\begin{proof}
For $c = \pm 2$, this follows because
$X^\cirr_\SLR(K) = \cL^ c \cap \beta^*\cL^c$ is a real semialgebraic
subset of the real semialgebraic set
$X^\cirr_\SLR(S_{2m}, S) = \cX^c_\bal(S_{2m}, S)$ and real
semialgebraic sets have finitely many connected components.  For
$c \in (-2, 2)$, consider the decomposition
\( \cX^c(K) = X^\cirr_\SU(K) \amalg X^c_{U_0}(K) \amalg
X^\cirr_\SLR(K) \) of $\cL^ c \cap \beta^*\cL^c$ from
Theorem~\ref{Prop:Lcap beta L}.  Each part of the decomposition is
real semialgebraic, and the inclusion of it into
$\cX^c_\bal(S_{2m}, S)$ is continuous; therefore
\(\cL^c \cap \beta^*\cL^c\) has finitely many connected components as
needed.
\end{proof}

\subsection{The total Lin invariant}

We orient \(\cX(F_m,T)\) and \(\cX(S_{2m},S)\) as in
Section~\ref{Sec:Orientations} and give \(\cL\) the orientation
inherited from \(\cX(F_m,T)\).  The slices
\(\cX^c(F_m,T), \cX^c_\bal(S_{2m},S)\), and $\cL^c$ are then oriented
using the submersions $\mtr$, where the range $\R$ has the reversed
orientation, see Theorem~\ref{thm: res for F_n orient} and
Corollary~\ref{cor: orient cX(S_n)}.  When $\beta$ is real
representation small, Lemma~\ref{lem: rrs gives compact} and
Theorem~\ref{thm:intersection} show that for each $c \in [-2, 2]$
there is a well-defined intersection number $\pairLbetaL$ and so we
set
\[
  h^c(\beta) := (-1)^{|\beta|} \pairLbetaL
\]
where \(|\beta|\) is the number of crossings in \(\beta\). Our setup
now allows us to easily show $h^c(\beta)$ does not depend on $c$:
\begin{theorem}
  \label{thm: h indep of c}
  If \(\beta\) is real representation small,
  \(h^c(\beta) = h^{c'}(\beta)\) for all \(c,c' \in [-2,2]\).
\end{theorem}

\begin{proof} We apply Theorem~\ref{Thm:InvarianceOfIntersection} with
\(M = \cX^{[-2,2]}(S_{2m})\), \(A= \cL^{[-2,2]}\),
\(B = \beta^*(\cL^{[-2,2]})\), and \(\pi=\mtr\).
\end{proof}

Given the theorem, if \(\beta \in B_{2m}\) is real representation
small, we define the \emph{total Lin invariant} \(h(\beta)\) as
\(h^c(\beta)\) for any \(c \in [-2,2]\). Next, we define the \(\SU\)
and \(\SLR\) Lin invariants. As in Section~\ref{sec: proof of Lcap
  beta L}, we write:
\[
  \cL^c = \cL^c_- \amalg \cL
  ^{c}_{0} \amalg \cL^c_+
\]
From Theorem~\ref{Prop:Lcap beta L}, we know that
\(\cL^c_+ \cap \beta^*(\cL^c_+) = X^c_\SU(K)\) and similarly for
\(\SLR\). Moreover, if $\cL^c \cap \beta^*(\cL^c)$ is compact and
\(c \not \in D_K\) (so \(X_{U_0}^c(K)\) is empty), then both
\(\cL^c_+ \cap \beta^*(\cL^c_+)\) and
\(\cL^c_- \cap \beta^*(\cL^c_-)\) will be compact. Hence, if \(\beta\)
is real representation small and
\(c \ \in [-2,2] \setminus D_\betahat\), we define
\begin{align*}
  h_\SLR^c(\beta) &= (-1)^{|\beta|} \pairLbetaLSLR \quad \text{and} \\
  h_\SU^c(\beta) & = (-1)^{|\beta|}\pairLbetaLSU.
\end{align*}

\begin{theorem}
\label{thm: hSL+hSU}
If \(\beta\) is real representation small, then
$h(\beta) = h_\SLR^c(\beta) + h_\SU^c(\beta)$ for all
\(c \in [-2,2] \setminus D_\betahat\). Moreover, $h_\SLR^c(\beta)$ and
$h_\SU^c(\beta)$ are constant on each connected component of
\( [-2,2] \setminus D_\betahat\).  Finally, $\pm 2 \not \in
D_\betahat$ and  $h^{\pm 2}_\SU(\beta) = 0$ and
$h^{\pm 2}_\SLR(\beta) = h(\beta)$.
\end{theorem}

\begin{proof}
To compute $h^c(\beta)$ as per Theorem~\ref{thm:intersection}, we
perturb $\cL^c$ to be transverse to $\beta^* \cL^c$ by an isotopy
supported in an arbitrarily small neighborhood of their intersection.
As \(c \notin D_\betahat\), we have \(X_{U_0}^c(\beta) = \emptyset\),
and by Theorem~\ref{Prop:Lcap beta L} the initial intersection of
$\cL^c$ with $\beta^* \cL^c$ is disjoint from $\cX^c_0(S_{2m}, S)$.
Thus we can perturb $\cL^c$ to a transverse $\cL'$ without changing
the intersection with $\cX^c_0(S_{2m}; S)$.  We then use $\cL'$,
$\cL'_{-}$, and $\cL_+'$ to compute $h^c(\beta), h^c_\SLR(\beta)$, and
$h^c_\SU(\beta)$ respectively.  Each point of
$\cL' \cap \beta^* \cL^c$ contributes $\pm 1$ to $h^c(\beta)$ and the
same to exactly one of $h^c_\SLR(\beta)$ and $h^c_\SU(\beta)$. This
proves the first claim.

For the second statement, suppose $[a, b]$ is disjoint from $D_\betahat$.
Then $\cL^{[a, b]}$ is disjoint from $\cX^{[a, b]}_0(S_{2m}, S)$ so we can
study $\cL^{[a, b]}_{-}$ and $\cL^{[a, b]}_{+}$ separately using the
proof of Theorem~\ref{thm: h indep of c} to see that
$h^c_\SU(\beta)$ and $h^c_\SLR(\beta)$ are constant on $[a, b]$.

Finally, we observed in Section~\ref{sec: alex and the reds} that
$\pm 2 \notin D_K$. By Theorem~\ref{thm:sX(S_n)}, we have
$\cX^{\pm 2}(S_{2m}, S) = X^{\pm 2, \irr}_\SLR(S_{2m}, S)$, so
$h^{\pm 2}_\SU(\beta) = 0$ and $h^{\pm 2}_\SLR(\beta) = h(\beta)$ by
definition.
\end{proof}

\section{Proof of invariance}
\label{sec: invariance}

The goal of this section is to prove:

\begin{theorem}
  \label{Thm:h is invariant}
  Suppose \(\beta \in B_{2m}\) and \(\beta' \in B_{2m'}\) with
  \(\betahat\) and \(\betahatprime\) isotopic to the same knot $K$.
  If $K$ is real representation small then \(h(\beta) = h(\beta')\).
  Moreover, for all $c $ in $[-2, 2] \setminus D_K$, we have
  \(h^c_\SLR(\beta) = h^c_\SLR(\beta')\) and
  \(h^c_\SU(\beta) = h^c_\SU(\beta')\).
\end{theorem}
Thus when \(K\) is real representation small, we define the
\emph{total Lin invariant} of $K$ by \(h(K) := h(\beta)\) where
\(\beta\) is any braid with \(\betahat = K\), and similarly
\(h^c_\SLR(K) := h^c_\SLR(\beta)\) and
\(h^c_\SU(K) := h^c_\SU(\beta)\).  Combined with Theorem~\ref{thm: hSL+hSU},
Theorem~\ref{Thm:h is invariant} gives Theorems~\ref{thm: hSLR intro}
and~\ref{Thm:MainTheorem} from the introduction.

To set up the proof, recall that
\(h^c(\beta) : = (-1)^{|\beta|} \pairLbetaL \). This set has finitely
many connected components by Lemma~\ref{lem: rrs gives finite}. Thus we
can write $h^c(\beta)$ in terms of local intersection numbers as in
(\ref{eq: component sum}):
\begin{equation}
\label{eq:local intersection}
  \pairLbetaL = \sum_{Z}  \pairLbetaLatZ
\end{equation}
where \(Z\) runs over the set of connected components of
\( \cL^c \cap \beta^* \cL^c\).  If \(Z\) is a connected component of
\(\cX^c(\beta)\), we define
\[
  n_{Z,\beta} = (-1)^{|\beta|} \pairLbetaLatZ.
\]

We will prove that the \(n_{Z,\beta}\) are invariants of the
underlying knot \(K\); to make this precise, we use the setup from
Section~\ref{sec: braids and plats} and Proposition~\ref{prop: Casson
  style}.  Given a knot $K = \betahat$ for some \(\beta \in B_{2m}\),
let $\EK$ be the exterior of $K$.  The inclusion $S_{2m} \to \EK$,
where $S_{2m} = \partial H_1$, induces an epimorphism
$\pi_1(S_{2m}) \to \pi_1(\EK)$. Concretely, using the standard
generators of
$\pi_1(S_{2m}) = \spandef{s_1,\ldots,s_{2m}}{w_{2m} = 1}$ where
\(w_{2m} = s_1\ldots s_{2m}\), we get a group presentation for
$\pi_1(\EK)$, namely:
\[
  \Pi_\beta = \PiPresentation{w_{2m}, H, \beta_*^{-1}(H)}
  \mtext{where \(H\) is the kernel of \(F_{2m} \to \pi_1(H_1)\).}
\]
By construction, \(\Pi_\beta\) comes
equipped with a preferred map \(i_\beta: F_{2m} \to
\Pi_\beta\), as well as an isomorphism $\Pi_\beta \to \pi_1(\EK)$ induced
by $D_{2m} \hookrightarrow S_{2m} \hookrightarrow \EK$.

Motivated by Theorem~\ref{Prop:Lcap beta L}, for $c \in [-2, 2]$ we
define
\[
  \cX^c(\Pi_\beta) := X^\cirr_\SU(\Pi_\beta, S)
  \amalg X^c_{U_0}(\Pi_\beta, S)
  \amalg X^\cirr_\SLR(\Pi_\beta, S).
\]
When $c = \pm 2$, since $\Pi_\beta$ is generated by $S$ then, by
Lemma~\ref{lem: empty char vars} and its analogue for $U_0$, the sets
$X^\cirr_\SU(\Pi_\beta, S)$ and $X^c_{U_0}(\Pi_\beta, S)$ are empty;
thus $\cX^c(\Pi_\beta) = X^\cirr_\SLR(\Pi_\beta, S)$ when $c = \pm 2$.

As in Section~\ref{sec: induced res}, we have an induced map
\(i_\beta^*: \cX^c(\Pi_\beta) \to \cX^c(F_{2m}, S)\) at the level of
sets only, not topological spaces.  For $c \in (-2, 2)$,
Theorem~\ref{Prop:Lcap beta L} implies $i_\beta^*$ is injective and its
image is \(\cX^c(\beta):= \cL^c \cap \beta^* \cL^c\); the same holds
for $c = \pm 2$ where $\cX^c_\bal(S_{2m}, S) = X^c_\SLR(S_{2m}, S)$ by
Proposition~\ref{prop: Casson style}.

If \(\phi: \Pi_{\beta'} \to \Pi_{\beta}\) is an isomorphism which
takes the conjugacy class of a meridian in \(\Pi_{\beta'}\) to the
conjugacy class of a meridian (or its inverse) in \(\Pi_{\beta}\), we
get a bijection \(\phi^*: \cX^c(\Pi_\beta) \to \cX^c(\Pi_{\beta'})\)
by considering the induced maps on character varieties. (In fact, any
isomorphism \(\phi\) must have this property by
\cite{GordonLuecke1989}; here, this will be obvious without appealing
to such deep results.)  We will prove:

\begin{theorem}
\label{Thm:n_Z is invariant}
Suppose \(\beta \in B_{2m}\) and \(\beta' \in B_{2m'}\) with
\(\betahat=\betahatprime = K\) and $c \in (-2, 2)$. Then there is an
isomorphism \(\phi: \Pi_{\beta'} \to \Pi_{\beta}\) and a homeomorphism
\(\psi: \cX^c(\beta) \to \cX^c(\beta')\) such that there is a
commuting square:
\[
  \begin{tikzcd}
    \cX^c(\Pi_{\beta}) \arrow{d}{\phi^*} \arrow{r}{i_\beta^*}  & \cX^c(\beta)  \arrow{d}{\psi} \\
      \cX^c(\Pi_{\beta'})  \arrow{r}{i_{\beta'}^*} &  \cX^c(\beta')
  \end{tikzcd}
\]
If in addition \(K\) is real representation
small and  \(Z\) is a connected component of \ \(\cX^c(\beta)\), then
\(n_{Z,\beta} = n_{\psi(Z),\beta'}\).
\end{theorem}

Assuming Theorem~\ref{Thm:n_Z is invariant}, we can show:

\begin{proof}[Proof of Theorem~\ref{Thm:h is invariant}]

For all $c$ in the open interval $(-2, 2)$, we combine the formula
\eqref{eq:local intersection} with Theorem~\ref{Thm:n_Z is invariant}
to get $h^c(\beta) = h^c(\beta')$, and then Theorem~\ref{thm: h indep
  of c} gives $h(\beta) = h(\beta')$.  Provided $c$ is further not in
$D_K$, we also get the corresponding statements for $h^c_{\SU}$ and
$h^c_{\SLR}$.  The case of $c = \pm 2$ then follows from the last
claim of Theorem~\ref{thm: hSL+hSU}.
\end{proof}

Theorem~\ref{Thm:n_Z is invariant} extends to $c = \pm 2$ as well:

\begin{theorem}
  \label{Thm:n_Z at ends}
  Suppose \(\beta \in B_{2m}\) and \(\beta' \in B_{2m'}\) with
  \(\betahat = \betahatprime = K\) and $c = \pm 2$.  Then the conclusions of
  Theorem~\ref{Thm:n_Z is invariant} hold.
\end{theorem}

\subsection{Plat moves}

Our proof of Theorem~\ref{Thm:n_Z is invariant} follows that of
Heusener \cite{Heusener2003}. (Lin's original proof \cite{Lin1992} is
similar but uses braid closures.)  The basic idea is this: first show
that if \(\betahat = \betahatprime\), then \(\beta\) and \(\beta'\)
are related by a sequence of moves, and then check invariance under
these moves. The first step was accomplished by Birman and Hilden:

\begin{theorem}[{\cite{Birman1976, Hilden1975}}]
\label{Thm:PlatMoves}

If \(\betahat\) and \(\betahatprime\) are isotopic, then \(\beta\)
and \(\beta'\) are related by a sequence of the following moves and
their inverses:
\begin{itemize}
\item (Type I) Replace \(\beta \in B_{2m}\) by \(\alpha \beta\) or
  \(\beta \alpha\), where \(\alpha \in B_{2m} \) is one of the
  following \emph{type I braids}:
  \(\sigma_1,\sigma_2\sigma_1^2 \sigma_2\), or
  \( \sigma_{2j} \sigma_{2j-1} \sigma_{2j+1} \sigma_{2j}\) for
  \(1\leq j \leq m-1\).
\item (Type II) Replace \(\beta \in B_{2m}\) by \(\beta \sigma_{2m} \in B_{2m+2}\).
\end{itemize}
\end{theorem}
These moves are illustrated in Figure 3 of \cite{Heusener2003}.

\subsection{Invariance under type I moves}

Suppose \(\beta' = \alpha \beta\) where \(\alpha\) is a type I
braid. Then \(\alpha_*(H) = H\), and so
\[
  \Pi_{\beta'} =
  \PiPresentation{w_{2m},H, \beta_*^{-1}\alpha_*^{-1} H}
  = \PiPresentation{w_{2m}, H, \beta_*^{-1} H} =
  \Pi_\beta
\]
and the isomorphism \(\phi: \Pi_{\beta'} \to \Pi_\beta\) is induced by
\(1_{F_{2m}}\). On the other hand, if \(\beta' = \beta \alpha \), then
\(\Pi_\beta' = F_{2m}/\pair{w_{2m}, H, \alpha_*^{-1} \beta_*^{-1}
  H}\), and the isomorphism \(\phi\) is induced by the map
\(\alpha_*: F_{2m} \to F_{2m}\) since $\alpha_*(w_{2m}) = w_{2m}$.

Next, we construct the homeomorphism \(\psi\) for a given
$c \in (-2, 2)$.  Our isomorphism \(\phi:\Pi_{\beta'} \to \Pi_\beta\)
is induced by an isomorphism \(\phihat: F_{2m'} \to F_{2m}\), where
each \(\phihat(s_i)\) is conjugate to some \(s_j\) or its inverse.
Thus by Section~\ref{sec: induced res}, there is a well-defined set
map \(\phihatstar: \cX^c(F_{2m}) \to \cX^c(F_{2m'})\) with a commutative
square
\begin{equation}
  \label{eq: pi beta sq}
  \begin{tikzcd}
    \cX^c(\Pi_{\beta}) \arrow{d}{\phi^*} \arrow{r}{i_\beta^*}  & \cX^c(F_{2m})  \arrow{d}{\phihatstar} \\
      \cX^c(\Pi_{\beta'})  \arrow{r}{i_{\beta'}^*} &  \cX^c(F_{2m'})
  \end{tikzcd}
\end{equation}
so \( \phihatstar\) maps \(\cX^c(\beta) = \im i_\beta^*\) bijectively
to \(\cX^c(\beta') = \im i_{\beta'}^*\).  Thus if we define
\(\psi = \phihatstar\), we get a commutative square as in
Theorem~\ref{Thm:n_Z is invariant}. When \(\beta' = \alpha \beta\),
the map $\phihatstar$ is the identity, and hence \(\psi\) is a
homeomorphism; when \(\beta' = \beta \alpha \), this follows from:

\begin{proposition}
\label{Prop:MoveI}
If \(\alpha \in B_{2m}\) is a braid of type I and $c \in (-2, 2)$, the
map \(\alpha^*\maps \cX^c(S_{2m},S) \to \cX^c(S_{2m},S)\) is a
diffeomorphism which fixes \(\cL^c\) setwise.  The map \(\alpha^*\)
reverses the orientation on \(\cL^c\) when $\alpha=\sigma_1$;
otherwise, it preserves it.
\end{proposition}

\begin{proof}

Recall from (\ref{eq: j'}) that
\(i'\maps \cX^c(F_m,T) \to \cX^c(S_{2m},S)\) is induced by
\(j'\maps \cRbar_{m} \to \cRbar_{2m}\), where
\(j'(v_1, \ldots v_m) = (v_1, -v_1, \ldots ,v_m,-v_m)\). Similarly,
\(\alpha^*\) on \(\cX^c(S_{2m},S)\) is induced by the map
\(\alpha^*\maps \cRbar_{2m} \to \cRbar_{2m}\) given
explicitly by the formulae in the proof of
Lemma~\ref{Lem:BraidActionOnF_n}. In each case, we check that
\(\alpha^*\circ j' = j' \circ \Phi_\alpha\) for some diffeomorphism
\(\Phi_\alpha\maps \cRbar_m \to \cRbar_m\).

To start, from (\ref{Eq:Braid Action}) we have
\[
  \sigma_1^*(v_1, -v_1, \ldots, v_m,-v_m) = (-v_1, v_1, v_2, -v_2,
  \ldots, v_m,-v_m),
\]
so we can take
\(\Phi_{\sigma_1}(v_1,\ldots, v_m) = (-v_1,v_2, v_3, \ldots, v_m)\).
The map \(\Qbar_t \to \Qbar_t \) given by \(v \mapsto -v\) is
orientation reversing, so \(\Phi_{\sigma_1}\) is orientation
reversing.  The map \(\Phi_{\sigma_1}\) commutes with the actions of
\(\R_{>0}\) and the \(G_t\), so it descends to an
orientation reversing diffeomorphism
\(\varphi_{\sigma_1}\maps \cL^c \to \cL^c\), which satisfies
\(\sigma_1^* \circ i' = i' \circ \varphi_{\sigma_1}\).

Similarly, we compute
\[
  (\sigma_2 \sigma_1^2 \sigma_2)^*(v_1, -v_1, v_2, -v_2, \ldots, v_m, -v_m)
  = (v_1',-v_1',v_2,-v_2, \ldots, v_m,-v_m)
\]
where
\(v_1' = \rho(-v_2) \cdot v_1 
\), and
\begin{multline}
(\sigma_{2j} \sigma_{2j-1} \sigma_{2j+1} \sigma_{2j})^*(v_1, -v_1, \ldots, v_m, -v_m) = \\
(v_1,-v_1, \ldots,v_{j+1},-v_{j+1},v_j,-v_j, \ldots, v_m,-v_m)
\end{multline}
so we take
\(\Phi_{\sigma_2 \sigma_1^2 \sigma_2} (v_1, \ldots, v_m) = (v_1',v_2,
\ldots, v_m)\) and
\[
  \Phi_{\sigma_{2j} \sigma_{2j-1} \sigma_{2j+1} \sigma_{2j}} (v_1, \ldots, v_m)
  = (v_1, \ldots, v_{j+1}, v_j, \ldots, v_m).
\]
To see these $\Phi$ preserve orientation, note that the first map is
effectively a rotation in the first $\Qbar_t$ coordinate (so homotopic
to the identity) and the second  interchanges two such
2-dimensional coordinates; see the proof of
Lemma~\ref{Lem:BraidActionOrientation} for complete details.
\end{proof}

\begin{corollary}
\label{Cor:MoveI}
Theorem~\ref{Thm:n_Z is invariant} holds when   \(\beta\)
and \(\beta'\) are related by a type I move.
\end{corollary}

\begin{proof}

We have already constructed \(\phi\) and \(\psi\) giving the
commutative square, so it remains to check the last statement.
If \(Z\) is a component of \(\cX^c(\beta)\), we write
\(Z' = \psi(Z)\).  Suppose \(\alpha = \sigma_1\), and let
\(\beta' = \alpha \beta\). Then we have
\[
  \pairLbraidLatZp{(\alpha \beta)^*} =
  \pairLbraidLatZp{\beta^* \alpha^*} =
  -\pairLbraidLatZ{\beta^*}
\]
by Proposition~\ref{Prop:MoveI}. Since
\(|\alpha \beta|= |\beta| + 1\), we see that
\(n_{Z',\beta'} = n_{Z,\beta}\).  Similarly, if \(\beta' = \beta \alpha\),
\begin{align*}
  \pairLbraidLatZp{(\beta \alpha)^*} &=
  \pairLbraidLatZp{\alpha^* \beta^*} \\
  &= \pairLbraidLatZ[(\alpha^*)^{-1}]{\beta^*} =
  -\pairLbraidLatZ{\beta^*}
\end{align*}
where in the second step we have used the fact that
\(\alpha^*\maps \cX^c(S_{2m},S) \to \cX^c(S_{2m},S)\) is an
orientation preserving diffeomorphism. Hence
\(n_{Z',\beta'} = n_{Z,\beta}\).

The proof for other \(\alpha\) of type I is exactly the same, except
that \(|\alpha \beta| \equiv |\beta| \bmod 2 \) and
\(\alpha^*\maps \cL^c \to \cL^c\) is orientation preserving, so rather
than a pair of canceling signs, we get no signs at all.
\end{proof}

\subsection{Invariance under the type II move} Given
\(\beta \in B_{2m}\), let \(\beta' = \beta \sigma_{2m} \in
B_{2m+2}\). If \(\Pi_\beta = \spandef{ s_1, \ldots, s_{2m}}{R}\), where
\(R\) is the set of relators, then
\[
  \Pi_{\beta'} = \spandef{ s_1, \ldots, s_{2m}, s_{2m+1}, s_{2m+2}}{R, s_{2m+1}s_{2m+2}, s_{2m}s_{2m+2}}.
\]
We define \(\phihat: F_{2m+2} \to F_{2m}\) by \(\phihat(s_i) = s_i\)
for \(i\leq 2m\), by \(\phihat(s_{2m+1}) = s_{2m}\), and by
\(\phihat(s_{2m+2}) = s_{2m}^{-1}\), so \(\phihat\) induces an
isomorphism \(\phi: \Pi_{\beta'} \to \Pi_\beta\). As in the previous
subsection, we define \(\psi = \phihatstar\) where $\phihatstar$ is as in
(\ref{eq: pi beta sq}); the restriction of $\psi$ to $\cX^c(\beta)$
will give the commutative square required by Theorem~\ref{Thm:n_Z is
  invariant}.

We now establish the properties of \(\psi\) needed to complete the proof of
Theorem~\ref{Thm:n_Z is invariant}.
For ease of notation, we write \(\cL_1 = \cL^c\) and
\(\cL_2 = \beta^*(\cL^c) \) in \( \cX^c(S_{2m})\) and similarly
\(\cL_1' = \cL^c\) and \(\cL_2' = (\beta')^*(\cL^c)\) in
\(\cX^c(S_{2m+2})\).

\begin{proposition}
  \label{Prop:MoveII}
  For each $c \in (-2, 2)$, the map
  \(\psi\maps \cX^c(F_{2m}) \to \cX^c(F_{2m+2})\) is a smooth
  embedding.  The image $X$ of \, $\cX^c(S_{2m})$ under $\psi$ is
  contained in \(\cX^c(S_{2m+2})\) and satisfies the following
  properties:
  \begin{enumerate}[label={(1)}]
  \item \(\cL_\ell' \cap X = \psi(\cL_\ell)\) for $\ell=1,2$.
  \item  \(\cL_1' \cap \cL_2' \subset X\).
  \item The normal bundle $\nu$ of \(X\) in $\cX^c(S_{2m + 2})$
    contains oriented 2-dimensional subbundles \(U_1,U_2\) such that
    \(\nu = -U_1 \oplus U_2\) and
    \(U_\ell|_{\psi(\cL_\ell)} = \nu_{\cL_\ell'/\psi(\cL_\ell)}\) as oriented
    bundles for $\ell=1,2$.
  \end{enumerate}
\end{proposition}

\begin{proof}

Let \(\Psi\maps \cRbar_{2m} \to \cRbar_{2m+2}\) be given by
\[
  \Psi(t, v_1, \ldots, v_{2m}) = (t, v_1,\ldots, v_{2m}, v_{2m}, -v_{2m}).
\]
Now \(\Psi\) is equivariant with respect to the actions of \(\R_{>0}\)
and the \(G_t\), so it descends to a smooth embedding
\(\cX^c(F_{2m}) \to \cX^c(F_{2m+2})\). It is easy to see that this
map agrees with the map \(\psi\) constructed at the beginning of the
subsection.  From Section~\ref{subsec:define res}, recall that
$\cRbar^c(S_{2m}) \subset \cRbar_{2m}$ is defined using
$F_{2m} \maps \cRbar_{2m} \to \SL{2}{\C}$ and similarly
$\cRbar^c(S_{2m+2}) \subset \cRbar_{2m+2}$ via the analogous
$F_{2m + 2}$. We next show $Y = \Psi\big(\cRbar^c(S_{2m})\big)$ is
contained in $\cRbar^c(S_{2m+2})$, and moreover
$Y = \Psi(\cRbar_{2m}) \cap \cRbar^c(S_{2m + 2})$.  Since
$F_{2m + 2} \circ \Psi = F_{2m}$, for $v \in \cRbar_{2m}$, if
\(F_{2m}(v) = I\) then \(F_{2m+2}\big(\Psi(v)\big) = I\) as well.  By
Lemma~\ref{lem:intersections}, this shows $\Psi$ takes
$\cRbar^c(S_{2m}) \setminus \cRbar_{2m}^{t \neq 0}$ into
$\cRbar^c(S_{2m+2}) \setminus \cRbar_{2m+2}^{t \neq 0}$.  The locus
where $t = 0$ is then covered by continuity of $\Psi$ and
Lemma~\ref{lem:M epsilon}, so $Y \subset \cRbar^c(S_{2m+2})$ and hence
$X \subset \cX^c(S_{2m + 2})$.  The stronger claim that
$Y = \Psi(\cRbar_{2m}) \cap \cRbar^c(S_{2m + 2})$ follows from
$F_{2m + 2} \circ \Psi = F_{2m}$ by examining the local functions
$f_\epsilon$ in Section~\ref{subsec:define res}.

Now we verify each of the three properties.

\noindent \textbf{Property 1:} Define \(\cN_\ell = \sigma^{-1}(\cL_\ell)\)
and \(\cN'_\ell= \sigma^{-1}(\cL'_\ell)\), and then let
\(i_1',i_2'\maps \cRbar_{m+1} \to \cRbar_{2m+2}\) be the inclusions
whose images are \(\cN_1',\cN_2'\). If \(u = (u_1,\ldots, u_{m})\) and
\(u' = (u_1,\ldots, u_{m+1})\), we have
\(i_1'(u') = (u_1,-u_1,\ldots, u_{m+1},-u_{m+1}).\) It is easy to see
that \(i_1'(u')\) is in \(Y\) if and only if \(u_{m+1} = -u_m\). It follows
that \(\cN_1' \cap Y = \Psi(\cN_1)\), and hence that
\(\cL'_1 \cap X = \psi(\cL_1)\).

Next, suppose that
\(i_2(u) = \beta^*(i_1(u)) = (v_1, \ldots, v_{2m})\).  Then
\begin{align*}
  i_2'(u') & = (\beta')^*(u_1,-u_1,\ldots, -u_m, u_{m+1},-u_{m+1}) \\
           & = \sigma_{2m}^*(v_1,\ldots, v_{2m},u_{m+1},-u_{m+1}) \\
           & = \big(v_1,\ldots, v_{2m-1}, u_{m+1},
             \rho(-u_{m+1}) \cdot v_{2m}, -u_{m+1}\big).
\end{align*}
Hence \( i_2'(u') \in Y\) if and only if
\(\rho(-u_{m+1}) \cdot v_{2m} = u_{m+1}\). This occurs if and only if
\(v_{2m} = u_{m+1}\), so
\(\cN_2' \cap Y= \setdefm{\big}{(v_1,\ldots,v_{2m},v_{2m},-v_{2m})}{\mbox{$v =
    i_2(u)$ with $u \in \cRbar_{2m}$}} = \Psi(\cN_2)\).
It follows that \(\cL_2' \cap X = \psi(\cL_2). \)

\noindent
\textbf{Property 2:} We use the descriptions of \(\cN_1'\) and
\(\cN_2'\) from Property 1. If
\[
  (v_1,\ldots, v_{2m-1},u_{m+1},\rho(-u_{m+1}) \cdot v_{2m}, -u_{m+1})
  \in \cN'_1,
\]
then \(\rho(-u_{m+1}) \cdot v_{2m} = u_{m+1}\). This can only happen
if \(v_{2m}=u_{m+1}\). Hence \(\cN_1' \cap \cN_2' \subset Y\), which
implies \(\cL_1' \cap \cL_2' \subset X\).

\noindent
\textbf{Property 3:} Consider the projection
\(p\maps \cRbar_1 \to \R\) that sends $(t, v)$ to $t$, and let
\(\Vtil = \ker dp\). Our standard orientation on \(\Qbar_t\) makes
\(\Vtil\) into an oriented 2-dimensional vector bundle over
\(\cRbar_1\). If \(\iota\maps \cRbar_1 \to \cRbar_1\) is given by
\(\iota(t,v) = (t,-v)\), then \(d\iota\maps \Vtil \to \Vtil\) is an
orientation reversing isomorphism. Note that the actions
of \(G_t\) and \(\R_{>0}\) on \(\cRbar_1\) extend to actions on
\(\Vtil\).

Define $\Yhat = \Psi(\cRbar_{2m})$ and let
\(\rho_i\maps \cRbar_{2m+2} \to \cRbar_1\) be given by
\(\rho_i(v) = (t(v),v_i)\).  Thus \(\Yhat\) is the locus where
\[
  \rho_{2m}(v) = \rho_{2m+1}(v) = \iota\big(\rho_{2m+2}(v)\big).
\]
Note that each \(\rho_i\) is equivariant with respect to the actions of the
\(G_t\) and \(\R_{>0}\).

Now let \(V = \rho_{2m}^*(\Vtil)|_{\Yhat}\), and consider the bundle maps
\[
  \alphatil\maps T\cRbar_{2m+2}|_{\Yhat} \to \Vtil \oplus \Vtil
  \mtext{and}
  \alpha\maps T\cRbar_{2m+2}|_{\Yhat} \to V \oplus V
\]
given by
\(\alphatil = \big(d\rho_{2m+1}-d\rho_{2m}, d(\iota \rho_{2m+2})
- d\rho_{2m}\big)\) and \(\alpha = \rho_{2m}^* \alphatil\).
Since the \(\rho_i\) are equivariant with respect to the actions of the
\(G_t\) and \(\R_{>0}\), so are \(\alphatil\) and \(\alpha\).

Let \(\nuhat\) be the normal bundle of \(\Yhat\) in \(\cRbar_{2m+2}\).
You can check that \(\alpha\) is surjective and
\(\ker \alpha = T\Yhat\), so ignoring orientations,
\(\nuhat \cong V \oplus V\). We now show this isomorphism is
orientation reversing. To see, this note that if
\(v' = \Psi(v) \in \Yhat\), then
\[
  T_{v'}\cRbar_{2m+2} =
  T_{v'}\Yhat \oplus T_{v_{2m}}\Qbar_{t} \oplus T_{-v_{2m}}\Qbar_{t}
\]
as an oriented vector space, where $\Yhat$ is oriented as the image of
$\cRbar_{2m}$. Since the dimension of \( T\Qbar_t \oplus T\Qbar_t \) is
even,
\(\nuhat|_{v'} \cong T_{v_{2m}}\Qbar_{t}\oplus T_{-v_{2m}}\Qbar_{t}\)
as oriented vector spaces. The map
\[
  \alpha\maps T_{v_{2m}}\Qbar_{t}\oplus T_{-v_{2m}}\Qbar_{t} \to V_{v'} \oplus V_{v'}
\]
is given by \((\bw_1, \bw_2) \mapsto \big(\bw_1,d\iota(\bw_2)\big)\), which is
orientation reversing. Finally, \(\alpha\) is equivariant with respect
to the actions of the \(G_t\) and \(\R_{>0}\), so our identification
\(\nuhat \cong -V\oplus V\) respects the actions of these groups.

We define \(V_1 \subset V\oplus V\) to be the image of the morphism
\(V \to V \oplus V\) given by \(\bv \mapsto (\bv,\bv)\), and let
\(V_2 = V \oplus 0 \subset V \oplus V\). Then \(V_1\) and \(V_2\) are
equivariant \2-dimensional subbundles of \(V \oplus V\)
with the orientation induced by the orientation on \(V\), and
\(V_1 \oplus V_2 \cong V \oplus V\) as oriented bundles.

We claim that
\(\nu_{\cN_\ell'/\cN_\ell} = V_\ell|_{\cN_\ell} \subset
\nuhat|_{\cN_\ell}\) as oriented equivariant vector bundles, where we
have identified $\cN_\ell$ with $\Psi(\cN_\ell)$ to reduce clutter.
To see this except for the orientation, first note that on \(\cN_1'\),
\(\rho_{2m+1} = \iota \rho_{2m+2}\), so \(\alpha(T\cN_1') = V_1\).
Similarly, \(\iota \rho_{2m+2} = \rho_{2m}\) on \(\cN_2'\), so
\(\alpha(T\cN_2') = V_2\) as unoriented vector spaces. For the
orientations, observe that for $u \in \cRbar_{m}$ and $u' \in
\cRbar_{m+1}$ as in the proof of Property 1, we have:
\[
  T_{u'}\cRbar_{m+1} \cong T_u\cRbar_m \oplus T_{u_{m+1}}\Qbar_t
\] as
oriented vector spaces, so for $v = i_\ell(u)$ and $v' = i'_\ell(u')$,
we have
\[
  T_{v'} \cN_\ell' = T_v \cN_\ell \oplus di'_\ell(T_{u_{m+1}}\Qbar_t).
\]
A computation shows that \(\alpha\big(di'_1(\bw)\big) = (\bw,\bw)\)
for $\bw \in T_{u_{m+1}}\Qbar_t$, which is an orientation preserving
map to \((V_1)_{v'}\). Similarly, you can check that
\(\alpha\big(di'_2(\bw)\big) = (-R(\bw), 0)\) where $R$ is rotation
clockwise by $2 \alpha$,
i.e.~$R =
\mysmallmatrix{\hphantom{-}\cos{2\alpha}}{\sin{2\alpha}}{-\sin{2\alpha}}{\cos{2\alpha}}$,
which is again orientation preserving.

Having studied \(\Yhat \subset \cRbar_{2m+2}\), we now turn our
attention to \(Y \subset \cRbar^c(S_{2m+2})\). The orientation on
\(\cRbar^c(S_{n})\) was determined using the map
\(\Ftil_n \maps \cRbar_{n} \to \Util\) from \eqref{Eq:Orientation 2},
which is a submersion for \(t\neq 0 \). Observe that
\(\Ftil_{2m+2} \circ \Psi = \Ftil_{2m}\), so we are in the situation
of Corollary~\ref{Cor:Normal bundles}, with $M' = \cRbar_{2m + 2}$,
$M = \Yhat$, $X' = \Util$, $X = \Itil$, $f = \Ftil_{2m+2}$,
$N' = \cRbar^c(S_{2m+2})$, and $N = Y$. Define
\(V_\ell' = V_\ell|_{Y}\). Applying Corollary~\ref{Cor:Normal bundles}
and noting that the dimension of $V_1'\oplus V_2'$ is even, we see
that
\[
  \nu_{\cRbar^c(S_{2m+2})/Y}= \nuhat|_Y = - V_1'\oplus V_2'
\]as oriented
bundles. Since \(\cN_\ell \subset Y\), we have
\(\nu_{\cN_\ell'/\cN_\ell} = V_\ell|_{\cN_\ell} = V_\ell'|_{\cN_\ell}\).

Finally, let \(\nu\) be the normal bundle of \(X\) in
\(\cX^c(S_{2m+2})\). The actions of \(\R_{>0}\) and the \(G_t\) on
\(Y\) and \(\cRbar^c(S_{2m+2})\) are free and extend to actions on
\( \nu_{\cRbar^c(S_{2m+2})/Y}\) which preserve the subbundles \(V_1'\)
and \(V_2'\). Passing to the quotient and using
Corollary~\ref{Cor:Normal Quotient}, we obtain subbundles
\(U_1,U_2 \subset \nu \) with \(U_1 \oplus U_2 = - \nu\) and
\(U_\ell|_{\cL_\ell} = \nu_{\cL_\ell'/\psi(\cL_\ell)}\). This
completes the proof of Property 3 and hence the proof of the
proposition.
\end{proof}

\begin{corollary}
  \label{Cor:MoveII}
  Theorem~\ref{Thm:n_Z is invariant} holds when \(\beta\) and
  \(\beta'\) are related by a move of type~II.
  \end{corollary}

\begin{proof}
Let \(\cL_\ell, \cL_\ell'\) be as above. By
Proposition~\ref{Prop:MoveII}, the hypotheses of
Lemma~\ref{Lem:Submanifold Intersection} are satisfied with the
opposite of the usual orientation on \(\cX^c(S_{2m+2})\). As
\(\dim U_1 =2\) is even, we get
\[
  \pairlocal{\cL_1}{\cL_2}{Z}
  = \pairlocal[\big]{\psi(\cL_1)}{\psi(\cL_2)}{\psi(Z)}
  = - \pairlocal[\big]{\cL_1'}{\cL_2'}{\psi(Z)}.
\]
Since \(|\beta'| = |\beta|+1\), it follows that
\[
  n_{Z,\beta} = (-1)^{|\beta|} \pairlocal{\cL_1}{\cL_2}{Z} =
  (-1)^{|\beta'|} \pairlocal[\big]{\cL_1'}{\cL_2'}{\psi(Z)} = n_{\psi(Z),\beta'}
\]
proving the corollary.
\end{proof}

We end this section with the proofs of its two main technical
theorems.

\begin{proof}[Proof of Theorem~\ref{Thm:n_Z is invariant}]
If the conclusion of Theorem~\ref{Thm:n_Z is invariant} holds for
the pairs \((\beta,\beta')\) and \((\beta', \beta'')\), it then
holds for \((\beta, \beta'')\) as well. Thus the general case of
Theorem~\ref{Thm:n_Z is invariant} follows from
Theorem~\ref{Thm:PlatMoves}, Corollary~\ref{Cor:MoveI}, and
Corollary~\ref{Cor:MoveII}.
\end{proof}

\begin{proof}[Proof of Theorem~\ref{Thm:n_Z at ends}]
For concreteness, we do the case $c = 2$.  Since
$\cX^2(\Pi_\beta) = X^{2, \irr}_\SLR(\Pi_\beta)$ and
$\cX^2(\beta) \subset \cX^2_\bal(S_{2m}, S) = X^{2, \irr}_\SLR(S_{2m},
S)$, the needed homeomorphism $\psi$ and commutative square exists just
from the maps on $\SLR$ character varieties.  So it remains to show
$n_{Z, \beta} = n_{\phi(Z), \beta'}$ when $K$ is real representation
small.  Rather than tackle this directly, it is more expedient to
prove it ``by continuity'' from Theorem~\ref{Thm:n_Z is invariant}
and Corollary~\ref{Cor:component sum}.

\begin{claim}
  For any knot $K$ there is an $a < 2$ such that $X^\cirr_\SU(K)$ is
  empty for all $c \in [a, 2]$.
\end{claim}
To see this, first note that $X_\SU(K)$ is compact as it is the image
of a closed subset of $(\SU)^{2m}$, which implies that the closure $C$
of $X^\irr_\SU(K)$ in $X_\SU(K)$ is compact.  Let $\chi$ be the
character in $C$ where $\tr$ has its maximum, say $a' = \tr(\chi)$. If
$\chi$ is reducible, then by Lemma~\ref{lemma: deform red} we have
$a \in D_K$ and so $a' < 2$ as $2 \notin D_K$.  If instead $\chi$ is
irreducible then $a' < 2$ by Lemma~\ref{lem: empty char vars}.  So any
$a$ with $a' < a < 2$ satisfies the claim.

Let $a$ satisfy the claim and further assume $a > \max D_K $.  Then
$\cX^{[a, 2]}(\Pi_\beta) = X^{[a, 2],\irr}_\SLR(\Pi_\beta)$ and hence
$\cX^{[a, 2]}(\beta)$ is contained in
$\cX_{-}^{[a, 2]}(S_{2m}, S) = X^{[a, 2], \irr}_\SLR(S_{2m}, S)$.
Hence $\cL^{[a, 2]} \cap \beta^* \cL^{[a, 2]}$ is a real semialgebraic
set with finitely many connected components, which are also its path
components.  Since we are working with character varieties only, we
also immediately have a ``parameterized'' version $\psi_c$ of $\psi$ in
Theorem~\ref{Thm:n_Z is invariant} which is still a homeomorphism.
Next, increase $a$ so that the components of the closed set
$\cL^{[a, 2]} \cap \beta^* \cL^{[a, 2]}$ correspond bijectively with
$\cL^2 \cap \beta^* \cL^2$.  For a connected component $Z$ of
$\cX^2(\beta)$, let $W$ be the component of
$\cL^{[a, 2]} \cap \beta^* \cL^{[a, 2]}$ that contains it. By
Corollary~\ref{Cor:component sum}, we have $n_{Z, \beta} =
\pair{\cL^c, \beta^* \cL^c}\big|_{W}$ for all $c \in [a, 2]$ and similarly
for $n_{\psi(Z), \beta'}$.  Applying Theorem~\ref{Thm:n_Z is
  invariant} for any $c$ in $[a, 2)$ proves the theorem.
\end{proof}


\section{Parabolics and the extended Lin invariant}
\label{sec: ext Lin}

In this section, we define the extended Lin invariant \(\eh(K)\), which
counts the heights of parabolic representations
\(\pi_1(S^3-K) \to \tSLR\). We begin by discussing parabolic
representations, and then review some facts about the \(\SLR\)
character varieties of knot complements and their boundaries. Finally,
we discuss the translation extension locus introduced in
\cite{CullerDunfield2018} and use it to define \(\eh(K)\).

\subsection{Parabolic representations}

Let  \(\rho \maps \pi_1(T^2) \to \SLR\) be a representation,  where
$T^2$ is the 2-torus.
Since
\(\pi_1(T^2) \cong \Z^2 \cong \pair{ \mu, \lambda} \) is abelian,
\(\rho(\mu)\) and \(\rho(\lambda)\) must belong to the same maximal
abelian subgroup of \(\SLR\). Up to conjugation by an element of
\(\SLR\), there are three such:
\begin{align*}
  T_e & = \setdef{\mysmallmatrix{\hphantom{-}\cos \theta}{\sin \theta }{-\sin \theta }{\cos \theta }}%
                 {\theta \in [0,2\pi)} & \text{(elliptic)}  \\
  T_p & = \setdef{\mysmallmatrix{\pm 1}{t}{0}{\pm 1}}{t \in \R}
                  & \text{(parabolic)} \\
  T_h & = \setdef{\mysmallmatrix{a}{0}{0}{1/{a}}}{a \in \R^\times}
                  & \text{(hyperbolic)}
\end{align*}
We call the representation \emph{elliptic}, \emph{parabolic}, or
\emph{hyperbolic} accordingly. Note that the center \(\{ \pm I\}\) of
\(\SLR\) is contained in all three subgroups; a representation whose
image in contained in \(\{ \pm I\}\) is called \emph{central}.

Now suppose \(K \subset S^3\) is a knot with
\(\EK = S^3 \setminus \nu(K)\) its exterior. We say that
\(\rho \maps \pi_1(\EK) \to \SLR\) is elliptic, parabolic, hyperbolic,
or central according to its restriction to \(\pi_1(\partial \EK)\).
We will always use generators $\mu, \lambda$ for
\(\pi_1(\partial \EK)\) where $\mu$ is a meridian for $K$ and $\lambda$ a
homological longitude.  We begin with a simple observation:

\begin{lemma} 
  \label{Lem:parabolics}
  If a nontrivial representation \(\rho \maps \pi_1(\EK) \to \SLR\)
  has \(\tr_ \mu \rho = \pm 2\), then \(\rho\) is parabolic.
\end{lemma}

\begin{proof}
Suppose \(\rho \maps \pi_1(\EK) \to \SLR\) with
\(\tr_ \mu \rho = \pm 2\). Since \(\rho(\mu)\) and \(\rho(\lambda)\)
belong to the same maximal abelian subgroup of \(\SLR\), either
\(\rho(\mu)\) is central or the representation $\rho$ is parabolic. If
\(\rho(\mu)\) is central and \(\rho(\lambda)\) is not, then \(\rho\)
descends to a nontrivial representation
\(\rho \maps \pi_1(\EK(\mu)) \to \PSLR\), a contradiction as
\(\EK(\mu) = S^3\).
\end{proof}

\begin{corollary}
\label{Cor:parabolics}
 If \(K\) is real representation small and \(h(K) \neq 0 \), then
  \(\pi_1(\EK)\) has an irreducible parabolic representation into
  $\SLR$.  Moreover, $2 h(K)$ is a signed count of conjugacy classes of
  such representations, where the signs and multiplicities come from
  the Casson-Lin picture.
\end{corollary}

\begin{proof}
By Theorem~\ref{thm: hSL+hSU}, we have $2 \not\in D_K$ and
$h^2_\SLR(K) = h(K)$.  So $h^2_\SLR(K)$ is nonzero, and in particular
\(X^{2,\irr}_\SLR(K)\) is nonempty by Theorem~\ref{Prop:Lcap beta L}.
By Lemma~\ref{Lem:parabolics}, each element of \(X^{2,\irr}_\SLR(K)\)
comes from an irreducible \emph{parabolic} representation, completing
the proof of the first claim.  As $K$ is real representation small,
the sets \(X^{2,\irr}_\SLR(K)\) and \(X^{-2,\irr}_\SLR(K)\) are finite,
so the second claim is immediate from the definitions and the fact
that $h^{-2}_\SLR(K)$ is also equal to $h(K)$.
\end{proof}

\subsection{The \(\SLR\) character variety of \(T^2\)}
\label{sec: SLR and T2}

We are principally interested in \(X_\SLR(K) = X_\SLR(\EK, \{\mu\})\),
where \(K\) is a knot in \(S^3\), and its image under the map
\(i^*: X_\SLR(K) \to X_\SLR(\partial \EK)\) coming from the inclusion
$i \maps \EK \to \partial \EK$.  To this end, we first describe
\(X_\SLR(T^2)\).

We define the \emph{elliptic locus} \(X^\elll_\SLR(T^2)\) to be the subset of
\(X_\SLR(T^2)\) in the image of the elliptic representations, and
similarly for the parabolic and hyperbolic loci.

\begin{figure}
  \begin{center}
    \begin{tikzoverlay}[scale=0.7]{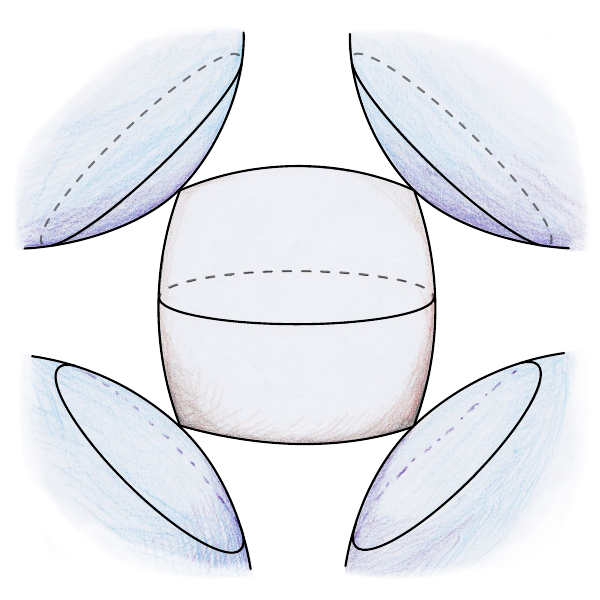}
      \begin{scope}[font=\footnotesize]
        \node[] at (49.8,37.1) {elliptic};
        \node[rotate=43] at (16.7,81.4) {hyperbolic};
        \begin{scope}[radius=0.4]
          \draw[fill=black] (30.08, 68.3) circle;
          \draw[fill=black] (68.95, 68.3) circle;
          \draw[fill=black] (68.95, 29.03) circle;
          \draw[fill=black] (30.08, 29.05) circle;
        \end{scope}
      \end{scope}
    \end{tikzoverlay}
  \end{center}

  \caption{The \(\SLR\) character variety $X_\SLR(T^2)$ of the
    \2-torus. The central pillowcase is the elliptic locus
    \(X^{\elll}_\SLR (T^2)\), with the rest the hyperbolic locus
    \(X^{\hypp}_\SLR (T^2)\); they meet at the four corners of the
    pillowcase, which is the parabolic locus.}
   \label{Fig:CV of T^2} 
\end{figure}
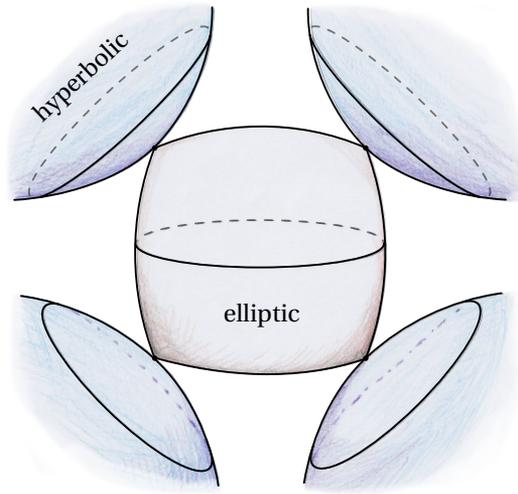

\begin{lemma}
  \label{lem: ell pillow}
  The locus \(X^{\elll}_\SLR (T^2)\) is the pillowcase orbifold
  \(T^2 /\{x\sim x^{-1}\}\), which is the flat orbifold with
  underlying space \(S^2\) and four orbifold points of order \(2\). It
  is shown in Figure~\ref{Fig:CV of T^2}.
\end{lemma}

\begin{proof}
By conjugation, every point in \(X^\elll_\SLR(T^2)\) is the image of a
representation \(\rho \maps \pi_1(T^2) \to T_e\). Such representations
are parametrized by \(T_e \times T_e \cong T^2\). Now if
\(A,A' \in T_e\), one has \(\tr A = \tr A'\) if and only if
\(A' = A^{\pm 1}\). It follows that
\(\rho, \rho' \maps \pi_1(T^2) \to T_e\) have the same character if
and only if \(\rho' = \rho^{\pm1}\), where
\(\rho^{-1}(x) = \rho(x)^{-1}\), proving the lemma.
\end{proof}

We will use the following explicit system of coordinates on
\(X^\elll_\SLR(T^2)\). Let
\(\mu^*,\lambda^* \in \Hom\big(\pi_1(T^2),\R\big) = H^1(T^2;\R)\) be
the basis algebraically dual to \(\{\mu, \lambda\}\), and let
\(\mubar^*, \lambdabarstar\) be their images in
\(\Hom\big(\pi_1(T^2),\R/2\Z\big)\). Identifying \(T_e\) with
\(\R/2\Z\) via $\theta \mapsto \theta/\pi$, we get coordinates on
\(\Hom\big(\pi_1(T^2),T_e\big)\), and hence on the pillowcase. We
usually normalize so that the $(\mubar^*, \lambdabarstar)$-coordinates
on the pillowcase are in $[0, 1] \times [-1, 1]$.  (The reason we
identify $T_e$ with $\R/2\Z$ rather than the seemingly more natural
$\R/\Z$ is to match the conventions of \cite{CullerDunfield2018},
where the focus is on $\PSL{2}{\R}$ rather than $\SL{2}{\R}$; note
also that the matrix
$\mysmallmatrix{\hphantom{-}\cos \theta}{\sin \theta }{-\sin \theta
}{\cos \theta }$ gives an anticlockwise rotation about the point $i$
in the upper-halfspace model of $\H^2$ through angle $2 \theta$.)

A similar argument shows that the \emph{hyperbolic locus}
\(X^\hypp_\SLR(T^2)\) consists of four distinct copies of
\(\R^2/ \{x \sim -x\}\), see also \cite{Gao2019}; we view the latter
space as \( \R^2\) with a single orbifold point of order \(2\). The
elliptic and hyperbolic loci intersect at the parabolic locus, which
consists of the four corners of the pillowcase, as illustrated in
Figure~\ref{Fig:CV of T^2}.

For the hyperbolic and elliptic loci of $X_\SLR(T^2)$, the preimage of
a character consists of either one or two conjugacy classes of
representations into \(\SLR\).  In contrast, the preimage of each
parabolic character contains a whole circle's worth of different
conjugacy classes, in addition to the conjugacy class of a central
representation.

\subsection{The character variety of a knot complement}
\label{subsec: X(K)}

We return to the setting of the exterior $\EK$ of \(K \subset S^3\),
where \(\mu\) and \(\ \lambda\) are its meridian and longitude.  Now
\(X_\SLR(K)\) decomposes as \(X_\SLR^\red(K) \cup X_\SLR^\irr(K)\),
and we discuss the two parts separately.

As in the case of $\SLC$ detailed in Section~\ref{subsec:reducibles},
the characters in \(X_\SLR^\red(K)\) are precisely those coming from
\(\rho\maps \pi_1(\EK) \to H\), where \(H\) is a maximal abelian subgroup
of \(\SLR\). Thus to understand \(X_\SLR^\red(K)\), it suffices to
study abelian representations. Any such representation factors through
\(H_1(\EK) \cong \Z\), so \(X_\SLR^\red(K) \cong X_\SLR^\red(U)\),
where \(U\) is the unknot.  Arguing as in Lemma~\ref{lem: ell
  pillow} shows that the set of elliptic reducible representations
$X^{\elll, \red}_\SLR(K)$ is homeomorphic to
\(S^1/\{x \sim x^{-1}\} \cong [0,1]\).  Also, for the inclusion
$i \maps \partial \EK \to \EK$, we have
\[
  i^*\big(X^{\elll, \red}_\SLR(K)\big) = \setdef{ \chi \in
    X_\SLR^\elll(\partial \EK)}{\lambdabarstar(\chi)=0}.
\]

Next, we consider the irreducible part of the character variety. For
the unknot $U$, we have $\pi_1(M_U) \cong \Z$, so every representation
is abelian and hence \(X^\irr_\SLR(U) = \emptyset\). For a more
interesting example, we consider the trefoil.

\begin{figure}
  \begin{center}
    \begin{tikzoverlay}[scale=0.7]{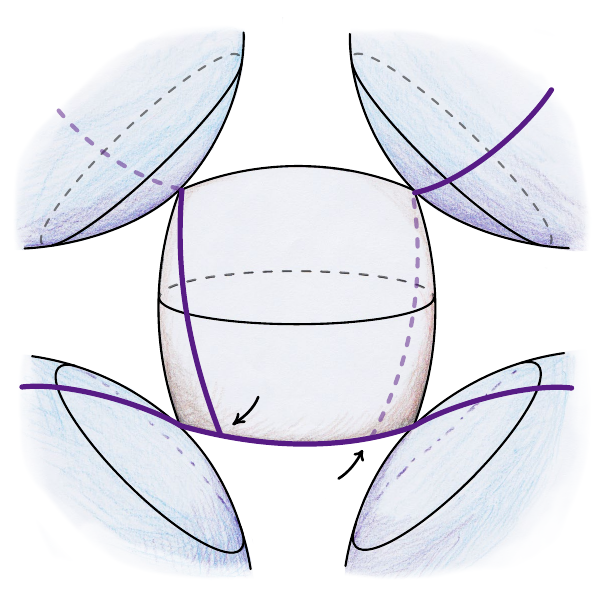}
      \begin{scope}[font=\footnotesize]
        \node[above=-1.5] at (43.6,33.9) {$\alpha = \frac{\pi}{6}$};
        \node[below left=-2 and -1] at (56.4,20.5) {$\alpha = \frac{5 \pi}{6}$};
        \node[right] at (30.4,59.5) {$i^*\big(\XSLR(M_K)\big)$};
      \end{scope}
    \end{tikzoverlay}
    \hspace{0.5cm}
    \begin{tikzpicture}[
  scale=2.8,
  nmdstd,
  line width=1.2pt,
  every node/.style={color=black},
  every circle/.style={radius=0.03}]
  
  \begin{scope}[color=axesgray]
    \draw (0, -1) rectangle (1, 1);
    \draw (0, -1) -- (1, -1)
    node[pos=0, below] {$0$}
    node[pos=0.5, below=0.2] {$\mubar^*$}
    node[pos=1, below] {$1$};
    
    \draw (0, -1) -- (0, 1);
    \node[left] at (0, 1) {$-1$};
    \node[left] at (0, 0) {$0$};
    \node[left=0.2] at (0, 0.3) {$\lambdabarstar$};
    \node[left] at (0, -1) {$1$};     

    \draw[->] (0, 0) -- (0, 0.5);
    \draw[->] (1, 0) -- (1, 0.5);

    \draw[->>] (0, 0) -- (0, -0.55);
    \draw[->>] (1, 0) -- (1, -0.55);

    \draw[->>>] (0, 1) -- (0.65, 1);
    \draw[->>>] (0, -1) -- (0.65, -1);

  \end{scope}

  \begin{scope}[line width=2pt]
    \draw[color=locus] (0, 0) -- (1, 0);
    \coordinate (A) at (0, 1);
    \coordinate (B) at (0 + 0.1666667, 0);
    \draw[color=locus] (A) -- (B);
    \draw[fill=otherparabolic, line width=0.5pt] (A) circle;
    \draw[fill=simplealex, line width=0.5pt] (B) circle node[below=0.1] {$\frac{1}{6}$};
    
    \coordinate (A) at (1, -1);
    \coordinate (B) at (1 - 0.1666667, 0);
    \draw[color=locus] (A) -- (B);
    \draw[fill=otherparabolic, line width=0.5pt] (A) circle;
    \draw[fill=simplealex, line width=0.5pt] (B) circle
    node[above=0.1] {$\frac{5}{6}$};
  \end{scope}
\end{tikzpicture}
      
    \end{center}
    \caption{ The image $\XSLR(M_K)$ inside $\XSLR(\partial M_K)$ for
      the trefoil knot $K$.  At right is the full picture, where
      $i^*(\XSLR^\red)$ is the approximately horizontal curve.  At
      left is $i^*\big(\XSLR^\elll(\EK)\big)$ in our
      $(\mubar^*, \lambdabarstar)$-coordinates on the pillowcase.
    }
    \label{Fig:Character of Trefoil}
\end{figure}

\subsection{The positive trefoil}
\label{Ex:Trefoil}
  
Let \(K=T(2,3)\) be the positive (right-handed) trefoil knot, which is
the plat closure of $\sigma_1^3 \in B_2$. Then \(\EK\) is Seifert
fibered, and we have a presentation
\[
  \pi_1(\EK) = \spandef{ x, y, f}{[x,f]=[y,f]=1, x^2 = f = y^3 }
\]
where \(f\) is the class of the Seifert fiber. Suppose
\(\rho \maps \pi_1(\EK) \to \SLR\) is irreducible. Since \(f\) is
central and \(\rho\) is irreducible, \(\rho(f)\) must be central. If
\(\rho(f) = I\), then \(\rho(x)^2 = I\), which implies
\(\rho(x) = \pm I\). This implies that \(\rho\) is reducible,
contradicting our assumption. Hence \(\rho(f) = -I\).

Since \(\rho(x)^2 = \rho(y)^3 = -I\), both \(\rho(x)\) and \(\rho(y)\)
are elliptic. The set of elliptic \1-parameter subgroups can be
identified with \(\H^2\) by taking a subgroup to its unique common
fixed point. Under this identification, the action of \(\SLR\) on such
subgroups by conjugation becomes its action on $\H^2$ by
isometries. Hence, after conjugation, we can assume that $\rho(x)$ and
$\rho(y)$ fix respectively $i$ and $t i$ in $\H^2$ for some $t > 1$;
equivalently
\[
\rho(x) = \twobytwomatrix[r]{\cos \theta}{ - \sin \theta}{\sin \theta}{\cos \theta}
\mtext{and}
\rho(y) = \twobytwomatrix[r]{\cos \phi}{ -t \sin \phi}{t^{-1} \sin \phi}{\cos \phi} 
\]
for some \(\theta, \phi \in (-\pi, \pi]\) and \(t \in [1, \infty)\).
As \(\rho(x)^2 = \rho(y)^3 = -I\), we must have
\(\theta = \pm \frac{\pi}{2}\) and \(\phi = \pm \frac{\pi}{3}\).
Moreover, if \(\theta = \pm \frac{\pi}{2}\),
\(\phi = \pm \frac{\pi}{3}\), and \(t \in [1, \infty)\), the above
formula determines a representation
\(\rho_{\theta, \phi, t} \maps \pi_1(\EK) \to \SLR\). The
representations \(\rho_{\theta,\phi,t}\) and
\(\rho_{-\theta,-\phi,t}\) are not conjugate in \(\SLR\), but they are
conjugate in \(\SLC\), so they have the same character. Hence 
$X_\SLR^{\irr}(K)$ consists of two arcs parametrized by
\(t \in (1,\infty)\).  Each arc limits on the reducible locus as
\(t \to 1\).

To describe their image in $X^\elll_\SLR(\partial \EK)$ under \(i^*\),
note there are curves on \(\partial \EK\) that are Seifert fibers, and
that \([f] = 6 \mu + \lambda \). Since an irreducible \(\rho\)
satisfies \(\rho(f) = -I\), the image \(i^*\big(X^\irr_\SLR(K)\big)\)
lies on the image of the line given by the equation
\(6 \mu^* + \lambda^* = 1\).  We have \(\mu = xyf^{-1}\), so if
\(\rho_t = \rho_{\pi/2,\pi/3,t}\), we compute
\(\tr \rho_t(\mu) = \frac{\sqrt3}{2} (t+t^{-1})\). At \(t=1\), the
representation \(\rho_1\) is reducible and
\(\tr_\mu \rho_1 = \sqrt{3} = 2 \cos \frac{\pi}{6}\). The trace
increases monotonically as \(t\) increases, until at \(t = \sqrt{3}\),
we have \(\tr \rho_{\sqrt{3}}(\mu) = 2\) and \(\rho_t\) is
parabolic. For larger values of \(t\), the representation \(\rho_t\)
is hyperbolic, and, as $t \to \infty$, the character of $\rho_t$ limits
to an ideal point of $X_\SLR(K)$ corresponding to the vertical
essential annulus in $\EK$.  The arc given by
\(\rho_t' =\rho_{\pi/2,-\pi/3,t}\) is very similar, except now
\(\tr \rho_t'(\mu) = - \frac{\sqrt3}{2} (t+t^{-1})\), so at \(t=1\) we
have \(\tr_\mu \rho_1 = -\sqrt{3} = 2 \cos \frac{5 \pi}{6}\) and the
trace decreases monotonically with \(t\). The character variety
\(X_\SLR(K)\) is shown in Figure~\ref{Fig:Character of Trefoil}.

\subsection{General picture}
\label{sec: gen pic}
  
In general, \(X^{\irr}_\SLR(K)\) has expected dimension \(1\)
although the actual dimension may be larger. If \(K\) is small,
Lemma~\ref{lem:small to few chars} implies that \(X^{\irr}_\SLR(K)\)
has dimension \(1\), and that \(K\) is real representation small.
As in the case of the trefoil, \(X^{\elll, \irr} _\SLR(K)\) is typically
noncompact and limits to both parabolic and reducible characters.  If
\(\chi\) is a limit point of the latter type with
\(\tr_\mu \chi = c\), then $c \in D_K$ by Lemma~\ref{lemma: deform
  red}. For the trefoil \(K=T(2,3)\), the roots of \(\Delta_K\) are
\(e^{\pm i\pi/3}\), and so $D_K = \{\sqrt{3}\}$, corresponding to the
reducible representations \(\rho_1, \rho_1'\) identified above.

We can interpret \(h^c_\SLR\) in terms of \(X_\SLR(K)\). For
\(\alpha \in [0, \pi]\), let \(V_\alpha\) be the vertical line/circle
in the pillowcase where the \(\mubar^*\)-coordinate is $\alpha/\pi$.
For $c = 2 \cos \alpha$ not in $D_K$, the invariant \(h^c_\SLR(K)\) is
a count (with signs and multiplicities) of the intersection of
$V_\alpha$ with $i^*\big(X^\irr_\SLR(K)\big)$. When \(X_\SLR(K)\)
is transversely cut out, one can show (as Heusener \cite{Heusener2003}
did for \(\SU\)) that the multiplicities can be used to orient the
arcs of \(X^\irr_\SLR(K)\) and make this a precise relationship.

For example, if \(K\) is the positive trefoil,
Figure~\ref{Fig:Character of Trefoil} shows that
\(X^{c, \irr}_\SLR(K)\) is empty for
\(\alpha \in (\frac{\pi}{6}, \frac{5\pi}{6})\), so
\(h^c_\SLR(K) = 0 \) in this range. For
\(\alpha \in [0, \frac{\pi}{6}) \cup (\frac{5 \pi}{6}, \pi]\),
\(X^{c, \irr}_\SLR(K)\) contains a single point, which suggests (but
does not prove) that \(h^c_\SLR(K) = 1 \) for such \(\alpha\). We will
show that this is indeed the case in Example~\ref{Ex:Trefoil2}.

\subsection{Symmetry}
\label{subsec: X_K sym}

The character variety in Figure~\ref{Fig:Character of Trefoil} is
symmetric under $\pi$ rotation around the vertical line through the
middle of the pillowcase.  This symmetry is present for any $K$ as we
now explain, following \cite[Lemma~6.1]{CullerDunfield2018}, and
extends to the resolved setting as we show in the below
Lemma~\ref{lem:b symmetry}.

For a space $Y$, a representation \(\eta\) of $\pi_1(Y)$ to the
center $\{\pm I\}$ of $\SLC$ is determined by its character, and the
resulting \emph{central characters} are classified by \(H^1(Y;\Z/2)\)
by viewing $\Z/2$ as $\{\pm I\}$. The group of central characters acts
on \(R_\C(Y)\) by \((\eta \cdot \rho) (g) = \eta(g) \rho(g)\), and
this descends to an action of \(H^1(Y;\Z/2)\) on \(X_\C(Y)\) and
$X(Y)$.

When \(Y = \EK\), the group \(H^1(\EK;\Z/2)\) has a unique nonzero
element \(\muhat^*\).
There is thus an involution
\(b \maps X(K) \to X(K)\) given by
\(b(\chi) = \muhat^*\cdot \chi\).  If
\(i^*(\chi) = (m, \ell)\), then
\[
  i^*\big(b(\chi)\big) = i^*\big(\muhat^* \cdot \chi\big) = i^*(\muhat^*)
  \cdot i^*(\chi) = \big(m+1, \ell\big) =
  \big(1-m, -\ell\big).
\]
This symmetry is clearly visible in Figure~\ref{Fig:Character of
  Trefoil}.

If \(\chi \in X^c(K)\), then \( b(\chi) \in X^{-c}(K)\), and $b$
restricts to maps $X^\cirr_\SLR(K) \to X^{-c,\irr}_\SLR(K)$ and
$X^\cirr_\SU(K) \to X^{-c,\irr}_\SU(K)$.  The same $\muhat^*$ action also
gives $b \maps X^c_{U_0}(K) \to X^{-c}_{U_0}(K)$.  Putting this
together, we get a map \(\varphi \maps \cX^c(K) \to \cX^{-c}(K)\),
which we use to investigate the effect of this symmetry on the local
intersection numbers \(n_{Z}\). Specifically, we next show the
following, which is needed for Corollary~\ref{cor: h of small mod 2}:
\begin{lemma}
  \label{lem:b symmetry}
  For $c \in [-2, 2]$, the map
  \(\varphi \maps \cX^c(K) \to \cX^{-c}(K)\) is a diffeomorphism.  If
  \(Z\) is a component of \(\cX^c(K)\), then
  \(n_{\varphi(Z)} = n_{Z}\).
\end{lemma}
Note that if we use the local intersection multiplicities to orient
the arcs of \(X^\irr_\SLR(K)\), as suggested in Section~\ref{sec: gen
  pic}, then the action of \(b\) on \(X^\irr_\SLR(K)\) is orientation
reversing because \(b\) preserves the orientation of \(X_\SLR(T^2)\)
but reverses the orientations of the vertical circles \(V_\alpha\).

\begin{proof}[Proof of Lemma~\ref{lem:b symmetry}]
First, assume that $c \in (-2, 2)$.  Given a plat diagram of \(K\)
corresponding to a braid \(\beta \in B_{2m}\), let
\(j \maps S_{2m} \to \EK\) be the inclusion of the splitting
2-sphere. Then \(j^*(\muhat^*)([s_i]) = -I\) for \(i=1,\ldots,2
m\). Consider the map \(\Phi \maps \cRbar_{2m} \to \cRbar_{2m}\) given
by \(\Phi(v) = -v\). From Section~\ref{sec: pi_c}, we have:
\begin{equation}
  \label{Eq:A beta}
  \pi'_{-c} \big(\Phi(v)\big)(s_i) = A\big(\pi-\alpha,t(v),-v_i\big) = - A(\alpha,t(v),v_i)
\end{equation}
so \(\pi'_{-c}\big(\Phi(v)\big) = j^*(\muhat^*) \cdot \pi'_c(v)\).
This identity also shows that
\(F_{-c}\big(\Phi(v)\big) = (-1)^{2m}F_c(v) = F_c(v)\), so
\(\Phi \maps \cRbar^c(S_{2m}) \to \cRbar^{-c}(S_{2m})\). Finally, it
is easy to check that \(\Phi\) commutes with the actions of \(G_t\)
and \(\R_{>0}\), so it descends to a smooth map
\(\cX^c(S_{2m}) \to \cX^{-c}(S_{2m})\) which you can check is $\phi$.
Note $\phi$ is a diffeomorphism since interchanging the roles of $c$
and $-c$ builds $\varphi^{-1}$.
  
Let \(\cL_0 ^c= \cL^c\) and \(\cL_1^c = \beta^*(\cL^c)\). It is
immediate from the definition that \(\varphi(\cL_0^c) = \cL_0^{-c}\)
as sets. Referring to equation~\eqref{Eq:Braid Action} and using
equation~\eqref{Eq:A beta} again, we see that \(\varphi\) is
equivariant with respect to the actions of \(B_{2m}\) on
\(\cX^c(S_{2m})\) and \(\cX^{-c}(S_{2m})\); the key point is that
\(\psi_t(-A) = \psi_t(A)\). Hence \(\varphi(\cL_1^c) = \cL_1^c\) as
sets. 

The local multiplicity
\(n_Z = (-1)^{|\beta|}\langle \cL^c_0, \cL^c_1\rangle|_Z\), so to
prove the lemma we must show
\(\langle \cL^c_0, \cL^c_1\rangle|_Z = \langle \cL^{-c}_0,
\cL^{-c}_1\rangle|_{\varphi(Z)}\). To do this, we examine the effect
of \(\varphi\) on the orientations of \(\cX^c(S_{2m})\), \(\cL_0^c\),
and \(\cL_1^c\). The antipodal map on \(\Qbar_t\) is orientation
reversing, so \(\Phi \maps \cRbar_{2m} \to \cRbar_{2m}\) multiplies
the orientation by a factor of \((-1)^{2m}=1\). As we observed above,
\(F_{-c}\big(\Phi(v)\big) = F_c(v)\). Referring to equation
(\ref{Eq:Orientation 2}), we see that
\(\Phi \maps \cRbar^c (S_{2m}) \to \cRbar^{-c}(S_{2m})\) also has
degree \(1\). Finally, \(\Phi\) commutes with the actions of \(G_t\)
and \(\R_{>0}\), so
\(\varphi \maps \cX^c (S_{2m}) \to \cX^{-c}(S_{2m})\) has degree \(1\)
as well.
 
A similar calculation shows that
\(\varphi \maps \cL_i^c \to \cL_i^{-c}\) has degree \((-1)^m\). So
\[
  \langle \cL_0^c,\cL_1^c\rangle|_Z
  = (-1)^m(-1)^m\langle \cL_0^{-c},\cL_1^{-c}\rangle_{\varphi(Z)}
  = \langle \cL_0^{-c},\cL_1^{-c}\rangle_{\varphi(Z)}
\]
as desired.
The remaining case of $c = \pm 2$ now follows ``by continuity'' as in
the proof of Theorem~\ref{Thm:n_Z at ends}.
\end{proof}

\subsection{The translation extension locus}
\label{sec: trans ex locus}

We next turn to studying representations to $\tSLR$ in the spirit of
\cite{CullerDunfield2018}.  Recall from our Section~\ref{sec: tSLR
  basics} the basic properties of $\tSLR$ and the space
$\XA_\tSLR(\Gamma)$ that encodes representations $\Gamma \to \tSLR$
modulo conjugation by $\tSLR$.  We also use from Section~\ref{sec:
  tSLR basics} the character variety
$\XA_\SLR(\Gamma) = R_\SLR(\Gamma) \sslash \SLR$, which is a branched
double cover of our usual
$X_\SLR(\Gamma) = R_\SLR(\Gamma)  \sslash \SLRpm$.  We start with some
simple examples.

\begin{example}
  \label{ex: torus chars}
  Suppose \(\Gamma = \pi_1(T^2) = \Z^2\).  Then
  \(\XA^{\elll}_\SLR(T^2)\) is a double cover of \(X^\elll_\SLR(T^2)\)
  branched over the set of central characters, which is the same as
  the set of parabolic points. In contrast, the map
  \(\XA^{\mathrm{hyp}}_\SLR(T^2) \to X^{\mathrm{hyp}}_\SLR(T^2)\) is
  one-to-one. We conclude that \(\XA_\SLR(T^2)\) consists of a torus
  (the branched double cover of the pillowcase) with four cones (the
  hyperbolic components) attached at the four parabolic points.
  Coordinates on \(\XA^{\elll}_\SLR(T^2)\) are given by pairs
  \(\mubar^*, \lambdabarstar\) with
  \(-1 \leq \mubar^*, \lambdabar ^* \leq 1\).  The generator
  \(a \maps \XA^{\elll}_\SLR(T^2) \to \XA^{\elll}_\SLR(T^2)\) of the
  branched covering group corresponds to flipping over the copy of
  $\H^2$ that $\SLR$ acts on, in particular reversing the rotation
  direction of each elliptic element; thus, in our coordinates we have
  \( a(\mubar^*, \lambdabarstar) = (-\mubar^*, -\lambdabarstar)\).
  There is also a map \(b \maps \XA_\SLR(T^2) \to \XA_\SLR(T^2)\)
  given by multiplication by the central character where
  \(\mu \mapsto -I\) and $\lambda \mapsto I$; in coordinates, we have
  \(b(\mubar^*, \lambdabarstar) = (\mubar^*+1,
  \lambdabarstar)\). (Caution: unlike $a$, the map $b$ is not a
  covering transformation for
  $\XA^{\elll}_\SLR(T^2) \to X^{\elll}_\SLR(T^2)$, but rather
  permutes the fibers of the map
  $\XA^{\elll}_\SLR(T^2) \to \XA^{\elll}_\PSLR(T^2)$.)

  Next, we consider \(\XA_\tSLR(T^2)\).  As \(T^2\) is a
  \(K(\Gamma,1)\), we have \(H^1(\Gamma) = H^1(T^2) = \Z^2\). The
  Euler class of the trivial representation is \(0\) and
  \(R_\SLR(\Gamma)\) is connected, so \(e(\rho) =0\) for all
  \(\rho \in R_\SLR(\Gamma)\).  It follows that \(\XA_\tSLR(T^2)\) is
  a \(\Z^2\) cover of \(\XA_\SLR(T^2)\), i.e.~its universal cover.  We
  will focus on the elliptic component of this space, which we view as
  \(\R^2\), with coordinates \(\mu^*, \lambda^*\).  Another viewpoint
  on these coordinates is described in \cite[\S
  3.5]{CullerDunfield2018}, namely one defines
  \[
    \trans \maps \XA_\tSLR(T^2) \to H^1(T^2; \R) \mtext{by} [\rho]
    \mapsto \trans \circ \rho.
  \]
  Concretely, in our $(\mu^*, \lambda^*)$-coordinates, this sends
  $[\rho]$ to
  $\big(\trans\big(\rho(\mu)\big),
  \trans\big(\rho(\lambda)\big)\big)$.  In particular, we have a
  homeomorphism
  $\XA^\elll_\tSLR(T^2) \xrightarrow{\trans} H^1(T^2; \R)$.
   
  The preimages of the parabolic points form a
  lattice \(\Z^2 \subset \R^2\), and the covering group of the
  previous paragraph corresponds to the subgroup
  \((2\Z)^2 \subset \Z^2\). There is a short exact sequence
  \[1 \to \Z \to Z(\tSLR) \to Z(\SLR) \to 1.\] The extension is
  nontrivial with \(Z(\tSLR) = \Z\). As before, the set of central
  characters \(H^1(\Z^2;Z(\tSLR)) \cong \Z^2\) acts on
  \(\XA_\tSLR(\Z^2)\); this is the action of the full lattice
  \(\Z^2 \subset \R^2\) by translation.
\end{example}

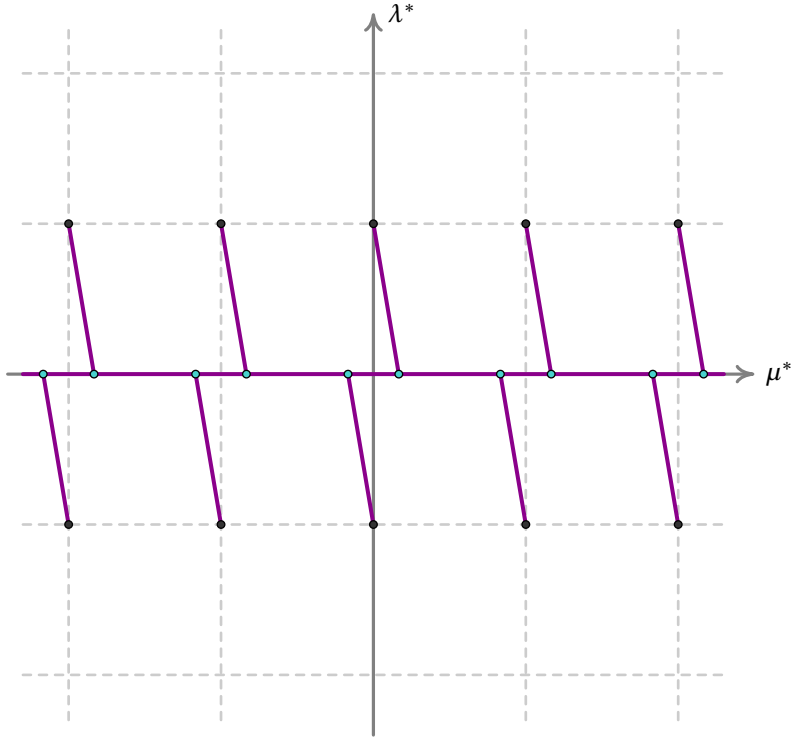
\begin{figure}
  \centering
  \begin{tikzpicture}[scale=2.0, nmdstd, line width=1.2pt,
  every circle/.style={radius=0.025}]
  \draw[color=black!20, dashed, line width=1 pt] (-2.3, -2.3) grid (2.3, 2.3);
  \begin{scope}[color=axesgray]
    \draw[->] (-2.4, 0) -- (2.5, 0) node[right, color=black] {$\mu^*$};
    \draw[->] (0, -2.4) -- (0, 2.4) node[right=0.04, color=black] {$\lambda^*$};
  \end{scope}
  \begin{scope}[line width=1.5pt]
    \draw[color=locus] (-2.3, 0) -- (2.3, 0);
    \foreach \x in {-2, -1, ..., 2}{
      \coordinate (A) at (\x, -1);
      \coordinate (B) at (\x - 0.1666667, 0);
      \draw[color=locus] (A) -- (B);
      \draw[fill=otherparabolic, line width=0.5pt] (A) circle;
      \draw[fill=simplealex, line width=0.5pt] (B) circle;
    }

    \foreach \x in {-2, -1, 0, ..., 2}{
      \coordinate (A) at (\x, 1);
      \coordinate (B) at (\x + 0.1666667, 0);
      \draw[color=locus] (A) -- (B);
      \draw[fill=otherparabolic, line width=0.5pt] (A) circle;
      \draw[fill=simplealex, line width=0.5pt] (B) circle;
    }
    
  \end{scope}
\end{tikzpicture}

  \vspace{0.5cm}
  
  \caption{ The translation extension locus
    \(i^*\big(\XA_\tSLR^\elll(K)\big)\) for the positive trefoil
    is shown above.  The related set
    \(i^*\big(\XA_\SLR^\elll(K)\big)\) in the torus
    $\XA^\elll_\SLR(\partial \EK)$ can be visualized by taking the
    four squares touching the origin and identifying the sides of the
    larger square they form by the action of $(2\Z)^2$.}
  \label{Fig:TrefoilCSLR}
\end{figure}

\begin{example}
  \label{ex: trans ext locus}

  Now suppose \(\Gamma = \pi_1(\EK)\).  Any nonparabolic point of
  \(X^{\elll}_\SLR(K)\) has two distinct preimages in \(\XA_\SLR(K)\),
  as does any irreducible parabolic by Lemma~\ref{lem: SLR branching}.
  However, the parabolic reducible characters have only a single
  preimage. We denote the covering action of \(\Z/2\) on
  \(\XA_\SLR(K)\) by \(a\) as well; it satisfies
  \(i^* \circ a = a \circ i^*\), so the image of
  \(i^*\big(\XA_\SLR(K)\big)\) in \(\XA_\SLR(\partial \EK)\) is
  determined by the image of \(i^*\big(X_\SLR(K)\big)\) in
  \(X_\SLR(\partial \EK)\).

  Now we consider lifts to \(\tSLR\). We have
  \(H^2(\Gamma) = H^2(K) = 0\), so every representation
  \(\Gamma \to \SLR\) lifts to \(\tSLR\), and \(\XA_\tSLR(\Gamma)\) is
  a covering space of \(\XA_\SLR(\Gamma)\) with deck group
  \(H^1(K) = \Z\) by Theorem~\ref{thm: tSLRchar}.  We are primarily
  interested in the elliptic part of \(\XA_\tSLR(K)\).  The closure of
  $i^*\big(\XA_\tSLR^\elll(K)\big)$ in
  $\XA_\tSLR^\elll(\partial \EK) \cong H^1(\partial \EK; \R)$ is the
  \emph{translation extension locus}
  $\ELG(\EK) \subset H^1(\partial \EK; \Z) $ studied in
  \cite{CullerDunfield2018}.  (Whenever $\EK$ is real representation
  small, one simply has
  $i^*\big(\XA_\tSLR^\elll(K)\big) = \ELG(\EK)$.)
  Figure~\ref{Fig:TrefoilCSLR} shows the translation extension locus
  of \(T(2,3)\); it is entirely determined by
  \(X_\SLR^\elll(\Gamma)\), since the latter space is connected.

  Lemma 6.1 of \cite{CullerDunfield2018} shows that the maps
  \(a\) and \(b\), which are automorphisms of \(\XA_\SLR(K)\), lift to
  automorphisms \(\atil, \btil\) of
  \(\XA_\tSLR(K)\). The maps \(\atil, \btil\) generate
  an action of the infinite dihedral group \(D_\infty\) on
  \(\XA_\tSLR(K)\). Their action on \(\XA_\tSLR(T^2)\) is given by
  \(\atil(\mu^*, \lambda^*) = (-\mu^*, -\lambda^*)\) and
  \(\btil(\mu^*, \lambda^*) = (\mu^*+1, \lambda^*)\).  This
  \(D_\infty\) action is evident in Figure~\ref{Fig:TrefoilCSLR}.
  For any knot, the translation extension locus is a finite union of
  analytic arcs and isolated points, whose quotient under $D_\infty$ is
  a finite graph \cite[Theorem 4.3]{CullerDunfield2018}.
\end{example}

\subsection{The extended Lin invariant}
\label{sec: extended Lin}

We now use the translation extension locus to define a refinement of
the invariant \(h(K)\) that also encodes the values of $\lambda^*$ for
lifts of $\chi \in X^2_\SLR(K)$ to $\tSLR$.  For this, we need:

\begin{lemma}
  \label{lem: para preimages}
  A nonparabolic \(\chi \in X^{\elll, \irr}_\SLR(K)\) has a unique
  preimage \(\chitil\in \XA_\tSLR(K)\) satisfying
  \(\mu^*(\chitil) = \mubar^*(\chi)\), where the latter is normalized
  to be in $[0, 1]$.  In contrast, if \(\chi\) is parabolic, there are
  two such preimages \(\chitil_\pm\), which satisfy
  \(\lambda^*(\chitil_-) = - \lambda^*(\chitil_+)\).
\end{lemma}

\begin{proof}
Throughout, you may find the illustration for the trefoil in
Figure~\ref{Fig:TrefoilCSLR} helpful.  By Lemma~\ref{lem: SLR
  branching}, there are two preimages $\chi_\pm$ of $\chi$ in
$\XA_\SLR(K)$.  Using $-1 \leq \mubar^*, \lambdabarstar \leq 1$ as our
fundamental domain for $\XA_\SLR(T^2)$, we have $\chi_- = a(\chi_+)$
and $\mubar^*(\chi_-) = -\mubar^*(\chi_+)$.

If $\chi$ is elliptic, take $\chi_+$ to be the one with $\mubar^*$ in
$(0, 1)$.  Then $\mubar^*(\chi_{-})$ is in $(-1, 0)$.  Now the lifts
of $\chitil_\pm$ to $\XA_\tSLR(K)$ have $\mu^*(\chitil_\pm)$ in
$\mubar^*(\chi_\pm) + 2 \Z$.  So there is a unique lift of $\chi$ with
\(\mu^* = \mubar^*(\chi)\), namely a lift of $\chi_+$,
proving the lemma in this case.

Suppose next that \(\chi\) is parabolic with $\chi(\mu) = 2$.  Then
there is one lift $\chitil_+$ of $\chi_+$ to $\XA_\tSLR(K)$ where
$\mu^*(\chitil_+) = 2k$ for each $k \in \Z$.  As the same is true for
$\chi_-$, we see there are exactly two lifts $\chitil_+$ and
$\chitil_-$ of $\chi$ to $\XA_\tSLR(K)$ with $\mu^* = 0$. We have
\(\atil \cdot \chitil_+ = \chitil_-\), so
\(\lambda^*(\chitil_-) = - \lambda^*(\chitil_+)\).  The final case of
a parabolic with $\chi(\mu) = -2$ is the same, except now we have
$\mu^*(\chitil_\pm) = 1$.
\end{proof}

Suppose \(X\) is a component of \(X^{2,\irr}_\SLR(K)\). By
Lemma~\ref{Lem:parabolics}, \(X\) is composed entirely of parabolic
representations. It follows that \(\mubar^*\) and \(\lambdabarstar\)
are integral (and hence constant) on \(X\); the same is true for
$\mu^*$ and $\lambda^*$ on each component of the preimage of $X$ in
\(\XA_\tSLR(K)\). Applying Lemma~\ref{lem: para preimages}, we see
that \(X\) has two preimages \(\Xtil_\pm\) in \(\XA_\tSLR(K)\) which
satisfy \(\mu^*|_{\Xtil_\pm} = \mubar^*|_X\), and that
\(\lambda^*|_{\Xtil_+} = - \lambda^*|_{\Xtil_-}\).

\begin{definition}
  \label{Def:extended invariant}
  If \(K\) is real representation small, the \emph{extended Lin
    invariant} of \(K\) is
  \[\eh(K) =
    \sum_{X} n_X \cdot \left(t^{{\lambda}^*(\Xtil_+)}+
      t^{{\lambda}^*(\Xtil_-)}\right) \in \Z[t, t^{-1}],
  \]
  where the sum runs over components of \(X^{2, \irr}_\SLR(K)\), and
  $n_X$ is well-defined by Theorem~\ref{Thm:n_Z at ends}.
\end{definition}

We will see in Example~\ref{Ex:Trefoil2} that
\(\eh(T(2,3)) = t+t^{-1}\), which is plausible from
Figure~\ref{Fig:TrefoilCSLR}.  Informally, the coefficient of \(t^i\)
in \(\eh(K)\) is a signed count of arcs exiting the translation
extension locus at the parabolic of height \(i\). The sign is
determined both by the orientation of the arc and which side of the
line \(\mu^*=0\) it exits from.


\begin{lemma}
\label{lem:eh properties}
For a real representation small knot $K$, the invariant
  \(\eh\) has the following properties.
  \begin{enumerate}

  \item 
    \label{item:htil to h}
    \(\eh(K)|_{t=1} = 2 h(K)\).
  \item
    \label{item: involution}
    \(\eh(K)\) is invariant under the involution of
    \(\Z[t^{\pm 1}]\) sending $t$ to $t^{-1}$.

  \item
    \label{item: Milnor-Wood}
    Define $\deg \eh(K)$ to be the maximum of $\abs{e}$ where $t^e$
    appears in $\eh(K)$.  Then \(\deg \eh(K) \leq 2g(K) - 1\) where
    \(g(K)\) is the Seifert genus of \(K\).
    
  \item
    \label{item: hyp sharp MW}

    If \(K\) is fibered and hyperbolic, the inequality in part
    (\ref{item: Milnor-Wood}) is strict.
  \end{enumerate}

\end{lemma}

\begin{proof}
Item (\ref{item:htil to h}) is immediate from the definition, and
(\ref{item: involution}) follows from the fact that
\(\lambda^*(\chitil_-) = - \lambda^*(\chitil_+)\) by Lemma~\ref{lem:
  para preimages}.  Next, (\ref{item: Milnor-Wood}) follows from the
Milnor-Wood inequality in the form of Proposition 6.5 from
\cite{CullerDunfield2018}.  Finally, for (\ref{item: hyp sharp MW}),
we use the argument of \cite[\S 3.5]{Calegari2006} as follows. If
$\rho \maps \pi_1\EK \to \PSLR$ corresponds to the putative
$t^{2g - 1}$ term of $\eh(K)$, then its restriction to the fiber $F$
must lie in the component of $X_{\PSLR}(F)$ corresponding to the
Teichm\"uller space of $F$.  Hence the bundle monodromy $\phi$ leaves
invariant that hyperbolic structure on $F$.  This forces $\phi$ to
have finite-order in the mapping class group of $F$, which contradicts
that $\EK$ is hyperbolic.
\end{proof}

\section{Properties of the invariant}
\label{sec: properties}

In this section, we collect some basic properties of $h(K)$. These
include its behavior under mirroring, its relationships with the
original Lin invariant and the Levine-Tristram signature, and its
parity.  The relationship with the signature allows us to compute
$h(K)$ for the first time in Example~\ref{Ex:Trefoil2}.

\subsection{Mirrors} The invariant $h(K)$ behaves nicely under taking
the mirror image:

\begin{proposition}
  \label{Prop:Mirror}
  If \(K\) is real representation small, then so is its mirror
  \(\Kbar\) and \(h(\Kbar) = - h(K)\).  Moreover,
  \(\eh(\Kbar) = - \eh(K)\).
\end{proposition}

\begin{proof}
First, the knot $\Kbar$ is real representation small as there is an
isomorphism from \(\pi_1(M_{\Kbar})\) to \(\pi_1(\EK)\) sending
meridians to meridians.  Second, note that if \(K = \betahat\), then
\(\Kbar = \betahat\vphantom{\beta}^{-1}\), and so for
$c \in [-2 , 2] \setminus D_K$ we have
\begin{align*}
  h(\Kbar)&  =  (-1)^{|\beta^{-1}|}\pair{\cL^c,\ {(\beta^{-1})^*
           \cL^c}}_{\cX^c_\bal(S_{2m})} && \\
          & =  (-1)^{|\beta|}\pairbetaLL
                     && \text{by functoriality and Lemma~\ref{Lem:BraidActionOrientation}}\\
          & = - (-1)^{|\beta|}\pairLbetaL
                     && \text{since \(\dim \cL^c = 2m-3\) is odd}   \\
          & = - h(K). &&
\end{align*}
The calculation above is underlaid by a
particular isomorphism from \(\pi_1(M_{\Kbar})\) to \(\pi_1(\EK)\).
While we have always identified $S_{2m}$ with $\partial H_1$ in
Figure~\ref{fig:plat}, we could have instead viewed it as
$\partial H_2$, and defined
\[
\cX^c_{!}(\beta) = \big((\beta^{-1})^* \cL^c \big) \cap \cL^c
\mtext{and}
  h_{!}(K) = (-1)^{|\beta|}\pair{(\beta^{-1})^* \cL^c, \ \cL^c}_{\cX^c_\bal(S_{2m})}
\]
\begin{figure}
\begin{center}
\begin{tikzoverlay}[width=9cm]{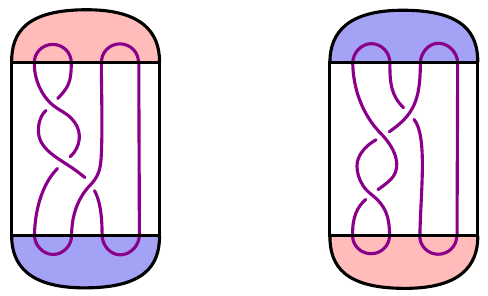}
  \begin{scope}[nmdstd, line width=1pt]
    \draw[<-] (63.2,31.6) -- (37.1,31.6) node[pos=0.5, above] {$\xi$};
    \draw[<-] (1, 47.65) -- +(-5, 0) node[left] {use for $h_!(K)$};
    \draw[<-] (1 ,12.4) -- +(-5, 0) node[left] {use for $h(K)$};
    \node at (-2, 30) {$\beta$};
    \node at (105, 30) {$\beta^{-1}$};
    \node[] at (17.8, 5.5) {$H_1$};
    \node[] at (17.8, -5) {\small $K = \betahat$};
    \node[] at (17.8,54.5) {$H_2$};
    \node[] at (83.1,5.5) {$H_1$};
    \node[] at (83.1, -5) {\small $\Kbar = \betahat\vphantom{\beta}^{-1}$};
    \node[] at (83.1,54.5) {$H_2$};
    \draw[<-] (100, 12.4) -- + (5, 0) node[right] {use for $h(\Kbar)$};
  \end{scope}
\end{tikzoverlay}
\end{center}

\caption{The homeomorphism \(\xi \maps \EK \to M_{\Kbar}\) is
  reflection in the plane of the page followed by $\pi$-rotation about
  the horizontal axis.  This takes the copy of $S_{2m}$ used to
  compute $h_!(K)$ at left to the one used to compute $h(\Kbar)$ at
  right.}
  \label{fig: flip and roll}
\end{figure}
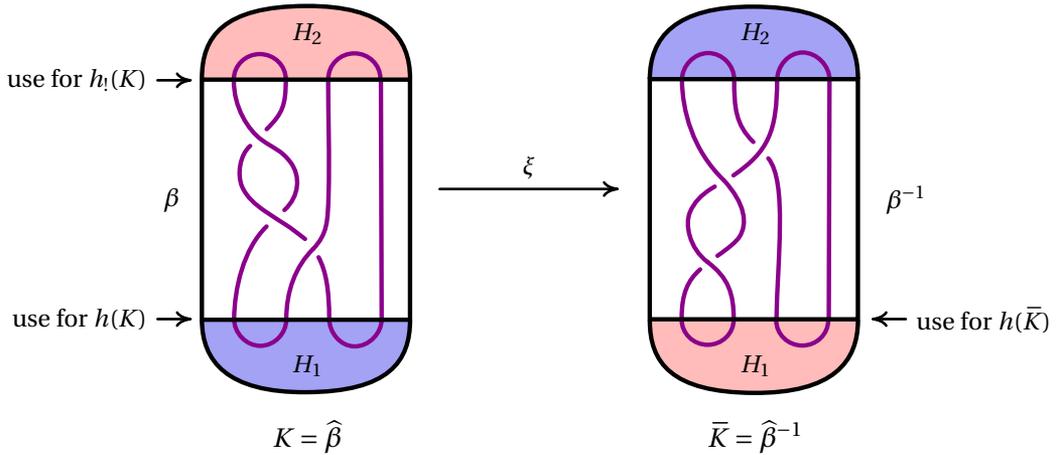
This is illustrated in the lefthand part of Figure~\ref{fig: flip and roll}.
As subsets of $\cX^c_\bal(S_{2m})$, note that
$\beta^*\big( \cX^c_{!}(\beta) \big) = \cX^c(\beta)$.  Moreover, the
previous calculation shows $h_{!}(K) = h(K)$.  Now consider the
homeomorphism $\xi$ from $M_K$ to $M_\Kbar$ shown in Figure~\ref{fig:
  flip and roll}.  This induces an isomorphism
\(\xi_* \maps \pi_1(\EK) \to \pi_1(M_{\Kbar})\), which takes $s_i$ in
$\partial H_2$ for $\EK$ to $s_i$ in $\partial H_1$ for $M_\Kbar$.
Note here that the orientations of the $s_i$ are preserved, not
reversed, so if we take $s_1$ as the meridian in both cases, we get
$\xi_*(\lambda) = \lambdabar\vphantom{\lambda}^{-1}$.  This shows that
$\cX^c_{!}(\beta) = \cX^c(\beta^{-1})$ as subsets of
$\cX^c_\bal(S_{2m})$, as in Proposition~\ref{prop: Casson style}.  We
then have $h_{!}(K) = -h(\Kbar)$ since both are computing the
intersection numbers of $\cL^c$ and $(\beta^{-1})^* \cL^c$, just in
the opposite order.  To extend this to $\htil_{!}(K) = -\htil(\Kbar)$,
given a connected component $X$ of $X^{2,\irr}_\SLR(K)$, we have to
show $\Xbar = \xi^*(X)$ contributes the same to $ -\htil(\Kbar)$ as
$X$ does to $\htil_{!}(K)$.  Now $X = \Xbar$ as subsets of
$\cX^2_\bal(S_{2m})$, so the contributions of the local intersection
numbers match.  Now consider the two preimages $\Xtil_\pm$ of $X$ in
$\XA_\tSLR(K)$. We want to compare $\lambda^*(\Xtil_\pm)$ and
$\lambdabarstar(\Xtil_\pm)$.  Since
$\lambda^*(\Xtil_+) = -\lambda^{*}(\Xtil_-)$ and
$\lambdabarstar(\Xtil_+) = -\lambdabarstar(\Xtil_-)$ by
Lemma~\ref{lem: para preimages}, as
$\xi_*(\lambda) = \lambdabar\vphantom{\lambda}^{-1}$ we see that the
corresponding terms of $-\htil(\Kbar)$ and $\htil_!(K)$ match,
completing the proof.
\end{proof}

\subsection{Relation with the signature}
\label{sec: rel with signature}

For \(\alpha \in (0, \pi)\) with $\Delta_K(e^{2 i\alpha}) \neq 0$, we
let \(h'_\alpha(K)\) be the Lin invariant as defined by Heusener in
\cite{Heusener2003}, where it is denoted $h^{(\alpha)}(K)$.  We will
show:

\begin{proposition}
  \label{lemma: us vs Heusener}
  When $K$ is real representation small and $c \notin D_K$, then 
  \(h^c_\SU(K) = h'_\alpha(K)\), where \(c = 2 \cos \alpha\).
\end{proposition}
This has the following important consequence:

\begin{corollary}
  \label{Cor: hSU=sig}
  For \(c \in [-2,2] \setminus D_K\), one has 
  \(h^c_\SU(K) = -\frac{1}{2} \sigma_K(e^{2i\alpha})\), where
  \(\sigma_K\) is the Levine-Tristram signature function.
\end{corollary}
This follows from the corresponding statement for \(h'_\alpha\), which
was proved by Herald \cite{Herald1997a} using gauge theory and by
Heusener and Kroll \cite{HeusenerKroll1998} following Lin's original
proof \cite{Lin1992} for \(\alpha = \frac{\pi}{2}\); see also \cite[\S
5]{Heusener2003}.  Here, our convention for the signature is that
positive knots have negative signature.  Note that \cite{Lin1992} and
\cite{HeusenerKroll1998} use the opposite convention, so their
statements have no minus sign.

\begin{figure}
  \begin{center}
    \begin{tikzpicture}[
  scale=2.6,
  nmdstd,
  line width=1.2pt,
  every node/.style={color=black},
  every circle/.style={radius=0.03}]

  \def\basicxscale{2.5}
  \def\basicyscale{1}
  
  \coordinate (X) at (\basicxscale, 0);
  \coordinate (Y) at (0, \basicyscale);

  \coordinate (U) at (${0.1666667}*(X)$);
  \coordinate (V) at (${1-0.1666667}*(X)$);

  \node[below right] at (U) {$\frac{1}{6}$};
  \node[above left] at (V) {$\frac{5}{6}$};

  \begin{scope}[dashed, color=simplealex, line width=1]
  \draw ($(U) + 1.2*(Y)$) -- ($(U) - 2.4*(Y)$);
  \draw ($(V) + 1.2*(Y)$) -- ($(V) - 2.4*(Y)$);
  \end{scope}
  
  \draw[coor grid, xstep=\basicxscale, ystep=\basicyscale] ($-0.3*(X)-1.2*(Y)$) grid ($1.3*(X)+ 1.2*(Y)$);
  \begin{scope}[color=axesgray]
    \draw[->] ($-0.3*(X)$) -- ($1.4*(X)$) node[right] {$\mu^*$};
    \draw[->] ($-1.2*(Y)$) -- ($1.3*(Y)$) node[right=0.04]{$\lambda^*$};
  \end{scope}

  \begin{scope}[line width=2pt]
    \draw[color=locus]  ($-0.25*(X)$) -- ($1.3*(X)$);
    \foreach \x in {0, 1}{
      \coordinate (A) at ($\x*(X) + (Y)$);
      \coordinate (B) at (${\x + 0.1666667}*(X)$);
      \coordinate (C) at ($\x*(X) - (Y)$);
      \coordinate (D) at (${\x - 0.1666667}*(X)$);
      \draw[color=locus, mid arrow=0.6] (A) -- (B); 
      \draw[color=locus, mid arrow=0.6] (D) -- (C);
      \draw[fill=otherparabolic, line width=0.5pt] (A) circle;
      \draw[fill=otherparabolic, line width=0.5pt] (C) circle;
      \draw[fill=simplealex, line width=0.5pt] (B) circle;
      \draw[fill=simplealex, line width=0.5pt] (D) circle;
    }
  \end{scope}
  \path (A) -- (B) node[above right, pos=0.5] {$\OTEL{K}$};


  \draw[color=Lline, mid arrow=0.7] ($0.95*(X) - 1.5*(Y)$) node[right] {$V_\alpha$} -- +($2.8*(Y)$);

  \coordinate (L1) at ($0.5*(U)$);
  \coordinate (L2) at ($0.5*(X)$);
  \coordinate (L3) at ($(V)!0.5!(X)$);
  \coordinate (L0) at  ($-1.8*(Y)$);
  
  \node[left] at (L0) {$\sigma_K$};
  \node at ($(L0)+(L1)$) {$0$};
  \node at ($(L0)+(L2)$) {$-2$};
  \node at ($(L0)+(L3)$) {$0$};
    
  \coordinate (L0) at  ($-2.0*(Y)$);
  
  \node[left] at (L0) {$h_\SU^c$};
  \node at ($(L0)+(L1)$) {$0$};
  \node at ($(L0)+(L2)$) {$1$};
  \node at ($(L0)+(L3)$) {$0$};

  \coordinate (L0) at  ($-2.2*(Y)$);

  \node[left] at (L0) {$h_\SLR^c$};
  \node at ($(L0)+(L1)$) {$1$};
  \node at ($(L0)+(L2)$) {$0$};
  \node at ($(L0)+(L3)$) {$1$};

\end{tikzpicture}
  \end{center}

  \vspace{-0.5cm}

  \caption{For the positive trefoil $K$, this figure shows a portion of
    $\OTEL{K}$ with its irreducible characters oriented according to
    Section~\ref{sec: gen pic}; in particular
    $\pair{\OTEL{K}, V_\alpha} = h^c_\SLR(K)$ for $c = 2 \cos \alpha$.
    The orientation is derived from Theorem~\ref{thm: hSL+hSU}
    and Corollary~\ref{Cor: hSU=sig}, which allow us to compute
    $h^c_\SLR(K)$ and $h^c_\SU(K)$ from $\sigma_K(e^{2 i \alpha})$.
    Note that this orientation is preserved by the translations
    in the $D_\infty$ action but reversed for rotations, as
    is consistent  with Lemma~\ref{lem:b symmetry}.
  }
  \label{fig: trefoil final}
\end{figure}
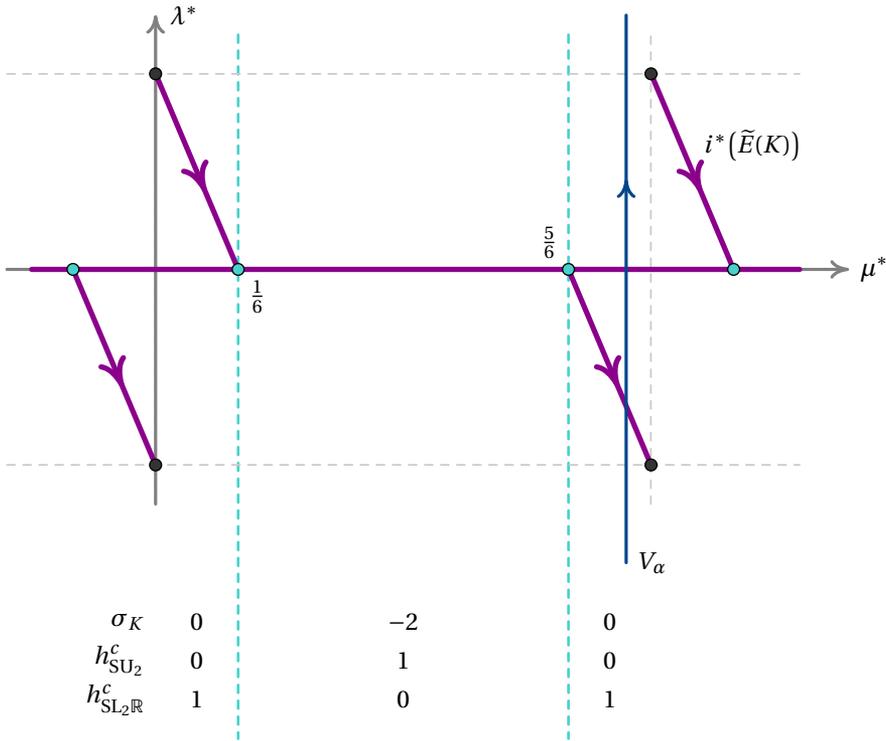

\begin{example} 
\label{Ex:Trefoil2}
Let \(K\) be the positive trefoil, so \(\sigma_K(e^{2i\alpha}) = -2\)
for \(\alpha \in (\frac{\pi}{6},\frac{5\pi}{6})\) and is \(0\) for
\(\alpha \in [0, \frac{\pi}{6}) \cup (\frac{5\pi}{6},1]\). We saw in
Example~\ref{Ex:Trefoil} that \(h^c_\SLR(K) = 0\) for
\(\alpha \in (\frac{\pi}{6},\frac{5\pi}{6})\), so \(h(K) = 0+1 =
1\). Applying the relation with the signature again, we see that
\(h^c_\SLR(K) = 1\) for
\(\alpha \in [0, \frac{\pi}{6}) \cup (\frac{5\pi}{6},1]\), as we
expected from looking at Figure~\ref{Fig:Character of Trefoil}.  We
summarize the full picture in Figure~\ref{fig: trefoil final}.
\end{example} 

\begin{proof}[Proof of Proposition~\ref{lemma: us vs Heusener}]
For \(\beta \in B_{2m}\), recall that 
\[
h^c_{\SU}(\betahat) 
 =  (-1)^{|\beta|} \langle L_1, L_2 \rangle_{X}, 
\]
where \(X = \cX^c_{+}(S_{2m}) = X^\cirr_\SU(S_{2m})\),
\(L_1 = \cL^c_+\), and \(L_2 = \beta^*(L_1)\).  Heusener's
$h'_\alpha(\betahat)$ is defined as the intersection number of the
same objects but with possibly different orientations and global sign.
If $X'$, $L_1'$, and $L_2'$ denote these manifolds with Heusener's
orientation, one has
\[
  h'_\alpha(\betahat) = (-1)^{m} \langle L'_1, L'_2 \rangle_{X'}.
\]
Thus to complete the proof, we must compare our orientations with
Heusener's.  We have already observed in Remark~\ref{rem: orient vs H
  for S_n} that \([X] = - [X']\).  The relation between \([L_k]\) and
\([L_k']\) is more complicated, and depends on $\beta$ as we now
detail.

The orientations on \([L_k]\) and \([L_k']\) are defined by
choosing a set of generators \(T_k\) for \(\pi_1(H_k)\) and using
\(T_k\) to identify \(X^\cirr_\SU(H_k,T_k)\) with
\(X^\cirr_\SU(F_m,T)\). (Our orientation on \(X^\cirr_\SU(F_m,T)\)
agrees with Heusener's, but this is a bit of a moot point. The
orientation on \(X^\cirr_\SU(F_m,T)\) appears twice (once for \(L_1\)
and once for \(L_2\)), so we would get the same answer even if the two
orientations on \(X^\cirr_\SU(F_m,T)\) were different.)

While our choice of generators for $\pi_1(H_1)$ is independent of
$\beta$, Heusener orients the knot $\betahat$ and uses generators of
$\pi_1(H_1)$ that are compatible with that orientation.  For related
reasons, while we always have $[L_2] = (\beta^*)_*[L_1]$, it turns out
that $[L'_2]$ is $(\beta^*)_* [L_1']$ or $- (\beta^*)_* [L_1']$
depending on $\beta$.

To make the comparison easy, we will work with particularly nice
$\beta$. Specifically, given $K$, we claim there is a $\beta$ with an
orientation on $\betahat$ so that:
\begin{enumerate}
\item
  \label{item: bridges}
  Every overbridge is oriented right to left, and every
  underbridge is oriented left to right.
  
\item
  \label{item: perm}
  The corresponding permutation $\betabar$ fixes all odd
  numbers and acts on the evens as the \emph{inverse} of the cycle
  $(2\ 4 \ 6\ \ 2m)$.
\end{enumerate}
Starting from any $\beta$ and any orientation on $\betahat$, part
(\ref{item: bridges}) is easy to ensure by adding a twist to reverse
each bridge with the wrong orientation.  To arrange (\ref{item:
  perm}), begin by walking along $\betahat$ in the prescribed
orientation.  As you go, number the overbridges $1, 2, \ldots, m$ and
do the same for the underbridges.  Now slide the overbridges in front
or behind each other to get a new braid $\beta'$ so that the
overbridge numbers appear in order from left to right; repeat this
process for the underbridges. You can now check that (\ref{item:
  perm}) must hold, so from now on we assume $\beta$ has both these
properties.

Condition (\ref{item: bridges}) means that in Section~4 of
\cite{Heusener2003} we have all $\epsilon^{(i)}_\ell = +1$, and in
this situation the construction of \cite{Heusener2003} matches ours on
the nose.  In particular, $[L_1] = [L'_1]$ and $[L_2] = [L'_2]$.
Since as noted $[X'] = -[X]$, to complete the proof we simply need to
show that $(-1)^{m + 1} = (-1)^{\abs{\beta}}$.  Since each braid
generator corresponds to a transposition in the symmetric group, the
latter is the sign of the permutation $\betabar$.  As $\betabar$ is an
$m$-cycle, it has sign $(-1)^{m+1}$ as needed.
\end{proof}

\subsection{Geometry at simple roots}

At a simple root of \(\Delta_K(t)\) on the unit cicle, the equivariant
signature \(\sigma_K\) must jump by \(\pm 2\).  The sign of this jump
is related to the local geometry of the character variety by part
(\ref{item: sig to side}) of the following lemma, which we thank Chris
Herald for explaining to us.  

\begin{lemma}
  \label{lem: herald power}
  Suppose a knot $K$ has a simple root $e^{i 2 \omega_0}$ of
  $\Delta_K(t)$ for $0 < \omega_0 < \pi$.  Then there is a smooth arc
  $\gamma \maps [0, 1] \to \Etil(K)$ such that
  \begin{enumerate}
  \item $i^*(\gamma(1)) = (\omega_0/\pi, 0)$ with $\gamma(1)$ a lift
    of $\chi^c \in X^\red_\SLR(K)$ where $c = 2 \cos \omega_0$.
    
  \item $i^* \circ \gamma \maps [0, 1] \to \Vtil(K)$ is a smooth
    embedding.
  \end{enumerate}
  If in addition the derivative
  \( (\tr_\lambda \circ \gamma)'(0) \neq 0\), then we have
  \begin{enumerate}[resume*]
  \item
    \label{item: sig to side}
    If $\sigma_K(e^{i 2 \omega})$ jumps by $-2$ as $\omega$
    increases past $\omega_0$, then $\gamma([0, 1))$ lies entirely
    above the horizontal axis.  If instead $\sigma_K(e^{i 2 \omega})$
    jumps by $+2$, then $\gamma([0, 1))$ lies entirely below the
    horizontal axis.
  \end{enumerate}
\end{lemma}

\begin{remark}
  \label{rem: long at roots}
  We do not know of an example of a simple root of $\Delta_K(t)$ on
  the unit circle where \( (\tr_\lambda \circ \gamma)'(0) = 0\).
  Indeed, it seems likely this never happens, and so (\ref{item: sig to
    side}) always holds.  (For multiple roots, it does arise, see the
  upper-right example in Figure~10 of \cite{CullerDunfield2018}.)  The
  question is subtle in that for other
  $\alpha \in \pi_1(\partial M_K)$ the derivative
  $(\tr_\alpha \circ \gamma)'(0)$ can vanish: for the trefoil,
  $\tr_\alpha \circ \gamma \equiv 2$ for $\alpha = 6 \mu + \lambda$ as
  noted in Section~\ref{Ex:Trefoil}.
\end{remark}

\begin{proof}[Proof of Lemma~\ref{lem: herald power}]
The idea follows \cite{HeraldZhang2019}, which adds \cite{Herald1997a,
  Herald1997b} to the results of \cite{HeusenerPortiPeiro2001} (see
also \cite{HeusenerPorti2005}) to strengthen
\cite[Lemma~7.3]{CullerDunfield2018}. To start, from \cite[Corollary
1.3]{HeusenerPortiPeiro2001}, one has a smooth arc
$\gammabar \maps [-1, 1] \to X(K)$ where $\gammabar(0) = \chi^c$ with
$\gammabar(t)$ in $X_\SLR^\irr(K)$ for $t < 0$ and in $X_\SU^\irr(K)$
for $t > 0$. As argued in \cite[page~2818]{HeraldZhang2019},
Corollary~6 of \cite{HeraldZhang2019} allows us to assume that
$\tr_\lambda\gammabar(t) = 2$ only for $t = 0$.  Additionally,
Lemma~7.3(4) of \cite{CullerDunfield2018}, or, more accurately, the
fact that
$\iota^* \maps H^1\big(\EK; \mathfrak{sl}_2(\C)_{\rho^+} \big)\to
H^1\big(\partial \EK; \mathfrak{sl}_2(\C)_{\rho^+}\big)$ is injective,
which follows from \cite[Equation 7.4]{CullerDunfield2018}, implies
there is some $\alpha \in \pi_1(\partial M_K)$ such that, after
possibly restricting the domain of $\gammabar$, we have
$\tr_\alpha \circ \gammabar$ is a smooth embedding.

Lifting to $\Etil(K)$ now gives all the claims except for the part of
(\ref{item: sig to side}) about which side of the horizontal axis
$\gamma([0, 1))$ lies on.  This follows from the corresponding
behavior on the $\SU$ side, which is understood by
\cite[Prop.~4.2]{Herald1997a} and \cite[Cor.~3]{Herald1997b}, see the
remark immediately after \cite[Prop.~4.2]{Herald1997a}.  The
hypothesis on \((\tr_\lambda \circ \gamma)'(0)\) ensures that if
\(\gammabar\) lies above the reducible line on the \(\SU\) side, it
lies below on the \(\SLR\) side, and vice versa. Rather than trace
through the various conventions, we can justify our claimed ``sign
rule'' that a negative jump in $\sigma_K(e^{i 2 \omega})$ means 
$\gamma$ is above the axis by looking at the positive trefoil in
Figure~\ref{fig: trefoil final}.
\end{proof}

\subsection{Parity of \(h\)}

The map \(b \maps X^c(K) \to X^{-c}(K)\) introduced in
Section~\ref{subsec: X_K sym} preserves \(X^0(K)\) setwise.  The fixed
points of $b$ on \(X^{0,\irr}(K)\) were classified by Nagasato and Yamaguchi in 
\cite[Corollary~1]{NagasatoYamaguchi2012}: they are all in
$X^{0, \irr}_\SU(K)$ and consist of characters of binary-dihedral
representations.  We now use this to show:

\begin{proposition}
  \label{cor: h of small mod 2}
  If \(K\) is a small knot, then
  \(h(K) \equiv \frac{1}{2} \sigma(K) \equiv \frac{1}{2}(\det(K) - 1) \bmod 2\).
\end{proposition}
 
\begin{proof}
By Lemma~\ref{lem:small to few chars}, as \(K\) is small, the set
\(X^{0, \irr}_\SLR(K)\) is finite.  As noted, $b$ has no fixed points
on \(X^{0, \irr}_\SLR(K)\), so that set can be partitioned into pairs
\(\{\chi, b(\chi)\}\) with \(\chi \neq b (\chi)\). Lemma~\ref{lem:b
  symmetry} tells us that the local intersection numbers \(n_{\chi}\)
and \(n_{b(\chi)}\) are equal. It follows that \(h^0_\SLR(K) \) is
even, and hence that $h(K) = h^0(K)$ is congruent to
$h_\SU^0(K) = -\frac{1}{2} \sigma(K)$ modulo 2, establishing $h(K)
\equiv \frac{1}{2}\sigma(K) \bmod 2$.

Continuing this idea, only the binary-dihedral characters can
contribute to $h_\SU^0(K)$, and there exactly are $(\det(K) - 1)/2$
such \cite{Klassen1991}.  If each binary-dihedral characters
contributes $\pm 1$ to $h_\SU^0(K)$, we would get that
$h_\SU^0(K) \equiv (\det(K) - 1)/2 \bmod 2$, implying our second
claim.  We do not know if this holds, but for any knot one has
\(\frac{1}{2} \sigma(K) \equiv \frac{1}{2}(\det(K) - 1) \bmod 2\)
by \cite[Theorem~5.6]{Murasugi1965}, completing the proof.
\end{proof}

Next, we relate the parity of \(h(K)\) to the \(\SLC\) character
variety, leading to Corollary~\ref{cor: really ideal} which shows
there are ideal points that are limits of $\SLR$ characters when
$\sigma(K)/2$ is odd.  Fix a braid \(\beta\) with \(\betahat = K\). As
in Section~\ref{Sec:Definition of h}, this determines a splitting
\(\EK = H_1 \cup_{S_{2m}} H_2\). We define \(X^{\irr}_\C(\beta) \) to
be the intersection of the two smooth varieties
\(L_{j,\C} = i_j^*\big(X^\irr_\C(H_j,T)\big)\) for \(j=1,2\)
inside the smooth variety \(X^\irr_\C(S_{2m},S)\). Note that this
intersection may be nonreduced, i.e. the components of
\(X^\irr_\C(\beta)\) may have multiplicities; if it is, we view it as
a nonreduced scheme rather than passing to the reduction. (As reduced
schemes, \(X^\irr_\C(\beta) = X^\irr_\C(K)\), but it is not clear that
the multiplicity of each component is an invariant of \(K\).)

When $K$ is small, every component $C$ of $X^\irr_\C(K)$ is a complex
curve and $\tr_\mu \maps C \to P^1(\C)$ is nonconstant (see
Lemma~\ref{lem:small to few chars}).  We define $d(K, \beta)$ to be
the (total) degree of the map
$\tr_\mu \maps X^\irr_\C(\beta) \to P^1(\C)$; this is a weighted sum of
all the Culler-Shalen seminorms of $\mu$ \cite{BoyerZhang1998} which
is related in turn to the A-polynomial
\cite[Section~8]{BoyerZhang2001}. We conjecture that $d(K, \beta)$
depends only on $K$ but do not prove this here.  We next show:

\begin{proposition}
  If \(K\) is a small knot, then \(d(K,\beta) \equiv h(K) \bmod 2\). 
\end{proposition}

\begin{proof}
For \(c \in \C\), let
\(L_{j,\C}^\cirr = L_{j,\C}^\irr \cap \tr_\mu ^{-1}(c)\). (This
should also be the scheme theoretic intersection, but
\(\tr_\mu \maps L_{j,\C}^\irr \to \C\) is a submersion by
Lemma~\ref{lem:c slice smooth}, so it is the same as the setwise
intersection.) Then
\[
  d(K,\beta) = \pair{ L_{1,\C}^\cirr , \ L_{2,\C}^\cirr}_{X^{c, \irr}_\C(S_{2n})}
\]
Now if \(A_{\C}, A'_{\C} \subset B_\C\) are the $\C$ points of
algebraic sets defined over \(\R\), by \cite[Chapter~13]{Fulton1998}
their real points \(A_{\R}\), \(A'_{\R}\), and \(B_\R\) satisfy
\[
  \pair{A_{\C}, A'_{ \C}}_{B_\C} \equiv
  \pair{A_{\R}, A'_{\R} }_{B_\R} \bmod 2.
\]
In our case, the real parts of \(L_{j,\C}^\cirr\) and
\(X^\cirr_\C(S_{2n})\) are the spaces we have been referring to as
\(L^\cirr_j\) and \(X^\cirr(S_{2n})\). Hence
\[
  d(K,\beta)  \equiv  \pair{L_1^\cirr , \ L_2^\cirr}_{X^{c,
      \irr}(S_{2n})} \equiv h(K) \bmod 2
\]
which proves the proposition.
\end{proof}

Recall that the affine algebraic set \(X^\irr_\R(K)\) and the map
\(\tr_\mu\) can be extended to a projective variety \(\Xhatirr_\R(K)\)
and a map \(\tr_\mu \maps \Xhatirr_\R(K) \to P^1(\R)\). The points in
\(\tr_\mu^{-1}(\infty)\) are called \emph{real ideal points}. We just
showed that the mod \(2\) degree of this map is given by \(h(K)\). If
the mod \(2\) degree is nonzero, \(\tr_\mu^{-1}(\infty)\) must be
nonempty. Hence we deduce:

\begin{corollary}
  \label{cor: really ideal}
  If \(K\) is a small knot with \(\sigma(K) \equiv 2 \bmod 4\), or equivalently 
  $\det(K) \equiv 3 \bmod 4$, then \(\Xhatirr_\R(K)\) contains a real
  ideal point.
\end{corollary}
Unsurprisingly, there are also knots with real ideal points where
\(\sigma(K) \not\equiv 2 \bmod 4\), for example the figure-8 knot
which has $\sigma(K) = 0$.  Moreover, real ideal points are extremely
common, even outside the context of exteriors of knots in $S^3$,
raising:
\begin{question}
  \label{qu: real ideal}
  Does every small knot in $S^3$ have a real ideal point?
\end{question}
An initial search using \cite{CullerAPoly} suggests that the answer
might well be yes.

%
 
\section{Applications to left orderings}
\label{sec: LO apps}

\subsection{Left Orderings}

Let \(G\) be a nontrivial group. We say that \(G\) is
\emph{left-orderable} (or for short, that \(G\) is LO), when there is
a total order on \(G\) which satisfies \(gx > gy\) whenever
\(x>y\). (By convention, the trivial group is not left-orderable.) The
study of left-orderable 3-manifold groups dates back to the work of
Boyer, Rolfsen, and Wiest \cite{BRW2005} and was given fresh impetus
by the L-space conjecture of Boyer-Gordon-Watson and Juhasz:

\begin{LspaceConj}[\cite{BGW2013, Juhasz2015}] If \(Y\) is a closed
  prime orientable 3-manifold then the following are equivalent:
\begin{enumerate}
\item \(Y\) is not a Heegaard Floer L-space, i.e.~$\HFhat_\red(Y) \neq 0$.
\item \(\pi_1(Y)\) is left-orderable.
\item \(Y\) admits a coorientable taut foliation. 
\end{enumerate}
\end{LspaceConj}
For brevity, when $\pi_1(Y)$ is LO we will call $Y$ itself LO.

In this section, we use the techniques of this paper to prove many $Y$
are LO, using the following theorem of Boyer-Rolfsen-Wiest: 
\begin{theorem}[\cite{BRW2005}]
  \label{Thm:SLR order}
  If \(Y\) is a prime 3-manifold with a nontrivial homomorphism
  \(\pi_1(Y) \to \tSLR\), then \(\pi_1(Y)\) acts faithfully on \(\R\),
  and hence \(Y\) is LO.
\end{theorem}
In particular, if \(\rho \maps \pi_1(Y) \to \SLR\) is a nontrivial
homomorphism whose Euler class \(e(\rho) \in H^2(Y; \Z)\) is zero,
then \(\rho\) lifts to \(\rhotilde \maps \pi_1(Y) \to \tSLR\) where
Theorem~\ref{Thm:SLR order} applies.  We will also consider the target
group $\PSLR$, where again any \(\rhohat \maps \pi_1(Y) \to \PSLR\)
has an Euler class \(e(\rhohat) \in H^2(Y; \Z)\) that obstructs
lifting $\rhohat$ to $\tSLR$.

We begin with results showing cyclic branched covers are LO
(Theorem~\ref{thm: branched LO}).  We then compute $h(K)$ for
alternating knots (Section~\ref{sec: determine h}), which in turn
strengthens the results on branched covers (Remark~\ref{rem: better
  branched bound}).  The rest of the section is devoted to showing
certain Dehn surgeries on $K$ are LO; key results there include
Theorems~\ref{thm: extended Lin to LO}, \ref{thm: 2-bridge LO}, and
\ref{thm:lens arcs}.

\subsection{Branched covers}
\label{sec: branched}

Let \(\Sigma_n(K)\) be the \(n\)-fold cyclic branched cover of a prime
knot \(K\).  The manifold \(\Sigma_n(K)\) is obtained by Dehn filling
the \(n\)-fold cyclic cover \(\Mtil_{K,n}\) of \(\EK\) along the curve
\(\mutil_n\) which is the \(n\)-fold cyclic cover of \(\mu\).
Moreover, \(\Sigma_n(K)\) is the \(n\)-fold cyclic orbifold cover of
the orbifold \((S^3,K_n)\) where the knot $K$ is labeled by
$\Z/n$. The orbifold fundamental group \(\pi_1(S^3,K_n)\) is
isomorphic to \(\pi_1(\EK)/\langle \mu^n \rangle\), and
\(\pi_1\big(\Sigma_n(K)\big)\) is the kernel of the homomorphism
\(\pi_1(S^3,K_n) \to \Z/n\) where $\mu \mapsto 1$.

\begin{remark}
  \label{Rem:PrimeCover}
  It follows from the Orbifold Sphere Theorem (see
  e.g.~\cite[Theorem~10.2]{Maillot2003}) that if \(\Sigma_n(K)\) is
  reducible, then the orbifold \((S^3,K_n)\) is reducible.
  Consequently, when $K$ is prime so are all \(\Sigma_n(K)\).
\end{remark}

Now suppose \(\rho \maps \pi_1(\EK) \to \SLR\) is a representation
with \(\tr\big(\rho(\mu)\big) = 2 \cos({k\pi}/{n})\) for some integer
$k$ with $0 < k < n$.  Then \(\rho\) is elliptic and satisfies
\(\rho(\mu^n) = \pm I\), so \(\rho\) descends to a representation
\(\rhohat \maps \pi_1(S^3, K_n) \to \PSLR\). We let \(\rhohat_n\) be
the restriction of \(\rhohat\) to \(\pi_1\big(\Sigma_n(K)\big)\).
The next result, which is Theorem~3.3 of \cite{Hu2015}, has been much
used to show branched covers are LO:

\begin{lemma}[\cite{Hu2015}] 
\label{Lem:Branched Cover}
If \(\rhohat\) is irreducible, then \(\rhohat_n\) is nontrivial and
\(e(\rhohat_n) = 0 \).
\end{lemma}






Let $A_K = \setdef{\alpha \in [0, \pi]}{2 \cos \alpha \in D_K}$ and
consider the corresponding partition
\[
  0 < \alpha_1 < \alpha_2 < \cdots < \alpha_k < \pi \mtext{of
    $[0, \pi]$.}
\]
The signature function $\sigma_K(e^{2i\alpha})$ is constant on the
interiors of the intervals of this partition.  Define $w_K$ to be the
minimum length of any of these intervals, including $[0, \alpha_1]$ and
$[\alpha_k, \pi]$.  Our main result on branched covers is:

\begin{theorem}
  \label{thm: branched LO}
  
  If a prime knot \(K\) is real representation small and
  \(\sigma_K(\omega)\) is not identically $0$, then \(\Sigma_n(K)\) is
  LO for all \(n> \frac{\pi}{w_K}\).
\end{theorem}

\begin{proof}
Taking $c = 2 \cos \alpha$, we can view $h_\SLR^c(K)$ as a function of
$\alpha$ in $[0, \pi] \setminus A_K$ rather than of $c$ in
$[-2, 2] \setminus D_K$; in the new setting, we will write
$h_\SLR^\alpha(K)$ and use other similar notation.  By
Theorem~\ref{thm: hSL+hSU} and Corollary~\ref{Cor: hSU=sig}, we
then have
\begin{equation}
  \label{eq: solve for SLR}
 h_\SLR^{\alpha}(K) = h(K) - h^{\alpha}_\SU(K) = h(K) + \frac{1}{2}\sigma_K(e^{2i\alpha}).
\end{equation}
As $\sigma_K$ is nonconstant, no matter what $h(K)$ is there will be
an interval $I$ in our partition where $h_\SLR^{\alpha}(K)$ is nonzero
on the interior of $I$.

As $n > \pi/w_K$, the interval $I$ has length greater than $\pi/n$.
Thus the interior of $I$ contains a point $\alpha_0 = \frac{\pi k}{n}$
where $k$ is an integer.  Thus $h^{\alpha_0}_\SLR(K)$ is nonzero,
which means that $X^{c_0, \irr}_\SLR(K)$ is nonempty for
$c_0 = 2 \cos \alpha_0$.  In particular, there is an irreducible
$\rho \maps \pi_1(\EK) \to \SLR$ which induces an irreducible
$\rhohat \maps \pi_1\big(S^3, K_n\big) \to \PSLR$.  By
Lemma~\ref{Lem:Branched Cover}, we have a nontrivial representation
\(\rhohat_n \maps \pi_1\big(\Sigma_n(K)\big) \to \PSLR\) with
\(e(\rhohat_n) = 0\).  As \(\Sigma_n(K)\) is prime by
Remark~\ref{Rem:PrimeCover}, it is LO by Theorem~\ref{Thm:SLR
  order}.
\end{proof}

\begin{remark}
  \label{rem: better branched bound}
  For any concrete signature function $\sigma_K$, one can often refine
  the result of Theorem~\ref{thm: branched LO}, especially if one
  knows $h(K)$.  For example, the alternating knot $K15a78855$ has
  $\Delta_K = 8t^6 - 21t^5 + 27t^4 - 27t^3 + 27 t^2 - 21 t + 8$, which
  is an irreducible polynomial all of whose roots are on the unit
  circle.  The resulting partition is roughly
  $[0, 0.139, 0.398, 0.963, 2.178, 2.743, 3.002, 3.141]$ and
  $\frac{1}{2}\sigma_K$ has values $0, -1, -2, -3, -2, -1, 0$ on the
  corresponding open intervals. We checked $K$ is small using the
  method of \cite{DunfieldGaroufalidisRubinstein2022}, so
  Theorem~\ref{thm: branched LO} applies for $n \geq 23$.  However, by
  Corollary~\ref{cor: h for alt} below, we have
  $h(K) = -\frac{1}{2}\sigma(K) = 3$.  Hence from (\ref{eq: solve for
    SLR}), we see that $h^\alpha_\SLR(K)$ is nonzero on every interval
  except the middle one $[0.963, 2.178]$.  As $\Delta_K$ has no
  cyclotomic factors, no point of the form $k \pi/n$ is in our
  partition. Thus we have $h^\alpha_\SLR(K)$ is valid and nonzero for
  $\alpha = \pi/n$ whenever $\alpha < 0.963$, that is for $n \geq 4$.
  As $\Sigma_2(K)$ is not LO, Theorem~\ref{thm: branched LO} settles
  the question of which $\Sigma_n(K)$ are LO except for $n = 3$.
\end{remark}

\subsection{Alternating knots}
\label{sec: determine h}

The results of the previous section can be used to compute \(h(K)\) for many knots. 

\begin{proposition}
  \label{prop:sigma h}
  Suppose a prime knot \(K\) is real representation small. If
  \(\Sigma_2(K)\) is not LO, then \(h(K) = -\frac{1}{2} \sigma(K)\).
\end{proposition}

\begin{proof} Suppose \([\rho] \in X^{0,\irr}_\SLR(K)\). By
Lemma~\ref{Lem:Branched Cover}, \(\rhohat_2 \in R_\PSLR(\Sigma_2(K))\)
is nontrivial and has \(e(\rhohat_2)=0\). Since \(\Sigma_2(K)\) is
prime and not LO, this contradicts Theorem~\ref{Thm:SLR order}. Hence
\(X^{0,\irr}_\SLR(K) = \emptyset\), and so \(h^0_\SLR(K) = 0 \). Since
\(0 \not \in D_K\) as noted in Section~\ref{sec: alex and the reds},
by Theorem~\ref{thm: hSL+hSU} and Corollary~\ref{Cor: hSU=sig} we
have
\(h(K) = h_\SU^0(K) = -\frac{1}{2}\sigma_K(-1) =
-\frac{1}{2}\sigma(K)\) as claimed.
\end{proof}

If the L-space conjecture holds, the hypothesis that \(\Sigma_2(K)\)
is not LO is equivalent to it being an L-space, a condition which has
been intensively studied. In particular, the Ozsv{\'a}th-Szab{\'o}
spectral sequence \cite{OzsvathSzabo2005} implies that if \(K\) has
thin Khovanov homology, then \(\Sigma_2(K)\) is an L-space. The
equivalence of conditions (i) and (ii) of the L-space conjecture have
been checked for many such knots, including alternating knots
\cite{BGW2013}. Hence we have:

\begin{corollary}
  \label{cor: h for alt}
  If a prime knot \(K\) is real representation small and alternating, then
  \( h(K) = -\frac{1}{2} \sigma(K)\).
\end{corollary}

\subsection{Dehn surgeries}
\label{sec: dehn}

Let \(\Ka\) be the Dehn filling of \(\EK\) along slope
\(\alpha = p \mu + q \lambda\).  In this subsection, we give criteria
for showing \(\Ka\) is LO.  The basic strategy follows
\cite{CullerDunfield2018}, using the translation extension locus which
is built from representations $\pi_1(\EK) \to \tSLR$ as discussed in
Section~\ref{sec: trans ex locus}.  Throughout, we freely use the
notation from Section~\ref{sec: ext Lin}, and for brevity we set
$\Etil(K) := \XA_\tSLR^\elll(K)$ and
$\Vtil(K) := \XA_\tSLR^\elll(\partial \EK)$.  Recall the translation
extension locus of $K$ is the image $\OTEL{K}$ inside $\Vtil(K)$,
where $i^*$ is induced by the inclusion
$i \maps \partial \EK \to \EK$.  (Provided $K$ is real representation
small, $\OTEL{K}$ is precisely the locus
$\mathit{EL}_\Gtil(\EK) \subset H^1(\partial \EK; \R)$ from
\cite{CullerDunfield2018} as discussed earlier in Section~\ref{sec:
  trans ex locus}.)  The new ingredient compared to
\cite{CullerDunfield2018} is using the tools in this paper to show
that $\OTEL{K}$ contains ``long arcs'' which give rise to explicit
ranges of LO fillings; in contrast, the main theorems of
\cite{CullerDunfield2018} only give orders on arbitrarily small
intervals of fillings. (For context, note Figure~\ref{fig: K12a380}
shows a case where a small interval is the best one can do.)

Define \(\Ltilde_{\alpha}\) in $\Vtil(K)$ to be those $\chitil$ where
$\rhotil(\alpha) = 1$ for some $\rhotil$ with $[\rhotil] = \chitil$;
in our $(\mu^*, \lambda^*)$-coordinates from Section~\ref{sec: trans
  ex locus}, this is the line of slope \(-\alpha\) passing though the
origin.  Let \(\Ltilde_{\alpha}^\circ \subset \Ltilde_{\alpha}\) be
the subset of nonparabolic points.  In this language, Lemma~4.4 of
\cite{CullerDunfield2018} becomes:

\begin{lemma}
\label{Lem:Dehn Fill}
If \(\Ltilde_{\alpha}^\circ\) intersects \(\OTEL{K}\), then \(\Ka\)
admits a nontrivial representation
\(\rho\maps \pi_1\big(\Ka\big) \to \tSLR\). Consequently, if \(\Ka\) is
prime then it is LO.
\end{lemma}

\begin{proof}
A nonparabolic $\chitil \in \Vtil(K)$ is in $\Ltil_\alpha$ if
and only if $\rhotil(\alpha) = 1$ for all
$\rhotil \in R_\tSLR(\partial \EK)$ representing $\chitil$. In
particular, if \(\Ltilde_{\alpha}^\circ\) intersects \(\OTEL{K}\), we
get a nontrivial representation $\pi_1\big(\Ka\big) \to \tSLR$; the
last conclusion now comes from Theorem~\ref{Thm:SLR order}.
\end{proof}

\begin{figure}
  \centering
  \begin{tikzoverlayabs}[width=0.9\textwidth]{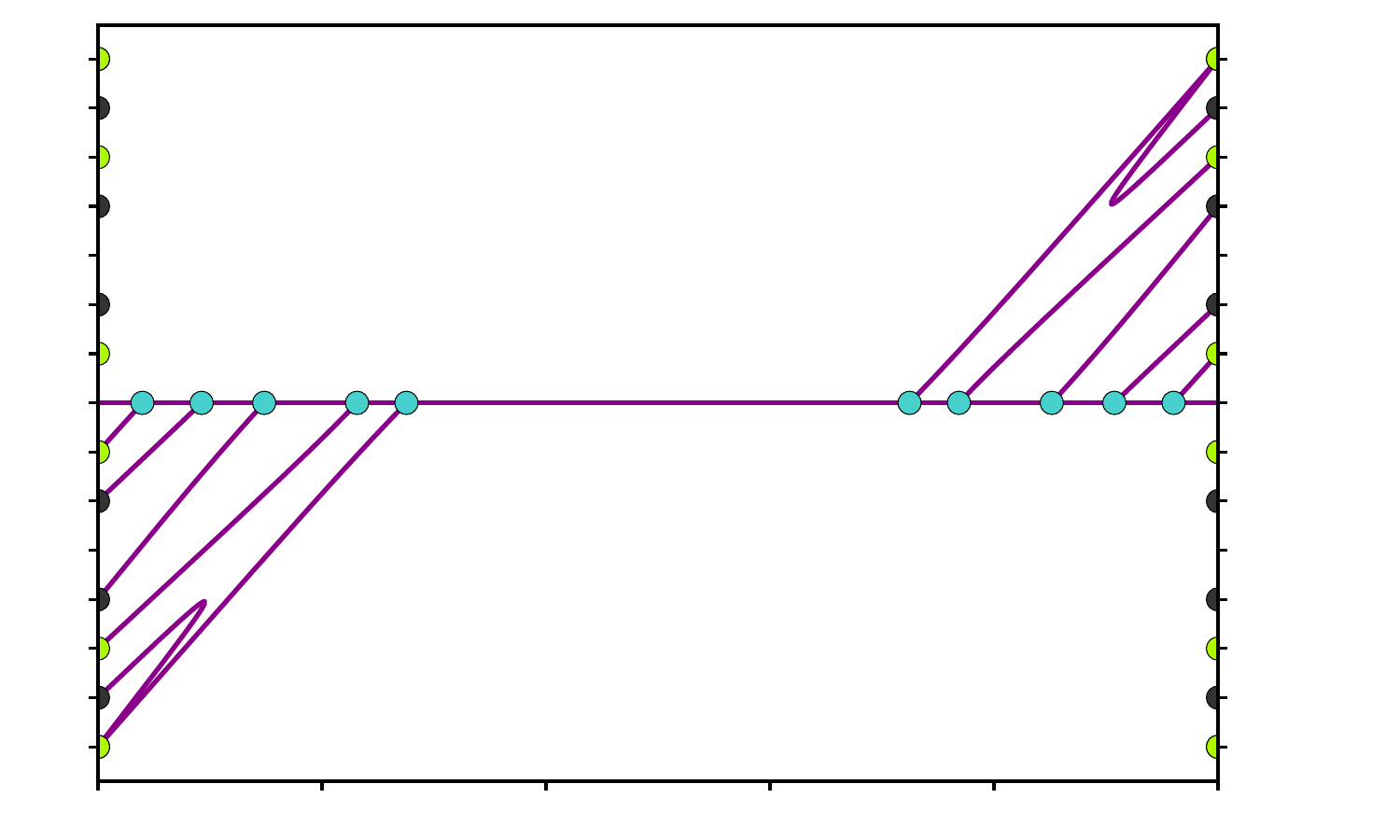}
  \begin{scope}[font=\footnotesize]
    \draw (0.086000, 0.952000) node[below right] {\small $s841$};
    \draw (0.070000, 0.049167) node[below] {$0.0$};
    \draw (0.230000, 0.049167) node[below] {$0.2$};
    \draw (0.390000, 0.049167) node[below] {$0.4$};
    \draw (0.550000, 0.049167) node[below] {$0.6$};
    \draw (0.710000, 0.049167) node[below] {$0.8$};
    \draw (0.870000, 0.049167) node[below] {$1.0$};
    \draw (0.057500, 0.110946) node[left] {$-7$};
    \draw (0.057500, 0.169441) node[left] {$-6$};
    \draw (0.057500, 0.227935) node[left] {$-5$};
    \draw (0.057500, 0.286429) node[left] {$-4$};
    \draw (0.057500, 0.344924) node[left] {$-3$};
    \draw (0.057500, 0.403418) node[left] {$-2$};
    \draw (0.057500, 0.461913) node[left] {$-1$};
    \draw (0.057500, 0.520407) node[left] {$0$};
    \draw (0.057500, 0.578901) node[left] {$1$};
    \draw (0.057500, 0.637396) node[left] {$2$};
    \draw (0.057500, 0.695890) node[left] {$3$};
    \draw (0.057500, 0.754385) node[left] {$4$};
    \draw (0.057500, 0.812879) node[left] {$5$};
    \draw (0.057500, 0.871374) node[left] {$6$};
    \draw (0.057500, 0.929868) node[left] {$7$};

    \draw (0.30500, 0.53) node[above] {$A_{-7}$};
    \draw (0.25500000, 0.53) node[above] {$A_{-5}$};
    \draw (0.20250000, 0.53) node[above] {$A_{-4}$};
    \draw (0.15500000, 0.53) node[above] {$A_{-2}$};
    \draw (0.1100000, 0.53) node[above] {$A_{-1}$};

    \draw (0.64250000, 0.51) node[below] {$A_{7}$};
    \draw (0.68500000, 0.51) node[below] {$A_{5}$};
    \draw (0.7500000, 0.51) node[below] {$A_{4}$};
    \draw (0.79500000, 0.51) node[below] {$A_{2}$};
    \draw (0.83500000, 0.51) node[below] {$A_{1}$};

    \draw (0.152, 0.31) node {$P$};
    \draw (0.79, 0.73) node {$Q$};
    \node[right] at (0.89, 0.07) {\small $\mu^*$};
    \node[left] at (0.026500, 0.96)  {\small $\lambda^*$};
  \end{scope}
\end{tikzoverlayabs}
  
  \caption{For the knot exterior $s841$, the set $\OTEL{K}$ consists
    entirely of images of good arcs.  There are 10 of type $A$,
    labeled near the endpoint of their reducible representation, and
    two of type $B$, specifically $P$ of type $B_{p_1, p_2}$ for
    $p_1 = (0, -7)$ and $p_2 = (0, -6)$ and $Q$ of type $B_{q_1, q_2}$
    for $q_1 = (1, 7)$ and $q_2 = (1, 6)$. Picture adapted from
    \cite[Figure~7]{CullerDunfield2018}.}
  \label{fig: good arcs 1}
\end{figure}
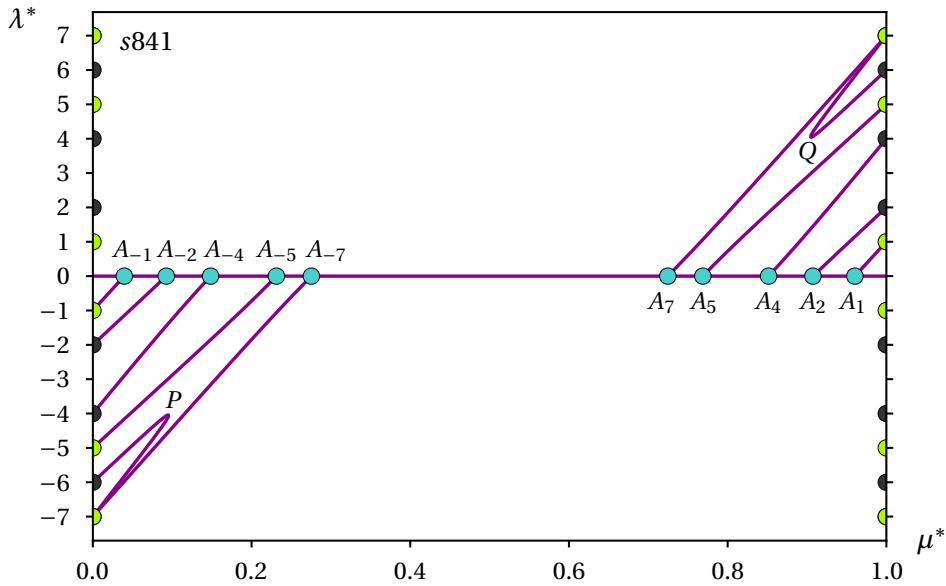

\begin{figure}
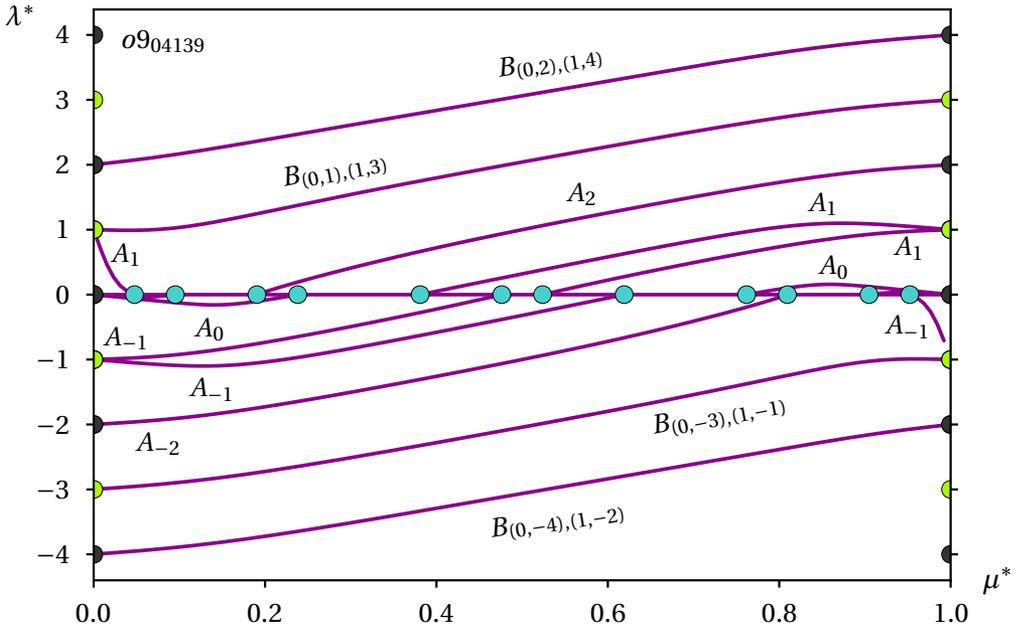

  \centering
  \input figures/o9_04139
  
  \caption{For the knot exterior $o9_{04139}$, the set $\OTEL{K}$
    again consists entirely of good arcs. The specific types are all
    labeled, with the caveat that there are actually four arcs of type
    $A_0$, only two of which are large enough to see clearly.  Picture
    adapted from \cite[Figure~8]{CullerDunfield2018}.}
  \label{fig: good arcs 2}
\end{figure}

\subsection{Good arcs}
\label{sec: good arc}

We say that \(\gamma:[0,1] \to \Etil(K)\) is a \emph{good arc} if the
endpoint \(\gamma(0)\) is parabolic, the other endpoint \(\gamma(1)\)
is either parabolic or reducible, the interior contains only elliptic
irreducibles, and finally
$i^*\big(\gamma(0)\big) \neq i^*\big(\gamma(1)\big)$.  A good arc
\(\gamma\) is of type \(A_k\) for \(k \in \Z\) if it starts at a
parabolic \(\chi\) with \(\lambda^*(\chi)=k\) and ends at a reducible.
It is of type \(B_{p_1,p_2}\) if its endpoints are parabolics
\(\chi_1, \chi_2\) with \(i^*(\chi_i)=p_i \in \Z^2\), where we have
identified the lattice of parabolic representations in \(\Vtil(K)\)
with \(\Z^2\). See Figures~\ref{fig: good arcs 1} and~\ref{fig: good
  arcs 2} for the images of some good arcs in the translation
extension locus in $\Vtil(K)$.  We show below in Lemmas~\ref{Lem:Type
  A} and~\ref{Lem:Type B} that many such good arcs lead to large
intervals where we can apply Lemma~\ref{Lem:Dehn Fill}.  Before
continuing, we note:

\begin{lemma}
  \label{lem: good don't touch}
  Suppose $\gamma$ is a good arc. Then the vertical lines where $\mu^* = k$ for
  $k \in \Z$ can meet $i^*(\gamma)$ only at its endpoints.
\end{lemma}

\begin{proof}
Since $\pi_1(S^3) = 1$, Lemma~\ref{Lem:Dehn Fill} gives the claim
for the line $\Ltil_\infty = \{\mu^* = 0\}$.  The claim for the other
vertical lines follows from the $D_\infty$ action discussed in
Example~\ref{ex: trans ext locus}.
\end{proof}
  
\begin{lemma}
  \label{Lem:Arcs}
  Suppose that \(K\) is real representation small.  If
  \(h(K) \neq 0\), then \(\Etil(K)\) contains either (a) an arc of
  type \(A_k\) or (b) an arc of type \(B_{p_1,p_2}\) where
  $p_1 = (0, y_1)$ and $p_2 = (1, y_2)$.  If \(\deg \eh(K) > 0 \),
  then \(\Etil(K)\) contains either (a) an arc of type \(A_k\) with
  \(k \neq 0\) or (b) an arc of type \(B_{p_1, p_2}\) with
  $p_1 = (0, y_1)$ with $y_1 \neq 0$.
\end{lemma}

\begin{proof}
If \(W\) is a path component of \(\cX^{[-2, 2]}(K)\), let
\(n_W^c = \sum_Z n_Z\), where \(Z\) runs over the set of components of
\(W \cap \cX^c(K)\). By Corollary~\ref{Cor:component sum},
\(n_W:=n_W^c\) is independent of \(c \in [-2,2]\). Since
\(h(K) \neq 0\), we can choose some component \(W\) with
\(n_W \neq 0 \). Since \(n_W^2 = n_W^{-2} \neq 0\), the component
\(W\) contains a parabolic \(\chi_1\) with \(\tr_\mu\chi_1 = 2\) and a
parabolic \(\chi_2\) with \(\tr_\mu\chi_2 = -2\).  Let \(\gamma\) be a
path in \(W\) between them. By restricting to a subarc if necessary,
we may assume that the interior of \(\gamma\) consists entirely of
nonparabolic points, retaining the property that $\tr_\mu\chi_1 = 2$
and $\tr_\mu\chi_2 = -2$.  If \(\gamma\) is contained in
\(X_\SLR(K)\), then \(\gamma\) lifts to an arc $\gammatil$ of type
\(B_{p_1, p_2}\).  Using the $D_\infty$ action, specifically the
index-two subgroup that preserves traces in $\SLR$, we can assume
$p_1 = (0, y_1)$.  Then $p_2 = (k, y_2)$ where $k$ must be odd.
Lemma~\ref{lem: good don't touch} forces $k = \pm 1$; if $k = -1$, act
by $\atil \in D_\infty$ to rotate $i^*(\gamma)$ about $(0,0)$ so that
$k = 1$ as desired.  If instead $\gamma$ is not contained in
\(X_\SLR(K)\), consider the first intersection point of \(\gamma\)
with \(\cX_{0}(K)\).  This initial segment of \(\gamma\) then lifts to
an arc of type \(A_k\). This proves the first claim.

For the second claim, let
\(\ntil_W = \sum_{Z} n_Z \big(t^{{\lambda}^*(\Ztil_+)}+
t^{{\lambda}^*(\Ztil_-)}\big),\) where the sum runs over components of
\(X^{2, \irr}_\SLR(K)\) contained in \(W\), as in
Definition~\ref{Def:extended invariant}.  For
$\chi \in X^{2, \irr}_\SLR(K)$, consider the two lifts $\chitil_{\pm}$
to $\XA_\tSLR(K)$ where $\mu^*(\chitil_\pm) = 0$.  These have
$\abs{ \lambda^*(\chitil_{+})} = \abs{\lambda^*(\chitil_{-})}$, and we
denote this common value $\abs{\lambda^*(\chitil)}$ without reference
to a particular lift.  As \(\sum_W \ntil_W = \eh(K)\), we can choose
some \(W\) with \(\deg \ntil_W \geq \deg \eh(K) > 0\).  Let
$\chi_1 \in W$ satisfy $\tr_\mu \chi_1 = 2$ and
$\abs{\lambda^*(\chitil_1)} = \deg \ntil_W$.

If \(n_W \neq 0\), we choose $\gamma$ in $W$ joining $\chi_1$ to
$\chi_2$ with $\tr_\lambda(\chi_2) = -2$.  Restricting to a subarc, we
may assume the interior of $\gamma$ consists entirely of nonparabolic
points, retaining $\tr_\mu(\chi_1) = 2$ and
$\abs{\lambda^*(\chitil_1)} = \deg \ntil_W$ and $p_1 \neq p_2$, but
perhaps not $\tr_\mu \chi_2 = -2$.  As before, either $\gamma$ is
entirely contained in \(X_\SLR^\irr(K)\), or there is a path from
\(\chi_1\) to a reducible. Lifting to \(\Etil(K)\) now proves the
second claim in the case $n_W \neq 0$.

Finally, if \(n_W=0\), note that \(\ntil_W|_{t=1} = 2 n_W = 0\), so
\(W\) must contain some other parabolic \(\chi_2\) with
$\tr_\mu(\chi_2) = 2$ and
\(\abs{\lambda^*(\chitil_2)} \neq \abs{\lambda^*(\chitil_1)}\). After
possibly changing the $\chi_i$, we can find an arc $\gamma$ in $W$
joining some $\chi_1$ and $\chi_2$ which satisfies the following:  1)
$\tr_\mu \chi_1 = 2$ and $\tr_\mu \chi_2 = \pm 2$;
2) $\abs{\lambda^*(\chitil_1)} = \deg \ntil_W$;  3) $p_1 \neq p_2$; and 4) 
 the interior of $\gamma$ consists of nonparabolic points.  The
 second claim now follows just as in the previous cases.
\end{proof}

The set of Dehn filling slopes \(S(K)\) of \(\EK\) is the projective space of
\(H_1(\partial \EK,\Z)\). It is  identified with
\(\Q \cup \{ \infty \} \) by our choice of basis
\(\pair{\mu, \lambda} \) for \(H_1(\partial \EK, \Z)\). We
define
\begin{align*}
  \Sirr(K) & = \setdefm{\big}{\alpha \in S(K)}%
             {\mbox{$\EK(\alpha)$ is irreducible}}, \\
  \SLO(K) & = \setdefm{\big}{\alpha \in S(K)}%
            {\mbox{$\pi_1\big(\EK(\alpha)\big)$ is left-orderable}},
            \mtext{and} \\
  \SSL(K) & = \setdefm{\big}{\alpha \in S(K)}%
            {\mbox{$\Ltilde^\circ_\alpha$ intersects $\OTEL{K}$}}. 
\end{align*} 
By Lemma~\ref{Lem:Dehn Fill}, we know that \(\SSL(K) \cap \Sirr(K)
\subset \SLO(K)\). As $\OTEL{K}$ contains the horizontal axis (as
the image of the reducible representations) and the $0$-surgery is
irreducible by \cite{Gabai1987}, the slope $0$ is in all three of
these sets.
  
Suppose \(\gamma\) is a good arc joining $\chi_1$ to $\chi_2$. Let
\(p_i=i^*(\chi_i) \in \Vtil(K)\). Suppose \(p_1, p_2 \neq (0,0)\) and at
most one $p_i$ is on the line $\mu^* = 0$.  In this case set
\[
  I^\circ(\gamma) = \setdef{\alpha \in S(K)}%
  {\mbox{$L_\alpha^\circ$ meets the interior of the line segment from $p_1$
      to $p_2$}}.
\]
When $p_i = (0, y_i)$ with $y_1$ and $y_2$ nonzero, we define
$I^\circ(\gamma) = S(K) \setminus \{\infty\}$ when the $y_i$ have
opposite signs and $I^\circ(\gamma)$ is empty otherwise.  Finally,
when either $p_1$ or $p_2$ is $(0,0)$, we take $I^\circ(\gamma)$ to be
empty.  We next show:

\begin{lemma}
  For a good arc $\gamma$, one has \( I^\circ(\gamma) \subset \SSL(K)\).
\end{lemma}

\begin{proof}
Whenever $\alpha$ is in $I^\circ(\gamma)$, the definitions give that
$p_1$ and $p_2$ are strictly separated from each other by the line
$L_\alpha$.  It follows that $L_\alpha$ meets $i^*(\gamma)$ at some
point in the latter's interior, which is not parabolic by the
definition of a good arc, and hence Lemma
\ref{Lem:Dehn Fill} applies.
\end{proof}

\begin{lemma} 
  \label{Lem:Type A}
  If \(\Etil(K)\) contains a good arc of type \(A_k\) for \(k \neq 0\),
  then either \((-\infty, |k|)\) or \((-|k|, \infty)\) is contained in
  \(\SSL(K)\).
\end{lemma}
Before proving this, we note Figure 5 of \cite{CullerDunfield2018}
shows that the condition $k \neq 0$ is crucial, as that example has
an $A_0$ arc which only yields $[-0.36, 3.6)$ in $\SSL(K)$.

\begin{proof}[Proof of Lemma~\ref{Lem:Type A}]
 Since \(\OTEL{K}\) is invariant under translation by
\((1,0)\), we can assume that the good arc $\gamma$ of type \(A_k\)
starts at the point \((0,k)\) and ends at \((x,0)\) for some
\(x \neq 0\) (since \(\pm 2 \not \in D_K\)).  Further, since
\(\OTEL{K}\) is invariant under the rotation $\atil \in D_\infty$
about the origin, we can assume \(k>0\). Suppose \(x>0\), so
\(I^\circ(\gamma) = (-\infty, 0)\). Next, consider the translate of
\(\gamma\) by \((-1, 0)\), which we denote by \(\gamma_{-1}\). Its
endpoints are \((-1,k)\) and \((x-1,0)\), so
\(I^\circ(\gamma_{-1}) = (0, k)\). The union of these two sets gives
the promised range since $0 \in \SSL(K)$ for all $K$. The argument
when \(x<0\) is symmetric.
\end{proof}

\begin{lemma}
  \label{Lem:Type B}
  Suppose \(\Etil(K)\) has a good arc of type \(B_{p_1,p_2}\) with
  \(p_i=(x_i,y_i)\).  Assume further that at least one $y_i$ is
  nonzero.  Then:
  \begin{enumerate}

  \item
    \label{item: flat}
    If \(y_1=y_2\neq 0\) then \(\Q \setminus
    \setdef{\frac{y_1}{n}}{n \in \Z} \subset \SSL(K).\)
    
  \item
    \label{item: near vert}
    If \(x_1=x_2\) then  \((-1,1)  \subset \SSL(K)\).

  \item
    \label{item: default}
    Otherwise, either \((-\infty, 1/2)\) or
    \((-1/2, \infty)\) is contained in \(\SSL(K)\).
\end{enumerate}
\end{lemma}

\begin{proof}

By Lemma~\ref{lem: good don't touch}, after translation we can assume
$i^*(\gamma)$ is contained in the strip $0 \leq \mu^* \leq 1$. In case
(\ref{item: flat}), we have \(p_1=(0,k)\) and \(p_2= (1,k)\) for some
\(k \neq 0\), and we can arrange that $k > 0$ using the $D_\infty$
action, resulting in the situation shown in Figure~\ref{fig: filling
  flat}. The translate of $\gamma$ by $(n, 0)$ then has \(I^\circ\)
equal to \(\big(\frac{k}{n}, \frac{k}{n+1}\big)\) or
\(\big(\frac{k}{n+1}, \frac{k}{n}\big)\) depending on the sign of $n$.
The union of these $I^\circ$ covers all of \(\Q\) except \(0\) and
those points of the form \(k/n\), \(n \in \Z\), completing the proof
of case (\ref{item: flat}).

\begin{figure}
  \begin{center}
    \begin{tikzpicture}[scale=1.8, TEL figure]

  \draw[coor grid] (-1.3, -0.3) grid (3.3, 3.3);
  \begin{scope}[color=axesgray]
    \draw[->] (-1.4, 0) -- (3.5,  0) node[right] {$\mu^*$};
    \draw[->] (0, -0.4) -- (0, 3.4) node[right=0.04]{$\lambda^*$};
  \end{scope}

  \draw[color=locus]
       (-1.3, 2.155) .. controls +(36.0:0.166) and +(107:0.22) .. (-1, 2);
  \foreach \x in {-1, 0, ..., 2}{
    \coordinate (A) at (\x, 2);
    \coordinate (B) at (\x + 1, 2);
    \draw[color=locus] (A) .. controls +(0.5, -0.35) and +(-0.2, 0.6) .. (B);
    \draw[parabolic] (A) circle;
  }
  \draw[color=locus]
       (3, 2) .. controls +(-32:0.143) and +(-172.5:0.1)  .. (3.3, 1.927);
  \draw[parabolic] (B) circle; 


  \begin{scope}[L line]
    \draw (0, 0) -- (1.3, 2.6) node[pos=1.05] {$\Ltil_{-y_1}$};
    \draw(0, 0) -- (-1.3, 2.6) node[pos=1.06] {$\Ltil_{y_1}$};
    \draw(0, 0) -- (2.6, 2.6) node[pos=1.05] {$\Ltil_{-y_1/2}$};
    \draw(0, 0) -- (3.4, 2.266666) node[pos=1.055] {$\Ltil_{-y_1/3}$};
  \end{scope}
  
  \node[below left=0 and -0.05] at (0, 2) {$p_1$};
  \node[below=0.1] at (1.05, 2) {$p_2$};
  \path (0, 2) .. controls +(0.5, -0.35) and +(-0.2, 0.6) .. (1, 2)
        node[pos=0.75, above left=0 and -0.2] {$i^*(\gamma)$};
\end{tikzpicture}
  \end{center}
  \caption{The situation of the proof of
    Lemma~\ref{Lem:Type B}(\ref{item: flat}), showing
    translates of $i^*(\gamma)$ for a good arc $\gamma$ of type
    $B_{p_1, p_2}$ with $y_1 = y_2 = 2$. The $\Ltil_\alpha$ drawn are
    among the few that are not guaranteed to meet \(\OTEL{K}\). }
  \label{fig: filling flat}
\end{figure}
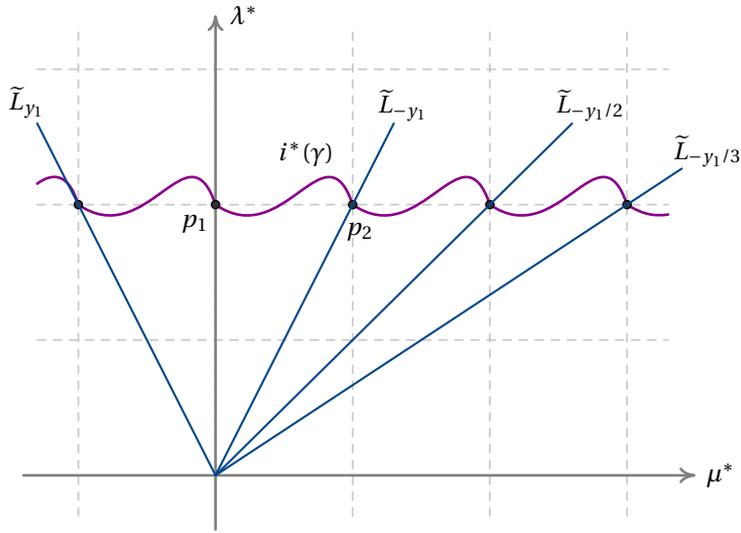

\begin{figure}
  \begin{center}
    \begin{tikzpicture}[scale=1.8, TEL figure]
  \draw[coor grid] (-0.1, -0.1) grid (4.3, 3.3);
  \begin{scope}[color=axesgray]
    \draw[->] (-0.2, 0) -- (4.5,  0) node[right] {$\mu^*$};
    \draw[->] (0, -0.2) -- (0, 3.4) node[right=0.04]{$\lambda^*$};
  \end{scope}

  \foreach \x in {0, ..., 4}{
    \coordinate (A) at (\x, 2);
    \coordinate (B) at (\x, 3);
    \draw[color=locus]
    (A) .. controls +(60:0.5) and +(-30:0.3) .. (B);
    \draw[parabolic] (A) circle;
    \draw[parabolic] (B) circle;
  }


  \coordinate (C) at (3.327, 2);
  \begin{scope}[L line]
    \filldraw (C) circle;
    \filldraw (4.186, 2.516) circle;
    \draw (0, 0) -- (3.3, 3.3) node[pos=1.045]{$\Ltil_{-1}$};
    \draw (0, 0) -- (31: 5.2) node[pos=1.04] {$\Ltil_{\alpha}$};
  \end{scope}
  
  \node[left=0.15] at (0, 2) {$y_1$};
  \node[left=0.15] at (0, 3) {$y_2$};
  \node[right] at (0.2, 2.55) {$i^*(\gamma)$};
  \node[below right] at (4, 2) {$(m, y_1)$};
  \node[above right] at (4, 3) {$(m, y_2)$};
  \node[below right=0 and -0.2] at (C) {$(x_3, y_1)$};
  
\end{tikzpicture}
  \end{center}
  \caption{The situation of the proof of Lemma~\ref{Lem:Type
      B}(\ref{item: near vert}), showing translates of
    $i^*(\gamma)$ for a good arc $\gamma$ of type $B_{p_1, p_2}$ with
    $x_1 = x_2$.}
  \label{fig: filling near vert}
\end{figure}
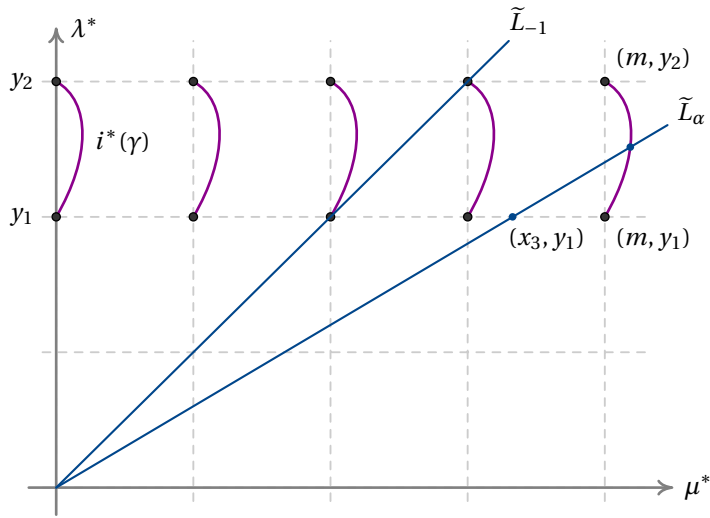

\begin{figure}
  \begin{center}
    \begin{tikzpicture}[scale=1.8, TEL figure]
  \draw[coor grid] (-2.3, -0.1) grid (3.3, 3.3);
  \begin{scope}[color=axesgray]
    \draw[->] (-2.3, 0) -- (3.5,  0) node[right] {$\mu^*$};
    \draw[->] (0, -0.2) -- (0, 3.4) node[right=0.04]{$\lambda^*$};
  \end{scope}

  \foreach \x in {-2, ..., 2}{
    \coordinate (A) at (\x, 1);
    \coordinate (B) at (\x + 1, 3);
    \draw[color=locus] (A) .. controls +(75:1) and +(-105:1) .. (B);
    \draw[parabolic] (A) circle;
    \draw[parabolic] (B) circle;
  }
  \draw[color=locus] (3, 1) .. controls +(75:0.3) and +(-120:0.3) .. (3.3, 1.707);
  \draw[color=locus] (-2.3, 2.295) .. controls +(57:0.25) and +(-105:0.25) .. (-2, 3);
  \draw[parabolic] (3, 1) circle;
  \draw[parabolic] (-2, 3) circle;


  \coordinate (C) at (2.355, 1);
  \begin{scope}[L line]
    \draw (0, 0) -- (3.3, 3.3) node[pos=1.045]{$\Ltil_{-n/2}$};
    \draw[name path=L alpha neg]
         (0, 0) -- (35: 3.9) node[pos=1.06] {$\Ltil_{\alpha < 0}$};
    \path[name path=bottom] (-2, 1) -- (2, 1);
    \filldraw[name intersections={of=L alpha neg and bottom, by=C}]
         (C) circle;
    \draw[name path=L alpha pos] (0, 0) -- (142: 3) node[pos=1.1] {$\Ltil_{\alpha > 0}$};
    \filldraw[name intersections={of=L alpha pos and bottom, by=D}]
         (D) circle;
    \path[name path=leftmost TEL] (-2, 1) .. controls +(75:1) and +(-105:1) .. (-1, 3);
    \path[name path=rightmost TEL] (2, 1) .. controls +(75:1) and +(-105:1) .. (3, 3);
    \filldraw[name intersections={of=L alpha pos and leftmost TEL, by=E}]
        (E) circle;
    \filldraw[name intersections={of=L alpha neg and rightmost TEL, by=F}]
        (F) circle;
  \end{scope}

  \node[right] at (0, 1) {$y_1$};
  \node[right] at (0, 3) {$y_2$};
  \node[above left=0 and -0.2] at (0.5, 2) {$i^*(\gamma)$};
  \node[below right] at (2, 1) {$(m, y_1)$};
  \node[below right] at (3, 3) {$(m + 1, y_2)$};
  \node[below right=0 and -0.3] at (C) {$(x_3, y_1)$};
  \node[above left] at (-2, 1) {$(\ell, y_1)$};
  \node[above left] at (-1, 3) {$(\ell + 1, y_2)$};
  \node[below left=0 and -0.2] at (D) {$(x_4, y_1)$};

\end{tikzpicture}
  \end{center}
  \caption{The situation of the proof of Lemma~\ref{Lem:Type B}%
    (\ref{item: default}) when $y_1 \geq 0$, showing that
    $\Ltilde^\circ_\alpha$ meets a translate of $i^*(\gamma)$ for any
    $\alpha \in (-n/2, \infty)$, where $n = y_2 - y_1$.}
\label{fig: filling default}
\end{figure}
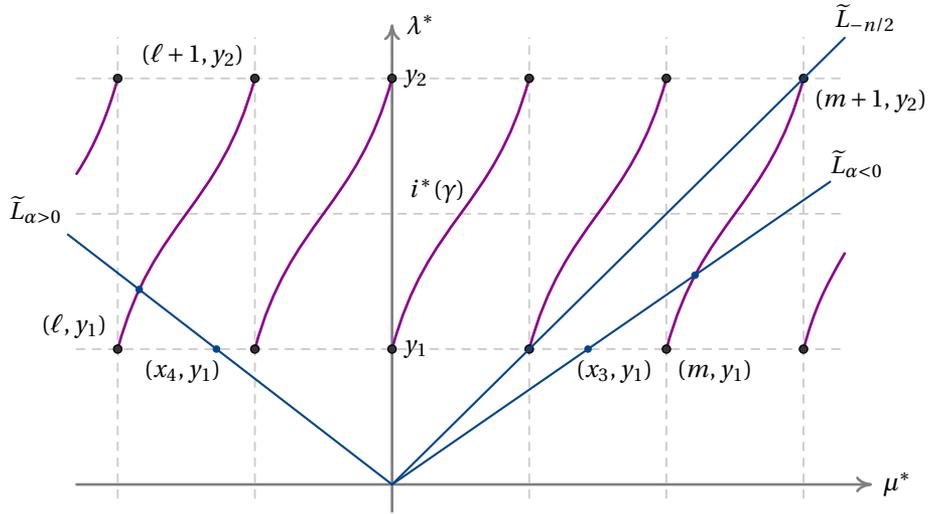

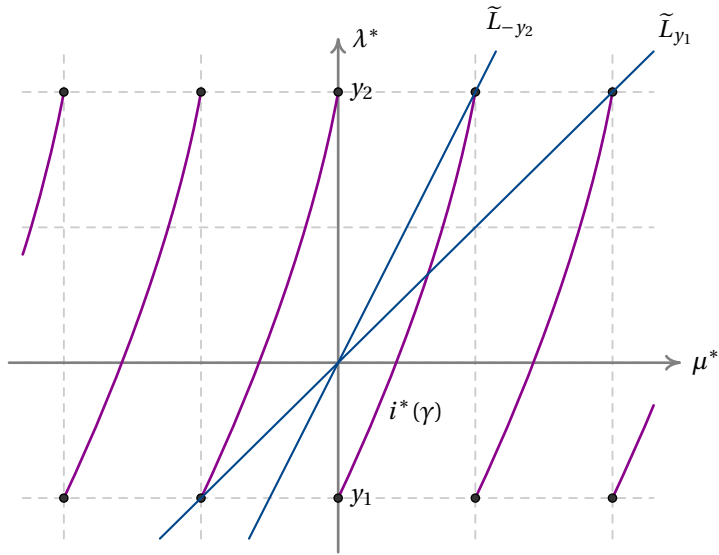
\begin{figure}
  \begin{center}
    \begin{tikzpicture}[scale=1.8, TEL figure]
  \draw[coor grid] (-2.3, -1.3) grid (2.3, 2.3);
  \begin{scope}[color=axesgray]
    \draw[->] (-2.4, 0) -- (2.5,  0) node[right] {$\mu^*$};
    \draw[->] (0, -1.4) -- (0, 2.4) node[right=0.04]{$\lambda^*$};
  \end{scope}

  \foreach \x in {-2, ..., 1}{
    \coordinate (A) at (\x, -1);
    \coordinate (B) at (\x + 1, 2);
    \draw[color=locus] (A) .. controls +(65:1.0) and +(-100:1.0) .. (B);
    \draw[parabolic] (A) circle;
    \draw[parabolic] (B) circle;
  }
  \draw[color=locus] (2, -1) .. controls +(65:0.2) and +(-112:0.2) .. (2.3, -0.315);
  \draw[color=locus] (-2.3, 0.8) .. controls +(72:0.35) and +(-100:0.35) .. (-2, 2);
  \draw[parabolic] (2, -1) circle;
  \draw[parabolic] (-2, 2) circle;
  

  \begin{scope}[L line]
    \draw (-0.65, -1.3) -- (1.15, 2.3) node[pos=1.06]{$\Ltil_{-y_2}$};
    \draw (-1.3, -1.3) -- (2.3, 2.3) node[pos=1.045]{$\Ltil_{y_1}$};
  \end{scope}

  \node[right] at (0, -1) {$y_1$};
  \node[right] at (0, 2) {$y_2$};

  \node[below right=0 and 0] at (0.28, -0.2) {$i^*(\gamma)$};
  
\end{tikzpicture}
  \end{center}
  \caption{The situation of the proof of Lemma~\ref{Lem:Type B}%
    (\ref{item: default}) when $y_1 < 0$, showing that
    $\Ltilde^\circ_\alpha$ meets a translate of $i^*(\gamma)$ for any
    $\alpha \in (-m, \infty)$, where $m = \max(-y_1, y_2)$.}
\label{fig: filling default alt}
\end{figure}

In case (\ref{item: near vert}), if $y_1$ and $y_2$ have opposite
signs, we are done as $I^\circ(\gamma)$ is all of $\SSL{K}$ except
$\infty$.  Otherwise, use the $D_\infty$ action and possibly
interchange $p_1$ and $p_2$ to assume $x_1 = x_2 = 0$ and
$0 \leq y_1 < y_2$, giving the situation in Figure~\ref{fig: filling
  near vert}.  If $0 < \abs{\alpha} < y_2 - y_1$, consider the point
where \(\Ltilde_{\alpha}\) intersects the line \(\lambda^*=y_1\).  For
concreteness, suppose this point is $(x_3, y_1)$ with $x_3 > 0$
(equivalently, assume $\alpha < 0$).  If $m$ is the smallest integer
with $x_3 < m$, then $\Ltilde^\circ_\alpha$ will have to meet the
translate of $\gamma$ joining $(m, y_1)$ to $(m, y_2)$.  As
$y_2 - y_1 \geq 1$, we have shown $(-1, 1) \subset \SSL(K)$.

Next, in case (\ref{item: default}), we can assume that
\(p_1 = (0, y_1)\) and \(p_2 = (1, y_2)\).  Now $y_1 \neq y_2$, and
let us assume $y_2 > y_1$ as the other case can be reduced to this one
by reflection across the vertical axis.  Moreover, we can also assume
$y_2 \geq 0$, since if not we apply the element of $D_\infty$ that
rotates around $(1/2, 0)$.

First suppose $y_1 \geq 0$ and set $n = y_2 - y_1$, giving the
situation of Figure~\ref{fig: filling default}.  If
$-n/2 < \alpha < 0$, then $\Ltil_\alpha$ has positive slope and
consider the point $(x_3, y_1)$ where $\Ltil_\alpha$ meets the line
$\lambda^*=y_1$.  If $m$ is the smallest integer where $x_3 < m$, then
$\Ltil^\circ_\alpha$ will meet the translate of $i^*(\gamma)$ joining
$(m, y_1)$ to $(m+1, y_2)$.  Hence $(-n/2, 0) \subset \SSL(K)$. If
$\alpha > 0$, let $(x_4, y_1)$ again be the point where $\Ltil_\alpha$
meets the line $\lambda^* = y_1$.  If $\ell$ is the largest integer where
$\ell < x_4$, then $\Ltil^{\circ}_\alpha$ will meet the translate of
$i^*(\gamma)$ joining $(\ell, y_1)$ to $(\ell + 1, y_2)$.  In
particular, we have $(-n/2, \infty) \subset \SSL(K)$ as needed.

The final subcase is when $y_1 < 0$, which is shown in
Figure~\ref{fig: filling default alt}.  There, you can see that
$\Ltilde^\circ_\alpha$ meets a translate of $i^*(\gamma)$ for any
$\alpha \in (-m, \infty)$, where $m = \max(-y_1 , y_2) \geq 1$.  This
completes the proof of case (\ref{item: default}), and hence the
lemma.
\end{proof}

\begin{theorem}
  \label{thm: extended Lin to LO}
  If \(\deg \eh(K) >0\), then \(A \cap \Sirr(K) \subset \SLO(K)\),
  where \(A\) is one of the following sets:
  \((\infty, 1/2), (-1/2, \infty), (-1,1)\), or
  \(\Q \setminus \setdef{\frac{1}{n}}{n \in \Z}\).
\end{theorem}

\begin{proof}
Just combine Lemmas~\ref{Lem:Dehn Fill}, \ref{Lem:Arcs},
\ref{Lem:Type A}, and \ref{Lem:Type B}.
\end{proof}

These techniques apply to many 2-bridge knots:
\begin{theorem}
  \label{thm: 2-bridge LO}
  If \(K\) is a 2-bridge knot with \(\sigma(K) \neq 0\), then
  \(\Etil(K)\) contains an arc of type \(A_k\) for some \(k \neq
  0\). Moreover, either $(-\infty, 1)$ or $(-1, \infty)$ is in
  $\SLO(K)$.
\end{theorem}
The positive trefoil knot $K$ from Section \ref{Ex:Trefoil} is a
\2-bridge knot and has $\sigma(K) = -2$.  There, the locus $\OTEL{K}$
consists of the orbit of a single $A_1$ arc under the $D_\infty$
action.  This shows that the final conclusion of Theorem~\ref{thm:
  2-bridge LO} cannot be strengthened to
$(-\infty, \infty) \subset \SLO(K)$ using the present methods. Of
course, the $-1$ Dehn surgery on $K$ is the Poincar\'e homology
sphere, which is not LO, but the same pattern holds for the knot $5_2$
where every Dehn surgery is expected to be LO.  

\begin{proof}[Proof of Theorem~\ref{thm: 2-bridge LO}]
As 2-bridge knots are small \cite{HatcherThurston1985}, 
$K$ is real representation small.  As $K$ is alternating, we have
\(h(K) = -\frac{1}{2} \sigma(K) \neq 0\) by Corollary~\ref{cor: h for
  alt}.  By Lemma~\ref{Lem:Arcs}, there is a good arc $\gamma$ of
either type $A_k$ or of type $B_{p_1, p_2}$ where $p_1 = (0, y_1)$ and
$p_2 = (1, y_2)$.  In the latter case, we would have
\(X^{0, \irr}_\SLR(K) \neq \emptyset\), but this is impossible since
\(\Sigma_2(K)\) is a lens space and so not LO (compare with the proof
of Proposition~\ref{prop:sigma h}). So we have a good arc of type
$A_k$.  By \cite[Theorem 2]{Riley1972},  if
\(\rho \maps \pi_1(\EK) \to \SLR\) is a parabolic representation then
\(\tr \rho(\lambda) = - 2\). This says that any arc in the translation
extension locus $\OTEL{K}$ ending on a parabolic must do so at
an odd height.  Thus $k \neq 0$ and Lemma~\ref{Lem:Type A} applies.
Finally, every Dehn surgery on $K$ is prime by \cite{Delman1995,
  Naimi1997}, allowing us to pass from $\SSL(K)$ to $\SLO(K)$ unimpeded.
\end{proof}

\begin{remark}
  \label{rem: extend LO}
  It is natural to ask whether Theorem~\ref{thm: 2-bridge LO} can be
  extended to other small Montesinos knots or to other small
  alternating knots.  In both cases one has
  $h(K) = -\frac{1}{2} \sigma(K)$ by Proposition~\ref{prop:
    Montesinos} and Corollary~\ref{cor: h for alt}.  Unlike for
  2-bridge knots, in these broader settings, parabolics with
  \(\tr_\lambda = 2\) and hence even height can occur.  For Montesinos
  knots, the $(-2, 3, 7)$ pretzel has several parabolics at even
  heights (see Figure~\ref{fig: pretzel constraints}).  For small
  alternating knots, Goerner's data \cite{PtolemyTables} includes nine
  such examples: $9_{34} = K9a28$, $10_{100} = K10a104$,
  $10_{108} = K10a119$, $10_{103} = K10a105$, $10_{104} = K10a118$,
  $10_{106} = K10a95$, $10_{102} = K10a97$, $10_{101} = K10a45$, and
  $11_{337} = K11a330$.  However, we have not seen any parabolics of
  height $0$ for alternating or Montesinos knots, so it is possible
  that the same conclusion as in Theorem~\ref{thm: 2-bridge LO} holds
  in these cases.
\end{remark}

\subsection{Knots with lens space surgeries}
\label{sec: lens space sur}

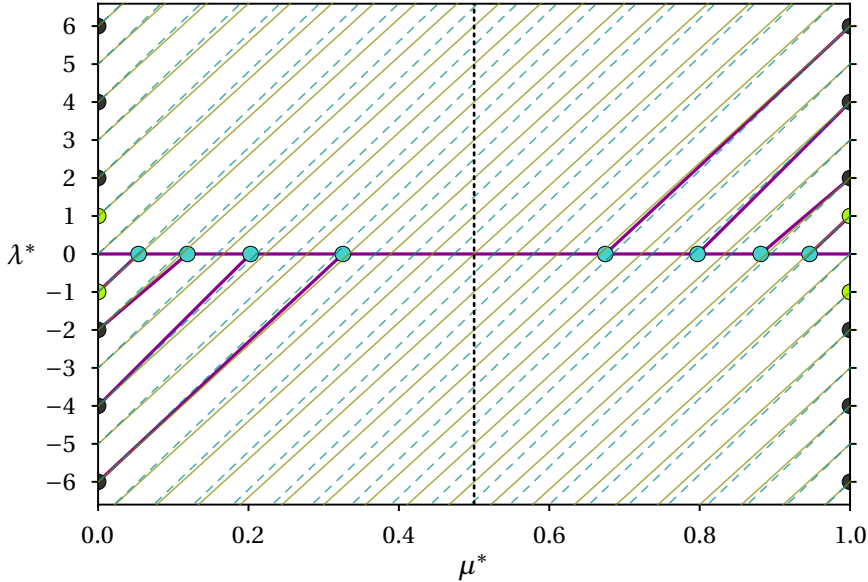
\begin{figure}
  \begin{center}
    \definecolor{Aslope}{HTML}{8B8B00}
\definecolor{Bslope}{HTML}{008B8B}

\begin{tikzoverlayabs}[width=351pt]{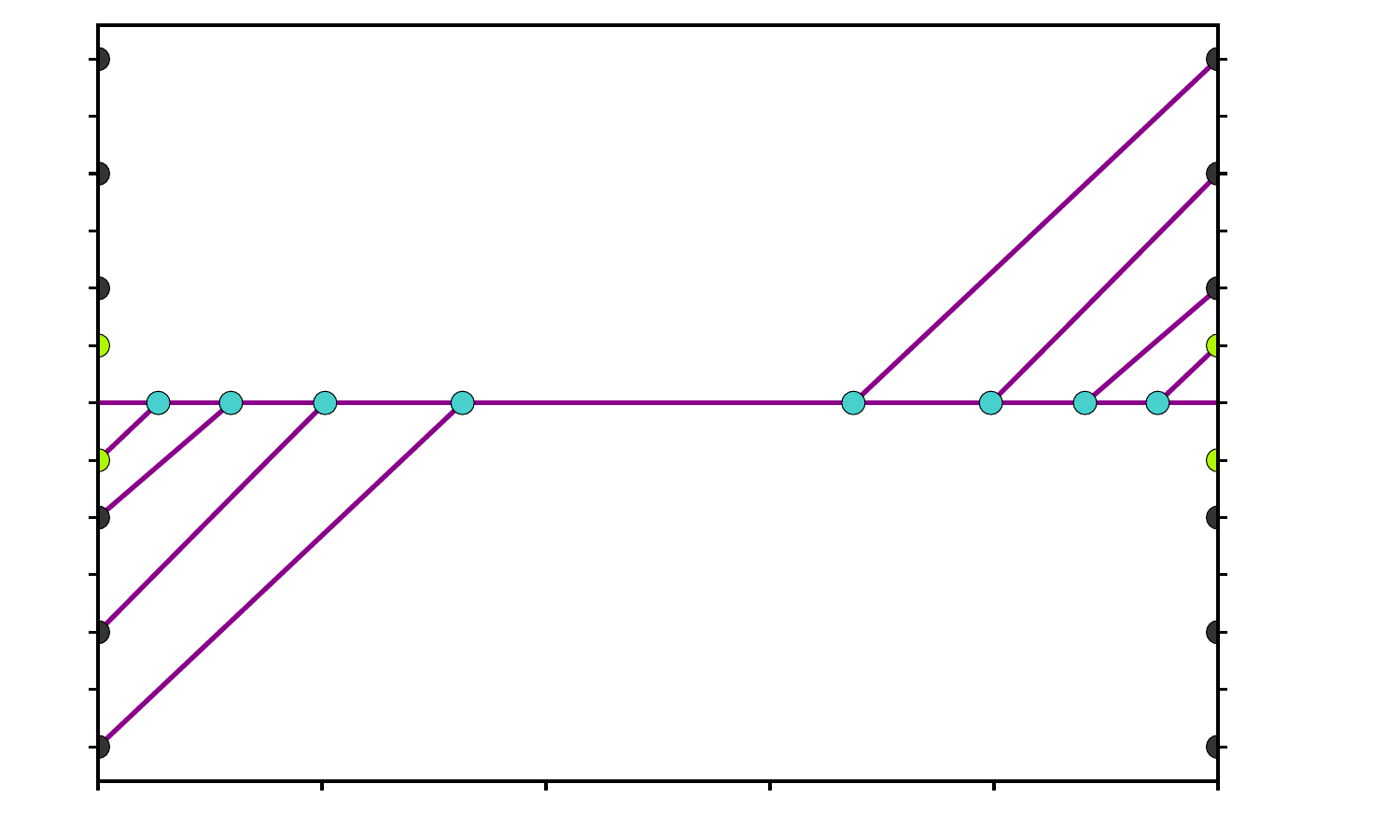}
  \begin{scope}[nmdstd]
    \begin{scope}[xshift=24.8, xscale=0.799,
      yshift=109.76, yscale=0.06812,
      line width=0.6pt]
       \path[clip] (0, -6.6) rectangle (1, 6.6);
       \begin{scope}[color=Aslope, opacity=0.7]
         \foreach \n in {-7, -6, ..., 25}{
           \draw (\n/18, 0) -- +(1, 18);
           \draw (\n/18, 0) -- +(-1, -18);
         }
       \end{scope}
       \begin{scope}[color=Bslope, opacity=0.7, dashed]
         \foreach \n in {-7, -6, ..., 25}{
           \draw (\n/19, 0) -- +(1, 19);
           \draw (\n/19, 0) -- +(-1, -19);
         }
       \end{scope}
       \draw[dotted, line width=1pt] (0.5, -7) -- (0.5, 7);
    \end{scope}
    
    \draw (0.070000, 0.049167) node[below] {$0.0$};
    \draw (0.230000, 0.049167) node[below] {$0.2$};
    \draw (0.390000, 0.049167) node[below] {$0.4$};
    \draw (0.550000, 0.049167) node[below] {$0.6$};
    \draw (0.710000, 0.049167) node[below] {$0.8$};
    \draw (0.870000, 0.049167) node[below] {$1.0$};
    \draw (0.057500, 0.110940) node[left] {$-6$};
    \draw (0.057500, 0.179173) node[left] {$-5$};
    \draw (0.057500, 0.247406) node[left] {$-4$};
    \draw (0.057500, 0.315639) node[left] {$-3$};
    \draw (0.057500, 0.383872) node[left] {$-2$};
    \draw (0.057500, 0.452105) node[left] {$-1$};
    \draw (0.057500, 0.520338) node[left] {$0$};
    \draw (0.057500, 0.588571) node[left] {$1$};
    \draw (0.057500, 0.656804) node[left] {$2$};
    \draw (0.057500, 0.725037) node[left] {$3$};
    \draw (0.057500, 0.793270) node[left] {$4$};
    \draw (0.057500, 0.861503) node[left] {$5$};
    \draw (0.057500, 0.929736) node[left] {$6$};
    \node[below] at (0.47, 0.0) {\small $\mu^*$};
    \node[left] at (0.017500, 0.520338)  {\small $\lambda^*$};
  \end{scope}

\end{tikzoverlayabs}
  
  \end{center}
  \caption{%
    This figure shows $\OTEL{K}$ for the $(-2, 3, 7)$ pretzel knot,
    whose exterior is small.  Surgery along the slopes $-18$ and $-19$
    both yield lens spaces, forcing $\OTEL{K}$ to avoid two families
    of lines per Lemma~\ref{lem:lens lines}, those of slope 18 (solid
    lines) and slope 19 (dashed lines).  The roots of $\Delta_K$ on
    $S^1$ are simple and are indicated by the blue dots.  Because of
    the locations of the roots, these lines rule out any arcs joining
    two distinct reducibles.  From Proposition~\ref{prop: Montesinos},
    we know $\OTEL{K}$ is disjoint from the vertical line
    $\mu^* = 1/2$.  Consequently, an arc starting at a particular
    reducible has only one or two possible parabolic ending points,
    basically dictating the picture shown (see the proof of
    Theorem~\ref{thm:lens arcs} for details). In particular, there
    must be an arc from the reducible at $\approx (0.6744, 0)$ to the
    parabolic at $(1, 6)$, and so $(-6, \infty) \subset \SSL(K)$.
    Here, $h(K) = -\frac{1}{2}\sigma(K) = -4$, and the signature
    function jumps by $+2$ at the first four roots of $\Delta_K$ and
    then by $-2$ at the remaining four. Picture adapted from
    \cite[Figure~3]{CullerDunfield2018}.}
  \label{fig: pretzel constraints}
\end{figure}

In this last subsection, we give strong constraints on $\OTEL{K}$ when
$K$ has a lens space surgery or more generally a Dehn surgery $\Ka$
with few $\PSLR$ representations.  We begin with an easy lemma that
nonetheless explains a great deal of the striking structure of
$\OTEL{K}$ for the knots shown in Figures 3 and 4 of
\cite{CullerDunfield2018}.  (The first of those figures is incorporated
here as Figure~\ref{fig: pretzel constraints}.)  To state it, for a
slope $\alpha = p \mu + q \lambda$ and an $n \in \Z$, we define
$\Ltil_{\alpha, n}$ to be the line in $\Vtil(K)$ of slope $-\alpha$ that
meets the horizontal axis at $(n/p, 0)$.  We also take
$\Ltil^\circ_{\alpha, n}$ to be its subset of nonparabolic points
where $\lambda^* \neq 0$, i.e.~remove $\Z^2$ and $(n/p, 0)$ from
$\Ltil_{\alpha, n}$.

\begin{lemma}
\label{lem:lens lines}
  Let $K$ be any knot where $X^\irr_\PSLR\big(\Ka\big)$ is empty. Then
  $\OTEL{K}$ is disjoint from
  $\bigcup \setdef{\Ltil^\circ_{\alpha, n}}{n \in \Z}$.
\end{lemma}
Figure~\ref{fig: pretzel constraints} shows just how constraining
Lemma~\ref{lem:lens lines} is when there are multiple lens space
surgeries.

\begin{proof}
The line $\Ltil_{\alpha, n}$ is exactly the locus where
$\trans_\alpha(\chitil) = n$, where $\trans_\alpha$ is the translation
number function of Corollary~\ref{cor: X of Gtil}.
Hence if $\chitil$ is in both $\Ltil^\circ_{\alpha, n}$ and
$\OTEL{K}$, then we get an irreducible elliptic representation
$\rhotil \maps \pi_1(\EK) \to \tSLR$ where $\rhotil(\alpha)$ is in the
center of $\tSLR$.  Quotienting $\tSLR$ by its center gives an
irreducible representation $\rho \maps \pi_1\big(\Ka\big) \to \PSLR$,
which is a contradiction.
\end{proof}

\begin{lemma} 
\label{Lem:no reducible}
Let \(W\) be a component of \(X^{[-2, 2],\irr}_\SLR(K)\), and suppose
that \(c, c' \in [-2,2] \setminus D_K\). If \(n_W^c \neq n_W^{c'}\),
then there is a path in $X_\SLR(K)$ which starts in either \(W^c\) or
\(W^{c'}\) and is contained in $W^{[c, c']}$ except for its final
endpoint, which is a reducible character $\chi$ with
$\tr_\mu(\chi) \in D_K \cap [c, c']$.
\end{lemma}

\begin{proof}
Since we can compute $n^c_W$ and $n^{c'}_W$ by summing over the
connected components of $Z$ of \(W^{[c,c'], \irr}(K)\), let \(Z\) be
such a component with \(n^c_{Z} \neq n^{c'}_{Z}\).  Let $\Zbar$ be the
component of $\cX^{[c, c']}(K)$ that contains $Z$.  As
$n^c_\Zbar = n^{c'}_\Zbar$, it follows that $\Zbar \neq Z$ and so
$\Zbar$ meets $\cX_0^{[c, c']}(K)$.  As \(n^c_{Z} \neq n^{c'}_{Z}\),
at least one of them must be nonzero, say $n^c_{Z} \neq 0$; in particular,
$Z^c$ is nonempty.  Thus there is a path in $\Zbar$ starting in $Z^c$
and ending in $\Zbar \cap \cX_0^{[c, c']}(K)$.  Taking an initial
segment of this path and using the map $\cX(K) \to X(K)$ gives the
desired path by Lemma~\ref{lemma: deform red}.
\end{proof}

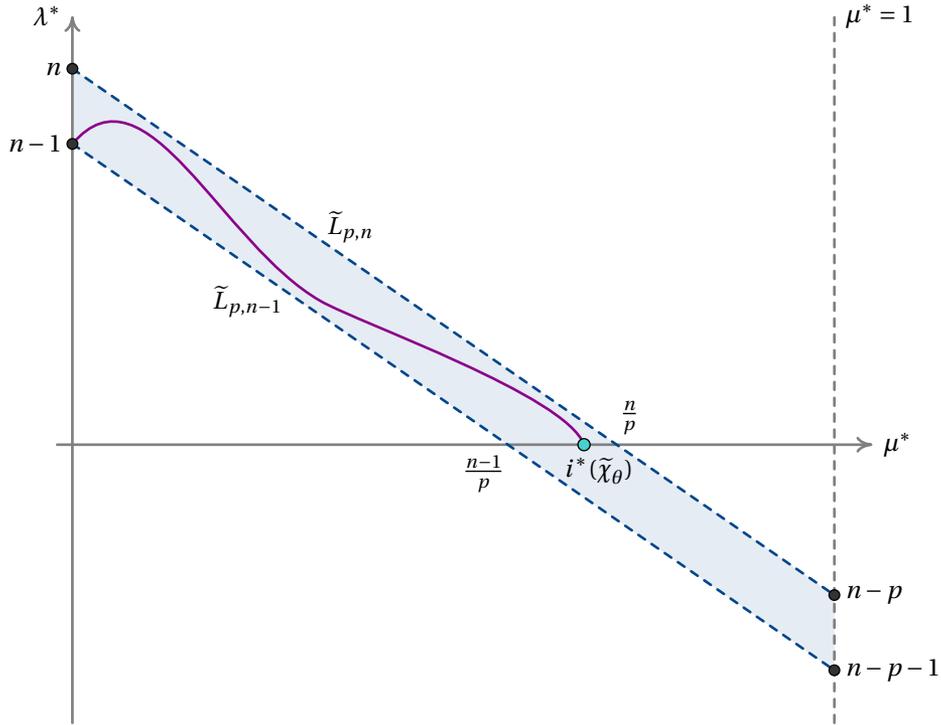
\begin{figure}
\begin{center}
  \begin{tikzpicture}[TEL figure]
  \def\n{5}
  \def\p{7}
  \def\s{0.7}
  \coordinate (X) at (10, 0);
  \coordinate (Y) at (0, 1);
  
  \coordinate (A) at (${\n - 1}*(Y)$);
  \coordinate (B) at ($\n*(Y)$);  
  \coordinate (C) at ($(X) + {\n - \p}*(Y)$);
  \coordinate (D) at ($(X) + {\n - \p - 1}*(Y)$);
  \coordinate (R) at (${\n - \p - 1 - \s}*(Y)$);
  \coordinate (S) at (${\n + \s}*(Y)$);
  \coordinate (E) at (${(\n - 1)/\p)}*(X)$);
  \coordinate (F) at (${(\n/\p)}*(X)$);
  \coordinate (G) at ($(E)!0.7!(F)$);
  \coordinate (H) at ($(A)!0.5!(G)+(-0.1, -0.09)$);

  \fill[color=Lline!10] (A) -- (B) -- (C) -- (D) -- cycle;

  \draw[color=locus]
      (A) .. controls +(50:1.5) and +(150:1.5) .. (H)
          .. controls +(-30:0.6) and +(110:0.7) .. (G);
  
  \begin{scope}[color=axesgray]
    \draw[->] (-0.2, 0) -- ($(X) + (0.5,  0)$) node[right] {$\mu^*$};
    \draw[->] (R) -- (S) node[left]{$\lambda^*$};
    \draw[dashed] ($(R) + (X)$) -- ($(S) + (X)$) node[right]{$\mu^*=1$};
  \end{scope}

  \begin{scope}[color=Lline, dashed]
    \draw (A) -- (D) node[pos=0.3, left=0.1] {$\Ltil_{p, n-1}$};
    \draw (B) -- (C) node[pos=0.3, right=0.2] {$\Ltil_{p, n}$};
  \end{scope}

  \foreach \P in {A, B, C, D}{
    \draw[parabolic] (\P) circle[radius=0.07];
  }

  \node[below left=0 and -0.1] at (E) {$\frac{n - 1}{p}$};
  \node[above right=0 and -0.1] at (F) {$\frac{n}{p}$};
  \draw[line width=0.5pt, fill=simplealex, radius=0.08] (G) circle;
  \node[left] at (B) {$n$}; 
  \node[left] at (A) {$n-1$};
  \node[right] at (C) {$n-p$}; 
  \node[right] at (D) {$n-p-1$};
  \node[below right=0.03 and -0.4] at (G) {$i^*(\chitil_\theta)$};
\end{tikzpicture}
\end{center}
\caption{A good arc ending at \(\chitil_\theta\) must have image lying
  in the shaded region.}
  \label{fig:Lens Surgery}
\end{figure}

\begin{theorem} 
\label{thm:lens arcs}
Suppose a prime knot $K$ is real representation small and that
\(\EK(p \mu + \lambda) \) is a lens space for some \(p>0\).  Suppose
that for some \(n=1,2,\ldots, p\) there is a unique
\(x_0 \in \big[\frac{n-1}{p}, \frac{n}{p}\big]\) with
\(\Delta_K(e^{2\pi i x_0})=0\), that the corresponding root of
\(\Delta_K\) is simple, and that
\(x_0 \neq \frac{n-1}{p},\frac{n}{p}\). Then \(\Etil(K)\) contains an
arc $\gamma$ of type \(A_k\), where \(k\) is one of \(n-1, n, n-1-p\),
or \(n-p\), and where the other endpoint of $i^*(\gamma)$ is
$(x_0, 0)$.  If $n \notin \{1, p\}$, then $(-\infty, 1) \subset \SSL(K)$.
\end{theorem}
A similar statement holds when \(\EK(-p \mu + \lambda)\) is a lens
space by taking the mirror.  Also, the hypothesis of a lens space
filling in Theorem~\ref{thm:lens arcs} can be replaced with the
requirement that $X^\irr_\PSLR\big(\EK(p \mu + \lambda)\big)$ is
empty.

\begin{remark}
  \label{rem: lens space conj}
  We conjecture that the hypothesis on the Alexander polynomial in
  Theorem~\ref{thm:lens arcs} always holds when
  \(\EK(p \mu + \lambda) \) is a lens space. By Greene's lens space
  realization theorem \cite{Greene2013}, the set of possible
  $\Delta_K$'s coincides with the set of Alexander polynomials of
  Berge knots.  For fixed $p$, there are only finitely many
  possibilities, and the hypothesis on their roots holds for
  $p \leq 1{,}000$ (which is some $5{,}265$ distinct $\Delta_K$).  In
  fact, with the exception of $p = 5$, there is always a root leading
  to the conclusion that $(-\infty, 1) \subset \SSL(K)$, indeed
  usually 10s or even 100s of such roots.  Here, it is unknown whether
  having a lens space surgery means $K$ is real representation small,
  and there are such $K$ that are not small \cite{Baker2005}.

\end{remark}

\begin{proof}[Proof of Theorem~\ref{thm:lens arcs}]
Let \(\omega=\frac{(n-1)\pi}{p}\) and \(\omega'=\frac{n\pi}{p}\), as
well as \(c = 2 \cos \omega\) and \(c' = 2 \cos \omega'\).  By
hypothesis, \(c, c' \not \in D_K\), and there is a unique root
$e^{2 i \theta}$ of $\Delta_K$ with $\theta \in [\omega, \omega']$.
Let \(\chitil_\theta\) be the point in \( X^\red_\tSLR(K)\) with
\(i^*(\chitil_\theta) = (\theta/\pi, 0)\).  Now if \(\gammatil\) is a
path in \(\Etil(K)\) starting at $\chitil_\theta$ whose interior
consists of nonparabolic characters, then \(i^*(\gammatil)\) is
contained in the shaded region \(R\) shown in Figure~\ref{fig:Lens
  Surgery} and meets $\partial R$ only if the other endpoint of
\(\gammatil\) is a parabolic.  In the latter case, the parabolic must
have height \(n-1, n, n-p\) or \(n-p-1\), and we can apply
Lemma~\ref{Lem:Type A} to the arc to conclude  $\big(-\infty, \min(n - 1,
n -p)\big) \subset
\SSL(K)$ provided $n \notin \{1, p\}$.  Thus to prove the
theorem it suffices to find such a $\gammatil$.

Since $e^{2 i \theta}$ is a simple root of $\Delta_K$, we have
\(\sigma_K(e^{2i\omega}) - \sigma_K(e^{2i\omega'}) = \pm 2\).  Using
Corollary~\ref{Cor: hSU=sig} and Theorem~\ref{thm: hSL+hSU}, we
compute
\begin{equation}
\label{Eq:simple root}
\pm 1 = h_\SU^c(K) - h_\SU^{c'}(K)
 = h_\SLR^{c'}(K) - h_\SLR^{c}(K)
 = \sum_{W} (n_W^{c'}-n_W^c)
\end{equation}
where the sum is over the path components $W$ of
$X^{[-2, 2], \irr}_\SLR(K)$.  Let $W$ be any such component with
$n_W^{c} \neq n_W^{c'}$, and consider the path given by
Lemma~\ref{Lem:no reducible}, which must end at the reducible
character $\chi_\theta$ where $\tr_\mu = 2 \cos \theta$.

If the closure $\Wbar$ of $W$ in $X^{[-2, 2]}_\SLR(K)$ contains a
reducible $\chi_{\theta'}$ other than $\chi_\theta$, let $\gamma$ be a
path in $W$ joining $\chi_\theta$ to $\chi_{\theta'}$.  Taking the
lift $\gammatil$ as before, note that $i^*(\gammatil)$ must exit $R$
to get to $i^*(\chitil_{\theta'})$.  As it can only do so at one of
the parabolic corners of $R$, an initial segment of $\gammatil$ is the
path we seek.  So now assume $\Wbar = W \cup \{\chi_\theta\}$, which
means $n_W^a$ can only change at $a = 2 \cos \theta$.  In particular,
$n_W^{2} = n_W^c$ and $n_W^{-2} = n_W^{c'}$.  As $n_W^c \neq n_W^{c'}$,
at least one of $n_W^2$ and $n_W^{-2}$ is nonzero. In particular, $W$
contains a parabolic character, so $\Wbar = W \cup \{\chi_\theta\}$
contains a path $\gamma$ from $\chi_\theta$ to a parabolic character
whose interior consists of irreducible nonparabolic characters.  The
lift of $\gamma$ starting at $\chitil_\theta$ is thus the path
$\gammatil$ we seek, proving the theorem.
\end{proof}

\section{Computations and conjectures}
\label{sec:computations} 

In this final section, we compute \(h\) for some additional classes of
small knots, including Montesinos knots, torus knots, and 98.7\% of
those with at most 11 crossings. We use \(h\) to give a new proof of
Riley's conjecture on parabolic representations of \2-bridge
knots. Finally, we make some conjectures about the behavior of the
extended Lin invariant \(\eh(K)\) for \2-bridge knots and knots
with lens space surgeries.

\subsection{Montesinos knots}

In Proposition~\ref{prop:sigma h}, we showed that
\(h(K) = -\frac{1}{2} \sigma(K)\) whenever \(K\) is real
representation small and \(\Sigma_2(K)\) is not LO. In fact, there are
are many knots which satisfy \(h(K) = -\frac{1}{2} \sigma(K)\) even
though \(\Sigma_2(K)\) is LO. The Montesinos knots provide a good
class of such examples.  Oertel showed a Montesinos knot is small if
and only if it has at most three rational tangles \cite{Oertel1984}. (A
Montesinos knot with one or two rational tangles is 2-bridge, hence
alternating.) We will show:

\begin{proposition}
  \label{prop: Montesinos}
  For any Montesinos knot $K = K(p_1/q_1, p_2/q_2, p_3/q_3)$, we have
  \(X^{0,\irr}_\SLR(K) = \emptyset\) and
  \(h(K) = -\frac{1}{2} \sigma(K)\).
\end{proposition}

The branched double covers of these knots are Seifert fibered spaces
over the disk with three exceptional fibres, many of which are LO and
admit nontrivial representations to \(\SLR\). As in Section~\ref{sec:
  branched}, we write \((S^3,K_n)\) to denote the orbifold where the
knot \(K\) is labeled by the cyclic group \(\Z/n\).

\begin{lemma}
  \label{lem: Ad-annoying}
  If $\rho \maps \pi_1(S^3, K_n) \to \PSLR$ is irreducible, then so is
  its restriction to $\pi_1\big(\Sigma_n(K)\big)$.
\end{lemma}

\begin{proof}
Set $\Gamma = \pi_1(S^3, K_n)$ and
$\Lambda = \pi_1\big(\Sigma_n(K)\big)$.  We will prove the
contrapositive: for any representation $\rho \maps \Gamma \to \PSLR$,
if $\rho(\Lambda)$ is reducible then so is $\rho(\Gamma)$.  We
consider the various possible restrictions of $\rho$ to $\Lambda$ from
most to least degenerate. To begin, if $\rho(\Lambda)$ is trivial,
then $\rho$ factors through $\Gamma/\Lambda = \Z/n$.  As every finite
subgroup of $\PSLR$ is reducible, we have $\rho(\Gamma)$ is reducible.

Suppose next that $\rho(\Lambda)$ is reducible but nontrivial.  The
fixed point set $F$ of $\rho(\Lambda)$ on $P^1(\C)$ must be either one
or two points, since if there were more then $\rho(\Lambda)$
would be trivial.  Since $\Lambda$ is normal in $\Gamma$, the full
group $\rho(\Gamma)$ leaves $F$ invariant.  Thus if $F$ is a single
point then $\rho(\Gamma)$ too is reducible.
So assume $F$ is two points, and let $L$ be an oriented geodesic
joining them, which may not lie in $\H^2$.  The stabilizer of $L$ in
$\PSLC$ is the abelian group $\C^\times$, so the restriction
$\rho|_\Lambda$ factors through
$\Lambda^{\mathrm{ab}} = H_1\big(\Sigma_2(K); \Z\big)$.  As $K$ is a
knot, the latter is finite of odd order.  In particular,
$\rho(\Lambda)$ is finite and moreover so is $\rho(\Gamma)$ since
$[\Gamma:\Lambda] = n$.  Thus again $\rho(\Gamma)$ must be reducible.
\end{proof}

\begin{proof}[Proof of Proposition~\ref{prop: Montesinos}]
Set $\Gamma = \pi_1(S^3, K_2)$ and
$\Lambda = \pi_1\big(\Sigma_2(K)\big)$ and suppose that some
$\rho \maps \Gamma \to \PSLR$ is irreducible.  As discussed in
\cite{Oertel1984}, the orbifold $(S^3, K_2)$ is Seifert fibered with
quotient orbifold $Q$ being a triangle with mirrored sides and
vertices labeled by dihedral groups of order $2 q_i$.  Hence if $f$
represents a regular fiber, we have:
\[
  1 \to \pair{f} \to \Gamma \to \pi_1(Q) \to 1.
\]
While $\pair{f}$ is not central in $\Gamma$, it is in $\Lambda$ since
$\Sigma_2(K)$ is Seifert fibered over the orbifold
$S^2(q_1, q_2, q_3)$ with \emph{orientable} fibers.  Thus $\rho(f)$
centralizes $\rho(\Lambda)$, which is irreducible by Lemma~\ref{lem:
  Ad-annoying}; by Lemma~\ref{lem: PSLR cent}, it follows that
$\rho(f) = 1$.  Thus on all of $\Gamma$, the representation \(\rho\)
factors through
\[
  \pi_1(Q) = \spandef{x, y, z}%
  {x^2 = y^2 = z^2 = (x y)^{q_1} = (x z)^{q_2} = (y z)^{q_3} = 1}.
\]
However, any two elements of $\PSLR$ of order two are either equal or
have product a nontrivial hyperbolic element.  This forces any
representation of $\pi_1(Q)$ to $\PSLR$ to be reducible, a
contradiction.  So \(X^{0,\irr}_\SLR(K) = \emptyset\) and
\(h(K) = -\frac{1}{2} \sigma(K)\) as claimed.
\end{proof}

\subsection{Knots with small crossing number}

There are 801 nontrivial prime knots with at most 11 crossings, of
which 601 are small \cite{BurtonCowardTillmann2013}.  All but possibly
8 of these 601 satisfy \(h = -\frac{1}{2} \sigma\). Indeed, some
85.9\% of these knots are either Montesinos or alternating
\cite{knotinfo, Castellano-MacasOwad2022}.  For the 85 knots that are
not in those classes, the question of whether \(\Sigma_2(K)\) is LO is
usually known from \cite{Dunfield2020} or can be determined using
those methods. If \(\Sigma_2(K)\) is not LO, then
\(h(K) = -\frac{1}{2} \sigma(K)\) by Proposition~\ref{prop:sigma
  h}. In this way, we can show that for small knots with \(\leq 11\)
crossings, \(h(K) = -\frac{1}{2} \sigma(K)\) for all but possibly the
following knots:
\(10_{161} = 10n31, 11n96, 11n111, 11n116, 11n135, 11n143, 11n145, \)
and \(11n183\). For these eight knots, \(\Sigma_2(K)\) is LO for all
but possibly \(\Sigma_2(11n143)\) and \(\Sigma_2(11n145)\) by
\cite{Dunfield2020}; we expect \(\Sigma_2(11n143)\) and
\(\Sigma_2(11n145)\) to be LO as they are not L-spaces.  (All eight
\(\Sigma_2(K)\) have co-orientable taut foliations by the techniques
of \cite{Dunfield2020}.)  

\subsection{Torus knots}\label{sec: ex torus}

The knots for which  we computed \(h\) in earlier sections all satisfy
\(h(K) = -\frac{1}{2} \sigma(K)\).  We now turn to torus knots and
show they exhibit more interesting behavior.  Let \(K = T(p,q)\) be
the positive \((p,q)\) torus knot. We can compute \(X_\SLR(M_K)\)
using the same method as for the trefoil in Example~\ref{Ex:Trefoil}.
The exterior of \(K=T(p,q)\) is Seifert-fibered and
\[
  \pi_1(M_K)  = \spandef{x,y,f}{x^p = f = y^q}
\]
where \(f\) is the class of the Seifert fiber. If
\(\rho\maps \pi_1(M_K) \to \SLR\) is irreducible, \(\rho(f)\) is central,
so \(\rho(x)^p = \rho(y)^q = \pm I\). As in Example~\ref{Ex:Trefoil},
we may assume that
\[
\rho(x) = \twobytwomatrix[r]{\cos \theta}{ - \sin \theta}{\sin \theta}{\cos \theta}
\mtext{and}
\rho(y) = \twobytwomatrix[r]{\cos \phi}{ -t \sin \phi}{t^{-1} \sin \phi}{\cos \phi} 
\]
For a fixed \(t\), there are \(2(p-1)(q-1)\) choices of \(\theta\) and
\(\phi\) which result in irreducible representations. The
representations for \((\theta, \phi)\) and \((-\theta, -\phi)\) have
the same character, so \(X_\SLR(M_K)\) consists of \((p-1)(q-1)\)
arcs, which are obtained by varying \(t\). Each arc tends to a
reducible character as \(t\to 1\), and its image in
\(X_\SLR(\partial \EK)\) lies on a line of slope \(-pq\), since
\(\rho(f) = \pm I\) for all \(\rho\) on a given arc and
$[f] = (pq) \mu + \lambda$ in $H_1(\partial \EK)$.

\begin{figure}
  \begin{center}
    \begin{tikzpicture}[
  scale=0.6,
  nmdstd,
  line width=1.2pt,
  every node/.style={color=black},
  every circle/.style={radius=0.12}]

  \def\basicxscale{20}
  \def\basicyscale{0.6}
  
  \coordinate (X) at (\basicxscale, 0);
  \coordinate (Y) at (0, \basicyscale);

  \draw[xstep=\basicxscale, ystep=\basicyscale, coor grid]
  ($11.5*(Y)$) grid ($(X) - 11.5*(Y)$);
  \begin{scope}[color=axesgray]
     \draw[->] (0, 0) -- ($1.04*(X)$) node[right] {$\mu^*$};
     \draw[->] ($-11.5*(Y)$) -- ($12.25*(Y)$) node[left=0.04]{$\lambda^*$};
  \end{scope}

  \begin{scope}
    \clip ($11.5*(Y)$) rectangle ($(X) - 11.5*(Y)$);
    \foreach \x in {1, 2, 3, 6, 7, 9, 11, 13, 14, 17, 18, 19}{
      \coordinate (A) at ($\x*(Y)$);
      \coordinate (B) at ($(X) + {\x - 20}*(Y)$);
      \draw[color=locus!30, dashed] (A) -- (B);
    }
  \end{scope}
  
  \foreach \x in {1, 2, 3, 6, 7, 11}{
    \coordinate (A) at ($\x*(Y)$);
    \coordinate (B) at (\x, 0);
    \draw[color=locus, line width=2pt] (A) -- (B);
    \draw[fill=otherparabolic, line width=0pt] (A) circle;
    \draw[fill=otherparabolic, line width=0pt] ($(A) + (X)$) circle;
    \draw[fill=simplealex, line width=0.5pt] (B) circle;
  }

  \foreach \x in {9, 13, 14, 17, 18, 19}{
    \coordinate (A) at ($(X) + {\x - 20}*(Y)$);
    \coordinate (B) at (\x, 0);
    \draw[color=locus, line width=2pt] (A) -- (B);
    \draw[fill=otherparabolic, line width=0pt] (A) circle;
    \draw[fill=otherparabolic, line width=0pt] ($(A) - (X)$) circle;
    \draw[fill=simplealex, line width=0.5pt] (B) circle;
  }

\end{tikzpicture}
  \end{center}
  
  \caption{This figure shows the basic structure of \(\OTEL{T(p,q)}\):
    it consists of arcs of slope $-pq$ joining a parabolic point to
    one of the $(p - 1)(q - 1)$ points of the form
    $\big(\frac{n}{pq}, 0\big)$ with neither $p$ nor $q$ dividing
    $n$. In particular, exactly one segment of the line of slope $-pq$
    through $\big(\frac{n}{pq}, 0\big)$ in the region
    $0 \leq \mu^* \leq 1$ is part of \(\OTEL{T(p,q)}\). The complete
    picture is thus determined by which parabolics of the form
    $(0, k)$ for $k > 0$ are ends of such arcs; these are
    characterized by Proposition~\ref{prop: Gamma(torus)}.  Here, the
    plot is that for $T(4, 5)$.  }
  \label{Fig:TorusKnot}
\end{figure}
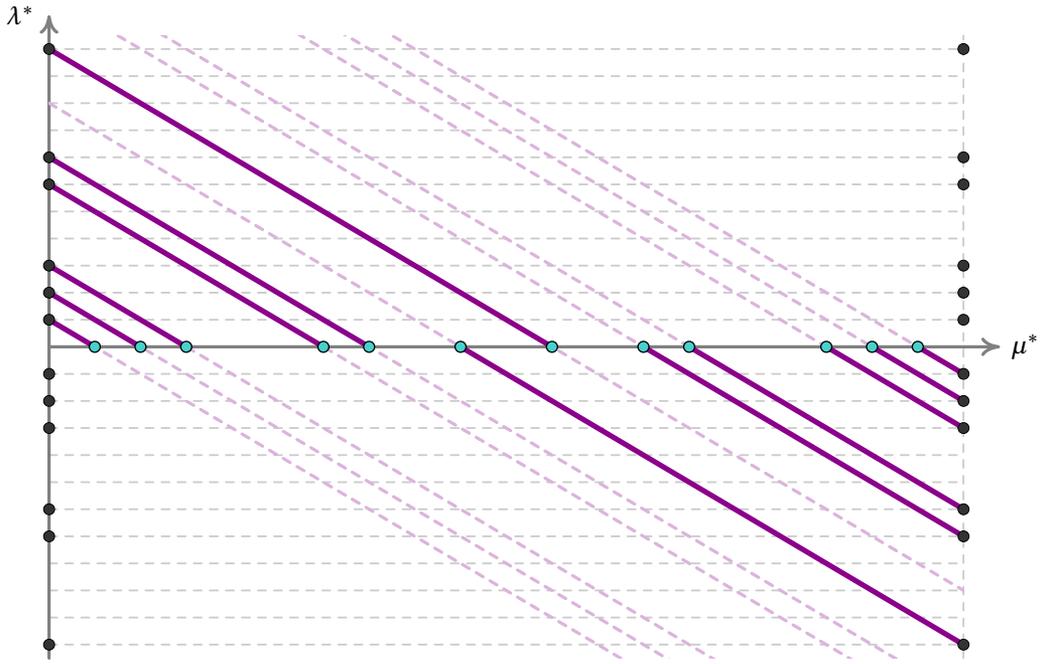

The Alexander polynomial of $K$ is
\begin{equation}
  \label{eq: torus alex}
  \Delta_K (t)  = \frac{(t^{pq}-1)(t-1)}{(t^p-1) (t^q-1)}
\end{equation}
which has \((p-1)(q-1)\) roots; specifically, for
$\zeta_{pq}= e^{\frac{2 \pi i}{pq}}$, the roots are $t = \zeta_{pq}^n$
for $1 \leq n < pq$ with $p$ and $q$ not dividing $n$.  Each root is
simple, and the number of arcs and roots is the same.  Hence arguing
as in Theorem~\ref{thm:lens arcs}, we see that the character variety
has a single arc ending at each point on the reducible line with
\(\mubar^* = \frac{n}{pq}\), for $n$ as before.  Lifting to
\(\Etil(K)\), we see that the translation extension locus has the form
shown in Figure~\ref{Fig:TorusKnot}. In order to describe \(\OTEL{K}\)
precisely, we must determine whether the arc ending at a given root of
the Alexander polynomial lies below the reducible line or above
it. Equivalently, it suffices to determine the set
\[
  \Gamma(K) = \setdef{ n \in \N}{%
    \mbox{an arc of $\Etil(K)$ ends at a parabolic of height $n$}}.
\]
The arc ending at $\mubar^* = \frac{n}{pq}$ will thus lie above the
reducible line if and only if \(n \in \Gamma(K)\).

\begin{proposition}
  \label{prop: Gamma(torus)}
  For the torus knot $K = T(p, q)$, one has
  \(\Gamma(K) = \N \setminus \Gamma_{p,q}\), where \(\Gamma_{p,q}\) is
  the additive semigroup generated by \(p\) and \(q\).
\end{proposition}

\begin{proof}
All the roots of \(\Delta_K\) are simple, so we can apply
Lemma~\ref{lem: herald power}. The arc \(\gamma\) considered in part
(a) of the lemma can be taken to have image a line segment of slope
\(-pq\) in the translation extension locus, so it satisfies
\((\tr_\lambda \circ \gamma)'(0) \neq 0\). We conclude that the
segment which limits to the reducible at \((\frac{n\pi}{pq},0)\) lies
above the reducible line if the jump in \(\sigma_K(\omega)\) at
\(\zeta_{pq}^n\) is \(-2\) and below it if the jump is
\(+2\).

Litherland \cite{Litherland1979} gave a simple formula for the sign of
the jump in \(\sigma_{K}\) at \(\zeta_{pq}^n\). It is
\begin{equation}
  \label{eq:Lit sign}
  (-1)^{\left\lfloor \frac{a}{q} \right\rfloor + \left\lfloor
      \frac{b}{p} \right\rfloor + \left \lfloor \frac{n}{pq} \right
    \rfloor}
  \mtext{where \(n = ap + bq\) for \(a, b \in \Z\).}
\end{equation}
If \(0 \leq n < pq\), there is a unique \(\nhat \in \Z\) with
\(\nhat \equiv n \bmod{pq}\) and \(\nhat = \ahat p + \bhat q\) with
\( 0 \leq \ahat < q\) and \(0\leq \bhat < p\).  Moreover, the sign in
(\ref{eq:Lit sign}) is unchanged if we replace $(a, b, n)$ by
$(\ahat, \bhat, \nhat)$.  Hence the jump at \(\zeta_{pq}^n\) is
positive if \(\nhat < pq\) and negative if \(\nhat >pq\). But
$0 \leq n < pq$ is in \(\Gamma_{p,q}\) if and only if \(n= \nhat \),
which is equivalent to \(\nhat <pq\). It follows that
\(\Gamma(K) = \setdef{n \in \N}{n \notin \Gamma_{p,q}}\).
\end{proof}

\begin{corollary}
  \label{Cor:T(p,q)}
  For $K = T(p, q)$, we have 
  \(h(K) = g(K) = (p-1)(q-1)/2\) and
  \[
    \eh(K) = \sum_{i \not \in \Gamma_{p,q}} \left(t^i + t^{-i}\right).
  \]
  In particular, \(\deg \eh(K) = 2g(K)-1\) and we have
  equality in Lemma~\ref{lem:eh properties}(\ref{item:
    Milnor-Wood}). \end{corollary}

\begin{proof} The formula for \(\eh\) follows from
Proposition~\ref{prop: Gamma(torus)} and the fact that all arcs of
\(\Etil(K)\) have the same left-to-right orientation by
Lemma~\ref{lem: herald power}; compare Figures~\ref{fig: trefoil
  final} and~\ref{Fig:TorusKnot}.  It is straightforward to show
\(\N \setminus \Gamma_{p,q}\) has \((p-1)(q-1)/2\) elements, the
largest of which is \(pq-p-q\), so the other two statements follow
from the formula for $\eh(K)$.
\end{proof}

This picture gives us a nice geometric interpretation of the signature
defect \(g(T(p,q)) + \frac{1}{2}\sigma(T(p,q))\): it is the number of
arcs in \(\OTEL{K}\) which cross the line \(\mu^*
=\frac{1}{2}\). Since \(\deg \eh = pq-p-q\) and the slope of each arc
is \(-pq\), the signature defect is \(0\) precisely when
\(1-\frac{1}{p}-\frac{1}{q} < \frac{1}{2}\); i.e.~when
\(\Sigma_2(K) = \Sigma(2,p,q)\) is Seifert fibered over a positively
curved orbifold. Recalling that \(h(K) = -\frac{1}{2}\sigma(K)\) when
\(K\) is Montesinos, we get an alternate proof of the fact that
\(T(p,q)\) is Montesinos exactly when \((p,q) = (2,2n+1), (3, 4)\), or
\((3,5)\).

\subsection{Berge knots}

Many more examples of knots with \(h(K) \neq -\frac{1}{2} \sigma(K)\)
are provided by Berge knots, and more generally, by L-space knots.  A
knot \(K \subset S^3\) is an \emph{L-space knot} when some positive
Dehn surgery on \(K\) is an L-space.  Positive torus knots, and, more
generally, any knot with a positive lens space surgery, such as the
Berge knots, are L-space knots.  For an L-space knot $K$, the surgery
\(K(p\mu + q\lambda)\) is an L-space if and only if
\(p/q \geq 2g(K) - 1\) by \cite{OzsvathSzaboLens}. The L-space
conjecture tells us to expect the same relation to hold
for left-orderings.

When \(K\) has a lens space surgery of slope \(p > 0\), the form of
\(\OTEL{K}\) is constrained by Lemma~\ref{lem:lens lines} and
Theorem~\ref{thm:lens arcs}. In particular, the lines of slope \(-p\)
through the parabolics intersect \(\OTEL{K}\) only at points that are
reducible or parabolic.  Moreover, whenever we have a unique root of
the Alexander polynomial in the interval
\(\left[\frac{k}{n},\frac{k+1}{n}\right]\), we get an arc of type
\(A_k\) in $\Etil(K)$.  For such knots, the pictures of
\cite{CullerDunfield2018} suggest that \(\OTEL{K}\) is similar to that
of a torus knot, in that it is composed entirely of arcs of type
\(A_k\). Numerical computations of \(\eh(K)\) for about 150 small
Berge knots exhibit some striking behavior, which we now describe.

The \emph{Milnor torsion} of a knot $K$ is the Laurent series
\(\tau_M(K) := \frac{\Delta_K(t)}{1-t}\) normalized so
\(\tau_M(K) = \sum_{n \geq 0} a_n t^n\) with \(a_n\neq 0\).  We say a
Laurent series \(\sum_{n \geq 0} a_n t^n\) is \emph{good} if
\(a_n \in \{0,1\}\) for all \(n\).  If \(K\) is an L-space knot, then
the Milnor torsion is good by \cite{OzsvathSzaboLens}.  We write
\(\eh_+(K)\) for the ``nonnegative exponent'' part of \(\eh(K)\);
i.e.~if \(\eh(K) = \sum_{i \geq 0 } c_i(t^i+t^{-i})\), then
\(\eh_+(K) = \sum_{i\geq 0} c_i t^i\). Based on the computations
discussed below, we posit:

\begin{conjecture} \label{conj:L-space torsion}
 If \(K\) is a small Berge knot, then
\begin{enumerate}
\item \label{item: h = r/2}
  \(h(K) = \frac{1}{2} r(K)\), where \(r(K)\) is the number of
  roots of \(\Delta_K(t)\) on the unit circle.
\item \label{item: h good}
  \(\eh_+(K)\) is good, and, moreover, 
\item \label{item: h + tau better}
  \(\eh_+(K) + \tau_M(K)\) is good. 
\end{enumerate}
\end{conjecture}

In other words, the coefficients of \(\eh_+(K)\) must all be $1$ and
fall into the ``gaps'' of the Milnor torsion
\(\tau_M(K) = \sum_{n\geq 0 } a_n t^n\) given by those \(n\) where
\(a_n = 0\). The symmetry of the Alexander polynomial implies that
\(a_n = 1-a_{2g-1-n}\), where \(g=g(K)\) is the Seifert genus of
\(K\). (Recall that Berge knots are fibered by \cite{Ni2007}, so
$\deg \Delta_K(t) = 2g$.)  Hence \(a_n=1\) for \(n \geq {2g}\), and, if
\(\tau_M(K)\) is good, precisely \(g\) of the \(2g\) coefficients
between \(a_{0}\) and \(a_{2g-1}\) will be
\(0\). Conjecture~\ref{conj:L-space torsion} predicts that
\(\frac{1}{2}r(K)\) of these \(g(K)\) gaps should be filled by the
nonzero coefficients of \(\eh_+(K)\).

When \(K=T(p,q)\) is a positive torus knot, then all the roots of
\(\Delta_K(t)\) lie on the unit circle, so \(r(K) = 2 g(K)\). In
addition, from (\ref{eq: torus alex}) we have 
\[
  \tau_M\big(T(p,q)\big) = \sum_{i \in \Gamma_{p,q}} t^i
  \mtext{and so}
  \eh_+(K) + \tau_M(K) = \sum_{i\geq 0 } t^i
\]
by Corollary~\ref{Cor:T(p,q)}. Thus Conjecture~\ref{conj:L-space
  torsion} holds in this case, and moreover all of the gaps in
$\tau_M$ are filled.

In contrast, if \(K\) is a hyperbolic Berge knot, some gaps in
$\tau_M$ must remain empty for the following reason. We have
\(a_{2g-1} = 1- a_0 = 0\), so the highest gap is always at height
\(2g(K)-1\). As every Berge knot is fibered \cite{Ni2007}, by
Lemma~\ref{lem:eh properties}(\ref{item: hyp sharp MW}) we have
\(\deg \eh(K) < 2g(K) -1 \), and the highest gap must remain unfilled.

As a consequence, if \(K\) is a hyperbolic L-space knot, the approach
of constructing left-orderings for Dehn surgeries \(K\) via \(\SLR\)
representations should not be completely effective at proving the
L-space conjecture for these manifolds. Indeed, the best expected
bound on the range of left-orderable surgery slopes coming from
Lemma~\ref{Lem:Type A} says that Dehn fillings of slope
\(< \deg \eh(K)\) should be LO.  In contrast, the L-space conjecture
says that we should expect all fillings of slope \(< 2g(K) - 1\) to be
LO.  This happens even for the $(-2, 3, 7)$ pretzel knot of
Figure~\ref{fig: pretzel constraints}, where $\SSL(K) = (-6, \infty)$
by \cite{Varvarezos2021} but we expect $\SLO(K) = (-9, \infty)$.

The slopes in this gap provide an interesting test case for the
L-space conjecture.
\begin{question}
\label{question:interval}
For some hyperbolic Berge knot \(K\), can we show that \(K(\alpha) \)
is LO for all rational \(\alpha \in (2g(K)-1-\epsilon, 2g(K)-1)\),
where \(\epsilon>0\)?
\end{question}
The best results in this direction are due to Krishna
\cite{Krishna2020} and Hu \cite{Hu2022}, whose combined work shows
that the $(2,-3,-7)$ pretzel knot has a sequence of such slopes
converging to 9.

We arrived at Conjecture~\ref{conj:L-space torsion} by considering the
158 small hyperbolic Berge knots whose exteriors can be triangulated
with at most 9 ideal tetrahedra \cite{DunfieldExceptional}.  The knots
in this sample have small hyperbolic volumes (all less than $7.1$), but
their genus is generally quite large: the median genus is $50.5$.
For these knots, a collection of parabolic representations has been
rigorously computed by Goerner \cite{PtolemyTables} using a method
based on ideal triangulations and Ptolemy coordinates; however, these
lists may not be complete due to limitations of this approach
\cite{GoernerZickert2018}.  We computed the heights of these
representations using \cite{pe} to find their possible contributions
to $\eh_+(K)$, e.g.~in Figure~\ref{fig: pretzel constraints} the
heights would be $\{1, 2, 4, 6\}$.  For all 158 knots, the parabolic
heights are all distinct and each one lies in a gap in $\tau_M(K)$,
which is consistent with parts (\ref{item: h good}) and (\ref{item: h
  + tau better}) of Conjecture~\ref{conj:L-space torsion}.

The fact that the parabolic heights fall into the gaps in
\(\tau_M(K)\) for even one of the larger of these knots is remarkable,
since on average 73\% of the gaps were filled by some parabolic
height.  For example, the manifold $o9_{03188}$, which is one of the
most complicated in the sample, has an $L(317,121)$ filling, Seifert
genus $139$, and $99$ parabolics. If we imagine that the parabolic
heights were randomly distributed in the interval \([0,2g(K)-1]\), the
probability that they all land on a gap would be
\(2^{-99} \approx 1.6 \times 10^{-30}\).  Even if we restrict to the
interval between \(1\) and \(219\) (the largest parabolic height),
where there are 122 gaps and 97 non-gaps, the probability of all the
heights landing in the gaps is less than \(10^{-23}\).

%
%
%
%

Turning to the count $r(K)$, the ratio $\frac{r(K)}{2g(K)}$ had mean
$0.74$ and was always in the interval $[0.60, 0.78]$; hence a majority
of the roots of $\Delta_K(t)$ were on the unit circle for all these
knots.  For comparison, define $h_G(K)$ to be the naive (unsigned)
count of parabolic $\SLR$ representations from
\cite{PtolemyTables}. As mentioned above, this may be an undercount of
the actual number.  In all cases, we have
$h_G(K) \leq \frac{1}{2} r(K)$ with equality in $133$ cases.  The sum
of all the $h_G(K)$ was $6{,}479$, whereas the sum of all
$\frac{1}{2} r(K)$ was $6{,}532$.  The largest difference
$\frac{1}{2} r(K) - h_G(K)$ observed was $4$.  The manifold
$o9_{03188}$ mentioned is one of two such examples as there
$h_G(K) = 99$ but $\frac{1}{2} r(K) = 103$.  The simplest example
where the counts differ is $v1392$ where $h_G(K) = 14$ but
$\frac{1}{2} r(K) = 16$.  Computations of $\OTEL{K}$ using \cite{pe}
are consistent with a pair of parabolic representations being missed.
Indeed, using the alternate triangulation
\texttt{kLvLAPQkaedgijhijijnxsvepnxxti\_abba} for $v0220$, we found
there are in fact $16$ distinct parabolic representations as predicted
by $r(K)$.

Finally, for 62 of the 158 manifolds, we have a plot generated by
\cite{pe} as part of \cite{CullerDunfield2018}; these all consist of
non-intersecting arcs of type \(A_k\) for $k \neq 0$.  For a more
complicated example than Figure~\ref{fig: pretzel constraints}, see
\cite[Figure~4]{CullerDunfield2018} for the plot for $v0220$ which has
$h_G(K) = \frac{1}{2}r(K) = 36$ with 47 gaps in $\tau_M$ and
$\OTEL{K}$ consisting of 72 arcs, making it a median-complexity
example in the overall sample.  Using Lemma~\ref{lem: herald power},
the pictures provide strong evidence that all parabolics from
\cite{PtolemyTables} contribute to $\eh_+(K)$ with positive signs,
further confirming Conjecture~\ref{conj:L-space torsion}.

\begin{figure}
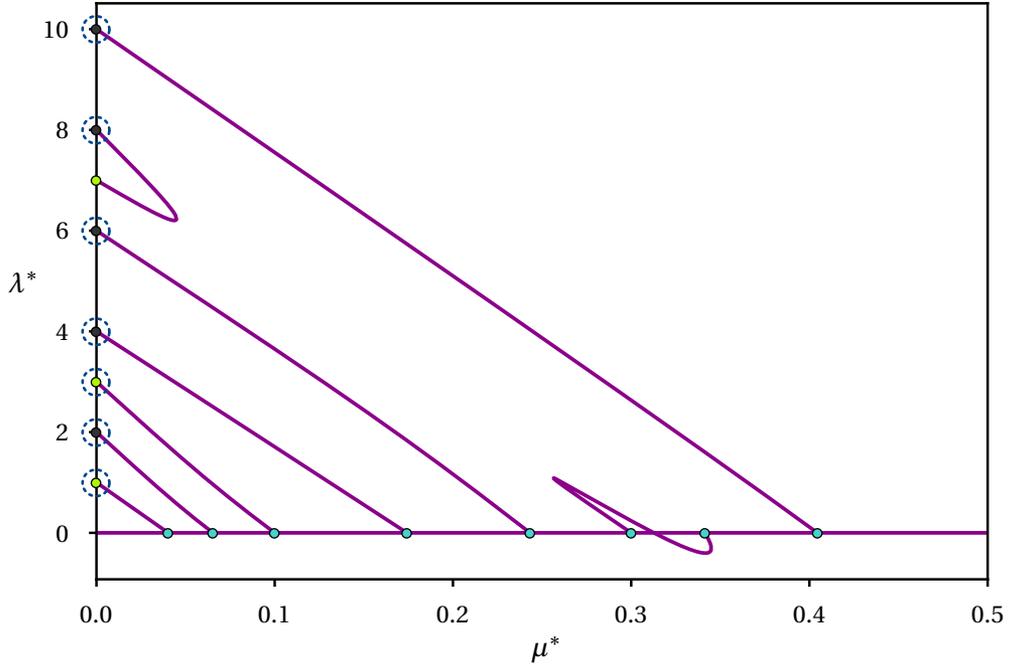

  \begin{center}
    \input figures/t09882
  \end{center}

  \vspace{-0.5cm}
  
  \caption{The manifold $t09882$ is the exterior of the mirror of an
    L-space knot $K$, half of whose $\OTEL{K}$ is shown here. The
    lighter dots on the vertical axis correspond to parabolics that
    are Galois conjugates of the holonomy representation of the
    hyperbolic structure on $M$; the darker dots on that axis are
    other parabolics.  A dotted circle about a parabolic point
    indicates a gap in $\tau_M$; further gaps at heights 12 and 17 are
    not shown. Here,
    $\Delta_K = (t^4 - t^3 + t^2 - t + 1) (t^{14} - t^{12} + t^7 - t^2
    + 1)$ has 16 roots on the unit circle.  The knot $K$ is not a
    Berge knot, having no lens space surgeries, and appears to satisfy
    none of the conclusions of Conjecture~\ref{conj:L-space torsion},
    having $h(K) = 6 < 8 = \frac{1}{2}r(K)$ and $\htil_+(K)$ is not
    good since the coefficients $t^{7}$ and $t^{8}$ should have
    opposite signs.}
  \label{fig: t09882}
\end{figure}

\begin{figure}
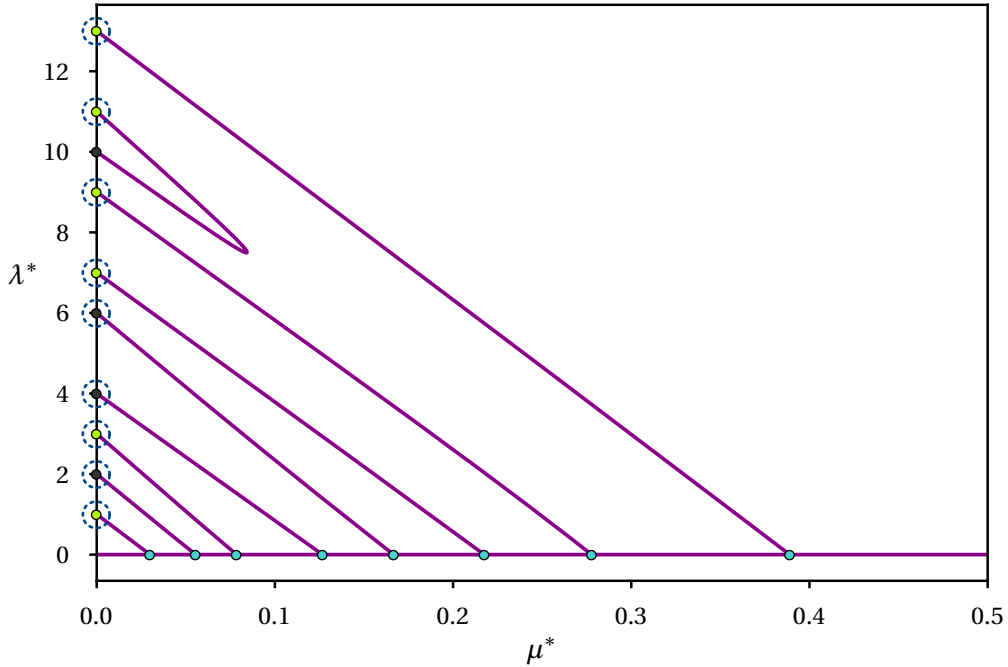

  \begin{center}
    \input figures/t08114
  \end{center}

  \vspace{-0.5cm}

  \caption{The manifold $t08114$ is the exterior of the mirror of an
    L-space knot $K$, half of whose $\OTEL{K}$ is shown here using the
    same conventions as Figure ~\ref{fig: t09882}. Here,
    $\Delta_K = (t^2 - t + 1)(t^6 - t^3 + 1)(t^{16} - t^{14} + t^{11}
    - t^9 + t^8 - t^7 + t^5 - t^2 + 1)$ has sixteen roots on the unit
    circle and $\tau_M$ has additional gaps at heights 15, 18, and
    23.  The knot $K$ appears to satisfy only part (\ref{item: h =
      r/2}) of Conjecture~\ref{conj:L-space torsion}, having
    $h(K) = 8 = \frac{1}{2}r(K)$ but $\htil_+(K)$ is not good since
    the coefficients $t^{10}$ and $t^{11}$ should have opposite
    signs.}
  \label{fig: t08114}
\end{figure}

\begin{remark}
  The pattern described in Conjecture~\ref{conj:L-space torsion} holds
  for many, but not all, L-space knots; see Figures~\ref{fig: t09882}
  and \ref{fig: t08114} for two examples where it appears not to hold.
  Caveat: The qualification in the last sentence is because pictures
  produced by \cite{pe} are not completely rigorous; see Section~5.1
  of \cite{CullerDunfield2018} for details.
\end{remark}

 \subsection{2-bridge knots}
\label{sec: almost riley}

If \(K\) is 2-bridge, it is alternating and small, so
\(h(K) = - \frac{1}{2} \sigma(K) \) by Corollary~\ref{cor: h for
  alt}. Using Corollary~\ref{Cor:parabolics}, we give a new proof of
the following result, which was conjectured by Riley almost 50 years
ago \cite{Riley1972, Riley1984} and proved recently by Gordon
\cite{Gordon2017}:

\begin{theorem}[(Riley Conjecture)] If \(K\) is a 2-bridge knot then
  \(\pi_1(\EK)\) admits at least $\abs{\sigma(K)}$ nonconjugate
  irreducible parabolic representations into \(\SLR\).
\end{theorem}
The \emph{Riley polynomial} \(p_K(y)\) introduced in \cite{Riley1972}
plays a central role in \cite{Riley1972,Gordon2017}. It is a
1-variable polynomial whose roots determine the parabolic
representations of a 2-bridge knot \(K\).  Riley proved that \(p_K\)
has only simple roots, which is a key ingredient in the proof
below.

\begin{proof}
We will show that \(\pi_1(M_K)\) admits at least
\(|h(K)| = \frac{1}{2}|\sigma(K)|\) conjugacy classes of representations
\(\rho\) with \(\tr_\mu(\rho) = 2\); the statement then follows from
the bijection $b \maps X_\SLR^{2,\irr}(K) \to X_\SLR^{-2,\irr}(K)$
discussed in Section~\ref{subsec: X_K sym}.
  
Since \(K\) is 2-bridge, we can write \(K = \widehat{\beta}\), where
\(K\) is a 4-strand braid.  We choose $\beta$ so that we can orient
\(K\) with all the underbridges and overbridges coherently
oriented from left to right. Our usual setup leads us to consider the
character varieties
\(L_i = X^{2,\irr}_\SLR(H_i) \subset X^{2,\irr}_\SLR(S_4)\), where
\(X^{2,\irr}_\SLR(S_4)\) is a real 2-dimensional surface (as in
Figure~\ref{fig: X(S_4)}(a)), and each \(L_i\) is a smooth curve in
it. Then \(X^{2,\irr}_\SLR(K) = L_1 \cap L_2\) and
\(h(K) = \pm \pair{L_1, L_2}\). To prove the theorem, it suffices to
show that \(L_1\) is transverse to \(L_2\), since if so each point of
$L_1 \cap L_2$ contributes $\pm 1$ to $h(K)$ and hence
$\abs{h(K)} \leq \#(L_1 \cap L_2)$.

Here is an outline of our strategy for proving that $L_1 \cap L_2$ is
transverse.  Consider the representation varieties
\(N_i =R^{2,\irr}_\SLR(H_i) \subset R^{2,\irr}_\SLR(S_4) \subset
R^{2,\irr}_\SLR(F_4) \).
Riley \cite{Riley1972} defines a map \(\gamma \maps \R_{>0} \to N_1 \)
such that 1)
\(\tau \circ \gamma \maps \R_{>0} \to X^{2,\irr}_\SLR(S_4)\) is an
embedding whose image is \(L_1' \subset L_1\) and 2) if
\(p \in L_1 \cap L_2\), then \(p \in L_1'\). Since \(L_1\) is
\1-dimensional, it is enough to show that for all
\(p = \tau\circ \gamma (u)\) in \(L_1 \cap L_2\), one has
\(T_p L_1 \not \subset T_p L_2\).  By the path lifting properties of
\cite[Lemma~2.11]{CullerDunfield2018}, this is equivalent to
\(\gamma'(u) \not \in T_{\gamma(u)} N_2\).  
To show the
latter, we will define a map \(f \maps R^{2,\irr}(F_4) \to \R\) such
that 3) \(f(N_2) = \{0\}\) and 4) \(f \circ \gamma(u) = up_K(u^2)\), where \(p_K\) is the Riley polynomial. Since all the roots of \(p_K\) 
are simple \cite[Theorem 3]{Riley1972}, it will follow that
\(\gamma'(u) \not \in T_qN_2\) whenever \(q = \gamma(u) \in N_2\).

It remains to define the maps \(\gamma\) and \(f\) and check claims
1)--4).  To start, let \(C \subset \SLR\) be the set of non-central
matrices of trace \(2\). It has two connected components corresponding
to whether a parabolic element rotates (really translates) its invariant
horocircles in $\H^2$ clockwise or anticlockwise.  Conjugating by
\(\SLR\) preserves each component of \(C\), but the two components are
exchanged by the action of the non-identity component of \(\SLRpm\).
We can identify \(Z = R^{2,\irr}_\SLR(F_2)\) and
\(Y = R^{2,\irr}_\SLR(F_4)\) with open subsets of \(C^2\) and \(C^4\),
respectively, so that \(N_1\) is the image of the map \(i:Z \to Y\)
given by \(i(A,B) = (A,A^{-1},B,B^{-1})\).

Following Riley, for $u \in \R_{>0}$ we define
\[
\gamma(u) = \big(A(u), \, A(u)^{-1}, \, B(u), \, B(u)^{-1}\big) \mtext{where}
A(u) = \twobytwomatrix{1}{u}{0}{1},
B(u) = \twobytwomatrix{1}{0}{-u}{1}.
\]
The image of \(\gamma\) is contained in \(N_1\), so
\(\tau \circ \gamma: \R_{>0} \to L_1\).  The map \(\tau \circ \gamma\)
is easily seen to be injective by considering
\(\tr\big( A(u) B(u) \big) = 2 - u^2\), and the same calculation shows
its derivative is injective as well. This establishes claim 1).

Next, we consider the image $L_1'$ of \(\tau \circ \gamma\). Recall
that \(C\) has two connected components that can be interchanged by
$\SLRpm$.  Consequently, \(Z\) has four connected components and
\(X^{2,\irr}_\SLR(F_2)\) has two, where the components of
\(X^{2,\irr}_\SLR (F_2)\) correspond to whether the invariant horoballs of
the two generators are rotated in the same or opposite directions.
The image $L_1'$ is contained in the former component since
conjugating by $\mysmallmatrix{0}{1}{-1}{0} \in \SLR$ interchanges
$A(u)$ and $B(u)$; moreover, it is all of that component by Lemma 2
of \cite{Riley1972}.  Turning to 2), recall \(\beta\) was chosen so
that the underbridges are oriented left to right, so the generators
\(t_1\) and \(t_2\) of $\pi_1(H_1)$ are both positively oriented
meridians for $K$; in particular, they are conjugate in
\(\pi_1(\EK)\).  Thus any $\rho \in N_1 \cap N_2$ gives rise to
$\pi_1(\EK) \to \SLR$ where $\rho(t_1)$ and $\rho(t_2)$ are conjugate
in $\SLR$, implying that the corresponding point in $L_1 \cap L_2$ is
in $L_1'$. This proves claim 2).
  
A plat diagram representing \(K\) determines a Heegaard splitting of
the complement, and hence a presentation of \(\pi_1(\EK)\) with one
generator for each underbridge and one relation for each overbridge.
The relations are redundant: any one relation is a consequence of the
others. If the plat diagram is chosen so that all of the underbridges
and overbridges are coherently oriented, this presentation is closely
related to the presentation obtained by starting with
\(\pi_1(S_{2n})\) and adding the relations \(s_{2i-1} = s_{2_i}^{-1}\)
and \(\beta^{-1}_*(s_{2i-1}) =\beta^{-1}_*(s_{2i}^{-1})\).
To be precise, if we
use the first set of relations to eliminate the generators \(s_{2i}\),
we get the presentation described above up to the operation of cyclically
permuting the elements in each relator.  In the case of a coherently
oriented 2-bridge knot \(K_{p/q}\) with \(q\) odd, the resulting
presentation is a 2-generator 1-relator presentation
\(\pi_1(\EK) = \spandef{ a, b }{r(a,b)} \) where
\(r(a,b)= waw^{-1}b^{-1}\) and \(w=w(a,b)\) is determined by \(p/q\).
 A precise formula for
$w(a,b)$ is given in Proposition 1 of \cite{Riley1972}. 

We now explain how this relates to \(N_1 \cap N_2\). Suppose
\(q = (A,A^{-1},B,B^{-1}) \in N_1\). Then
\((\beta^*)^{-1}(q) = (C_1,C_2,C_3,C_4)\), where each \(C_i\) is a
word in \(A\) and \(B\) determined by \(\beta\). The condition that
\(q\in N_2\) is equivalent to the relations \(C_1 = C_2^{-1}\) and
\(C_3 = C_4^{-1}\), and by the discussion above, each of these are
equivalent to the relation \(r(A,B)\). It follows that there are words
\(v_1 = v_1(A,B)\) and \(v_2 = v_2(A,B)\) such that \(v_1C_1v_2 = wA\)
and \(v_1C_2^{-1}v_2 = Bw\).

Suppose \(q = (A_1,A_2,A_3,A_4) \in Y\). If
\((\beta^*)^{-1}(q) = (c_1,c_2,c_2,c_4)\), where the $c_i$ are words
in $A_1,A_2,A_3,A_4$, we define
\(F:C^4 \to M_{2\times 2}(\R)\) by
\[F(q) = v_1(A_1,A_3) (c_1 - c_2^{-1}) v_2(A_1,A_3).\] Then if
\(q \in N_2\), we have \(c_1 = c_2^{-1}\), so \(F(q) = 0\). We define \(f(q)\)
to be the upper right entry of the matrix \(F(q)\); claim 3) follows
immediately.  Finally, we compute
\[F(\gamma(u)) = v_1(u)(C_1(u) - C_2^{-1}(u)) v_2(u) = w(u)A(u) -
  B(u) w(u) ,\] where we write \(v_1(u) = v_1(A(u),B(u))\) etc. In
\cite[Theorem 2]{Riley1972}, Riley computes the quantity on the right
and shows it has the form
\(\mysmallmatrix{0}{up_K(u^2)}{up_K(u^2)}{0}\), where \(p_K(y)\) is the
Riley polynomial.This proves claim 4) and hence completes the proof
of the theorem.
\end{proof}

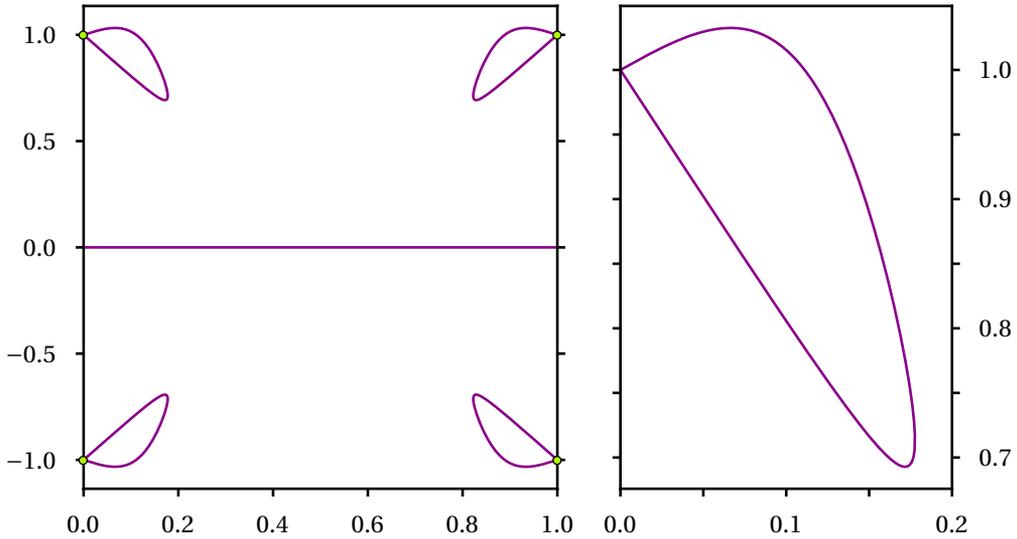
\begin{figure}
  \begin{center}
    \begin{tikzoverlay}[width=0.95\textwidth]{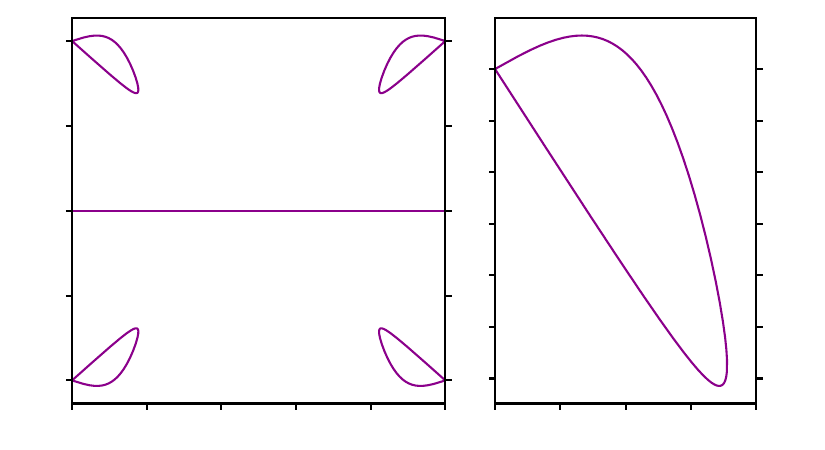}
  \begin{scope}[nmdstd]
    \draw (8.732323, 4.000000) node[below] {$0.0$};
    \draw (17.781937, 4.000000) node[below] {$0.2$};
    \draw (26.831551, 4.000000) node[below] {$0.4$};
    \draw (35.881165, 4.000000) node[below] {$0.6$};
    \draw (44.930778, 4.000000) node[below] {$0.8$};
    \draw (53.980392, 4.000000) node[below] {$1.0$};
    \draw (7.090909, 8.433681) node[left] {$-1.0$};
    \draw (7.090909, 18.718103) node[left] {$-0.5$};
    \draw (7.090909, 29.002525) node[left] {$0.0$};
    \draw (7.090909, 39.286947) node[left] {$0.5$};
    \draw (7.090909, 49.571370) node[left] {$1.0$};
    \draw (60.003119, 4.000000) node[below] {$0.0$};
    \draw (75.839944, 4.000000) node[below] {$0.1$};
    \draw (91.676768, 4.000000) node[below] {$0.2$};
    \draw (93.318182, 8.664514) node[right] {$0.7$};
    \draw (93.318182, 21.168313) node[right] {$0.8$};
    \draw (93.318182, 33.672111) node[right] {$0.9$};
    \draw (93.318182, 46.175910) node[right] {$1.0$};
    \begin{scope}[shift={(8.73232323, 29.00252525)}]
      \coordinate (X) at (45.24806892, 0);
      \coordinate (Y) at (0, 20.56884445);
      \foreach \x in {0, 1}{
        \foreach \y in {-1, 1}{
          \draw[line width=0.5pt, radius=1.5pt, fill=galoisgeom]
          ($\x*(X) + \y*(Y)$) circle;
        }
      }
    \end{scope}
    \begin{scope}[shift={(60.00311943, -78.86207689)},
      xscale=158.36824124, yscale=125.03798694]
    \end{scope}
  \end{scope}
\end{tikzoverlay}
  \end{center}

  \vspace{-0.25cm}

  \caption{The 2-bridge knot $K_{41/9} = K10a107 = 10_{17}$ is an
    example where Conjecture~\ref{conj: enhanced riley} requires
    cancellation of distinct parabolics when computing $\eh(K)$. Here
    $\sigma_K = 0$, but $X^{2, \irr}_\SLR(K)$ consists of four points,
    all at height $\pm 1$. Assuming everything is transversely cut out,
    either orientation on the lobe of $\OTEL{K}$ shown at right
    results in $\eh(K) = 0$ as predicted.  Plot made using \cite{pe}.
  }
  \label{fig: K10a107}
\end{figure}

The torus knot \(T(p,q)\) is 2-bridge only if \(p=2\). If
\(K = T(2,2n+1)\), Corollary~\ref{Cor:T(p,q)} tells us that
\(\eh(K) = t^{-2n+1} + t^{-2n+3}+ \cdots + t^{2n-3}+ t^{2n-1}\). In
general, we conjecture that if \(K\) is 2-bridge, then
\(\eh(K)=\eh\big(T(2,2n+1)\big)\), where
\(-2n = \sigma(K) = \sigma\big(T(2,2n+1)\big)\). In other words:

\begin{EnhancedRiley}
  \label{conj: enhanced riley}
  If \(K\) is a 2-bridge knot, we conjecture that 
  \[
    \eh(K) =- \frac{t^{\sigma(K)} - t^{-\sigma(K)}}{t-t^{-1}}.
  \]
\end{EnhancedRiley}

Note that this is compatible with Riley's theorem \cite{Riley1972}
that if \(K\) is 2-bridge, any parabolic representation
\(\rho\maps \EK \to \SLR\) satisfies \(\tr \rho(\lambda) = -2\);
equivalently, any parabolic must have an odd height in the translation
extension locus. By using the Riley polynomial to find all the
parabolic $\SLR$ representations of a 2-bridge knot \(K_{p/q}\), we
checked numerically that the conjecture holds ``modulo \(2\)'' for
all 2-bridge knots \(K_{p/q}\) with \(p<500\). That is, the number of
parabolics at an odd height \(i\) is odd if \(|i| < |\sigma(K)|\) and
is even otherwise.  There are $12{,}929$ knots in this range after
accounting for $(p, q) \sim (p, \pm q^{\pm 1})$.  More than 20\% have
strictly more than $\sigma(K)$ parabolics, and so would require
cancellation when computing $\eh(K)$ in order for
Conjecture~\ref{conj: enhanced riley} to hold.  A simple example of
this is shown in Figure~\ref{fig: K10a107}.  A complicated example is
the 32-crossing knot $K = K_{479/29}$ where $\sigma = 2$ and so we
expect $\eh_+ = t$. In fact, there actually 55 parabolics that
contribute to $\eh_+$; the \emph{unsigned} count of these parabolics
is $21 t + 16 t^3 + 10 t^5 + 6 t^7 + 2 t^9$, which is indeed equal to
$t$ modulo 2.

\begin{figure}
  \begin{center}
    \begin{tikzoverlay}[width=0.9\textwidth]{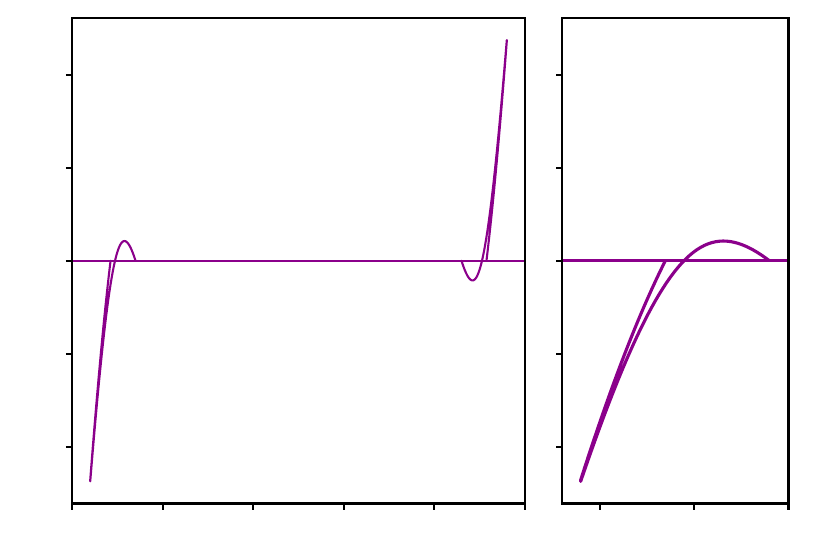}
  \begin{scope}[TEL figure]
    \draw (8.755051, 4.000000) node[below] {$0.0$};
    \draw (19.726094, 4.000000) node[below] {$0.2$};
    \draw (30.697138, 4.000000) node[below] {$0.4$};
    \draw (41.668182, 4.000000) node[below] {$0.6$};
    \draw (52.639226, 4.000000) node[below] {$0.8$};
    \draw (63.610269, 4.000000) node[below] {$1.0$};
    \draw (7.113636, 12.522526) node[left] {$-0.4$};
    \draw (7.113636, 23.792828) node[left] {$-0.2$};
    \draw (7.113636, 35.063131) node[left] {$0.0$};
    \draw (7.113636, 46.333434) node[left] {$0.2$};
    \draw (7.113636, 57.603737) node[left] {$0.4$};
    \draw (72.723204, 4.000000) node[below] {$0.05$};
    \draw (84.151375, 4.000000) node[below] {$0.10$};
    \draw (95.579545, 4.000000) node[below] {$0.15$};
    \foreach \x in {13.399508495990055, 16.426023262355727, 55.939296617644274, 58.96581138400995}{
      \draw[line width=0.2pt, radius=1pt, fill=simplealex] (\x, 35.06313132) circle;
    }

    \foreach \x in {80.64694194354745, 93.2574201358182}{
      \draw[line width=0.5pt, radius=2pt, fill=simplealex] (\x, 35.06313132) circle;
    }
    \begin{scope}[shift={(8.75505051, 35.06313132)}, xscale=54.85521886, yscale=56.35151439]
    \end{scope}
    \begin{scope}[shift={(61.29503367, 35.06313132)}, xscale=228.56341190, yscale=56.35151439]
    \end{scope}
  \end{scope}
\end{tikzoverlay}
  \end{center}

  \vspace{-0.5cm}

  \caption{The 2-bridge knot $K_{61/21} = K14a2459$ is an example
    where $\OTEL{K}$ contains an arc that starts and ends at
    reducibles but whose interior consists of irreducibles.  While
    such arcs are extremely common in $X_\SU(K)$, this seems to be the
    first such observed in $X_\SLR(K)$. Here,
    $\Delta_K = 5 t^4 - 15 t^3 + 21 t^2 - 15 t + 5$ which has all
    roots on the unit circle.  Moreover, $h(K) = 0$ as
    $X^{2,\irr}_\SLR(K)$ is empty, and $\sigma_K = 0$ except in
    intervals between each pair of close roots where it is $2$.
    Lemma~\ref{lem: herald power} forces the arc to cross the
    horizontal axis. In the closeup at right, the two regions enclosed
    by $\OTEL{K}$ must have equal area: the derivative of the
    Chern-Simons invariant/Seifert volume/Godbillon-Vey invariant is
    essentially $\eta = x \dy - y \dx$, see \cite{Khoi2003}, so as
    $d\eta = 2 \dx \, \dy$, the loop created from the small segment of
    the axis between the two roots and the curved parts of $\OTEL{K}$
    must have signed area $0$.  Plot made using \cite{pe}.}
  \label{fig: K14a2459}
\end{figure}

\begin{figure}
  \begin{center}
    \begin{tikzoverlay}[width=0.95\textwidth]{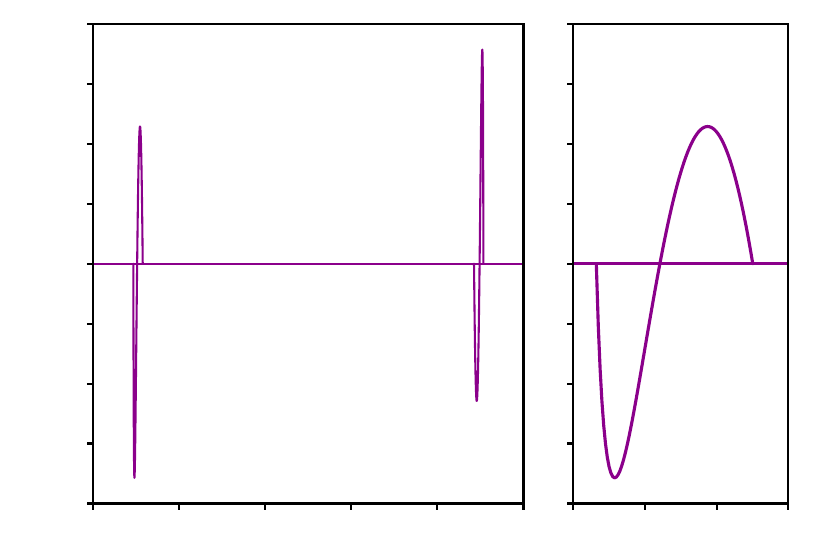}
  \begin{scope}[TEL figure]
    \draw (11.300505, 4.000000) node[below] {$0.0$};
    \draw (21.731650, 4.000000) node[below] {$0.2$};
    \draw (32.162795, 4.000000) node[below] {$0.4$};
    \draw (42.593939, 4.000000) node[below] {$0.6$};
    \draw (53.025084, 4.000000) node[below] {$0.8$};
    \draw (63.456229, 4.000000) node[below] {$1.0$};
    \draw (9.659091, 5.641414) node[left] {$-0.004$};
    \draw (9.659091, 12.905934) node[left] {$-0.003$};
    \draw (9.659091, 20.170455) node[left] {$-0.002$};
    \draw (9.659091, 27.434975) node[left] {$-0.001$};
    \draw (9.659091, 34.699495) node[left] {$0.000$};
    \draw (9.659091, 41.964015) node[left] {$0.001$};
    \draw (9.659091, 49.228535) node[left] {$0.002$};
    \draw (9.659091, 56.493056) node[left] {$0.003$};
    \draw (9.659091, 63.757576) node[left] {$0.004$};
    \draw (69.490320, 4.000000) node[below] {$0.09$};
    \draw (78.182941, 4.000000) node[below] {$0.10$};
    \draw (86.875561, 4.000000) node[below] {$0.11$};
    \draw (95.568182, 4.000000) node[below] {$0.12$};
    \foreach \x in {16.1621877808475, 17.2998073220793, 57.4569266879207, 58.5945462291525}{
      \draw[line width=0.2pt, radius=1pt, fill=simplealex] (\x, 34.69949495) circle;
    }
    \foreach \x in {72.2847795179100, 91.2451052036539}{
      \draw[line width=0.5pt, radius=2.5pt, fill=simplealex] (\x, 34.69949495) circle;
    }
    \begin{scope}[shift={(11.30050505, 34.69949495)},
      xscale=52.15572391, yscale=7264.52020278]
    \end{scope}
    \begin{scope}[shift={(-8.74326599, 34.69949495)},
      xscale=869.26206510, yscale=7264.52020278]
    \end{scope}
  \end{scope}
\end{tikzoverlay}
  \end{center}

  \vspace{-0.5cm}

  \caption{The 2-bridge knot $K_{77/20} = K12a380$ is another example
    where $\OTEL{K}$ contains an arc that starts and ends at
    reducibles, but whose interior consists of irreducibles (compare
    Figures~\ref{fig: t09882} and~\ref{fig: K14a2459}).  It is
    remarkable how small $\SSL(K)$ is here, contained in
    $[-0.022, 0.038]$.  Here,
    $\Delta_K = (2t^2 - 3t + 2) (3t^2 - 5t + 3)$ which has all 
    roots on the unit circle.  Moreover, $h(K) = 0$ as $X^2_\SLR(K)$
    is empty, and $\sigma_K = 0$ except in intervals between each pair
    of close roots where it is $2$. Plot made using \cite{pe}.  }
  \label{fig: K12a380}
\end{figure}

Although alternating and Montesinos knots also satisfy
\(h(K) = - \frac{1}{2} \sigma(K)\), Conjecture~\ref{conj: enhanced
  riley} does not extend to either class of knots: Goerner's tables of
parabolic representations \cite{PtolemyTables} contain knots in each
of these classes that admit parabolic representations with
\(\tr \rho(\lambda) = 2\).  Finally, we note that although for small
values of \(p\) the translation extension locus of \(K_{p/q}\)
contains only arcs of type \(A_k\) (joining reducibles to parabolics),
there are 2-bridge knots for which the translation
extension locus provably contains arcs joining two reducibles. Two
such examples are shown in Figures~\ref{fig: K14a2459} and~\ref{fig:
  K12a380}; we found them by looking for knots where the number of
real roots of the Riley polynomial is smaller than the number of roots
of the Alexander polynomial on the unit circle.

{\RaggedRight
\bibliographystyle{nmd/math} 
\bibliography{signature}
}

\pagebreak

\end{document}